\documentclass[twoside]{amsbook}
 \pdfoutput=1
\usepackage[utf8]{inputenc}
\usepackage[T1]{fontenc}
\usepackage[english]{babel}
\usepackage{amsmath,amssymb,amscd}
\usepackage{bbm}

\expandafter\let\expandafter\oldproof\csname\string\proof\endcsname
\let\oldendproof\endproof
\renewenvironment{proof}[1][\proofname]{%
  \oldproof[\it\hspace*{-2.2em} #1]%
}{\oldendproof}

\newenvironment{miabstract}%
{\cleardoublepage \null \vfill \begin{center}%
\bfseries \abstractname \end{center}}%
{\vfill \null }

\usepackage[pdftex,svgnames]{xcolor}

\usepackage[bookmarks=true,hyperindex,pdftex,colorlinks, citecolor=blue,linkcolor=blue,urlcolor=blue,linktocpage=true]{hyperref}

\usepackage{fancybox}

\makeindex

\theoremstyle{plain}
\newtheorem{Theo}{Theorem}[chapter]
\newtheorem{Prop}[Theo]{Proposition}
\newtheorem{Coro}[Theo]{Corollary}
\newtheorem{Lemm}[Theo]{Lemma}

\theoremstyle{definition}
\newtheorem{Defi}[Theo]{Definition}
\newtheorem{Obse}[Theo]{Observation}
\newtheorem{Exam}[Theo]{Example}
\newtheorem{Prob}[Theo]{Problem}

\theoremstyle{remark}
\newtheorem{Rema}[Theo]{Remark}

\newcommand{\conv}{\operatorname{conv}}
\newcommand{\aconv}{\operatorname{aconv}}
\newcommand{\Lip}{\operatorname{Lip}}
\newcommand{\Id}{\operatorname{Id}}
\newcommand{\ext}{\operatorname{ext}}
\newcommand{\eps}{\varepsilon}
\newcommand{\Real}{\operatorname{Re}}
\newcommand{\re}{\operatorname{Re}}
\newcommand{\Imag}{\operatorname{Im}}

\newcommand{\supp}{\operatorname{supp}}
\newcommand{\dist}{\operatorname{dist}}

\newcommand{\spn}{\operatorname{span}}

\newcommand{\diam}{\operatorname{diam}}

\newcommand{\GS}{\operatorname{gSlice}}
\newcommand{\GF}{\operatorname{gFace}}

\newcommand{\Slice}{\operatorname{Slice}}
\newcommand{\Face}{\operatorname{Face}}

\newcommand{\Spear}{\operatorname{Spear}}

\newcommand{\interior}{\operatorname{int}}
\newcommand{\slope}{\operatorname{slope}}

\newcommand{\AAA}{\mathcal{A}}

\newcommand{\FF}{\mathcal{F}}

\newcommand{\N}{\mathbb{N}}
\newcommand{\R}{\mathbb{R}}
\newcommand{\C}{\mathbb{C}}
\newcommand{\T}{\mathbb{T}}
\newcommand{\K}{\mathbb{K}}
\newcommand{\weakstar}{w^\ast}

\numberwithin{section}{chapter}
\numberwithin{equation}{chapter}

\title[Spear operators between Banach spaces]{Spear operators between Banach spaces}

\author[Kadets]{Vladimir Kadets}
\address[Kadets]{\newline School of Mathematics and Informatics,\newline V.~N.~Karazin Kharkiv National University,\newline
pl. Svobody 4, 61022 Kharkiv, Ukraine.\newline
\href{http://orcid.org/0000-0002-5606-2679}{ORCID: \texttt{0000-0002-5606-2679}}} \email{v.kadets@karazin.ua}

\author[Mart\'{\i}n]{Miguel Mart\'{\i}n}
\address[Mart\'{\i}n]{\newline Departamento de An\'{a}lisis Matem\'{a}tico,\newline Facultad de
 Ciencias,\newline Universidad de Granada,\newline 18071 Granada, Spain.\newline
\href{http://orcid.org/0000-0003-4502-798X}{ORCID: \texttt{0000-0003-4502-798X} }
 }
\email{mmartins@ugr.es}

\author[Mer\'{\i}]{Javier Mer\'{\i}}
\address[Mer\'{\i}]{\newline Departamento de An\'{a}lisis Matem\'{a}tico,\newline Facultad de
 Ciencias,\newline Universidad de Granada,\newline 18071 Granada, Spain.\newline
\href{http://orcid.org/0000-0002-0625-5552}{ORCID: \texttt{0000-0002-0625-5552} }
 }
\email{jmeri@ugr.es}

\author[P\'erez]{Antonio P\'erez}
\address[P\'erez]{\newline Departamento de Matem\'{a}ticas,\newline Universidad de Murcia,\newline
30100 Espinardo (Murcia), Spain.\newline
\href{http://orcid.org/0000-0001-8600-7083}{ORCID: \texttt{0000-0001-8600-7083} }} \email{antonio.perez7@um.es}

\thanks{The research of the first author is done in frames of Ukrainian Ministry of Science and Education Research Program 0115U000481, and it was partially done during his stay in Murcia under the support of MINECO/FEDER projects MTM2014-57838-C2-1-P and partially during his visits to the University of Granada which were supported by the Spanish MINECO/FEDER project MTM2015-65020-P. The second and the third authors were partially supported by Spanish MINECO/FEDER project  MTM2015-65020-P, and by Junta de Andaluc\'{\i}a and FEDER grant FQM-185. The fourth author was partially supported by the MINECO/FEDER project MTM2014-57838-C2-1-P and a Ph.D.\ fellowship of ``La Caixa Foundation''.}



\date{January 2017}

\begin{document}

\begin{titlepage}
\setcounter{page}{1}
\newlength{\centeroffset}
\setlength{\centeroffset}{-0.5\oddsidemargin}
\addtolength{\centeroffset}{0.5\evensidemargin}
\thispagestyle{empty} \vspace*{\stretch{1}}
\noindent\hspace*{\centeroffset}\makebox[0pt][l]{\begin{minipage}{\textwidth}
\flushright
{\Huge\bfseries \textcolor[rgb]{0.00,0.00,.50}{Spear operators\\ between Banach spaces}

}
\noindent\textcolor[rgb]{0.00,0.50,1.00}{\rule[-1ex]{\textwidth}{5pt}}\\[2.5ex]
\hfill\emph{\Large \textcolor[rgb]{0.00,0.00,0.50}{Monograph}}
\end{minipage}}

\vspace{\stretch{1}}
\noindent\hspace*{\centeroffset}\makebox[0pt][l]{\begin{minipage}{\textwidth}
\flushright {\bfseries \Large
\textcolor[rgb]{0.00,0.00,1.00}{Vladimir Kadets\linebreak Miguel Mart\'{\i}n\linebreak Javier Mer\'{\i}\linebreak Antonio P\'{e}rez} \\[3ex]}
January 2017
\end{minipage}}

\vspace{\stretch{2}}
\newpage
$ $ \thispagestyle{empty}

\end{titlepage}

\maketitle

\begin{miabstract}
The aim of this manuscript is to study \emph{spear operators}: bounded linear operators $G$ between Banach spaces $X$ and $Y$ satisfying that for every other bounded linear operator $T:X\longrightarrow Y$ there exists a modulus-one scalar $\omega$ such that
$$
\|G + \omega\,T\|=1+ \|T\|.
$$
To this end, we introduce two related properties, one weaker called the alternative Daugavet property (if rank-one operators $T$ satisfy the requirements), and one stronger called lushness, and we develop a complete theory about the relations between these three properties. To do this, the concepts of spear vector and spear set play an important role. Further, we provide with many examples among classical spaces, being one of them the lushness of the Fourier transform on $L_1$. We also study the relation of these properties with the Radon-Nikod\'{y}m property, with Asplund spaces, with the duality, and we provide some stability results. Further, we present some isometric and isomorphic consequences of these properties as, for instance, that $\ell_1$ is contained in the dual of the domain of every real operator with infinite rank and the alternative Daugavet property, and that these three concepts behave badly with smoothness and rotundity. Finally, we study Lipschitz spear operators (that is, those Lipschitz operators satisfying the Lipschitz version of the equation above) and prove that (linear) lush operators are Lipschitz spear operators.
\end{miabstract}

\tableofcontents

\chapter{Introduction}

\section{Motivation}
Let $X$ and $Y$ be Banach spaces. The main goal of this manuscript is to study bounded linear operators $G:X\longrightarrow Y$ such that for every other bounded linear operator $T:X\longrightarrow Y$ there is a modulus-one scalar $\omega$ such that the norm equality
\begin{equation}\label{eq:spear-operator}
\|G + \omega\,T\|=1+ \|T\|
\end{equation}
holds. In this case, we say that $G$ is a \emph{spear operator}. When $X=Y$ and $G=\Id$ is the identity operator, the study of this equation goes back to the 1970 paper \cite{D-Mc-P-W}. It was proved there that for $T:X\longrightarrow X$ bounded and linear, the existence of a modulus-one scalar $\omega$ such that the norm equality
\begin{equation}\label{eq_0.1}
\| \Id + \omega\,T\|=1+\|T\|
\end{equation}
holds, is equivalent to the equality between the numerical radius of $T$ and its norm. The list of spaces for which the identity is a spear operator (they are called spaces with numerical index one) contains all $C(K)$ spaces and $L_1(\mu)$ spaces, as well as some spaces of analytic functions and vector-valued functions, which motivated the intensive study of this class of spaces in the past decades. Let us say that, as we will see here, the extension to general operators produces other important (and very different from the identity) examples of operators which are spear. One of the most striking one is the Fourier transform on the $L_1$ space on a locally compact Abelian group (see subsection \ref{ssecFourier}).

The concepts of numerical range and numerical radius of operators, and the one of numerical index of Banach spaces, played an important role in operator theory, particularly in the study and classification of operator algebras (see the fundamental Bonsall and Duncan books \cite{B-D1,B-D2}, the survey paper \cite{KaMaPa}, and the sections \S2.1 and \S2.9 of the recent book \cite{Cabrera-Rodriguez}). For general operators $G$, a concept of numerical range of operators with respect to $G$ has been recently introduced \cite{Ardalani}, and there is a relation between  \eqref{eq:spear-operator} and numerical ranges, analogous to the case of \eqref{eq_0.1}. For instance, a spear operator is geometrically unitary in the strongest possible form. These concepts provide also a natural motivation for the study of spear operators; we will give a short account on this in section~\ref{subsect:numericalranges}.

The property of $G$ being a spear operator is formulated in terms of all bounded linear operators between two Banach spaces, which leads to several difficulties for its study in abstract spaces, and also  in some concrete spaces. It would be much more convenient to have a geometric definition of this property (in terms of $G$), but unfortunately until now a description of this property in pure geometrical terms has not been discovered, even for the case when $G=\Id$. In order to manage this difficulty for this case, two other Banach space properties were introduced: the alternative Daugavet Property (aDP for short) and lushness (see \cite{martinOikhberg} and \cite{NumIndexDuality}, respectively). These two properties are of geometric nature, the aDP is weaker and lushness is stronger than the fact that the identity is a spear operator. On the other hand, in some classes of Banach spaces these properties are equivalent (say, in Asplund spaces, in spaces with the Radon-Nikod\'{y}m property and, more generally, in SCD spaces introduced in \cite{SCDsets}). The study of these two properties has been crucial in the development of the theory in the case when $G=\Id$. We will give a short account on this in section~\ref{subsec:thecaseoftheId}, where we also present the main examples and results about SCD spaces.

So, now naturally appears the task of re-constructing the theory of the aDP and lushness in such a way that it could be applied to spear operators. We will do so here. On this way, we not only transfer the known results to the new setting, but in fact make much more. Namely, we introduce a unified approach to a huge number of previously known results, substantially simplify the system of notations and, in many cases, present the general results for operators $G$ in a more clear way than it was done earlier for the identity operator. To do so it has been crucial to study the new concept of target operator (see section \ref{sec:StrongTargetOperators}), which plays the same role here that the concept of strong Daugavet operator plays for the study of the Daugavet property \cite{KadSW2} and of Daugavet centers \cite{BosGnar}. Let us also mention that the aDP, to be a target operator and lushness are separably determined properties, while we do not know whether the concept of spear operator is. This separable determination allows us to use the full power of the theory of SCD sets and operators, a task which will be crucial in the development of the subject.

Finally, another motivation to do this job was the potential applicability of the extended theory to the study of non-linear Lipschitz maps, which would allow to extend some recent results from \cite{Wang-Huang-Tan} and \cite{LipschitzSlices}. In particular, this will allow to reprove the results of the latter papers in a more reasonable (i.e.\ linear) way. Namely, the standard technique of Lipschitz-free spaces (see \cite{GodKalt} or \cite{Weaver}) reduces equation \eqref{eq_0.1} for a non-linear Lipschitz map $T: X \longrightarrow X$ to an analogous equation for the linearization of $T$, but this linearization acts from the Lipschitz-free space $\FF(X)$ to $X$. Hence, in order to use this technique, we are in need of studying equation \eqref{eq:spear-operator} instead of \eqref{eq_0.1}, that is, to study spear operators between two different spaces. See section \ref{sec:Lipschitz} for more details.

\section{Notation, terminology, and preliminaries} By $\mathbb{K}$ we denote the scalar field ($\mathbb{R}$ or $\mathbb{C}$), and we use the standard notation $\T:=\{ \lambda \in \mathbb{K} \colon |\lambda|=1 \}$ for its unit sphere. We use the letters $X$, $Y$, $Z$ for Banach spaces over $\mathbb{K}$ and by subspace we always mean closed subspace. In some cases, we have to distinguish between the real and the complex case, but for most results this difference is insignificant. The closed unit ball and unit sphere of $X$ are denoted respectively by
$B_{X}$ and $S_X$. We denote the Banach space of all bounded linear operators from $X$ to $Y$ by $L(X,Y)$, and write $L(X)$ for $ L(X,X)$. The identity operator is denoted by $\Id$, or $\Id_X$ if it is necessary to precise the space. The dual space of $X$ is denoted by $X^\ast$, and $J_X:X\longrightarrow X^{\ast\ast}$ denotes the natural isometric inclusion of $X$ into its bidual $X^{\ast\ast}$. For $x_0\in X^\ast$ and $y_0\in Y$, we write $x_0^\ast\otimes y_0$ to denote the rank-one operator given by $[x_0^\ast\otimes y_0](x)=x_0^\ast(x) y_0$ for every $x\in X$. We write $X=Y\oplus_1 Z$ and $X=Y\oplus_\infty Z$ to mean that $X$ is, respectively, the \emph{$\ell_1$-sum} and the \emph{$\ell_\infty$-sum} of $Y$ and $Z$. In the first case, we say that $Y$ is an \emph{$L$-summand} of $X$; in the second case, we say that $Y$ is an \emph{$M$-summand} of $X$.
\index{l1sum@$\ell_1$-sum}%
\index{linfinitysum@$\ell_\infty$-sum}%
\index{lsummand@$L$-summand}%
\index{msummand@$M$-summand}%

For a subset $A \subset X$ we write
$\| A \| := \sup{\{ \| x\| \colon x \in A \}}$ if $A$ is bounded and $\|A\|=\infty$ if it is unbounded.
Observe that this function has the following properties:
\begin{equation*}
\| \lambda A\| = |\lambda| \, \| A\|,\quad
\| A + B\| \leqslant \| A - C \| +  \|C + B\|, \quad
\| A-B\| \geqslant \bigl| \| A- C\| - \| B - C \| \bigr|,
\end{equation*}
for every $\lambda\in \K$ and every subsets $A$, $B$, $C$ of $X$. The diameter of a (bounded) set $A \subset X$ can be calculated as $\diam(A)= \| A- A\|$. We will also use the notation $$\| F \pm x\| := \max{\{ \| F+x\|, \| F-x\| \}}.$$ For a subset $A \subset X$ and for $x \in X$ we write
$$
\T A:=\{ \omega a \colon \omega \in \T,\, a \in A \} \quad \text{and}\quad  \T x := \{ \omega x \colon \omega \in \T \}.
$$
A subset $A$ of $X$ is said to be \emph{rounded} if $\T A=A$.

Given $A \subset X$, we denote by $\conv{(A)}$ the \emph{convex hull} of $A$, and by $\aconv{(A)}$ its \emph{absolutely convex hull}, i.e.\  $\aconv{(A)}=\conv{(\T A)}$. We say that $A \subset B_{X}$ is \emph{norming} for $Z \subset X^{\ast}$ if for every $f \in Z$ we have that $\| f\|=\sup_{x \in A}{|f(x)|}$ or, equivalently, $B_{X}=\overline{\aconv}^{\sigma(X,Z)}(A)$.
\index{rounded}%
\index{convex hull}%
\index{convex hull@$\conv(\cdot)$}%
\index{absolutely convex hull}%
\index{absolutely convex hull@$\aconv(\cdot)$}%
\index{norming}%

Given $x^{\ast} \in X^{\ast}$ and $\alpha > 0$, we put
\[
\Slice(A,x^{\ast}, \alpha):= \left\{ x \in A\colon \Real{x^{\ast}(x)} > \sup\nolimits_A \bigl(\Real x^{\ast}\bigr) - \alpha\right\},
\]
and we say that this is a \emph{slice} of $A$. A \emph{face} of $A$ is a (non-empty) subset of the form
\index{slice}%
\index{slice1@$\Slice(\cdot,\cdot,\cdot)$}%
\index{face1@$\Face(\cdot,\cdot)$}%
\index{face}%
$$
\Face(A,x^\ast):=\left\{ x \in A\colon \Real{x^{\ast}(x)} = \sup\nolimits_A \bigl(\Real x^{\ast}\bigr) \right\},
$$
where $x^\ast\in X^\ast$ is such that its real part attains its supremum on $A$. If $A \subset X^{\ast}$ and the functional defining the slice or the face is taken in the predual, i.e.\  $x^\ast=x \in X\equiv J_X(X) \subset X^{\ast \ast}$, then $\Slice(A, x, \alpha)$ is called a \emph{$\weakstar$-slice} of $A$ and $\Face(A,x)$ is called a \emph{$\weakstar$-face} of $A$.
\index{wstarslice@$\weakstar$-slice}%
\index{wstarface@$\weakstar$-face}%

Given $B \subset B_{X}$, $F \subset B_{X^\ast}$ and $\eps > 0$, we define
\[
\GS(B,F,\eps):= \left\{ x \in B \colon \sup_{x^\ast \in F}{\Real{x^\ast}(x)} > 1 - \eps \right\}
\]
and we call it a \emph{generalized slice} of $B$ (observe that it is a union of slices when non-empty). We also define
\index{generalized slice}%
\index{generalized slice@$\GS(\cdot,\cdot,\cdot)$}%
\index{generalized face}%
\index{generalized face@$\GF(\cdot,\cdot)$}%
\[
\GF(B,F):=\left\{ x \in B \colon \sup_{x^\ast \in F}{\Real{x^\ast}(x)}= 1 \right\},
\]
and call it a \emph{generalized face} of $B$. If $F = \{ z^\ast\}$, then we will simply write $\GS(B,F,\eps) = \GS(B,z^\ast, \eps)$ (which is a slice of $B$ when non-empty). The following easy results about generalized slices will be frequently used.

\begin{Rema}
Let $X$ be a Banach space, let $B\subset B_X$ be a rounded set, let $A\subset B_{X^\ast}$ be a set, and let $z^\ast\in S_{X^\ast}$. Then:
\begin{enumerate}
\item[(a)] $\aconv{\: \GS(B,A, \eps)} = \conv{\: \GS(B,\T A, \eps)}$;
\item[(b)] if $\re z^\ast$ attains its supremum on $B$, then the set $$\GF(B,\T\,z^\ast)=\bigl\{x\in B\colon |z^\ast(x)|=1\bigr\}$$ coincides with $\T \Face(B,z^\ast)$.
\end{enumerate}
\end{Rema}

Given $A\subset X$, we denote by $\ext A$ the set of extreme points of $A$. When $A$ is convex and it is compact in a locally convex topology, then the set of its extreme points has many good topological properties. Here, we will repeatedly use the following ones.

\index{extreme@$\ext$}%
\begin{Lemm}\label{Kreinlemma}
Let $X$ be a Banach space and let $A\subset X^\ast$ convex and weak$^\ast$-compact.
\begin{enumerate}
\item[(a)] (Choquet's Lemma) If $x^\ast \in \ext A$, then for every weak$^\ast$-neighborhood $U$ of $x^\ast$ in $A$ there is a weak$^\ast$-slice $S$ of $A$ such that $x^\ast \in S \subset U$. In other words, the weak$^\ast$-slices of $A$ containing $x^\ast$ form a base of the relative weak$^\ast$-neighborhoods of $x^\ast$ in $A$.
\item[(b)] (Milman's Theorem) If $D \subseteq A$ satisfies that $\overline{\conv}^{w^{\ast}}(D) \supseteq A$, then  $\overline{D}^{w^{\ast}}\supseteq  \ext A$.
\item[(c)] $(\ext A, \weakstar)$ is a Baire space, so the intersection of every sequence of G$_\delta$ dense subsets of $\ext A$ is again (G$_\delta$) dense.
\end{enumerate}
\end{Lemm}

Assertion (a) can be found in \cite[p.~107]{ChoquetLecturesV2}; (b) is an immediate consequence of (a) and Hahn-Banach separation Theorem; (c) appears in \cite[p.~146, Theorem 27.9]{ChoquetLecturesV2}.

Recall that a Banach space $X$ is said to be \emph{strictly convex} if $\ext B_{X} = S_{X}$, and \emph{smooth} if the mapping $x\longmapsto \|x\|$ is G\^{a}teaux differentiable at every point of $X\setminus \{0\}$ (equivalently, for each $0 \neq x \in X$ there is a unique $x^{\ast} \in S_{X^{\ast}}$ with $x^{\ast}(x) = \|x\|$).
\index{strictly convex}%
\index{smooth}%
If moreover, the mapping $x\longmapsto \|x\|$ is Fr\'{e}chet differentiable at every point of $X\setminus \{0\}$, then $X$ is said to be \emph{Fr\'{e}chet smooth}.
\index{Fr\'{e}chet smooth}%
\index{strongly extreme}%
A point $x\in B_X$ is said to be \emph{strongly extreme} if given a sequence $(y_n)$ in $X$ such that $\left\|x\pm y_n\right\|\longrightarrow 1$, we have that $\lim y_n = 0$. A point $x\in B_X$ is \emph{denting} if it belongs to slices of $B_X$ of arbitrarily small diameter. If $X$ is a dual space and the slices can be taken to be weak$^\ast$-open, then the point is called \emph{weak$^\ast$-denting}. Observe that denting points are strongly extreme and strongly extreme points are extreme points, and none of the implications reverses in general (see \cite{KunenRosenthal}, for instance).
\index{denting}%
\index{weak$^\ast$-denting}%
\index{CKX@$C(K,X)$}%
\index{L1muX@$L_1(\mu,X)$}%
\index{LinftymuXinfty@$L_\infty(\mu,X)$}%
Finally, we recall some common notation for spaces of vector-valued function spaces. Given a compact Hausdorff topological space $K$ and a Banach space $X$, $C(K,X)$ is the Banach space of all continuous functions from $K$ into $X$ endowed with the supremum norm. Given a positive measure space $(\Omega,\Sigma,\mu)$ and a Banach space $X$, $L_\infty(\mu,X)$ is the Banach space of all (clases of) measurable functions from $\Omega$ into $X$ which are essentially bounded, endowed with the essential supremum norm; $L_1(\mu,X)$ is the Banach space of all (clases of) Bochner-integrable functions from $\Omega$ into $X$, endowed with the integral norm.

\section{A short account on the results for the identity}\label{subsec:thecaseoftheId}
Our goal here is to briefly present some results about Banach spaces with numerical index one and the related properties aDP and lushness. For more information and background, we refer the reader to \cite{SCDsets, LushNumerical, NumIndexDuality, ChoiGarKimMaestre_JMAA2014, D-Mc-P-W,  KadThes, KaMaPa,KSSW, martinOikhberg} and references therein. Let us start with the main definitions.

\begin{Defi}
Let $X$ be a Banach space.
\begin{enumerate}
\item[(a)] $X$ has \emph{numerical index one} \cite{D-Mc-P-W} if $\Id_X$ is a spear operator, i.e.\
    $
    \|\Id + \T\,T\|=1 + \|T\|
    $
    for every $T\in L(X)$.
\item[(b)] $X$ has the \emph{Daugavet property} (\emph{DPr} in short) \cite{KSSW} if
    $
    \|\Id + T\|=1 + \|T\|
    $
    for every rank-one $T\in L(X)$.
\item[(c)] $X$ has the \emph{alternative Daugavet property} (\emph{aDP} in short) \cite{martinOikhberg} if the equality
    $
    \|\Id + \T\,T\|=1 + \|T\|
    $ holds for every rank-one $T\in L(X)$.
\item[(d)] $X$ is \emph{lush} \cite{NumIndexDuality} if for every $x_0\in B_X$, every $\eps>0$, and every $x\in S_X$, there exists $x^{\ast} \in \Slice(B_{X^\ast},x, \eps)$ such that
$
\dist\bigl(x_{0}, \aconv(\GS(S_{X},x^{\ast}, \eps))\bigr) < \eps.
$
\end{enumerate}
\end{Defi}
\index{numerical index one}%
\index{Daugavet property}%
\index{DPr}%
\index{alternative Daugavet property}%
\index{aDP}%
\index{lush}%

These four properties are related as follows:
$$
\begin{CD}
\Ovalbox{lush} @>>>  \Ovalbox{numerical index one} @>>> \Ovalbox{aDP} @>>> \Ovalbox{DPr}.
\end{CD}
$$
The first implication appears in \cite[Proposition~2.2]{NumIndexDuality}, while the second and the third ones are obvious. None of them reverses in general (see \cite[Remarks 4.2.(a)]{LushNumOneDual}, \cite[Remark 2.4]{martinOikhberg}, and \cite[Example 3.2]{martinOikhberg}), and the DPr and numerical index one are not related. Examples of lush spaces include, among others, $L_1(\mu)$ spaces and their isometric preduals (so, in particular, $C(K)$ spaces), uniform algebras, and finite codimensional subspaces of $C[0,1]$ (see \cite{LushNumerical, ChoiGarKimMaestre_JMAA2014,KaMaPa}). They also are the main examples of Banach spaces with numerical index one. On the other hand, $C([0,1],X)$ and $L_1([0,1],X)$ have the DPr (and so the aDP) no matter the range space $X$ \cite{KSSW}, while they have numerical index one if and only $X$ does (see e.g.\ \cite[\S 1]{KaMaPa}). This provides many examples of spaces having the aDP but failing to have numerical index one as, for instance, $C([0,1],\ell_2)$ and $L_1([0,1],\ell_2)$. Finally, $C([0,1],\ell_2)\oplus_\infty c_0$ is an example of a space with the aDP, but failing the DPr and not having numerical index one.

One of the milestones of the theory was reached in the 2010 paper \cite{SCDsets}, where a general condition was given to make properties (a), (c), and (d) equivalent, the SCD property.

\begin{Defi}[\mbox{\cite{SCDsets}}]
Let $X$, $Y$ be Banach spaces. A bounded subset $A\subset X$ is said to be \emph{slicely countably determined} (\emph{SCD} in short) if there exists a countable family $\{S_n\colon n\in \N\}$ of slices of $A$ such that $A \subset \overline{\conv}{B}$ whenever $B \cap S_{n} \neq \emptyset$ for every $n \in \N$. The space $X$ is said to be \emph{SCD} if every convex bounded subset of $X$ is SCD. Finally, a bounded linear operator $T\in L(X,Y)$ is an \emph{SCD operator} if $T(B_X)$ is an SCD  subset of $Y$.
\end{Defi}
\index{slicely countably determined}%
\index{SCD}%
\index{SCD space}%
\index{SCD operator}%

While the definition of SCD set is valid for arbitrary bounded sets, it was only introduced for convex bounded sets in \cite{SCDsets} since it is mainly used for the image of the unit ball by bounded linear operators. The next result contains the main examples of (convex) SCD sets, spaces, and operators.

\begin{Exam}[\mbox{\cite{SCDsets}}]
\label{Exam:SCD-set-spaces-operators} Let $X$, $Y$ be Banach spaces.
\begin{enumerate}
\item[(a)] A convex bounded subset $A$ of $X$ is SCD provided:
    \begin{enumerate}
    \item[(a.1)] $A$ is separable and has the convex point of continuity property; in particular, $A$ is separable and has the Radon-Nikod\'{y}m property.
    \item[(a.2)] $A$ is separable and it does not contain $\ell_1$-sequences; in particular, $A$ is separable and Asplund.
    \end{enumerate}
\item[(b)] The following conditions on $X$ imply that every separable subspace of $X$ is SCD:
    \begin{enumerate}
    \item[(b.1)] $X$ has the convex point of continuity property; in particular, $X$ has the Radon-Nikod\'{y}m Property.
    \item[(b.2)] $X$ does not contain copies of $\ell_1$; in particular, $X$ is Asplund.
    \end{enumerate}
\item[(c)] If a Banach space has the DPr, then its unit ball is not an SCD set. Actually it is shown in \cite[Example 2.13]{SCDsets} that given $x_{0} \in S_{X}$ and a sequence of slices $(S_{n})_{n \in N}$ of $S_{X}$ we can find $x_{n} \in S_{n}$ for each $n \in \N$ so that $x_{0}$ does not belong to the closed linear hull of $\{ x_{n} \colon n \in \N\}$. In particular, the unit balls of $C[0,1]$ and $L_1[0,1]$ are not SCD sets.
\item[(d)] The following conditions on an operator $T\in L(X,Y)$ guarantee that its restriction to every separable subspace of $X$ is SCD:
    \begin{enumerate}
    \item[(d.1)] $T$ does not fix any copy of $\ell_1$.
    \item[(d.2)] $T(B_X)$ has the convex point of continuity property; in particular, $T(B_X)$ has the Radon-Nikod\'{y}m Property.
    \end{enumerate}
\end{enumerate}
\end{Exam}

The main applications of the SCD property in our context are the following.

\begin{Prop}[\mbox{\cite[\S 4]{SCDsets}}]
Let $X$ be a Banach space. If $X$ has the aDP and $B_Y$ is SCD for every separable subspace $Y$ of $X$, then $X$ is lush. In particular, if $X$ has the aDP and it has the convex point of continuity property, the Radon-Nikod\'{y}m Property, or it does not contain copies of $\ell_1$, then $X$ is lush.
\end{Prop}

\begin{Prop}[\mbox{\cite[\S 5]{SCDsets}}]
Let $X$ be a Banach space with the aDP. For every $T\in L(X)$ such that $T(B_Y)$ is SCD for every separable subspace $Y$ of $X$, one has $
\|\Id +\T\,T\|=1+\|T\|.
$
In particular, this happens if $T(B_X)$ has the convex point of continuity property, or the Radon-Nikod\'{y}m Property, or if $T$ does not fix copies of $\ell_1$.
\end{Prop}

Observe that SCD sets are separable, and we do not know whether numerical index one is separably determined. However, the aDP and lushness are, and this is crucial in the way the results above were proved.

Let us finally say that lushness is surprisingly related to the study of Tingley's problem about extensions of surjective isometries between unit spheres of Banach spaces \cite{TanHuangLiu-MU} and to the study of norm attaining operators \cite{ChoiKim-lush, KimLeeMartin-Studia}.

\section{A brushstroke on numerical ranges and numerical indices}\label{subsect:numericalranges} Our aim here is to give a short account on numerical ranges which will be very useful for the motivation and better understanding of the concepts of spear vector and spear operator. The study of the numerical range of an operator started with O.~Toeplitz's field of values of a matrix of 1918, a concept which quickly extended to bounded linear operators on Hilbert spaces. An extension of it to elements of unital Banach algebras was used in the 1950's to relate geometrical and algebraic properties of the unit (that the unit is a vertex of the unit ball of every complex unital Banach algebra), and in the developing of Vidav's characterization of $C^\ast$-algebras. Later on, in the 1960's, F.~Bauer and G.~Lumer gave independent but related extensions of Toeplitz's numerical range to bounded linear operators on arbitrary Banach spaces, which do not use the algebraic structure of the space of operators. All these notions are essential to define and study when an operator on a general Banach space is hermitian, skew-hermitian, dissipative\ldots We refer the reader to the monographs by F.~Bonsall and J.~Duncan \cite{B-D1,B-D2} for background. In the 1985 paper \cite{MarMenaPayaRod1985}, an abstract notion of numerical range, which had already appeared implicitly in the 1950's paper \cite{Bohn-Karlin}, was developed. We refer to sections \S2.1 and \S2.9 of the very recent book \cite{Cabrera-Rodriguez} for more information and background.

Given a Banach space $Z$ and a distinguished element $u\in S_Z$, we define the \emph{numerical range} and the \emph{numerical radius} of $z\in Z$ \emph{with respect to} $(Z,u)$ as, respectively,
$$
V(Z,u,z):=\bigl\{z^\ast(z)\colon z^\ast\in \Face(S_{Z^\ast},u)\bigr\}, \quad v(Z,u,z):=\sup\bigl\{|\lambda|\colon \lambda \in V(Z,u,z)\bigr\}.
$$
\index{numerical range}%
\index{V@$V(\cdot,\cdot,\cdot)$}%
\index{numerical radius}%
\index{v@$v(\cdot,\cdot,\cdot)$}%
Then, the \emph{numerical index} of $(X,u)$ is
\index{numerical index}%
\index{numerical index@$N(\cdot,\cdot)$}%
$$
N(X,u):=\inf\bigl\{v(Z,u,z)\colon z\in S_Z\bigr\} = \max\bigl\{k\geqslant 0\colon k\|z\|\leqslant v(Z,u,z)\ \forall z\in Z\bigr\}.
$$
With this notation, $u\in S_Z$ is a \emph{vertex} if $v(Z,u,z)\neq 0$ for every $z\in Z\setminus\{0\}$ (that is, $\Face(S_{Z^\ast},u)$ separates the points of $Z$), $u\in S_Z$ is a \emph{geometrically unitary} element if the linear hull of $\Face(S_{Z^\ast},u)$ equals the whole space $Z^\ast$ or, equivalently, if $N(X,u)>0$ (see e.g.\ \cite[Theorem~2.1.17]{Cabrera-Rodriguez}). By Definition \ref{def-spear}, $u\in S_Z$ is a \emph{spear vector} if $\|u + \T\,z\|=1+\|z\|$ for every $z\in Z$, and this is equivalent, by Hahn-Banach Theorem, to $N(Z,u)=1$. So, spear vectors are geometrically unitary elements (in the strongest possible way!), geometrically unitary elements are vertices, and vertices are extreme points. None of this implications reverses (see \cite[\S 2.1]{Cabrera-Rodriguez}). Let us also comment that the celebrated
Bohnenblust-Karlin Theorem \cite{Bohn-Karlin} states that algebraic unitary elements of a unital complex Banach algebra (i.e.\ elements $u$ such that $u$ and $u^{-1}$ have norm one) are geometrically unitary, see \cite{Rodriguez-JMAA-unitaries} for a detailed account on this.
\index{vertex}%
\index{geometrically unitary}%

We then have an easy way to define the numerical range of an operator: given a Banach space $X$ and $T\in L(X)$, the \emph{algebra numerical range} or \emph{intrinsic numerical range} of $T$ is just $V(L(X),\Id,T)$. As in order to study this concept we have to deal with the (wild) dual of $L(X)$, there are other concepts of numerical range which simplify such a task. The (Bauer) \emph{spatial numerical range} of $T$ is defined as
\index{algebra numerical range}%
\index{intrinsic numerical range}%
\index{spatial numerical range}%
$$
W(T):=\bigcup\nolimits_{x\in S_X} V(X,x,Tx)\,=\,\bigl\{x^\ast(Tx)\colon x\in S_X,\, x^\ast\in S_{X^\ast},\, x^\ast(x)=1\bigr\}.
$$
It is a classical result that the two numerical ranges are related as follows:
$$
\overline{\conv}\, W(T) = V(L(X),\Id,T)
$$
(see e.g.\ \cite[Proposition~2.1.31]{Cabrera-Rodriguez}), so they produce the same numerical radius of operators.

Let now $X$, $Y$ be Banach spaces and let us deal with numerical ranges with respect to $G\in L(X,Y)$ with $\|G\|=1$. First, the \emph{intrinsic numerical range} of $T\in L(X,Y)$ \emph{with respect to} $G$ is easy to define: just consider
\index{intrinsic numerical range with repect to an operator}%
$$
V(L(X,Y),G,T)=\bigl\{\Phi(T)\colon \Phi\in L(X,Y)^\ast,\, \|\Phi\|=\Phi(G)=1\bigr\},
$$
and so we have the corresponding numerical radius $v(L(X,Y),G,T)$ and numerical index $N(L(X,Y),G)$. But this definition forces us to deal with the dual of $L(X,Y)$, which is not a nice task. On the other hand, the possible extension of the definition of spatial numerical range has many problems as, for instance, it is empty if  $G$ does not attain its norm; moreover, even in the case when $G$ is an isometric embedding, it does not have a good behaviour, see \cite{MarMerPay}. Very recently, a new notion has appeared \cite{Ardalani}: the \emph{approximated spatial numerical range} of $T\in L(X,Y)$ \emph{with respect to} $G$ is defined by
\index{approximated spatial numerical range with respect to an operator}%
\index{W@$\widetilde{W}_G(\cdot)$}%
$$
\widetilde{W}_G(T):= \bigcap\nolimits_{\eps>0} \overline{\bigl\{y^\ast(Tx)\colon y^\ast\in S_{Y^\ast},\, x\in S_X,\, \re y^\ast(Gx)>1-\eps\bigr\}}.
$$
We then have the corresponding numerical radius and numerical index:
$$
v_G(T)=\sup\bigl\{|\lambda|\colon \lambda\in \widetilde{W}_G(T)\bigr\}, \quad n_G(X,Y)=\inf\bigl\{v_G(T)\colon T\in L(X,Y),\,\|T\|=1\bigr\}.
$$
The relationship between these two numerical ranges is analogous to the one for the identity operator \cite[Theorem~2.1]{Mar-numrange-JMAA2016}:
$$
\conv\,\widetilde{W}_G(T) = V(L(X,Y),G,T)
$$
for every norm-one $G\in L(X,Y)$ and every $T\in L(X,Y)$.
Therefore, both concepts produce the same numerical radius of operators and so, the same numerical index of $G$, the same concepts of vertex and geometrically unitary elements. In particular, both numerical ranges produce the same concept of spear operator: $G\in L(X,Y)$ is a spear operator if and only if $N\bigl(L(X,Y),G\bigr)=1$ if and only if $n_G(X,Y)=1$.

\section{The structure of the manuscript}\label{subsection:structure}
The main part of this manuscript is divided in seven chapters.  In chapter \ref{sec:spears} we recall the concept of spear vector and introduce the new concept of spear set. These concepts are used here as ``leitmotiv'' to give a unified presentation of the concepts of spear operator, lush operator, alternative Daugavet property, and other notions that we will introduce here for operators. We collect some properties of spear sets and vectors, together with some examples of spear vectors.

Chapter \ref{chapter:spear-aDP-target-lush} includes the main definitions of the manuscript for operators: spearness, the alternative Daugavet property and lushness. We start presenting some preliminary results and easy examples of spear operators in section \ref{section:spear-operators}. Next, in section \ref{sec:aDP} we study operators $G \in L(X,Y)$ with the alternative Daugavet Property (aDP). These are operators satisfying $\| G + \T\,T\| = 1 + \| T\|$
for every rank-one operator $T \in L(X,Y)$. This definition is a  generalization of the aDP, which is analogous to the generalization of the DPr given by Daugavet centers \cite{Bosenko-Kadets}: $G \in L(X,Y)$ is a \emph{Daugavet center} if $\| G + T\| = 1 + \| T\|$
for every rank-one operator $T \in L(X,Y)$.
\index{Daugavet center}%
Of course, Daugavet centers have the aDP, but the converse result is not true. There is a clear parallelism between the study of the aDP and of Daugavet centers. We give several characterizations of operators with the aDP (some of them in terms of spear sets) and prove that this is a separably determined property. Section \ref{sec:StrongTargetOperators} starts with the definition of target operator for $G \in L(X,Y)$. This property guarantees that an operator $T \in L(X,Y)$ that has it satisfies $\|G + \T\,T\| = 1 + \| T\|$. Interestingly, if $G$ has the aDP and the operator $T$ is SCD, then $T$ is a target for $G$, and this will be frequently used to deduce important results. Our new concept of target operator naturally plays an analogous role that the one played by strong Daugavet operators in the study of the DPr \cite{KadSW2} and in the study of Daugavet centers \cite{BosGnar}. Let us say that even for the case $G=\Id_X$, this concept is new and provides with non trivial new results. We characterize target operators for a given operator $G$, show that this property is separably determined, and prove that if $G$ has the aDP, then every operator whose restriction to separable subspaces is SCD is a target for $G$. In section \ref{sec:Lushness}, we introduce the notion of lush operator, which generalizes the concept of lush space. This generalization is closely connected with target operators from the previous section, which, on the one hand, reduces some results about lush spaces to results from the previous section, and on the other hand, gives more motivation for the study of target operators. We give several characterizations of lush operators, prove that this property is separably determined, show that the aDP and lushness are equivalent when every separable subspace of the domain space is SCD (so, for instance, when the domain is Asplund, has the Radon-Nikod\'{y}m Property, or does not contain copies of $\ell_1$), and present some sufficient conditions for lushness which will be used in the chapter about examples and applications. Besides, we prove that lush operators with separable domain fulfill a stronger version of lushness which has to do with spear functionals.

Chapter \ref{sect:classical-examples} is devoted to present some examples in classical Banach spaces. Among other results, we show that the Fourier transform is lush, we characterize operators from $L_1(\mu)$ spaces which have the aDP, and we study lushness, spearness and the aDP for operators which arrive to spaces of continuous functions. In particular, we show that every uniform algebra isometrically embeds by a lush operator into the space of bounded continuous functions on a completely regular Hausdorff topological space (for unital algebras, this space is just its Choquet boundary).

Next, we devote chapter \ref{sec.examples-appl} to provide further results on our properties. We characterize lush operators when the domain space has the Radon-Nikod\'{y}m Property or the codomain space is Asplund, and we get better results when the domain or the codomain is finite-dimensional or when the operator has rank-one. Further, we study the behaviour of lushness, spearness and the aDP with respect to the operation of taking adjoint operators; in particular, we show that these properties pass from an operator to its adjoint if the domain has the Radon-Nikod\'{y}m Property or the codomain is $M$-embedded; we also show that the aDP and spearness pass from an operator to its adjoint when the codomain is $L$-embedded.

In chapter \ref{sect:consequences} we provide with some isomorphic and isometric consequences of the properties as, among others, that the dual of the domain of an operator with the aDP and infinite rank contains $\ell_1$ in the real case. Many results showing that the aDP, spearness and lushness do not combine well with rotundity or smoothness properties are also presented.

We study Lipschitz spear operators in chapter \ref{sec:Lipschitz}. These are just the spear vectors of the space of Lipschitz operators between two Banach spaces endowed with the Lipschitz norm. The main result here is that every (linear) lush operator is a Lipschitz spear operator, a result which can be applied, for instance, to the Fourier transform. We also provide with an analogous result for aDP operators and for Daugavet centers.

Finally, a collection of stability results for our properties is given in chapter \ref{sec:stability}. We include results for various operations like absolute sums, vector-valued function spaces, and ultraproducts. The results we got are in most cases extensions of previously known results for the case of the identity.

We complement the manuscript with a collection of open problems in chapter \ref{sec:OpenProblems}.

\vspace{1em}

We finish the introduction with a diagram about the relationships between the properties of a norm-one operator $G\in L(X,Y)$ that we have presented in this introduction:
$$
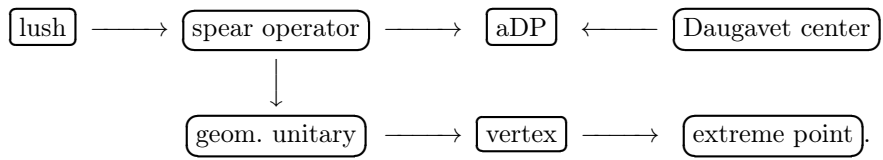
\begin{CD}
\Ovalbox{lush}@>>> \Ovalbox{spear operator} @>>> \Ovalbox{aDP} @<<< \Ovalbox{Daugavet center} \\ @. @VVV  \\
@. \Ovalbox{geom.~unitary}  @>>> \Ovalbox{vertex}  @>>> \Ovalbox{extreme point}.
\end{CD}
$$
None of the implications above reverses, and Daugavet centers and spear operators do not imply each other.

\chapter{Spear vectors and spear sets}
\label{sec:spears}

The following definition will be crucial in our further discussion.

\begin{Defi}[\mbox{\cite[Definition~4.1]{Ardalani}}]\label{def-spear} Let $X$ be a Banach space. An element $z \in S_X$ is a \emph{spear} (or \emph{spear vector}) if $ \| z + \T\,x \| = 1 + \| x\| $ for every $x \in X$. We write $\Spear(X)$ to denote the set of all elements of a Banach space $X$ which are spear.
\end{Defi}
\index{spear}%
\index{spear vector}%
\index{spear(X)@$\Spear(\cdot)$}%

As we commented in section \ref{subsect:numericalranges}, this is equivalent to the fact that $N(X,z)=1$. In particular, the definition was motivated in \cite{Ardalani} by the fact that  $\Id_{X}$ is a spear element of $L(X)$ if and only if $X$ has numerical index one. Let us also comment that the concept of spear vector appeared, without name, in the paper \cite{LimaIntersectionBalls} by \AA.~Lima about intersection properties of balls. It had also appeared tangentially in the monograph \cite{Lindenstrauss-MemAMS1964} by J.~Lindenstrauss about extension of compact operators.

\begin{Rema}
Observe that using a standard convexity argument, $z\in S_X$ is a spear vector if and only if $\|z + \T\,x\|=2$ for every $x\in S_X$.
\end{Rema}

The next notion extends the definition of spear from vectors to sets.

\begin{Defi} Let $X$ be a Banach space.
$F \subset B_{X}$ is called a \emph{spear set} if $\|F + \T\,x\| = 1 + \| x\|$ for every $x \in X$.
\end{Defi}
\index{spear set}%

Observe that if $F\subset B_X$ is a spear set, then every subset of $B_X$ containing $F$ is also a spear set. In particular, if a subset $F$ of $B_X$ contains a spear vector, then $F$ is a spear set. On the other hand, it is not true that every spear set contains a spear vector: in every Banach space $X$, $F=S_X$ is obviously a spear set, but there are Banach spaces containing no spear vectors at all (for instance, a two-dimensional Hilbert space, see Example \ref{Exam:basicsSpears}.(h)).

We start the exposition with the following fundamental result.

\begin{Theo}\label{Theo:strongSpearCharacterization}
Let $X$ be a Banach space, let $F$ be a subset of $B_{X}$ and let $\mathcal{A} \subset B_{X^{\ast}}$ with $B_{X^{\ast}} = \overline{\conv}^{\weakstar}{(\mathcal{A})}$. The following statements are equivalent:
\begin{enumerate}
\item[(i)] $\| F + x\| = 1 + \| x\|$ for every $x \in X$.
\item[(ii)] $B_{X^{\ast}} = \overline{\conv}^{w^{\ast}}{(\GS(\mathcal{A},F,\eps))}$ for every $\eps > 0$.
\item[(iii)] $B_{X^{\ast}} = \overline{\conv}^{w^{\ast}}{(\GS(\ext B_{X^{\ast}},F,\eps))}$ for every $\eps > 0$.
\item[(iv)] $\GF(\ext B_{X^{\ast}},F)$ is a dense subset of $(\ext B_{X^{\ast}}, \weakstar)$.
\end{enumerate}
If $X = Y^{\ast}$ is a dual Banach space, this is also equivalent to \begin{enumerate}
\item[(v)] $B_{Y}=\overline{\conv}{(\GS(S_{Y},F,\eps))}$ for every $\eps > 0$.
\end{enumerate}
Moreover, remark that the set $\GF(\ext B_{X^{\ast}},F)$ is G$_\delta$ which makes item (iv) more applicable.
\end{Theo}

\begin{proof}
(i) $\Rightarrow$ (ii): Given  $\eps > 0$, we just have to check that every $\weakstar$-slice $S$ of $B_{X^{\ast}}$ intersects $\GS(\mathcal{A},F, \eps)$. We can assume that $S = \Slice(B_{X^\ast}, x_{0}, \delta)$ for some $x_{0} \in S_{X}$ and $\eps > \delta > 0$. Using (i) and the condition on $\mathcal{A}$, we can find $x_{0}^{\ast} \in \mathcal{A}$ and $z_{0} \in F$ such that
$$
\Real{x_{0}^{\ast}(z_{0})} + \Real{ x_{0}^{\ast}(x_{0})} > 2 - \delta.
$$
In particular, $\Real{x_{0}^{\ast}(z_{0})}> 1 - \delta$ and $\Real{x_{0}^{\ast}(x_{0})} > 1 - \delta$, so $x_{0}^{\ast} \in S \cap \GS(\mathcal{A},F, \eps)$.

(ii) $\Rightarrow$ (i):
Given $x \in S_X$ and $\eps>0$, the hypothesis allows us to find $x^\ast\in \GS(\mathcal{A},F,\eps)$ such that $\re x^\ast(x)>1-\eps$. Also, by definition of $\GS(\mathcal{A},F,\eps)$, there is $z\in F$ such that $\re x^\ast(z)>1-\eps$. Now,
\begin{align*} \| F + x \|  & \geqslant \|z+x\| \geqslant \re x^\ast(z) + \re x^\ast(x)>2-2\eps,
\end{align*}
and the arbitrariness of $\eps$ gives the result.

The equivalence between (i) and (iii) is just a particular case of the already proved equivalence between (i) and (ii) since $\mathcal{A}=\ext B_{X^{\ast}}$ satisfies the condition above by the Krein-Milman Theorem.

(iii) $\Rightarrow$ (iv): For each $\eps > 0$, $\GS(\ext B_{X^{\ast}},F, \eps)$ is a relatively weak$^\ast$-open subset of $\ext B_{X^{\ast}}$, as it can be written as union of weak$^\ast$-slices. Moreover, condition (iv) together with Milman's Theorem (see Lemma \ref{Kreinlemma}.(b)) yields that the set $\GS(\ext B_{X^{\ast}},F, \eps)$ is weak$^\ast$-dense in $\ext B_{X^{\ast}}$. Using that the set $(\ext B_{X^{\ast}}, \weakstar)$ is a Baire space (see Lemma~\ref{Kreinlemma}.(c)), we conclude that
\begin{equation}\label{eq:Theo:strongSpearCharacterization-G-delta}
\GF(\ext B_{X^{\ast}},F)=\bigcap\nolimits_{n\in \N} \GS(\ext B_{X^{\ast}},F,1/n)
\end{equation}
satisfies the properties above.

(iv) $\Rightarrow$ (iii): Given $\eps > 0$, since $\GF(\ext B_{X^{\ast}},F)\subset\GS(\ext B_{X^{\ast}},F, \eps)$ we deduce that this last set is also dense in $(\ext{B_{X^{\ast}}}, \omega^{\ast})$, and so using the Krein-Milman Theorem we conclude that
$$ B_{X^{\ast}} = \overline{\conv}^{w^{\ast}}{(\ext{B_{X^\ast}})} \subset \overline{\conv}^{w^{\ast}}{(\GS(\ext B_{X^{\ast}},F,\eps))} $$

Finally, if $X=Y^\ast$, (v) is a particular case of (ii) with $\mathcal{A} = B_{Y}$ by Goldstine's Theorem.

The ``moreover'' part follows from equation~\eqref{eq:Theo:strongSpearCharacterization-G-delta}.
\end{proof}

Now we may present a characterization of spear sets which is an easy consequence of the above theorem and the fact that $\|F+\T\,x\|=\|\T\,F + x\|$ for every set $F$ and every vector $x$.

\begin{Coro}
\label{Coroll:spearSet} Let $X$ be a Banach space and let $\mathcal{A} \subset B_{X^{\ast}}$ with $\overline{\conv}^{\omega^{\ast}}{(\mathcal{A})}=B_{X^{\ast}}$. For $F \subset B_{X}$, the following assertions are equivalent:
\begin{enumerate}
\item[(i)] $F$ is a spear set, i.e.\ $\|F + \T\,x \| = 1 + \| x\|$ for each $x \in X$.
\item[(ii)] $B_{X^{\ast}} = \overline{\aconv}^{w^{\ast}}{(\GS(\mathcal{A},F,\eps))}$ for every $\eps > 0$.
\item[(iii)] $\GF(\ext B_{X^{\ast}},\T F)$ is a dense G$_\delta$ subset of $(\ext B_{X^{\ast}},\weakstar)$.
\end{enumerate}
If $X=Y^{\ast}$ is a dual Banach space, this is also equivalent to
\begin{enumerate}
\item[(iv)] $B_{Y} = \overline{\aconv}{(\GS(S_{Y},F,\eps))}$ for every $\eps > 0$.
\end{enumerate}
\end{Coro}

The following result is of interest in the complex case.

\begin{Prop}
\label{Prop:MLURset} Let $X$ be a Banach space.
If $F \subset B_{X}$ is a spear set, then $\| F \pm x\|^{2} \geqslant 1 + \| x\|^{2}$ for every $x \in X$.
\end{Prop}

\begin{proof}
Let $x \in X$ and $\eps > 0$. Using Corollary \ref{Coroll:spearSet}.(ii), we get that
\begin{align*}
\| F \pm x\|^{2} & \geqslant \sup\bigl\{|x^{\ast}(z) \pm x^{\ast}(x)|^{2} \colon z \in F,\, x^{\ast} \in \GS(S_{X^{\ast}},F,\eps)\bigr\} \\ & \geqslant \sup\bigl\{|x^{\ast}(z)|^{2} + |x^{\ast}(x)|^{2}
\colon z \in F,\, x^{\ast} \in \GS(S_{X^{\ast}},F,\eps)\bigr\}
\\
& \geqslant \sup\bigl\{(1-\eps)^{2} + |x^{\ast}(x)|^{2}\colon x^{\ast} \in \GS(S_{X^{\ast}},F,\eps)\bigr\} \\ & = (1 - \eps)^2 + \| x\|^{2}.\qedhere
\end{align*}
\end{proof}

The case in which a spear set is a singleton coincides, of course, with the concept of spear vector of Definition~\ref{def-spear}. Most of the assertions of the next corollary follow from Corollary \ref{Coroll:spearSet} and the fact that $\GF(\ext B_{X^{\ast}},\T z) = \{ x^{\ast} \in \ext B_{X^{\ast}} \colon |x^{\ast}(z)| = 1 \}$ is $\weakstar$-closed. The other ones are consequences of the general theory of numerical range spaces (see section \ref{subsect:numericalranges}) and can be found in \cite[\S2.1]{Cabrera-Rodriguez}.

\begin{Coro}
\label{Defi:spearVector} Let $X$ be a Banach space and let $\mathcal{A} \subset B_{X^{\ast}}$ with $B_{X^\ast}=\overline{\conv}^{\omega^{\ast}}{(\mathcal{A})}$. The following assertions are equivalent for $z\in S_X$:
\begin{enumerate}
\item[(i)] $z\in \Spear(X)$ (that is, $\|z + \T\,x\| = 1 + \| x\|$ for every $x \in X$).
\item[(ii)] $B_{X^{\ast}} = \overline{\aconv}^{\weakstar}{(\Slice(\mathcal{A},z, \eps))}$ for each $\eps > 0$.
\item[(iii)$_\R$] If $X$ is a real space, $B_{X^\ast}=\conv \bigl(\Face(S_{X^\ast},z)\cup - \Face(S_{X^\ast},z)\bigr)$.
\item[(iii)$_\C$] If $X$ is a complex space, $\interior(B_{X^\ast}) \subset \aconv \bigl(\Face(S_{X^\ast},z)\bigr)$ so, in particular, $B_{X^\ast} = \overline{\aconv} \bigl(\Face(S_{X^\ast},z)\bigr)$.
\item[(iv)] $|x^{\ast}(z)| = 1$ for every $x^{\ast} \in \ext B_{X^{\ast}}$.
\end{enumerate}
If $X=Y^{\ast}$ is a dual Banach space and $z=y^\ast\in S_{Y^\ast}$, this is also equivalent to:
\begin{enumerate}
\item[(v)] $B_{Y} = \overline{\aconv}{\bigl(\Slice(S_{Y},y^\ast,\eps)\bigr)}$ for every $\eps > 0$.
\end{enumerate}
\end{Coro}

\begin{proof}
The equivalence between (i), (ii) and (iv) is just particular case of Corollary \ref{Coroll:spearSet}, as it is the equivalence with (v) when $X$ is a dual space.

(i) $\Rightarrow$ (iii) is contained in \cite[Theorem 2.1.17]{Cabrera-Rodriguez} (both in the real and in the complex case), but we give the easy argument here. Recall that (i) is equivalent to the fact that $v(X,z,x)=\|x\|$ for every $x\in X$ or, equivalently, $\Face(S_{X^\ast},z)$ is norming for $X$ or, equivalently,
\begin{equation}\label{eq:spearvector-01}
B_{X^\ast}=\overline{\aconv}^{w^\ast}(\Face(S_{X^\ast},z)).
\end{equation}
In the real case, we have that the set
$$\aconv(\Face(S_{X^\ast},z))=\conv\bigl(\Face(S_{X^\ast},z)\cup -\Face(S_{X^\ast},z)\bigr)$$
is weak$^\ast$-compact as so is $\Face(S_{X^\ast},z)$, and the result follows from \eqref{eq:spearvector-01}. In the complex case, for $0<\rho<1$ we take $n\in \N$ such that $(1-\rho)B_\C\subseteq \conv \{z_1,\ldots,z_n\}$, where $\{z_1,\ldots,z_n\}$ are the $n$th roots of $1$ in $\C$. Then we have
\begin{align*}
(1-\rho) \aconv(\Face(S_{X^\ast},z))&= (1-\rho)\conv(B_\C\,\Face(S_{X^\ast},z))\\
&\subseteq \conv\left({\textstyle \bigcup\nolimits_{k=1}^n z_k \Face(S_{X^\ast},z)}\right).
\end{align*}
Since $\conv\left(\bigcup\nolimits_{k=1}^n z_k \Face(S_{X^\ast},z)\right)$ is weak$^\ast$-compact and is contained in the set $\aconv(\Face(S_{X^\ast},z))$, it follows from \eqref{eq:spearvector-01} that $(1-\rho)B_{X^\ast}\subset\aconv(\Face(S_{X^\ast},z))$, and this gives the result moving $\rho\downarrow 0$.

The implication (iii) $\Rightarrow$ (i) follows immediately from (ii) $\Rightarrow$ (i) for the particular case $\mathcal{A} = B_{X^{\ast}}$.
\end{proof}

The next surprising result about spear vectors of a dual space appeared literally in \cite[Corollary~3.5]{AcoBecRod} and it is also consequence of the earlier \cite[Theorem~2.3]{Godefroy-Indumathi} using Corollary \ref{Defi:spearVector}.(iii). In both cases, the main tool is the use of norm-to-weak upper semicontinuity of the duality and pre-duality mappings. We include here an adaptation of the proof of \cite[Theorem~2.3]{Godefroy-Indumathi} to our particular situation which avoids the use of semicontinuities (which, on the other hand, are automatic in our context, see \cite[Fact 2.9.3 and Theorem 2.9.18]{Cabrera-Rodriguez}).

\begin{Theo}[\mbox{\cite[Theorem~2.3]{Godefroy-Indumathi},
\cite[Corollary~3.5]{AcoBecRod}}]\label{Theo:spears_of_the_dual}
Let $X$ be a Banach space and let $z^\ast\in S_{X^\ast}$.
Then, $z^\ast \in \Spear(X^{\ast})$ if and only if $ B_{X} = \overline{\aconv}{\,\Face(S_{X}, z^{\ast})} $.
\end{Theo}

\begin{proof}
The ``if'' part follows immediately from Corollary \ref{Defi:spearVector}.(v), so we just have to prove the ``only if'' part. To simplify, we will denote $F = \GF(S_{X}, \{z^{\ast}\})$ and $F^{\ast \ast} = \Face(S_{X^{\ast \ast}}, z^{\ast})$.

\emph{Claim 1}. {\slshape If $A \subset B_{X^{\ast \ast}}$ satisfies $\sup_{A}{\Real{z^{\ast}}} = 1$, then $\dist(A, F^{\ast \ast}) = 0$.}\  Indeed, let $\eps > 0$ and $a^{\ast \ast} \in A$ with $\Real{z^{\ast}(a^{\ast \ast})} > 1 - \eps$. Since $B_{X^{\ast \ast}} = \overline{\aconv}{\,F^{\ast \ast}}$ by Corollary \ref{Defi:spearVector}.(iii), we can find $x_{1}^{\ast \ast}, \ldots, x_{m}^{\ast\ast}\in F^{\ast \ast}$ and $\lambda_{1}, \ldots, \lambda_{m}\in [0,1]$ with $\sum_{k=1}^m \lambda_k=1$, and $\theta_{1}, \ldots, \theta_{m}\in \mathbb{T}$ such that
$$
\left\| a^{\ast \ast} - \sum_{n=1}^{m}{\lambda_{n} \theta_{n} x_{n}^{\ast \ast}} \right\| < \eps \hspace{4mm} \mbox{ and } \hspace{4mm} \Real{z^{\ast}\left( \sum_{n=1}^{m}{\lambda_{n} \theta_{n} x_{n}^{\ast \ast}} \right)} = \sum_{n=1}^{m}{\lambda_{n} \Real{\theta_{n}}} > 1 - \eps. $$
Therefore,
\begin{align*}
\dist{(A, F^{\ast \ast})} & \leqslant \left\| a^{\ast \ast} - \sum_{n=1}^{m}{\lambda_{n} x_{n}^{\ast \ast}} \right\| \leqslant \eps + \left\| \sum_{n=1}^{m}{\lambda_{n} \theta_{n} x_{n}^{\ast \ast}} - \sum_{n=1}^{m}{\lambda_{n} x_{n}^{\ast \ast}} \right\|\\ & \leqslant \eps + \sum_{n=1}^{m}{\lambda_{n} |1 - \theta_{n}|}
 \leqslant \eps + \sum_{n=1}^{m}{\lambda_{n} 2 \sqrt{1 - \Real{\theta_{n}}}}\\ & \leqslant \eps + 2 \sqrt{\sum_{n=1}^{m}{\lambda_{n} (1 - \Real{\theta_{n}})}} \leqslant \eps + 2 \sqrt{\eps}.
\end{align*}
Since $\eps > 0$ was arbitrary, the claim is proved.

\emph{Claim 2}. {\slshape Given $x_{0} \in B_{X}$ with $\dist(x_{0}, F^{\ast \ast}) < \eps$ we have that the set $A=B_{X} \cap (x_{0} + \eps B_{X})$ satisfies $\sup_{A}{\Real{z^{\ast}}} = 1$.}\ Indeed, let $x_{0}^{\ast \ast} \in F^{\ast \ast}$ with $\| x_{0}^{\ast \ast} - x_{0}\| < \eps$. Using the Principle of Local Reflexivity \cite[Theorem 6.3]{basisLinear}, we have that for each $\delta > 0$ we can find an element $x_{\delta} \in X$ such that
$\| x_{\delta}\| \leqslant 1$, $\Real{z^{\ast}(x_{\delta})} > 1 - \delta$, and $\| x_{\delta} - x_{0} \|< \eps$.

\emph{Claim 3}. {\slshape $F\neq \emptyset$ (in particular, $F=\Face(S_X,z^\ast)$) and if $A \subset B_{X}$ satisfies $\sup_{A}{\Real{z^{\ast}}} = 1$, then $\dist(A, F) = 0$.}\  Indeed, let $A_{0}=A$ and fix $\eps > 0$. By Claim 1, we have that $\dist(A_{0}, F^{\ast \ast})=0$, so taking $x_{0} \in A_{0}$ with $\dist(x_{0}, F^{\ast \ast}) < \frac{\eps}{2}$ we have that $A_{1} = B_{X} \cap (x_{0} + \frac{\eps}{2} B_{X})$ satisfies $\sup_{A_1}\Real{z^{\ast}} = 1$. Repeating the same process with $A_{1}$, we can find $x_{1} \in A_{1}$ such that $A_{2} := B_{X} \cap (x_{1} + \frac{\eps}{4} B_{X})$ satisfies $\sup_{A_2} \Real{z^{\ast}}= 1$. Iterating this process, we will have a Cauchy sequence $(x_{n})$ whose limit $z \in B_{X}$ satisfies that $\dist{(z, A)} \leqslant \eps$ and that $z^{\ast}(z) = 1$, so $z\in F$.

\emph{Claim 4}. {\slshape $B_{X} = \overline{\aconv}\,{F}$.}\  Indeed, given a slice $S$ of $B_{X}$, Corollary \ref{Defi:spearVector}.(v) shows that  $\sup_{\mathbb{T}S}{\Real z^{\ast}} = 1$, so Claim 3 provides that   $\dist(\mathbb{T}S, F) = 0$. It is now routine to show that $S \cap \mathbb{T}F \neq \emptyset$ for every slice $S$ of $B_X$.
\end{proof}

The following proposition collects all the properties of spear vectors we know. In order to prove it we need the following technical lemma.

\begin{Lemm}
	\label{Lemm:convergenceSpears}Let $X$ be a Banach space and let $(F_{n})_{n \in \N}$ be a decreasing sequence of spear sets of $X$ such that $\diam(F_{n})$ tends to zero. If $z \in \bigcap_{n \in \N}{F_{n}}$, then $z$ is a spear element.
\end{Lemm}

\begin{proof}
	For every $x \in X$ we can write
	\[ \|z + \T\,x \|  \geqslant \| F_{n} + \T\,x\| - \| F_{n} - z\| = 1 + \|x\| - \| F_{n} - z\|. \]
	But the hypothesis implies that $\lim_{n}{\| F_{n} - z \|} = 0$.
\end{proof}

\begin{Prop}
\label{Prop:spearVectorsProperties}
Let $X$ be a Banach space. Then:
\begin{enumerate}
\item[(a)] $\| z \pm x\|^{2} \geqslant 1 + \| x\|^{2}$ for each $z \in \Spear(X)$ and every $x \in X$.
\item[(b)]Every $z\in \Spear(X)$ is a strongly extreme point of $B_X$. In particular, $\Spear(X) \subset \ext B_{X}$.
\item[(c)] $J_X\bigl(\Spear(X)\bigr) \subset \Spear(X^{\ast\ast})$. In particular, $J_X\bigl(\Spear(X)\bigr) \subset \ext B_{X^{\ast\ast}}$.
\item[(d)] $\Spear(X)$ is norm-closed.
\item[(e)] If $\dim(X) \geqslant 2$, then $\Spear(X)$ is nowhere-dense in $(S_{X},\| \cdot\|)$.
\item[(f)] If $B_{X} = \overline{\conv}{(\Spear(X))}$ then $\Spear(X^{\ast}) = \ext B_{X^{\ast}}$.
\item[(g)] If $B_{X^\ast} = \overline{\conv}^{\weakstar}{(\Spear(X^\ast))}$, then
    $x \in S_{X}$ is a spear element if and only if $x$ is a strongly extreme point of $B_X$ if and only if $x \in \ext B_{X^{\ast \ast}}$.
\item[(h)]  If $X$ is strictly convex and $\dim(X) \geqslant 2$, then $\Spear(X)=\emptyset=\Spear(X^{\ast})$.
\item[(i)] If $X$ is smooth and $\dim(X) \geqslant 2$, then $\Spear(X)=\emptyset$.
\end{enumerate}
If $X$ is a \textbf{real} space, we can add:
\begin{enumerate}
\item[(j)] If $\Spear(X)$ is infinite, then $X$ contains a copy of $c_{0}$ or $\ell_1$.
\item[(k)] If $z^{\ast} \in \Spear(X^\ast)$ and $x^{\ast} \in X^{\ast}$ is norm-attaining with $\| z^{\ast} - x^{\ast}\| < 1 + \| x^{\ast}\|$, then $z^{\ast} + x^{\ast}$ is norm-attaining and $\| z^{\ast} + x^{\ast}\| = 1 + \| x^{\ast}\|$.
\item[(l)] If $X$ is smooth and $\dim(X) \geqslant 2$, then $\Spear(X^{\ast}) = \emptyset$.
\end{enumerate}
\end{Prop}

\begin{proof}[Proof of Proposition~\ref{Prop:spearVectorsProperties}]
Statement (a) is given by Proposition \ref{Prop:MLURset}, and (b) is an obvious consequence of it.

(c). Fixed $z\in \Spear(X)$, we have that $B_{X^\ast}= \overline{\aconv}^{\|\cdot\|}\bigl(\Face(S_{X^\ast},z)\bigr)$ by Corollary~\ref{Defi:spearVector}.(iii)  so, using Goldstine's Theorem for $X^\ast$, we obtain that
$$
B_{X^{\ast\ast\ast}}= \overline{\aconv}^{\sigma(X^{\ast\ast\ast},X^{\ast\ast})} J_{X^\ast}\bigl(\Face(S_{X^\ast},z)\bigr)$$ and, a fortiori,
$$
B_{X^{\ast\ast\ast}}=\overline{\aconv}^{\sigma(X^{\ast\ast\ast}, X^{\ast\ast})} \bigl(\Face(S_{X^{\ast\ast\ast}},J_X(z))\bigr).
$$
Now, Milman's Theorem (see Lemma \ref{Kreinlemma}.(b)) gives that
$$
\ext(B_{X^{\ast\ast\ast}})\subset \overline{\Face(S_{X^{\ast\ast\ast}},\T\,J_X(z))}^{\sigma(X^{\ast\ast\ast}, X^{\ast\ast})}.
$$
Then, Corollary \ref{Defi:spearVector}.(iv) gives that $J_X(z)\in \Spear(X^{\ast\ast})$. An alternative proof is the following: for $z\in \Spear(X)$ we have that $|x^\ast(z)|=1$ for every $x^\ast\in \ext(B_{X^\ast})$ by Corollary \ref{Defi:spearVector}.(iv), and so $|x^{\ast\ast\ast}(J_X(z))|=1$ for every $x^{\ast\ast\ast}\in \ext(B_{X^{\ast\ast\ast}})$ by \cite[Proposition 3.5]{Sharir}, so $J_X(x)\in \Spear(X^{\ast\ast})$ by using again Corollary \ref{Defi:spearVector}.(iv).

(d). Given a norm-convergent sequence $(x_{n})_{n \in \N}$ in $\Spear(X)$, apply Lemma \ref{Lemm:convergenceSpears} to the family of sets $F_{n}:= \overline{\{ x_{m}\colon m \geqslant n \}}$.

(e). Fixed $e_{0}^{\ast} \in \ext B_{X^\ast}$, take an element $x_{0} \in \ker{e_{0}^{\ast}} \cap S_{X}$. Given $z \in \Spear(X)$, there exists $\theta_{0} \in \T$ such that $\|z + \theta_0 x_0\|=2$ and so, by convexity, $\| z + \delta \theta_{0}x_{0} \| = 1 + \delta$ for every $\delta>0$. Then, for each $\delta>0$, $v = (z + \delta \theta_{0}x_{0})/(1 + \delta)$ belongs to $S_{X}$, satisfies $\| v - z\|  \leqslant 2 \delta$, and $|e_{0}^{\ast}(v)| = (1 + \delta)^{-1} \neq 1$, so $v$ is not a spear by Corollary \ref{Defi:spearVector}.(iv).

(f). Fix $e^{\ast} \in \ext B_{X^{\ast}}$ and let $x^{\ast} \in X^{\ast}$. The hypothesis implies that for every $\eps > 0$ we can find $z \in \Spear(X)$ with $|x^{\ast}(z)| > \| x^{\ast}\| - \eps$. Since $|e^{\ast}(z)| = 1$ (as $z$ is a spear) we conclude that $\| \T\,e^{\ast} + x^{\ast} \| \geqslant |e^{\ast}(z)| + |x^{\ast}(z)| \geqslant 1 + \| x^{\ast}\| - \eps$. This gives that $\ext B_{X^\ast}\subset \Spear(X^\ast)$ and the equality follows from (b).

(g). If we assume that $x \in \ext B_{X^{\ast \ast}}$, then $|x^{\ast}(x)| = 1$ for each $x^{\ast} \in \Spear(X^\ast)$ by Corollary~\ref{Defi:spearVector}.(iv). But Milman's theorem (see Lemma \ref{Kreinlemma}.(c)) applied to our assumption gives that $\Spear(X^{\ast})$ is weak$^\ast$-dense in $\ext B_{X^{\ast}}$, so we conclude that $|x^{\ast}(x)| = 1$ for every $x^{\ast} \in \ext B_{X^{\ast}}$, which implies that $x$ is a spear by Corollary~\ref{Defi:spearVector}.(iv).

(h). If $x_0\in \Spear(X)$, we have that $\|x_0+\T\,x\|=2$ for every $x\in S_X$. If $X$ is strictly convex, this implies that $S_X\subset \T\,x_0$ and so $\dim(X)=1$. Next, suppose that there exists $z^{\ast} \in \Spear(X^{\ast})$, so $B_{X} = \overline{\aconv}{(\Face(S_{X},z^{\ast}))}$ by Theorem \ref{Theo:spears_of_the_dual}. If $\dim(X) \geqslant 2$ then $\Face(S_{X},z^{\ast})$ contains at least two points, so $X$ is not strictly convex.

(i). If $X$ is smooth, the set $\Face(S_{X^\ast},z)$ is a singleton for every $z\in S_X$. If there is $z\in \Spear(X)$, then the above observation and Corollary \ref{Defi:spearVector}.(iii) imply that $X^\ast$ is one-dimensional, so $X$ is one-dimensional as well.

(j). Is just a reformulation of \cite[Proposition~2]{RealNumIndexOne}, but we include the short argument for completeness. Suppose that $X$ does not contain $\ell_1$. Then, by Rosenthal's $\ell_1$-Theorem
\cite[Chapter XI]{DiestelSeq} there is a weakly Cauchy sequence $(x_n)_{n\in\N}$ of distinct members of $\Spear(X)$. Write $Y$ for the closed linear span of $\{x_n\colon n\in \N\}$ and observe that, obviously, $x_n\in \Spear(Y)$ for every $n\in \N$. Therefore, by Corollary \ref{Defi:spearVector}.(iv), the fact that the sequence $(x_n)$ is weakly Cauchy and that we are in the real case, we have that $\ext(B_{Y^\ast})= \bigcup\nolimits_{n\in\N} \bigl(E_n\cup-E_n\bigr)$ where
$$
E_n:=\bigl\{y^\ast\in \ext(B_{Y^\ast})\colon y^\ast(x_k)=1\ \text{for } k\geqslant n\bigr\} \qquad (n\in \N).
$$
As $\{x_n\}$ separates the points of $Y^\ast$, each $E_n$ is finite, so $\ext(B_{Y^\ast})$ must be countable. Fonf's Theorem \cite{Fonf} gives us that $X\supseteq c_0$, finishing the proof.

The proof of (k) is based on ideas of \cite{LimaIntersectionBalls}. Let $F=\Face(S_{X^{\ast \ast}},z^\ast)$. By Corollary \ref{Defi:spearVector}.(iii), we have that $B_{X^{\ast \ast}} = \conv{(F \cup -F)}$. Let $x^{\ast} \in X^\ast$ attain its norm at $x \in S_{X}$, i.e.\  $x^{\ast}(x) = \| x^{\ast}\|$, and suppose that $\| z^{\ast} - x^{\ast}\|< 1 + \| x^{\ast}\|$. We can write $x = (1-\lambda) x_{1}^{\ast \ast} - \lambda x_{2}^{\ast \ast}$ for some $0 \leqslant \lambda \leqslant 1$ and $x_{1}^{\ast \ast}, x_{2}^{\ast \ast} \in F$. If we assume that $0 < \lambda \leqslant 1$, then $x_{2}^{\ast \ast}(x^{\ast}) = -\| x^{\ast}\|$ necessarily, which is not possible as $|x_{2}^{\ast \ast}(x^{\ast}) - x_{2}^{\ast \ast}(z^{\ast})| \leqslant \| x^{\ast} - z^{\ast}\| < 1 + \| x^{\ast}\|$. Therefore, $\lambda = 0$ and we get that $x = x_{1}^{\ast \ast} \in F \cap S_{X}$, so $[z^{\ast} + x^{\ast}](x) = 1 + \| x^{\ast}\|$.

(l). Suppose that there exists $z^{\ast} \in \Spear(X^{\ast})$. By (j), we can take a norm-attaining functional $x_{0}^{\ast} \in S_{X^{\ast}}$ with $0 < \| z^{\ast} - x_{0}^{\ast}\| < 2$, such that $z^{\ast} + x_{0}^{\ast}$ is norm-attaining and $\|z^{\ast} + x_{0}^{\ast}\|=2$. Hence there is $x_{0} \in S_{X}$ with $x_{0}^{\ast}(x_{0}) = z^{\ast}(x_{0}) = 1$, which means that $X$ is not smooth.
\end{proof}

Below we list the known examples of spear vectors from \cite{Ardalani} and some easy-to-check generalizations of those examples.

\begin{Exam}
\label{Exam:basicsSpears}
{$\,$}
\begin{enumerate}
\item[(a)] \cite[p.~170]{Ardalani} If $\Gamma$ is an arbitrary set, then
    $$
    \Spear\bigl(\ell_1(\Gamma)\bigr)=\bigl\{\theta e_\gamma\colon \theta\in \T,\, \gamma\in \Gamma\bigr\},
    $$
    where $e_{\gamma}$ is the function on $\Gamma$ with value one at $\gamma$ and zero on the rest. The proof is straightforward.
\item[(b)] \cite[p.~170]{Ardalani} Let $(\Omega, \Sigma, \mu)$ be a measure space. The spear vectors of the space $L_{1}(\Omega, \Sigma, \mu)$ are functions of the form $\theta \mathbbm{1}_{A}/\mu(A)$ where $\theta \in \T$ and $A \in \Sigma$ is an atom. That is, spear vectors coincide with extreme points of the unit ball. The proof of this result is straightforward.
\item[(c)] \cite[p.~170]{Ardalani} If $K$ is a Hausdorff compact space, then
    $$
    \Spear\bigl(C(K)\bigr)=\bigl\{f\in C(K)\colon |f(t)|=1 \ \forall t\in K\bigr\}.
    $$
    That is, again spear vectors coincide with extreme points of the unit ball. Again, the proof is elementary.
\item[(d)] With the same ideas of the example above, one can even show that if $X$ is a Banach space, then for $f\in C(K,X)$ one has that
$$
f\in \Spear\bigl(C(K,X)\bigr) \ \ \Longleftrightarrow \ \ f(t)\in \Spear(X) \text{ for every $t\in K$}.
$$
\item[(e)] Let $(\Omega, \Sigma, \mu)$ be a measure space and let $g\in L_\infty(\mu)$. Then,
$$
g\in \Spear(L_\infty(\mu)) \ \ \Longleftrightarrow \ \ |g(t)|=1 \text{ for $\mu$-almost every $t$}.
$$
Indeed, suppose that $g\in \Spear(L_\infty(\mu))$
    and there exists a measurable subset $A$ with $\mu(A)>0$ such that $|g(t)|<1$ for every $t\in A$. Then we have that
    $$
    A=\bigcup\nolimits_{n\in \N} \{t\in A\colon |g(t)|\leqslant 1-1/n\},
    $$
    so there exists $n\in \N$ such that the measurable set $B:=\{t\in A\colon |g(t)|\leqslant 1-1/n\}$ has positive measure. Now,
    $$
    \|g + \T\, \mathbbm{1}_{B}\|_\infty  \leqslant 2-1/n < 2 =1 + \|\mathbbm{1}_B\|_\infty
    $$
    and thus $g$ is not a spear vector, a contradiction. Conversely, suppose that $|g(t)|=1$ for $\mu$-almost every $t$. For $f\in L_\infty(\mu,Y)$ and $\eps>0$, there is $A\in \Sigma$ with $\mu(A)>0$ such that $|f(t)|\geqslant \|f\|_\infty-\eps$ for every $t\in A$. By the hypothesis, there is $A'\in \Sigma$, $A'\subset A$, $\mu(A')>0$ such that $|g(t)|=1$ for every $t\in A'$. Now, using the compactness of $\T$ we can give a lower bound for $\|g + \T\,f\|_\infty$. Indeed, fixed an $\eps$-net $\mathbb{T}_\eps$ of $\T$ we can find an element $\theta_1 \in \mathbb{T}_\eps$ and a subset $A''$ of $A'$ with positive measure such that $|g(t)+\theta_1 f(t)|\geqslant 1+|f(t)|(1-\eps)$ for every $t\in A''$. Therefore, we can write
    \begin{align*}
    \|g + \T\,f\|_\infty&\geqslant \inf_{t\in A''} |g(t)+ \theta_1\,f(t)|\\
    &\geqslant \inf_{t\in A''} 1 + |f(t)|(1-\eps) \geqslant 1 + (\|f\| - \eps)(1-\eps)
    \end{align*}
    and the arbitrariness of $\eps$ gives $\|g + \T\,f\|_\infty\geqslant 1+\|f\|$.
\item[(f)] It will be proved in Corollary \ref{Coro:spearLinfinity} that the vector valued case of (e) is also valid. Let $(\Omega, \Sigma, \mu)$ be a measure space, let $X$ be a Banach space and let $g\in L_\infty(\mu,X)$. Then $g\in \Spear(L_\infty(\mu,X))$ if and only if $g(t)\in \Spear(X)$ for $\mu$-almost every $t$. Actually, the proof of the ``if'' part is just an straightforward adaptation of the corresponding one for (e).
\item[(g)] It is straightforward to show that if the linear span of $z\in S_X$ is an $L$-summand of $X$, then $z\in \Spear(X)$. Observe that this is what happens in examples (a) and (b) above. On the other hand, the converse result does not hold. Indeed, just note that $C(\Delta)$, where $\Delta$ is the Cantor set, contains no proper $L$-summand by Behrends $L$-$M$ Theorem (see \cite[Theorem~I.1.8]{HWW} for instance).
\item[(h)] The space $L_p(\mu)$ contains no spear vector if $1<p<\infty$ and $\dim(L_p(\mu))\geqslant 2$ (use Proposition \ref{Prop:spearVectorsProperties}.(h), for instance).
\item[(i)] Let $X_1$, $X_2$ be Banach spaces and let $X=X_1\oplus_\infty X_2$. Then, $(z_1,z_2)\in\Spear(X)$ if, and only if, $z_1\in\Spear(X_1)$ and  $z_2\in\Spear(X_2)$. The proof is straightforward.
\item[(j)] Let $X_1$, $X_2$ be Banach spaces and let $X=X_1\oplus_1 X_2$. Then, $(z_1,z_2)\in\Spear(X)$ if, and only if, either $z_1\in\Spear(X_1)$ and $z_2=0$ or $z_1=0$ and $z_2\in\Spear(X_2)$. The proof is again straightforward.
\end{enumerate}
\end{Exam}

As an application of examples (b) and (e) above and Theorem \ref{Theo:spears_of_the_dual}, we get easily the following well-known old result (see \cite{MartinPaya-CLspaces} for an exposition which also covers the complex case).

\begin{Coro}\label{example:XorX*=L1=>casi-almost-CL}
Let $X$ be a Banach space such that either $X$ or $X^\ast$ is isometrically isomorphic to an $L_1(\mu)$ space. Then, $B_X=\overline{\aconv}\, \Face(S_X, x^\ast)$ for every $x^\ast \in \ext (B_{X^\ast})$.
\end{Coro}

\chapter[Spearness, the aDP and lushness]{Three definitions for operators: spearness, the alternative Daugavet property, and lushness}\label{chapter:spear-aDP-target-lush}

This is the main chapter of our manuscript, as we introduce and deeply study the main definitions: the one of spear operator, the weaker of operator with the alternative Daugavet property and the stronger lush operator.

\section{A first contact with spear operators}\label{section:spear-operators}

Even though it has been given in the introduction, we formally state the definition of spear operator as it is the main concept of the manuscript.

\begin{Defi}
Let $X$, $Y$ be Banach spaces and let $G\in L(X,Y)$ be a norm-one operator. We say that $G$ is a \emph{spear operator} if the norm equality
$$
\|G + \T\,T\|=1 + \|T\|
$$
holds for every $T\in L(X,Y)$, that is, if $G\in \Spear\bigl(L(X,Y)\bigr)$.
\end{Defi}
\index{spear operator}%

We would like now to list some of the equivalent reformulations of the concept of spear operator which one can get particularizing the results of the previous chapter. We will also include a characterization in terms of numerical ranges that comes from section \ref{subsect:numericalranges}.

\begin{Prop}
Let $X$, $Y$ be Banach spaces and let $G\in L(X,Y)$ be a norm-one operator. The following assertions are equivalent:
\begin{enumerate}
\item[(i)] $G$ is a spear operator, i.e.\ $\|G + \T\,T\|=1 + \|T\|$ for every $T\in L(X,Y)$.
\item[(ii)] $|\zeta(G)|=1$ for every $\zeta\in \ext \bigl(B_{L(X,Y)^\ast}\bigr)$.
\item[(iii)] Given $\eps>0$,
$$
\|T\|=\sup \bigl\{|y^\ast(T x)|\colon y^\ast\in S_{Y^\ast},\, x\in S_X,\, \re y^\ast(Gx)>1-\eps\bigr\}
$$
for every $T\in L(X,Y)$.
\end{enumerate}
\end{Prop}

\begin{proof}
(i) and (ii) are equivalent by Corollary \ref{Defi:spearVector}, and (iii) $\Rightarrow$ (i) is immediate, so only (i) $\Rightarrow$ (iii) needs an explanation. For $T\in L(X,Y)$ and $\eps>0$, write
$$
W_\eps(T):=\bigl\{|y^\ast(T x)|\colon y^\ast\in S_{Y^\ast},\, x\in S_X,\, \re y^\ast(Gx)>1-\eps\bigr\}.
$$
Fix now $\eps_0>0$ and observe that
$$
v_G(T)= \sup \bigcap\nolimits_{\eps>0} W_\eps(T) \leqslant \sup W_{\eps_0}(T) \leqslant \|T\|.
$$
On the other hand, $G$ is a spear operator if and only if $N(L(X,Y),G)=1$ or, equivalently, if $n_G(X,Y)=1$, that is, if $v_G(T)=\|T\|$ for every $T\in L(X,Y)$, and so all the inequalities above become equalities.
\end{proof}

We next present some examples of spear operators which may help to better understand the definition and see how far from the Identity a spear operator can be. The first family appeared in \cite[Theorem 4.2]{Ardalani}.

\begin{Prop}\label{Prop:examplesell1gamma-c0gamma}
Let $\Gamma$ be an arbitrary set, let $X$, $Y$ be Banach spaces, and let $\{e_\gamma\}_{\gamma \in \Gamma}$ be the canonical basis of $\ell_{1}(\Gamma)$ (as defined in Example \ref{Exam:basicsSpears}.(a)).
\begin{enumerate}
\item[(a)] $G \in L(\ell_{1}(\Gamma),Y)$ is a spear operator if and only if $G(e_{\gamma})\in \Spear(Y)$ for each $\gamma \in \Gamma$.
\item[(b)] $G \in L(X, c_{0}(\Gamma))$ is a spear operator if and only if $G^{\ast}(e_\gamma) \in \Spear(X^\ast)$ for every $\gamma\in \Gamma$.
\end{enumerate}
\end{Prop}

\begin{proof}
As indicated, (a) appears in \cite[Theorem 4.2]{Ardalani} and the proof is elementary. Let us prove (b). The sufficiency of the condition is given by the obvious fact that $G$ is a spear operator when $G^\ast$ is (as taking adjoint preserves the norm) and the result in (a). For the necessity, suppose that there is $\xi\in \Gamma$ such that $G^\ast(e_\xi)\notin \Spear(X^\ast)$ and find $x_0^\ast \in S_{X^\ast}$ such that $\|G^\ast(e_\xi) + \T\,x_0^\ast\|<2$. We then consider the norm-one operator $T\in L(X,c_0(\Gamma))$ given by $[Tx](\xi)=x_0^\ast(x)$ and $[Tx](\gamma)=0$ if $\gamma\neq \xi$ for every $x\in X$, and observe that $T^\ast(e_\xi)=x_0^\ast$ and $T^\ast(e_\gamma)=0$ for $\gamma\neq \xi$. Therefore,
$$
\|G+\T\,T\|=\|G^\ast + \T\,T^\ast\|=\sup_{\gamma\in \Gamma}\|G^\ast(e_\gamma) + \T\,T^\ast(e_\gamma)\|<2,
$$
so $G$ is not a spear operator.
\end{proof}

This result will be improved in Example \ref{Exam:ell_1-gamma-all-equivalent}. More involved examples of spear operators, lush operators, operators with the aDP\ldots will appear in chapters \ref{sect:classical-examples}, \ref{sec.examples-appl}, and \ref{sec:Lipschitz}.

The following observations follow straightforwardly from the definition of spear operator.

\begin{Rema}\label{Rema:SpearOperators}
Let $X,Y$ be Banach spaces and let $G \in L(X,Y)$.
\begin{enumerate}
\item[(i)] Composing with isometric isomorphisms preserves spearness: Let $X_1$, $Y_1$ be Banach spaces, and let $\Phi_1\in L(X_1,X)$ and $\Phi_2\in L(Y,Y_2)$ be isometric isomorphisms. Then $G\in \Spear\bigl(L(X,Y)\bigr)$ if and only if $\Phi_2 G \Phi_1\in \Spear\bigl(L(X_1,Y_1)\bigr)$.
\item[(ii)] We may restrict the codomain of a spear operator keeping the property of being spear operator: If $G$ is a spear operator and  $Z$ is a subspace of $Y$ containing $G(X)$, then $G:X\longrightarrow Z$ is a spear operator. On the other hand, the extension of the codomain does not preserve spears: the map $j: \mathbb{K} \longrightarrow \mathbb{K} \oplus_{\infty} \mathbb{K}, \, j(x) = (x,0)$ is not a spear operator.
\item[(iii)] As an easy consequence of (i) and (ii), we get that the following assertions are equivalent: (a) $X$ has numerical index one (i.e.\ $\Id_{X}$ is a spear), (b) there exists a Banach space $Z$ and an isometric isomorphism which is a spear in $L(X,Z)$ or $L(Z,X)$, (c) there exists a Banach space $W$ and an isometric embedding of $X$ into $W$ which is a spear operator.
\end{enumerate}
\end{Rema}

\section{Alternative Daugavet Property}
\label{sec:aDP}

We start presenting the definition of the alternative Daugavet property for an operator, which extends the analogous definition for a Banach space (through the Identity).

\begin{Defi}
Let $X$, $Y$ be Banach spaces. We say that $G\in L(X,Y)$ has the \emph{alternative Daugavet property} (\emph{aDP} in short), if the norm equality
\begin{equation}
\label{equa:aDPequation}\tag{\textrm{aDE}}
\| G + \T\, T \| = 1 + \| T\|
\end{equation}
holds for every rank-one operator $T\in L(X,Y)$.
\end{Defi}
\index{alternative Daugavet property for an operator}%
\index{aDP for an operator}%

Substituting $T=0$ in \eqref{equa:aDPequation} we deduce that if $G$ has the aDP then $\|G\|=1$.

The following fundamental result characterizes the aDP of an operator in terms of the behaviour of the operator with respect to slices, spear sets\ldots

\begin{Theo}
\label{Theo:aDPCharacterization}
Let $G\in L(X,Y)$ be a norm-one operator between two Banach spaces $X$, $Y$, let $\mathcal{B} \subset B_{X}$ with $B_{X} = \overline{\conv}{(\mathcal{B})}$ and let $\mathcal{A} \subset B_{Y^{\ast}}$ with $\overline{\conv}^{\omega^{\ast}}{(\mathcal{A})}=B_{Y^{\ast}}$. The following assertions are equivalent:
\begin{enumerate}
\item[(i)] $G$ has the aDP.
\item[(ii)] $G(S)$ is a spear set for every slice $S$ of $\mathcal{B}$.
\item[(ii$^{\ast}$)] $G^{\ast}(S^{\ast})$ is a spear set for every weak$^{\ast}$-slice $S^{\ast}$ of $\mathcal{A}$.
\item[(iii)] For every $y_{0} \in S_{Y}$ and $\eps > 0$
\[ B_{X} = \overline{\conv}{\bigl( \{ x \in \mathcal{B} \colon \|  Gx + \T\,y_{0}\| > 2 - \eps \}\bigr)}.\]
\item[(iv)] For every $x_{0}^{\ast} \in X^{\ast}$, the set
\[ \bigl\{ y^{\ast} \in \ext B_{Y^{\ast}} \colon \| G^{\ast} y^{\ast} + \T\,x_{0}^{\ast}\| = 1 + \| x_{0}^{\ast}\| \bigr\}
\]
is a dense G$_{\delta}$ set in $(\ext B_{Y^{\ast}}, \weakstar)$.
\end{enumerate}
\end{Theo}

\begin{proof}
(i) $\Rightarrow$ (ii): Let $S = \Slice(\mathcal{B},x_{0}^{\ast}, \eps)$ be a slice of $\mathcal{B}$ where $x_{0}^{\ast} \in S_{X^{\ast}}$ and $0 < \eps < 1$. Given any $0 \neq y_{0} \in Y$ consider the rank-one operator $T= x_{0}^{\ast} \otimes y_{0}\in L(X,Y)$ (i.e.\ $T(x) = x_{0}^{\ast}(x) y_{0}$ for every $x \in X$) which satisfies that $\| T\|=\|y_{0}\|$. Since $\| G + \mathbb{T}\, T\| = \| G^{\ast} + \mathbb{T}\, T^{\ast}\|$, for every $0 < \delta < 1$ there exists $y_{0}^{\ast} \in S_{Y^{\ast}}$ such that
\begin{equation}\label{equa:aDPCharacterizationAux}
\| G^{\ast} y_{0}^{\ast}  + \mathbb{T} x_{0}^{\ast} y_{0}^{\ast}(y_{0}) \| > 1 + \| y_{0}\|(1 - \eps \delta) \end{equation}
Making a rotation of $y_{0}^{\ast}$ if necessary, we can assume that $ 0 \leqslant y_{0}^{\ast}(y_{0}) \leqslant \| y_{0}\|$. Using the hypothesis on $\mathcal{B}$, we deduce from \eqref{equa:aDPCharacterizationAux}  the existence of some $x_{0} \in \mathcal{B}$ satisfying
$$ |y_{0}^{\ast}(G x_{0})| + y_{0}^{\ast}(y_{0}) \re x_{0}^{\ast}(x_{0}) > 1 + \| y_{0}\|(1 - \eps \delta).  $$
Since $\| G x_{0}\| \leqslant 1$, we deduce that $\re x_{0}^{\ast}(x_{0}) > 1 - \eps \delta > 1 - \eps$. Hence $x_{0} \in S$ and so
\[
\begin{split}
\| G(S) + \T\, y_{0}\| & \geqslant \| Gx_{0} + \T\, y_{0} \| > \| Gx_{0} + \T\, x_{0}^{\ast}(x_{0}) y_{0} \| - \eps \delta\\
& > 1 + \| y_{0}\| (1 - \eps \delta) - \| y_{0}\|\eps \delta > 1 + \| y_{0}\|(1 - 2 \delta).
\end{split}
\]

(ii) $\Rightarrow$ (iii): Given $y_{0} \in S_{Y}$, $\eps > 0$ and a slice $S$ of $\mathcal{B}$, since $G(S)$ is a spear set, we can find $x \in \T S$ with $\| Gx + y_{0}\| > 2 - \eps$, which means that every slice $S$ of $\mathcal{B}$ intersects
\[
\T\{ x \in \mathcal{B} \colon \| Gx + y_{0}\| > 2- \eps \} = \{ x \in \mathcal{B} \colon \| Gx + \T\,y_{0}\| > 2- \eps \}.
\]
Therefore
\[
B_{X} = \overline{\conv}{(\mathcal{B})} \subset  \overline{\conv}{\bigl(\{ x \in \mathcal{B} \colon \| Gx + \T\,y_{0}\| > 2- \eps \}\bigl)}.
\]

(iii) $\Rightarrow$ (i): Let $T\in L(X,Y)$ be a rank-one operator. By a convexity argument, we may and do suppose that $\|T\|=1$. Then, it is of the form $T = x_{0}^{\ast} \otimes y_{0}$ for some $y_{0} \in S_Y$ and $x_{0}^{\ast} \in S_{X^{\ast}}$. Given $\eps > 0$, the hypothesis implies that $\Slice(\mathcal{B},x_{0}^{\ast}, \eps)$ intersects $\{ x \in \mathcal{B} \colon \|G x  + \T\, y_{0} \| > 2 - \eps \}$, so there exists $x_{0} \in \Slice(\mathcal{B},x_{0}^{\ast}, \eps)$ and $\theta_{0} \in \T$ such that $\| G(x_{0}) + \theta _{0} y_{0} \| > 2 - \eps$. Hence
\begin{align*}
\| G + \T\,T \| & \geqslant \|G x_{0} +  \theta _{0}  x_{0}^{\ast}(x_{0})y_{0} \| \\ & \geqslant \| G x_{0} + \theta_{0} y_{0} \| - |x_{0}^{\ast}(x_{0}) - 1| > 2 - 2\eps .
\end{align*}

(ii) $\Rightarrow$ (iv): We may and do suppose that $\|x_0^\ast\|=1$. Given $\eps>0$, the set $G\bigl(\Slice(\mathcal{B},x_{0}^{\ast}, \eps)\bigr)$ is a spear of $B_Y$ by hypothesis, which by Corollary \ref{Coroll:spearSet}.(iii) means that the set $\GF\bigl(\ext B_{Y^{\ast}},\T G\bigl(\Slice(\mathcal{B},x_{0}^{\ast}, \eps)\bigr)\bigr)$ is a dense G$_\delta$ set in the Baire space $(\ext B_{Y^{\ast}}, \weakstar)$. Hence,
$$
\bigcap_{m \in \N}\GF\bigl(\ext B_{Y^{\ast}},\T G\bigl(\Slice(\mathcal{B},x_{0}^{\ast}, 1/m)\bigr)\bigr)
$$
is also dense in $\ext B_{Y^{\ast}}$, and it is easy to check that
\[
\bigcap_{m \in \N}\GF\bigl(\ext B_{Y^{\ast}},\T G\bigl(\Slice(\mathcal{B},x_{0}^{\ast}, 1/m)\bigr)\bigr) = \bigl\{ y^{\ast} \in \ext B_{Y^\ast}\colon \| G^{\ast} y^{\ast} + \T\,x_{0}^{\ast} \| = 2\bigr\}.
\]

(iv) $\Rightarrow$ (ii$^\ast$): Fix $x_0^\ast\in X^\ast$. If $S^\ast$ is any weak$^\ast$-slice of $B_{Y^\ast}$, then $S^\ast \cap \ext B_{Y^{\ast}}$ is a non-empty open subset of $(\ext B_{Y^{\ast}}, \weakstar)$. The hypothesis implies that the slice $S^{\ast}$ contains an element of $\{ y^{\ast} \in \ext B_{Y^{\ast}}\colon \| G^{\ast}y^{\ast} + \T\,x_{0}^{\ast} \| = 1 + \| x_{0}^{\ast}\|\}$, so $\|G^{\ast}(S^{\ast}) + \T\, x_{0}^{\ast}\| = 1 + \| x_{0}^{\ast}\|$.

(ii$^\ast$) $\Rightarrow$ (i): Let $T = x_{0}^{\ast} \otimes y_{0}$, where $x_{0}^{\ast} \in X^{\ast}$ and $y_{0} \in S_Y$, be an arbitrary rank-one operator. Given any $\eps > 0$ put $S^{\ast} = \Slice(\mathcal{A},y_{0}, \eps)$. Notice that for every $y^{\ast} \in S^\ast$,
$$
\| T^\ast y^{\ast} - x_{0}^\ast\| = \|y^{\ast}(y_{0}) x_{0}^\ast - x_{0}^\ast\| < \|x_0^\ast\|\eps,
$$
so using that $G^{\ast}(S^{\ast})$ is a spear we deduce that
\begin{equation*}
\| G + \T\, T \| = \| G^{\ast} + \T\, T^{\ast}\| \geqslant \| G^{\ast}(S^{\ast}) + \T\,x_{0}^{\ast} \| - \|x_0^\ast\|\eps = 1 + (1-\eps)\| x_{0}^{\ast}\|.\qedhere
\end{equation*}
\end{proof}

The next result shows that the aDP is separably determined, and will be very useful in the next section where we deal with SCD operators.

\begin{Prop}
\label{Prop:separablyDetermined}
Let $X$, $Y$ be Banach spaces and let $G\in L(X,Y)$. Then, $G$ has the aDP if and only if for every separable subspaces $X_{0} \subset X$ and $Y_{0} \subset Y$, there exist separable subspaces $X_{\infty}$, $Y_{\infty}$ satisfying $X_{0} \subset X_{\infty} \subset X$ and $Y_{0} \subset Y_{\infty} \subset Y$ and such that $G(X_{\infty}) \subset Y_{\infty}$ and $G|_{X_{\infty}} : X_{\infty} \longrightarrow Y_{\infty}$ has  the aDP.
\end{Prop}

\begin{proof}
Suppose first that $G$ has the aDP. Pick a sequence $(x_n)_{n\in \N}$ of $S_X$ with $\sup_n \|G x_n\|=1$ and consider $X_{1} = \overline{\spn}(X_{0}\cup \{x_n\colon n\in \N\})$ and $Y_{1} = \overline{Y_{0} + G(X_{1})}$, both separable subspaces. By Theorem \ref{Theo:aDPCharacterization}.(iii) we have that
\[ B_{X_{1}} \subset \overline{\conv}{\left( \{ x \in S_{X} \colon \| Gx + \T\,y_{1}\| > 2 - \eps\} \right)} \mbox{ \: for every $y_{1} \in S_{Y_{1}}$ and $\eps > 0$.} \]
But since $B_{X_{1}}$ and $S_{Y_{1}}$ are separable, it is easy to deduce the existence of a countable set $A_{1} \subset S_{X}$ such that
\[ B_{X_{1}} \subset \overline{\conv}{\left(  \{ x \in A_{1} \colon \|G x + \T\,y_{1} \| > 2 - \eps\} \right)} \mbox{ for every $y_{1} \in S_{Y_{1}}$ and $\eps > 0$}. \]
Define then $X_{2} = \overline{\spn}{(X_{1} \cup A_{1})}$ and $Y_{2} = \overline{Y_{1} + G(X_{2})}$, which are again separable. Repeating the same process as above, we can construct an increasing sequence of closed separable subspaces $X_{n} \subset X$ and $\overline{G(X_{n})} \subset Y_{n} \subset Y$ such that
\[
B_{X_{n}} \subset \overline{\conv}{\left(\{ x \in S_{X_{n+1}} \colon \| Gx +  \T\,y_{n} \| > 2 - \eps \} \right)} \mbox{ \: for every $y_{n} \in S_{Y_{n}}$ and $\eps > 0$.}
\]
This implies that $X_{\infty} := \overline{\bigcup_{n \in \N}{X_{n}}}$ and $Y_{\infty}:= \overline{\bigcup_{n \in \N}{Y_{n}}}$, satisfy that
\[
B_{X_{\infty}} \subset \overline{\conv}{\left(\{ x \in S_{X_{\infty}}\colon \| Gx + \T\, y\| > 2 - \eps \} \right)} \mbox{ \: for every $y \in S_{Y_{\infty}}$ and $\eps > 0$,}
\]
which means that $G:X_{\infty} \longrightarrow Y_{\infty}$ has the aDP by using again Theorem \ref{Theo:aDPCharacterization}.(iii).

Conversely, take a non-null rank-one operator $T\in L(X,Y)$, consider a separable subspace $X_0\subset X$ such that $\|G|_{X_0}\|=\|G\|=1$ and $\|T|_{X_0}\|=\|T\|$, and write $Y_0=\overline{G(X_0) + T(X_0)}$. By hypothesis, there are separable subspaces $X_0\subset X_\infty \subset X$ and $Y_0\subset Y_\infty \subset Y$ such that $G|_{X_\infty}:X_\infty \longrightarrow Y_\infty$ has norm one and has the aDP. As $T$ is rank-one and $T|_{X_0}\neq 0$, it follows that $T(X)\subset T(X_0)\subset Y_0 \subset Y_\infty$ and $\|T|_{X_\infty}\|=\|T\|$. Then we may apply that $G|_{X_\infty}$ has the aDP to get that $\|G|_{X_\infty} + \T\, T|_{X_\infty}\|=1+\|T\|$. But, clearly,
\[
\|G + \T\, T\|\geqslant \|G|_{X_\infty} + \T\, T|_{X_\infty}\|=1+\|T\|,
\]
and the reverse inequality is always true, so $G$ has the aDP.
\end{proof}

As we did for spear operators, we may directly deduce from the definition the following three elementary results about operators with the aDP.

\begin{Rema}\label{Rema:aDPOperators}
Let $X,Y$ be Banach spaces and $G\in L(X,Y)$.
\begin{enumerate}
\item[(i)] The composition with isometric isomorphisms preserves the aDP: If $X_1$, $Y_1$ are Banach spaces and $\Phi_1\in L(X_1,X)$, $\Phi_2\in L(Y,Y_2)$ are isometric isomorphisms, then, $G\in L(X,Y)$ has the aDP if and only if $\Phi_2 G \Phi_1\in L(X_1,Y_1)$ has the aDP.
\item[(ii)] If $G$ has the aDP and $Z$ is a subspace of $Y$ containing $G(X)$, then $G:X\longrightarrow Z$ has the aDP. However, the property of aDP is not preserved by extending the codomain of the operator, as the same example of Remark \ref{Rema:SpearOperators} shows.
\item[(iii)] As an easy consequence of (i) and (ii), we have that the following statements are equivalent: (a) $X$ has the aDP, (b) there exist a Banach space $Z$ and an isometric isomorphism in $L(X,Z)$ or in $L(Z,X)$ which has the aDP, (d) there exist a Banach space $W$ and an isometric embedding $G\in L(X,W)$ which has the aDP.
\end{enumerate}
\end{Rema}

\section{Target operators}
\label{sec:StrongTargetOperators}

Our goal in this section is to present and study the concept of target operator, which will be the key in the next section to relate the aDP and lushness so, in particular, to relate the aDP and spear operators. As far as we know, this is a new concept even in the particular case in which $G$ is the identity operator of a Banach space.

\begin{Defi}\label{Defi:target}
Let $X$, $Y$, $Z$ be Banach spaces and let $G\in L(X,Y)$ be a norm-one operator. We say that $T \in L(X,Z)$ is a \emph{target} for $G$ if each $x_{0} \in B_{X}$ has the following property:
\begin{equation}\label{eq:diamond}\tag{\ensuremath{\Diamond}}
\begin{aligned}
 &\text{For every $\eps > 0$ and every $y \in S_{Y}$, there is $F \subset B_{X}$ such that} \\
 &\conv{F} \subset \{ x \in B_{X} \colon \| G x + y\| > 2 - \eps \} \: \mbox{ and }\: \dist\bigl(Tx_{0}, T\bigl(\aconv(F)\bigr)\bigr) < \eps.
\end{aligned}
\end{equation}
\end{Defi}
\index{target}%
\index{00diamond@$\Diamond$}%
Remark that if $F$ satisfies \eqref{eq:diamond} then there is a finite subset of $F$ satisfying the same condition.

At the end of the section we will include a result characterizing spear vectors in terms of target operators which will allow to better understand this definition, see Proposition \ref{Prop:char-spearpoints-with-diamond}. Next, we provide with several characterizations of this kind of operators which will be very useful in the sequel.

\begin{Prop}\label{Prop:vectorsPropertyDiamond}
Let $X$, $Y$, $Z$ be Banach spaces, let $G\in L(X,Y)$ with $\|G\|=1$, let $T\in L(X,Z)$, and let $\mathcal{A} \subset B_{Y^{\ast}}$ with $\overline{\conv}^{\omega^{\ast}}{(\mathcal{A})} = B_{Y^{\ast}}$. Given $x_0\in B_X$, the following assertions are equivalent:
\begin{enumerate}
\item[(i)] $x_{0}$ satisfies \eqref{eq:diamond}.
\item[(ii)] For every $\eps > 0$ and $y \in S_{Y}$ there is $y^{\ast} \in \Slice(\mathcal{A}, y, \eps)$ such that
$$\dist\bigl(T x_{0}, T\bigl(\aconv{\GS(S_{X}, y^{\ast}, \eps)}\bigr)\bigr) < \eps.$$
\item[(iii)]  For every $\eps > 0$, the set
$$
\mathcal{D}_{T}^{\varepsilon}(\mathcal{A}, x_{0}) = \bigl\{ y^{\ast} \in \mathcal{A} \colon \dist{\bigl(Tx_{0}, T(\aconv{\GS(S_{X}, G^{\ast}y^{\ast}, \eps)})\bigr)} < \eps\bigr\}
$$
intersects every $\omega^{\ast}$-slice of $\mathcal{A}$.
\item[(iv)] The set
$$
\mathcal{D}_{T}(x_{0}) = \left\{ y^{\ast} \in \ext{B_{Y^\ast}} \colon Tx_{0} \in \overline{T\bigl(\aconv{ \GS(S_{X}, G^{\ast}y^{\ast}, \eps)}\bigr)} \mbox{ for every $\eps > 0$}\right\}
$$
 is a dense (G$_{\delta}$) subset of $(\ext B_{Y^\ast}, w^{\ast})$.
\end{enumerate}
\end{Prop}

\begin{proof}
(i) $\Rightarrow$ (ii): Let $\eps > 0$ and $y \in S_{Y}$. Fixed $0<\delta<1$ such that $\delta^2+\delta+\frac{\delta}{1-\delta} <\eps$, by \eqref{eq:diamond} in Definition \ref{Defi:target}, we can find $F = \{ x_{1}, \ldots, x_{n}\}\subset B_X$, $\lambda_{1}, \ldots, \lambda_{n} \geqslant 0$ with $\sum{\lambda_{k}} = 1$, and $\theta_{1}, \ldots, \theta_{n}\in \mathbb{T}$ such that
\begin{equation}\label{eq:propDiamondAux1}
\left\| \sum_{k=1}^{n}{\lambda_{k} G (x_{k})} + y \right\|>2-\delta^2 \qquad \mbox{ and } \qquad \left\| T x_{0} - \sum_{k=1}^{n}{\lambda_{k} \theta_{k} T (x_{k})} \right\| < \delta^{2}. \end{equation}
Let $a^{\ast} \in \mathcal{A}$ be such that
$$ \Real{a^{\ast}}\left(\sum_{k=1}^{n}{\lambda_{k} G (x_{k})}+ y\right) > 2 - \delta^2. $$
Then $\Real{a^{\ast}(y)} >1-\delta^2> 1 - \eps$ and, moreover, $J = \{ k \colon \Real{a^{\ast}}(G x_{k}) > 1 - \delta \}$ satisfies
$$ 1 - \delta \sum_{k \notin J}{\lambda_{k}} = \sum_{k \in J}{\lambda_{k}} +  (1 - \delta) \sum_{k \notin J}{\lambda_{k} } \geqslant \sum_{k=1}^{n}{\lambda_k\Real{a^{\ast}}(G x_{k})} > 1 - \delta^{2}.  $$
Hence, we get that
\begin{equation}\label{eq:propDiamondAux2}
\sum_{k \notin J}{\lambda_{k}}<\delta.
\end{equation}
Using the right-hand side inequality of \eqref{eq:propDiamondAux1} and \eqref{eq:propDiamondAux2}, we can write
\begin{align*}
\left\| T x_{0} - \sum_{k \in J}{\frac{\lambda_{k}}{\sum_{j \in J}{\lambda_{j}}} \theta_{k} T(x_{k})} \right\| & \leqslant \delta^{2}  + \left\| \sum_{k=1}^{n}{ \lambda_{k} \theta_{k} T x_{k} } - \sum_{k \in J}{\frac{\lambda_{k}}{\sum_{j \in J}{\lambda_{j}}} \theta_{k} T x_{k}} \right\|\\
& \leqslant \delta^{2} + \delta + \left\| \sum_{k \in J}{ \lambda_{k} \theta_{k} T x_{k} } - \sum_{k \in J}{\frac{\lambda_{k}}{\sum_{j \in J}{\lambda_{j}}} \theta_{k} T x_{k}} \right\|\\
& \leqslant \delta^{2} + \delta + \left| 1 - \frac{1}{\sum_{k \in J}{\lambda_{k}}} \right|\\
& \leqslant \delta^{2} + \delta + \frac{\delta}{1 - \delta}<\eps.
\end{align*}
(ii) $\Rightarrow$ (iii): statement (ii) claims that $\mathcal{D}_T^\delta(\mathcal{A},x_0)$ intersects $\Slice(\mathcal{A},y,\delta)$ for every $y\in S_Y$ and every $\delta>0$. Since for every $\varepsilon>0$, $\mathcal{D}_T^\varepsilon(\mathcal{A},x_0)$ contains $\mathcal{D}_T^\delta(\mathcal{A},x_0)$ for every $0<\delta<\varepsilon$, we conclude the result.

\noindent(iii) $\Rightarrow$ (iv): Using Lemma \ref{Kreinlemma}.(a), we have that $\mathcal{D}_{T}^{\eps}(\ext{B_{Y^{\ast}}},x_{0})$ intersects every $\omega^{\ast}$-open subset of $\ext{B_{Y^{\ast}}}$. In other words, $\mathcal{D}_{T}^{\eps}(\ext{B_{Y^{\ast}}},x_{0})$ is a dense subset of $(\ext{B_{Y^{\ast}}}, \omega^{\ast})$. Since
$$ \mathcal{D}_{T}(x_{0}) = \bigcap_{m \in \N}{\mathcal{D}_{T}^{1/m}(\ext{B_{Y^{\ast}}},x_{0})}, $$
we just have to show that $\mathcal{D}_{T}^{\eps}(\ext{B_{Y^{\ast}}},x_{0})$ is open and then apply Lemma \ref{Kreinlemma} (c). Indeed, notice that for every $a_{0}^{\ast} \in \mathcal{D}_{T}^{1/n}(\mathcal{A},x_{0})$ we can find a finite subset $F$ of $\GS(S_{X},G^{\ast}a_{0}^{\ast}, 1/n)$ such that $\dist(Tx_{0}, T(\aconv{F})) < 1/n$. The set
$U:=\bigcap_{x \in F} \GS(\mathcal{A},G x, 1/n)$ is a relatively weak$^{\ast}$-open subset of $\mathcal{A}$, which contains $y_{0}$ by definition. Also $U \subset \mathcal{D}_{T}^{1/n}(\mathcal{A}, x_{0})$, since $F \subset \GS(S_{X},G^{\ast}y^{\ast}, 1/n)$ for every $y^{\ast} \in U$. So, $\mathcal{D}_{T}^{\eps}(\ext{B_{Y^{\ast}}},x_{0})$ is open as desired.

(iv) $\Rightarrow$ (i): Given $y \in S_{Y}$ and $\eps > 0$, $\mathcal{D}_{T}(x_{0})$ intersects $\Slice(\ext{B_{Y^{\ast}}}, y, \eps)$. Taking an element $y^{\ast}$ in such intersection,  by definition of $\mathcal{D}_{T}(x_{0})$, we can find a finite set $F$ in $\GS(S_{X}, G^{\ast} y^{\ast}, \eps)$ such that
$$ \dist{(T x_{0}, T(\aconv{F}))} < \eps. $$
But the condition $F \subset \GS(S_{X}, G^{\ast} y^{\ast}, \eps)$ yields that every $x \in \conv{F}$ satisfies
\[\| Gx + y \| \geqslant \Real{y^{\ast}}(G x) + \Real{y^{\ast}(y)} > 2 - 2 \eps.\qedhere \]
\end{proof}

One of the main applications of target operators is the following result.

\begin{Prop}
\label{Prop:aDpPlusTarget}
Let $X$, $Y$ be Banach spaces and let $G\in L(X,Y)$ be a norm-one operator. If $T \in L(X,Y)$ is a target for $G$, then $$\|G + \T\,T\| = 1 + \| T\|.$$
\end{Prop}

\begin{proof}
We can assume that $\| T\| = 1$. For $0 < \eps < 1$, take $x_{0} \in B_X$ with $\| T x_{0}\| > 1- \eps$ and write $y_{0} := Tx_{0} / \| T x_{0}\|$. By Proposition \ref{Prop:vectorsPropertyDiamond}, there exists $y_{0}^{\ast} \in \Slice(S_{Y^{\ast}},y_{0}, \eps)$ with
$$
\dist\bigl(Tx_{0}, T\bigl(\aconv\GS(S_{X},G^{\ast}y_0^{\ast}, \eps)\bigr)\bigr) < \eps.
$$
We can then find $n \in \N$ and elements $x_{k} \in \GS(S_{X},G^{\ast}y_{0}^{\ast}, \eps)$, $\theta_{k} \in \T$, $\lambda_{k} \geqslant 0$ for $k=1,\ldots,n$ such that
\[ \sum_{k=1}^{n}{\lambda_{k}}=1  \qquad \text{and} \qquad \left\| T x_{0} - T \left( \sum_{k=1}^{n}{\lambda_{k} \theta_{k} x_{k}} \right) \right\| < \eps. \]
Since $|y_{0}^{\ast}(T x_{0})| > 1 - \eps$, a standard convexity argument leads to the existence of some $k \in \{ 1, \ldots, n\}$ with $|y_{0}^{\ast}(T x_{k})| \geqslant 1 - 2 \eps$. Therefore
\begin{equation*}
\| G + \T\, T\| \geqslant |y_{0}^{\ast}(Gx_{k})| + |y_{0}^{\ast}(T x_{k})| > 2 - 3 \eps.\qedhere
\end{equation*}
\end{proof}

The following observations can be proved directly from the definition of target, but they follow easier from Proposition \ref{Prop:vectorsPropertyDiamond}.

\begin{Rema}\label{Rema:target-norm-one}
Let $X$, $Y$, $Z$, $Z_1$, $Z_2$ be Banach spaces and let $G\in L(X,Y)$ be a norm-one operator.
\begin{enumerate}
\item[(a)] If $T\in L(X,Z)$ is a target for $G$, then $\lambda T$ is a target for $G$ for every $\lambda\in \mathbb{K}$.
\item[(b)] Let $T_1\in L(X,Z_1)$ and $T_2\in L(X,Z_2)$ be operators such that $\|T_2 x\|\leqslant \|T_1 x\|$ for every $x\in X$. If $T_1$ is a target for $G$, then so is $T_2$.
\end{enumerate}
\end{Rema}

Target operators are separably determined, and this fact will be crucial in the study of the relationship with the aDP and SCD.

\begin{Theo}
\label{Theo:targetSeparablyDetermined}
Let $X$, $Y$, $Z$ be Banach spaces, let $G\in L(X,Y)$ be a norm-one operator and consider $T \in L(X,Z)$. Then, $T$ is a target for $G$ if and only if for every separable subspaces $X_{0} \subset X$ and $Y_{0} \subset Y$, there exist separable subspaces $X_{\infty}$, $Y_{\infty}$ satisfying  $X_{0} \subset X_{\infty} \subset X$ and $Y_{0} \subset Y_{\infty} \subset Y$ and such that $G|_{X_{\infty}} \in L(X_{\infty}, Y_{\infty})$ has norm one and $T|_{X_{\infty}} \in L(X_{\infty}, Z)$ is a target for $G|_{X_{\infty}}$.
\end{Theo}

\begin{proof}
Let us assume first that $T$ is a target for $G$.

\noindent\emph{Claim:}  given separable subspaces $\widetilde{X}\subset X$ and $\widetilde{Y}\subset Y$, we can find a countable set $B \subset B_{X}$ with the following property (P): \emph{given $x_{0} \in B_{\widetilde{X}}$, $y_{0} \in S_{\widetilde{Y}}$ and $\eps > 0$, there exists $F \subset B$ with $\conv{(F)} \subset \{ x \in B_{X} \colon \| Gx + y \| > 2 - \eps \}$ and $\dist(T x_{0}, T(\aconv{F})) < \eps$}.

Indeed, fixing $C_{0}$ and $D_{0}$ countable dense subsets of $B_{\widetilde{X}}$ and $B_{\widetilde{Y}}$ respectively, we can apply the definition of target operator to construct a countable set $B \subset B_{X}$ satisfying the property (P) for all $x_{0} \in C_{0}$, $y_{0} \in D_{0}$ and $\eps \in \mathbb{Q}^{+}$. But using the density of $C_{0}$ and $D_{0}$, it turns out that $B$ has the same property for each $x_{0} \in B_{\widetilde{X}}$, $y_{0} \in B_{\widetilde{Y}}$ and $\eps > 0$.

Now we prove the theorem: We can assume that $\| G|_{X_{0}}\| = \| G\| = 1$. Put $X_{1} = X_{0}$ and $Y_{1} = \overline{\spn}{(Y_{0} \cup G(X_{0}))}$, both separable Banach spaces. Using the claim, we deduce the existence of a countable set $B_{1} \subset B_{X}$ with the property (P) for $X_{1}$ and $Y_{1}$. Define $X_{2}=\overline{\spn}{(X_{1} \cup B_{1})}$ and $Y_{2} = \overline{\spn}{(Y_{1} \cup G(X_{2}))}$. Repeating this process inductively, we construct increasing sequences of closed separable subspaces $X_{n} \subset X$ and $G(X_{n}) \subset Y_{n} \subset Y$ such that $B_{X_{n+1}}$ has the property (P) above for $X_{n}$ and $Y_{n}$. Taking $X_{\infty}:=\overline{\bigcup_{n \in \N}{X_{n}}}$ and $Y_{\infty}:=\overline{\bigcup_{n \in \N}{Y_{n}}}$, we conclude the result, as $G|_{X_{\infty}}$ and $T|_{X_{\infty}}$ satisfy the definition of target operator by construction.

Let us check the converse implication. Given $x_{0} \in B_{X}$ and $y_{0} \in S_{Y}$ we can find separable subspaces $x_{0} \in X_{\infty} \subset X$ and $y_{0} \in Y_{\infty} \subset Y$ with the properties above, so applying the definition of target for $T|_{X_{\infty}}$, $G|_{X_{\infty}}$ and the previous elements we get the result.
\end{proof}

We need one more ingredient to be able to present the main result about the relationship between target operators and SCD operators.

\begin{Prop}
\label{Prop:target-W*T-dense-sufficient}
Let $X$, $Y$, $Z$ be Banach spaces and let $G\in L(X,Y)$ be a norm-one operator. Let $T \in L(X,Z)$ be an operator such that the set
\begin{equation*}
\mathcal{D}_T:=\left\{ y^{\ast} \in \ext B_{Y^{\ast}} \colon T(B_{X}) \subset \overline{T(\aconv{ \GS(B_{X}, G^{\ast}y^{\ast}, \eps)})} \mbox{\, for every $\eps > 0$} \right\}
\end{equation*}
is dense in $(\ext B_{Y^{\ast}}, \weakstar)$.
Then, $T$ is a target for $G$, and in the case of $Z=Y$, we have $\|G +\T\,T\| = 1 + \| T\|$.
\end{Prop}

\begin{proof}
In the notation of Proposition \ref{Prop:vectorsPropertyDiamond}, the inclusion $\mathcal{D}_T \subset \mathcal{D}_{T}(x_{0})$ holds for every $x_{0} \in B_{X}$, so by Proposition \ref{Prop:vectorsPropertyDiamond}.(iv), this means that $T$ is a target for $G$. If $Z = Y$, an application of Proposition \ref{Prop:aDpPlusTarget} implies that $\|G + \T\,T\| = 1 + \| T\|$.
\end{proof}

The converse of the above result holds for operators with separable image.

\begin{Prop}
\label{Prop:target-W*T-dense-necessary-separable}
Let $X$, $Y$, $Z$ be Banach spaces and let $G\in L(X,Y)$ be a norm-one operator. Suppose that $T \in L(X,Z)$ is a target for $G$ such that $T(X)$ is separable. Then,
\begin{equation*}
\mathcal{D}_T:=\left\{ y^{\ast} \in \ext B_{Y^{\ast}} \colon T(B_{X}) \subset \overline{T(\aconv{ \GS(B_{X}, G^{\ast}y^{\ast}, \eps)})} \mbox{\, for every $\eps > 0$} \right\}
\end{equation*}
is a dense G$_\delta$ subset of $(\ext B_{Y^{\ast}}, \weakstar)$.
\end{Prop}

\begin{proof}
Since $T$ is a target for $G$, for every $x_0\in B_X$ we have by Proposition~\ref{Prop:vectorsPropertyDiamond} that the set
$$
\mathcal{D}_{T}(x_{0}):= \left\{ y^{\ast} \in \ext B_{Y^{\ast}} \colon Tx_{0} \in \overline{T\bigl(\aconv (\GS(S_{X},G^{\ast}y^{\ast}, \eps))\bigr)} \mbox{ for every $\eps > 0$} \right\}
$$
is a dense G$_\delta$ subset of $(\ext B_{Y^{\ast}}, \weakstar)$. If we choose a sequence $(x_n)_{n\in \N}$ in $B_X$ so that $(Tx_n)_{n \in \N}$ is dense in $T(B_{X})$, then  $\mathcal{D}_T = \bigcap_{n \in \N}\mathcal{D}_{T}(x_{n})$. Since all $\mathcal{D}_{T}(x_{n}, \ext{B_{Y^\ast}})$ are dense $G_{\delta}$ subsets of $(\ext B_{Y^{\ast}}, \weakstar)$, so is $\mathcal{D}_T$ (see Lemma~\ref{Kreinlemma}.(c)).
\end{proof}

\begin{Theo}
\label{Theo:SCDOperators}
Let $X$, $Y$, $Z$, $Z_1$ be Banach spaces and let $G\in L(X,Y)$ be an operator with the aDP. If for $T \in L(X,Z)$ there is an SCD operator $T_1\in L(X,Z_1)$ such that $\|Tx\|\leqslant \|T_1x\|$ for every $x\in X$, then $T$ is a target for $G$. In the case of $Z=Y$, we have $\|G +\T\,T\| = 1 + \| T\|$.
\end{Theo}

\begin{proof}
By Remark~\ref{Rema:target-norm-one} we may assume that $T$ is an SCD operator and that $\| T\| \leqslant 1$. Let $\bigl\{ \widehat{S}_{n} \colon n \in \N\bigr\}$ be a determining family of slices for $T(B_{X})$. Then,  $S_{n}:=T^{-1}(\widehat{S}_{n}) \cap S_{X}$ is a slice of $S_{X}$ for each $n\in\N$. Since $G$ has the aDP, Theorem~\ref{Theo:aDPCharacterization}.(ii) tells us that $G(S_{n})$ is a spear set for every $n \in \N$, which implies that $\GF(\ext B_{Y^{\ast}},\T G(S_{n}))$ is a dense $G_{\delta}$ set in $(\ext B_{Y^{\ast}}, \weakstar)$ by Corollary~\ref{Coroll:spearSet}.(iii). As $(\ext B_{Y^{\ast}}, \weakstar)$ is a Baire space (see Lemma~\ref{Kreinlemma}.(c)), we deduce that the intersection $\bigcap_{n \in \N}{\GF(\ext B_{Y^{\ast}},\T G(S_{n}))}$ is weak$^{\ast}$-dense in $\ext B_{Y^{\ast}}$. Observe that, by Proposition~\ref{Prop:target-W*T-dense-sufficient}, it suffices to show that this intersection is contained in $\mathcal{D}_T$. Given $y_{0}^\ast$ belonging to this intersection, we have that for every $n \in \N$ and $\eps > 0$, $G(S_{n}) \cap \T \Slice(B_{Y},y_{0}^\ast, \eps) \neq \emptyset$. Therefore, $S_{n} \cap \T \GS(B_{X},G^{\ast}y_{0}^{\ast}, \eps) \neq \emptyset$, and so $T(\T \GS(B_{X}, G^{\ast}y_{0}^{\ast}, \eps)) \cap \widehat{S}_{n} \neq \emptyset$. Using that the family $\bigl\{\widehat{S}_{n} \colon n \in \N\bigr\}$ is determining for $T(B_{X})$, we conclude that
\begin{equation*}
T(B_{X}) \subset \overline{\conv{(T(\T \GS(B_{X}, G^{\ast}y_{0}^{\ast}, \eps)))}} = \overline{T(\aconv{ \GS(B_{X}, G^{\ast}y_{0}^{\ast}, \eps)})}
 \end{equation*}
and, therefore, $y^\ast_0\in \mathcal{D}_T$.
\end{proof}

As a consequence of the previous results we may present a class of operators which is a two-sided operator ideal consisting of operators $T$ satisfying the condition $\|G+\T\,T\|=1+\|T\|$ whenever $G$ has the aDP. Let us recall the needed definitions which we borrow from \cite{SCDsets} and \cite{Ka-She}. Let $X$, $Y$ be  Banach spaces, an operator $T\in L(X,Y)$ is \emph{hereditarily SCD}  if $T(B_X)$ is a hereditarily SCD set, that is, if every convex subset $B$ of $T(B_X)$ is SCD. Obviously, hereditarily SCD operators are SCD. The operator $T$ is \emph{HSCD-majorized} if there is a Banach space $Z$ and a hereditarily SCD operator $\widetilde{T}\in L(X,Z)$ such that $\|Tx\|\leqslant \|\widetilde{T}x\|$ for every $x\in X$. It is shown in \cite[Theorem~3.1]{Ka-She} that the class of HSCD-majorized operators is a two sided operator ideal. By Example \ref{Exam:SCD-set-spaces-operators}, this ideal contains those operators with separable range such that the image of the unit ball has the Radon-Nikod\'{y}m Property, or the convex point of continuity property, or it is an Asplund set, and those operators with separable rank which do not fix copies of $\ell_1$.
\index{hereditarily SCD}%
\index{HSCD-majorized}%

\begin{Coro}
\label{Coro:HSCD-majorized-Operators}
Let $X$, $Y$, $Z$ be Banach spaces and let $G\in L(X,Y)$ be an operator with the aDP. If $T \in L(X,Z)$ is a HSCD-majorized operator then $T$ is a target for $G$. In the case of $Z=Y$, we have $\|G +\T\,T\| = 1 + \| T\|$.
\end{Coro}

\begin{Coro}
\label{Coro:Operator-ideal-aDP}
Let $X$, $Y$ be Banach spaces and let $G\in L(X,Y)$ be an operator with the aDP. The class of operators $T \in L(X,Y)$ satisfying $\|G +\T\,T\| = 1 + \| T\|$ contains the component in $L(X,Y)$ of the two-sided operator ideal formed by HSCD-majorized operators.
\end{Coro}

Even in the case when $G=\Id$, the result above was unknown.

We can extend Theorem~\ref{Theo:SCDOperators} to the non-separable setting in the following way.

\begin{Prop}
\label{Prop:SCDaDPseparablesubspaces}
Let $X$, $Y$, $Z$ be Banach spaces and let $G\in L(X,Y)$. If $G$ has the aDP and $T \in L(X,Z)$ satisfies that $T(B_{X_{0}})$ is an SCD set for every separable subspace $X_{0}$ of $X$, then $T$ is a target for $G$. Therefore, if $Z=Y$ then $\|G + \T\,T \| = 1 + \| T\|$.
\end{Prop}

\begin{proof}
Our aim is to use Theorem \ref{Theo:targetSeparablyDetermined} to deduce that $T$ is a target for $G$. To do so let us fix separable subspaces $X_{0} \subset X$ and $Y_{0} \subset Y$. By Proposition \ref{Prop:separablyDetermined} we can find separable subspaces $X_{0} \subset X_{\infty} \subset X$ and $Y_{0} \subset Y_{\infty} \subset Y$ such that $G(X_{\infty}) \subset Y_{\infty}$ and $G|_{X_{\infty}}: X_{\infty} \longrightarrow Y_{\infty}$ has norm one and the aDP. Now, as $T|_{X_\infty}:X_\infty\longrightarrow Z$ is SCD, Theorem~\ref{Theo:SCDOperators} tells us that $T|_{X_\infty}$ is a target for $G|_{X_\infty}$ and we can apply Theorem~\ref{Theo:targetSeparablyDetermined} to get that $T$ is a target for $G$.
\end{proof}

Using the known results about SCD sets (which were commented in section \ref{subsec:thecaseoftheId}), we may provide with the following consequence.

\begin{Coro}\label{Coro:aDP+SCD=target}
Let $X$, $Y$, $Z$ be Banach spaces and let $G\in L(X,Y)$. Suppose that $G$ has the aDP and $T \in L(X,Z)$ satisfies that $T(B_{X})$ has one of the following properties: Radon-Nikod\'{y}m Property, Asplund Property, convex point of continuity property or absence of $\ell_1$-sequences. Then, $T$ is a target for $G$. Therefore, if $Z=Y$ then $\|G + \T\,T \| = 1 + \| T\|$.
\end{Coro}

\begin{proof}
In \cite[\S 5]{SCDsets} (see Examples~\ref{Exam:SCD-set-spaces-operators}) it is shown that any of the previous properties implies that the requirements of Proposition \ref{Prop:SCDaDPseparablesubspaces} are satisfied.
\end{proof}

To finish the discussion about which operators are targets for a given aDP operator, we may also extend Corollary \ref{Coro:Operator-ideal-aDP} to operators with non separable range as follows. Given two Banach spaces $X$ and $Y$, consider the class of those operators $T\in L(X,Y)$ such that for every separable subspace $X_0$ of $X$, $T|_{X_0}$ is HSCD-majorized. As a consequence of the cited result \cite[Theorem~3.1]{Ka-She}, this class is a two sided operator ideal. Therefore, extending straightforwardly the proof of Proposition \ref{Prop:SCDaDPseparablesubspaces}, we get the following result.

\begin{Coro}
\label{Coro:Operator-ideal-aDP-nonseparable}
Let $X$, $Y$ be Banach spaces and let $G\in L(X,Y)$ be an operator with the aDP. The class of operators $T \in L(X,Y)$ satisfying $\|G +\T\,T\| = 1 + \| T\|$ contains the component in $L(X,Y)$ of the two-sided operator ideal of those operators such that their restrictions to separable subspaces are HSCD-majorized. Moreover, this ideal contains those operators for which the image of the unit ball has one of the following properties: Radon-Nikod\'{y}m Property, Asplund Property, convex point of continuity property, or absence of $\ell_1$-sequences.
\end{Coro}

We finish this section with two results concerning elements satisfying property \eqref{eq:diamond} in Definition \ref{Defi:target}.

\begin{Prop}\label{Prop:char-spearpoints-with-diamond}
Let $X$ be a Banach space and let $x_{0} \in B_{X}$. Then, $x_{0}$ is a spear vector if and only if $x_{0}$ belongs to $\ext B_{X^{\ast \ast}}$ and satisfies \eqref{eq:diamond} with $G=T=\Id_{X}$.
\end{Prop}

\begin{proof}
Using Proposition \ref{Prop:vectorsPropertyDiamond}, $x_{0}$ has property \eqref{eq:diamond} for $G=T = \Id_{X}$ if and only if $$\mathcal{D}(x_{0}) := \left\{ x^{\ast} \in \ext B_{X^{\ast}} \colon x_{0} \in \overline{\aconv}\bigl(\GS(S_{X},x^{\ast}, \eps)\bigr) \mbox{ for every $\eps > 0$} \right\}$$ is dense in $(\ext B_{X^{\ast}}, \weakstar)$. If $x_{0}$ is a spear, then $x_{0} \in \ext B_{X^{\ast \ast}}$ by Proposition \ref{Prop:spearVectorsProperties}.(b). Moreover, the definition of spear yields that $|x^{\ast}(x_{0})| = 1$ for each $x^{\ast} \in \ext B_{X^{\ast}}$, and so $ \mathcal{D}(x_{0}) = \ext B_{X^{\ast}}$. Let us see the converse. If $x_{0} \in \ext B_{X^{\ast \ast}}$ then for each $x^{\ast} \in \mathcal{D}(x_{0})$ and $\eps > 0$, we have that $x_{0}$ is an extreme point of $\overline{\aconv}^{\sigma(X^{\ast \ast}, X^{\ast})}\bigl(\GS(S_{X},x^{\ast}, \eps)\bigr)$. By Milman's Theorem (see Lemma \ref{Kreinlemma}.(b)), we deduce that $x_{0} \in \overline{\GS(S_{X},\T x^{\ast}, \eps)}^{\sigma(X^{\ast \ast}, X^{\ast})}$, and so $|x^{\ast}(x_{0})| \geqslant 1 - \eps$. Since $\eps > 0$ is arbitrary and $\mathcal{D}(x_{0})$ is weak$^\ast$-dense in $\ext B_{X^{\ast}}$, we conclude that $|x^{\ast}(x_{0})|=1$ for each $x^{\ast} \in \ext B_{X^{\ast}}$, which means that $x_{0}$ is a spear by Corollary~\ref{Defi:spearVector}.(iv).
\end{proof}

\begin{Prop}
\label{Prop:convexityDiamondPoints}
Let $X$, $Y$, $Z$ be Banach spaces and let $G\in L(X,Y)$ be a norm-one operator. Given $T \in L(X,Z)$, the set of points $x_0\in B_X$ satisfying \eqref{eq:diamond} in Definition \ref{Defi:target} is absolutely convex and closed.
\end{Prop}

\begin{proof}
Denote by $B$ the set of points satisfying \eqref{eq:diamond} for $G$ and $T$. Fixed $x_{0} \in \overline{\aconv}{\:B}$, $\eps > 0$,and $y \in S_{Y}$, there is a finite subset $F \subset B$ with $\dist(x_{0}, \aconv{\: F}) < \eps$. By Proposition \ref{Prop:vectorsPropertyDiamond}.(iv), we have that $\bigcap_{b \in F}{\mathcal{D}_{T}(b, \ext{B_{Y^\ast}})}$ is a dense G$_ \delta$ subset of $(\ext B_{Y^{\ast}}, w^{\ast})$. Take any $y^{\ast} \in \left[\bigcap_{b \in F}{\mathcal{D}_{T}(b, \ext{B_{Y^\ast}})} \right] \cap \GS(B_{Y^{\ast}},y, \eps)$.
Then, $T b$ belongs to $\overline{T\bigl(\aconv{\GS(S_{X},G^{\ast}y^{\ast}, \eps)}\bigr)}$ for each $b \in F$, and so
\begin{equation*}
\dist{\bigl(T x_{0}, T\bigl(\aconv{\GS(S_{X},G^{\ast}y^{\ast}, \eps)}\bigr)\bigr)} \leqslant \| T\| \dist{(x_{0}, \aconv{F})} < \| T\| \eps.
\end{equation*}
A straightforward normalization gives that $x_0$ satisfies \eqref{eq:diamond} for $G$ and $T$, so $B=\overline{\aconv}{\:B}$, as desired.
\end{proof}

\section{Lush operators}
 \label{sec:Lushness}

We start with the definition of lush operator, which generalizes the concept of lush space when applied to the Identity.

\begin{Defi}\label{Defi:lushOperator}
Let $X$, $Y$ be Banach spaces and let $G\in L(X,Y)$ be a norm-one operator. We say that $G$ is \emph{lush} if $\Id_{X}$ is a target for $G$.
\end{Defi}
\index{lush operator}%

From the definition of target (or better from Remark \ref{Rema:target-norm-one}.(b)), it follows immediately the following observation.

\begin{Rema}
$G$ is lush if and only if every operator whose domain is $X$ is a target for $G$. In particular, every lush $G$ is a spear operator, that is,
$$
\|G + \T\,T\|=1 + \|T\|
$$
for every $T\in L(X,Y)$.
\end{Rema}

Let us summarize the results of the previous section when applied to lushness.

\begin{Prop}
\label{Prop:characterization-lushness}
Let $X$, $Y$ be Banach spaces, let $\mathcal{A} \subset B_{Y^{\ast}}$ with $\overline{\conv}^{\omega^{\ast}}{(\mathcal{A})} = B_{Y^{\ast}}$ and let $\mathcal{B} \subset B_{X}$ with $\overline{\aconv}\,\mathcal{B}=B_X$. Then the following assertions are equivalent for a norm-one operator $G\in L(X,Y)$:
\begin{enumerate}
\item[(i)] $G$ is lush.
\item[(ii)] For every $x_{0} \in \mathcal{B}$, $y \in S_{Y}$ and $\eps > 0$ there is $F \subset B_{X}$ such that
$$ \conv{(F)} \subset \{ x \in B_{X} \colon \| Gx + y\| > 2 - \eps \} \: \mbox{ and } \: \dist{(x_{0}, \aconv{(F)})} < \eps. $$
\item[(iii)] For every $x_0\in \mathcal{B}$, $y\in S_Y$ and $\eps>0$ there exists $y^{\ast} \in \Slice(\mathcal{A},y, \eps)$ such that
 $$
 \dist\bigl(x_{0}, \aconv(\GS(S_{X},G^{\ast}y^{\ast}, \eps))\bigr) < \eps.
 $$
\item[(iv)]
For every $x_0\in \mathcal{B}$, the set
$$
\mathcal{D}(x_{0})=\left\{y^{\ast} \in \ext B_{Y^{\ast}} \colon x_{0} \in \overline{\aconv} \bigl(\GS(S_{X},G^{\ast}y^{\ast}, \eps)\bigr) \mbox{ for every $\eps > 0$} \right\}
$$
is a dense (G$_\delta$) subset of $(\ext B_{Y^{\ast}},w^\ast)$.
\item[(v)] For every $x_0\in \mathcal{B}$, every $y\in S_Y$ and every $\eps>0$, there exists $y^\ast\in \ext(B_{Y^\ast})$ such that
$$
y\in \Slice(S_{Y},y^\ast,\eps)\quad \text{ and } \quad x_{0} \in \overline{\aconv}\bigl(\GS(S_{X},G^{\ast}y^{\ast}, \eps)\bigr).
$$
\item[(vi)] For every separable subspaces $X_{0} \subset X$ and $Y_{0} \subset Y$, we can find separable subspaces $X_{0} \subset X_{\infty} \subset X$ and $Y_{0} \subset Y_{\infty} \subset Y$ such that $G(X_{\infty}) \subset Y_{\infty}$, $\|G|_{X_{\infty}}\|=1$ and $G|_{X_{\infty}}: X_{\infty} \longrightarrow Y_{\infty}$ is lush.
\end{enumerate}
\end{Prop}

\begin{proof}
The equivalences are consequence of Theorem \ref{Theo:targetSeparablyDetermined}, together with Proposition \ref{Prop:convexityDiamondPoints} to pass from $\mathcal{B}$ to $B_{X} = \overline{\aconv}{\,\mathcal{B}}$.
\end{proof}

Next, we get from the previous sections some conditions for an operator having the aDP to be lush. The main result in this line is the next one, which follows from Theorem \ref{Theo:SCDOperators} applied to $T = \Id_{X}$.

\begin{Theo}\label{Theo:aDP+SCD=lush}
Let $X$, $Y$ be Banach spaces and let $G\in L(X,Y)$ be a norm-one operator. Suppose that $B_{X}$ is SCD. Then, $G$ has the aDP if and only if $G$ is lush.
\end{Theo}

As all the properties involved in the above result are separably determined, we have the following generalization.

\begin{Coro}
Let $X$, $Y$ be Banach spaces and let $G\in L(X,Y)$ be a norm-one operator. Suppose that $B_{X_{0}}$ is SCD for every separable subspace $X_{0} \subset X$. Then, $G$ has the aDP if and only if $G$ is lush.
\end{Coro}

\begin{proof}
Since every lush operator is a spear, it has in particular the aDP. The converse is consequence of Proposition \ref{Prop:SCDaDPseparablesubspaces} applied to $T = \Id_{X}$.
\end{proof}

The most interesting particular cases of the above results are summarized in the next corollary, which uses the examples of SCD spaces provided in Example \ref{Exam:SCD-set-spaces-operators}.

\begin{Coro}\label{Coro:aDP+RNP-CPCP..=lush}
Let $X$, $Y$ be Banach spaces and let $G\in L(X,Y)$ be a norm-one operator. Suppose that $X$ has one of the following properties: Radon-Nikod\'{y}m Property, Asplund Property, convex point of continuity property or absence of isomorphic copies of $\ell_1$. Then, $G$ has the aDP if and only if $G$ is lush.
\end{Coro}

A result of this kind for the codomain space will be given in Proposition \ref{Prop:characLushAsplund}: if $G:X\longrightarrow Y$ has the aDP and $Y$ is Asplund, then $G$ is lush.

Our next aim is to provide the following sufficient conditions for an operator to be lush which will be used in the next chapters.

\begin{Prop}
\label{Prop:lushSufficient}
Let $X$, $Y$ be Banach spaces and let $G\in L(X,Y)$ be a norm-one operator. Then, each of the following conditions ensures $G$ to be lush.
\begin{enumerate}
\item[(a)] The set $\bigl\{ y^{\ast}  \in B_{Y^{\ast}} \colon G^{\ast}y^{\ast}\in \Spear(X^\ast) \bigr\}$ is norming for $Y$.
\item[(b)] The set
$
\bigl\{ y^{\ast} \in \ext B_{Y^{\ast}} \colon G^{\ast}y^{\ast}\in \Spear(X^\ast)\bigr\}
$
is dense in $(\ext B_{Y^{\ast}},\weakstar)$.
\item[(c)] $B_{X} = \overline{\conv} \bigl\{ x \in B_{X} \colon Gx\in \Spear(Y)\bigr\}$.
\end{enumerate}
\end{Prop}

\begin{proof}
The fact that (a) implies lushness follows from Proposition \ref{Prop:characterization-lushness}.(v), as Theorem \ref{Theo:spears_of_the_dual} gives that $B_X=\overline{\aconv}\bigl(\Face(S_X,G^\ast y^\ast)\bigr)$ for every $y^\ast \in B_{Y^\ast}$ such that $G^\ast y^\ast \in \Spear(X^\ast)$. Condition (b) is a particular case of condition (a). Finally, by using Corollary \ref{Defi:spearVector}.(iv), condition (c) implies that every $y^\ast \in \ext (B_{Y^\ast})$ satisfies $\displaystyle B_X = \overline{\conv} \T\,\Face(S_X,G^\ast y^\ast)$, so $G^\ast y^\ast$ is a spear vector by (the easy part of) Theorem \ref{Theo:spears_of_the_dual}.
\end{proof}

We do not know whether the conditions (a) or (b) above are necessary for lushness in general, but they are when the domain space is separable as the following deep result shows. We will see later that they are also necessary when the codomain is an Asplund space (see Proposition \ref{Prop:characLushAsplund}).

\begin{Theo}\label{Theo-lush=casi-almost-CL}
Let $X$ be a separable Banach space and let $Y$ be a Banach space. If $G\in L(X,Y)$ is lush, then the set
$
\Omega=\bigl\{ y^{\ast} \in \ext B_{Y^{\ast}} \colon G^{\ast}y^{\ast}\in \Spear(X^\ast)\bigr\}
$
is a G$_\delta$ dense subset of $(\ext B_{Y^{\ast}},\weakstar)$. In other words, if $G$ is lush, there exists a G$_\delta$ dense subset $\Omega$ of $(\ext B_{Y^{\ast}},\weakstar)$ such that
$$
B_X = \overline{\aconv}\bigl(\Face(S_X,G^\ast y^\ast)\bigr)
$$
for every $y^\ast\in \Omega$.
\end{Theo}

\begin{proof}
This is consequence of Proposition \ref{Prop:target-W*T-dense-necessary-separable} and the characterization of spear vectors given in Corollary \ref{Defi:spearVector}. The last part is a consequence of Theorem \ref{Theo:spears_of_the_dual}.
\end{proof}

On the other hand, condition (c) of Proposition \ref{Prop:lushSufficient} is not in general necessary for lushness: consider $X=Y=c_0$ and $G=\Id$, which is lush as $c_0$ is lush, but $\Spear(Y)$ is empty as $B_{c_0}$ contains no extreme points. We will see later that condition (c) is necessary when the domain space has the Radon-Nikod\'{y}m Property (see Proposition \ref{prop:RNP-equivalence}).

We finish the section with some elementary observations analogous to the ones given for spear operators and for operators with the aDP.

\begin{Rema}\label{Rema:lushOperators-elementary}
Let $X,Y$ be Banach spaces and $G\in L(X,Y)$.
\begin{enumerate}
\item[(i)] The composition with isometric isomorphisms preserves lushness: If $X_1$, $Y_1$ are Banach spaces and $\Phi_1\in L(X_1,X)$, $\Phi_2\in L(Y,Y_1)$ are isometric isomorphisms, then, $G\in L(X,Y)$ is lush if and only if $\Phi_2 G \Phi_1\in L(X_1,Y_1)$ is lush.
\item[(ii)] If $G$ is lush and $Z$ is a subspace of $Y$ containing $G(X)$, then $G:X\longrightarrow Z$ is lush. However, lushness is not preserved by extending the codomain of the operator, as the same example of Remark \ref{Rema:SpearOperators} shows.
\item[(iii)] As an easy consequence of (i) and (ii), we have that the following statements are equivalent: (a) $X$ is lush, (b) there exist a Banach space $Z$ and an isometric isomorphism in $L(X,Z)$ or in $L(Z,X)$ which is lush, (d) there exist a Banach space $W$ and an isometric embedding $G\in L(X,W)$ which is lush.
\end{enumerate}
\end{Rema}

\chapter{Some examples in classical Banach spaces}\label{sect:classical-examples}
Our aim here is to present examples of operators which are lush, spear, or have the aDP, defined in some classical Banach spaces. One of the most intriguing examples is the Fourier transform on $L_1$, which we prove that is lush. Next, we study a number of examples of operators arriving to spaces of continuous functions. In particular, it is shown that every uniform algebra is lush-embedded into a space of bounded continuous functions. Finally, examples of operators acting from spaces of integrable functions are studied.

\section{Fourier transform} \label{ssecFourier}

Let $H$ be a locally compact Abelian group and let $\sigma$ be the Haar measure on $H$. The dual group $\Gamma$ of $H$ is the set of all continuous homomorphisms $\gamma: H \longrightarrow \T$ endowed with a topology that makes it a locally compact group (see \cite[\S 1.2]{RudinGroups} for the details). If $L_{1}(H)$ is the space of $\sigma$-integrable functions over $H$, and $C_{0}(\Gamma)$ is the space of continuous functions on $\Gamma$ which vanish at infinity, then the \emph{Fourier transform} $\mathcal{F}: L_{1}(H) \longrightarrow C_{0}(\Gamma)$ is defined as
\index{Fourier transform}%
\[
\mathcal{F}(f) : \Gamma \longrightarrow \mathbb{C}, \qquad \bigl[\mathcal{F}(f)\bigr](\gamma) = \int_{H}{f(x) \gamma(x^{-1}) \: d\sigma(x)}.
\]

\begin{Theo}\label{Theo:Fourier-transform}
Let $H$ be a locally compact Abelian group and let $\Gamma$ be its dual group. Then, the Fourier transform $\mathcal{F}: L_{1}(H) \longrightarrow C_{0}(\Gamma)$ is lush. In particular, $\mathcal{F}$ is a spear operator, that is,
$$
\|\mathcal{F}+\T\,T\|=1 + \|T\|
$$
for every $T\in L(L_{1}(H), C_{0}(\Gamma))$.
\end{Theo}

\begin{proof}
For each $\gamma \in \Gamma$, $\mathcal{F}^{\ast}(\delta_{\gamma})$ corresponds to the function $g \in L_{\infty}(H) \equiv L_{1}(H)^{\ast}$ given by $g(x) = \gamma(x^{-1})$ for every $x \in H$. Hence, $|g(x)|=1$ for every $x \in H$, and so $\mathcal{F}^{\ast}(\delta_{\gamma})$ is a spear of $L_\infty(H)\equiv L_{1}(H)^{\ast}$ by Example \ref{Exam:basicsSpears}.(e). As $\T\,\{\delta_\gamma\colon \gamma\in \Gamma\}$ is the set of extreme points of $B_{C_0(\Gamma)^\ast}$, Proposition~\ref{Prop:lushSufficient}.(b) shows that $\mathcal{F}$ is lush.
\end{proof}

\section{Operators arriving to sup-normed spaces}

Our goal here is to study various families of operators arriving to spaces of continuous functions. We start with a general result.

\begin{Prop}\label{Prop-ejem-X-into-CK}
Let $X$ be a Banach space, let $L$ be a locally compact Hausdorff topological space and let $G \in L(X, C_0(L))$ be a norm-one operator. Consider the following statements:
\begin{enumerate}
\item[(i)] The set $\bigl\{t\in L\colon G^\ast\delta_t\in \Spear(X^\ast)\bigr\}$ is dense in $L$.
\item[(ii)] $G$ is lush.
\item[(iii)] $G$ is a spear operator.
\item[(iv)] $G$ has the aDP.
\item[(v)] $\{ G^{\ast} \delta_{t} \colon t \in U\}$ is a spear set of $B_{X^\ast}$ for every open subset $U \subset L$.
\end{enumerate}
Then $(i) \Rightarrow (ii) \Rightarrow (iii) \Rightarrow (iv) \Leftrightarrow (v)$. Besides, we have the following:
\begin{enumerate}
\item[(a)] If $L$ is scattered, then all of the statements are equivalent.
\item[(b)] If $X$ is separable, then (i) $\Leftrightarrow$ (ii).
\end{enumerate}
\end{Prop}

\begin{proof}
The implications (ii) $\Rightarrow$ (iii) $\Rightarrow$ (iv) are clear. Using the fact that we can identify homeomorphically $L \equiv \{ \delta_{t} \colon t \in L \}\subset C_0(L)^\ast$ and the fact that $\ext{B_{C_0(L)^\ast}} = \T \{ \delta_{t} \colon t \in L \}$, we conclude easily (i) $\Rightarrow$ (ii) from Proposition \ref{Prop:lushSufficient}, while the equivalence (iv) $\Leftrightarrow$ (v) follows from Theorem \ref{Theo:aDPCharacterization}.

(a). $L$ is scattered if and only if $C_0(L)$ is Asplund (see \cite[Comment after Corollary 2.6]{BPBAsplund}, for instance). Notice that $G^{\ast}: C_0(L)^{\ast} \longrightarrow X^{\ast}$ is weak$^{\ast}$-weak$^{\ast}$-continuous, so we just have to prove that given an open set $U \subset L$, there exists $t \in U$ such that $G^{\ast} \delta_{t}$ is a spear. Since $G^{\ast}(U)$ is weak$^{\ast}$-fragmentable (see \cite[Theorem 11.8]{basisLinear} for the relation and defintion between Asplundness and fragmentability), for each $\eps > 0$ there exists a weak$^\ast$-open set $V$ satisfying that $V \cap G^{\ast}(U)$ has diameter less than $\eps$ . Now, since $G^{\ast}$ is weak$^\ast$-continuous, we can find an open set $W \subset L$ with $G^{\ast}(\overline{W}) \subset V \cap G^{\ast}(U)$. Because of the local compactness of $L$ we can select $W$ in such a way that $\overline{W}$ is compact. In particular $G^{\ast}(\overline{W})$ is a closed spear set with diameter less than $\eps$. It is clear that we can iterate this process to construct a decreasing sequence $\bigl(W_{n}\bigr)_{n\in \N}$ of open subsets of $U$ such that $\diam\bigl(G^{\ast}(\overline{W_{n}})\bigr)$ tends to zero. Since $G^{\ast}(\overline{W_{n}})$ is a spear set by (v), it follows then from Lemma~\ref{Lemm:convergenceSpears} that an element $t \in \bigcap_{n}{\overline{W_{n}}} \subset U$ must satisfy that $G^{\ast}\delta_{t}$ is a spear.

(b). If $X$ is separable, the result follows from Theorem \ref{Theo-lush=casi-almost-CL}.
\end{proof}

Let us mention that (a) is also a particular case of Proposition \ref{Prop:characLushAsplund} in the next chapter, using the stated above equivalence between $L$ being scattered and $C_0(L)$ being Asplund. In the more restrictive case in which $L$ has the discrete topology (which is trivially scattered), the result will be also proved, with a different approach, in Example \ref{Exam:ell_1-gamma-all-equivalent}.

The next result characterizes lush spaces. We need some notation. Let $\Omega$ be a completely regular Hausdorff topological space and denote by $C_{b}(\Omega)$ the Banach space of all scalar bounded and continuous functions on $\Omega$ endowed with the supremum norm.

\begin{Theo}\label{Theo:lushnessCharacterizationCb} Let $X$ be a Banach space. Then, $X$ is lush if and only if the canonical inclusion $J: X \longrightarrow C_{b}(\ext{B_{X^{\ast}}}) $ is lush.
\end{Theo}

Recall that for a completely regular Hausdorff topological space $\Omega$ the corresponding $C_b(\Omega)$ can be canonically seen as a $C(K)$ space by taking $K = \beta \Omega$, the Stone-Cech compactification of $\Omega$. Since $\Omega$ is a dense subset of $\beta \Omega$, we have that the set $\mathbb{T} \{ \delta_{t} \colon t \in \Omega \}$ is dense in $(\ext{C_{b}(\Omega)^\ast}, \omega^{\ast})$.

\begin{proof}
The space $\Omega := \ext{B_{X^{\ast}}}$ endowed with the weak$^\ast$-topology is completely regular. As we mentioned before, the set
$$ \mathcal{A} = \{ \delta_{x^{\ast}} \colon x^{\ast} \in \ext{B_{X^{\ast}}} \} $$
is dense in $(\ext{C_{b}(\Omega)^{\ast}}, \omega^{\ast})$ being moreover a Baire space with the induced topology (Lemma \ref{Kreinlemma}.(c)). By Proposition \ref{Prop:vectorsPropertyDiamond}, $J$ is lush if and only if for every $x_{0} \in B_{X}$ and $\eps > 0$ the set $\mathcal{D}_{J}^\eps(\mathcal{A}, x_{0})$ intersects every $\omega^{\ast}$-slice of $\mathcal{A}$. Actually, the properties of $\mathcal{A}$ given in Lemma \ref{Kreinlemma} allow us to repeat the arguments in the proof of Proposition \ref{Prop:vectorsPropertyDiamond} to deduce that $J$ is lush if and only if given an arbitrary $x_{0} \in B_{X}$ the set
$$
\{ y^{\ast} \in \mathcal{A}\colon x_{0} \in \overline{\aconv}{\, \GS(B_{X}, J^{\ast} y^{\ast}, \eps)}  \mbox{ for every $\eps > 0$} \}
$$
is dense in $(\mathcal{A}, \omega^{\ast})$. But using the natural homeomorphism between  $(\ext{B_{X^{\ast}}},w^\ast)$ and $\mathcal{A}$, we can reformulate the previous condition as
$$
\bigl\{ x^{\ast} \in \ext{B_{X^{\ast}}} \colon x_{0} \in \overline{\aconv}{\GS(B_{X}, x^{\ast}, \eps)}  \mbox{ for every $\eps > 0$} \bigr\}
$$
being dense in $(\ext{B_{X^{\ast}}}, \omega^{\ast})$ for every $x_0\in B_X$, which, by Proposition \ref{Prop:characterization-lushness}, is equivalent to say that $\Id_X$ is lush, that is, $X$ is lush.
\end{proof}

The following very general result will allow us to deduce many other interesting examples.

\begin{Theo}\label{Theo:supNormLushInclusion}
Let $\Gamma$ be a non empty set and let $ \Gamma \in \mathcal{A} \subset \mathcal{P}(\Gamma)$. Let $X \subset Y \subset \ell_{\infty}(\Gamma)$ be Banach spaces satisfying the following properties:
\begin{enumerate}
\item[(i)] For every $y \in S_{Y}$, $\eps > 0$ and $A \in \mathcal{A}$ there are $b \in \mathbb{K}$ and $U \in \mathcal{A}$  such that  $A \supset U$,
$$|b| = \sup_{t \in A}{|y(t)|} \hspace{3mm} \mbox{ and } \hspace{3mm} |y(t) - b| < \eps \hspace{3mm} \mbox{ whenever $t \in U$}. $$
\item[(ii)] For each $A \in \mathcal{A}$ there is $h \in \ell_{\infty}(\Gamma)$ such that
$$h(\Gamma)\subseteq [0,1],\quad \supp{(h)} \subset A,\quad \| h\|_{\infty} = 1 \quad \mbox{ and } \quad \dist(h, X) < \eps. $$
\end{enumerate}
Then, the inclusion $J:X \longrightarrow Y$ is lush.
\end{Theo}

\begin{proof}
Fix $x \in S_{X}$, $y \in S_{Y}$ and $0 < \eps < 1$.
Since $\Gamma \in \mathcal{A}$, we can find $A \in \mathcal{A}$ and $b \in \mathbb{T}$ such that
$$
|y(t) - b| < \frac{\eps}{9} \hspace{3mm} \mbox{ for each $t \in A$}.
$$
Again by (i) we can find $U \in \mathcal{A}$ such that $A \supset U$,  and $a \in \overline{\mathbb{D}}$ with $|a| = \sup_{t \in A}{|x(t)|}$ such that
$$
|x(t) - a| < \frac{\eps}{9} \quad \mbox{ for each $t \in U$}.
$$
By (ii), there is $h \in \ell_{\infty}(\Gamma)$, with $h(\Gamma)\subseteq [0,1]$, $\supp{(h)} \subset U$, $\| h\|_{\infty} = 1$ and $\dist(h, X) < \eps/9$. Let us fix $x_{0} \in X$ with $\| h - x_{0}\|_{\infty} < \eps/9$.

We claim that for each $\gamma \in \K$ with $|a + \gamma b| = 1$ one has that
\begin{equation*}
\| x + \gamma b x_{0} \|_{\infty} \leqslant  1 + \frac{\eps}{3}.
\end{equation*}
Indeed, the conditions $|a + \gamma b| = 1$, $|b| = 1$ and $|a| \leqslant 1$  imply that $|\gamma| \leqslant 2$. We distinguish two cases: if $t\notin U$ then $h(t)=0$ which gives $|x_0(t)|<\eps/9$, so we get that
$$
|x(t) + \gamma b x_{0}(t)| \leqslant 1 + \frac{\eps}{3}.
$$
On the other hand, if $t \in U$ then
\begin{align*}
|x(t) + \gamma b x_{0}(t)| &\leqslant |x(t) - a| + |a (1 - h(t)) + h(t) (a + \gamma b)| + |b| |\gamma| |x_{0}(t) - h(t)| \\ & < \frac{\eps}{9} + 1 + \frac{2 \eps}{9} = 1 + \frac{\eps}{3}.
\end{align*}
This finishes the proof of the claim.

Observe that $0$ belongs to the convex hull of the set $\{\gamma \in \K\colon |a + \gamma b| = 1\}$ as $|a| \leqslant 1 = |b|$. We can then find $\gamma_{1}, \gamma_{2} \in \{\gamma \in \K\colon |a + \gamma b| = 1\}$ and $ 0 \leqslant \lambda \leqslant 1$ such that $\lambda \gamma_{1} + (1 - \lambda) \gamma_{2} =0$. Take $t_{0} \in U$ with $h(t_{0}) > 1 - \eps/9$, pick $\theta_0, \theta_1, \theta_2 \in \T$ satisfying
$\theta_1\big(x(t_0) + \gamma_{1} b x_{0}(t_0)\big) \geqslant 0$, $\theta_2\big(x(t_0) + \gamma_{2} b x_{0}(t_0)\big) \geqslant 0$ and $\theta_0 y(t_0)  \geqslant 0$.
Define
$$
x_{1} = \theta_1\frac{x + \gamma_{1} b x_{0}}{1 + \eps/3} \qquad \text{and} \qquad x_{2} = \theta_2\frac{x + \gamma_{2} b x_{0}}{1 + \eps/3}.
$$
By the claim above, we have that $x_{j} \in B_{X}$. Besides, we can write
\begin{align*}
\left(1 + \frac{\eps}{3}\right) |x_{j}(t_{0})| & = |x(t_{0}) + \gamma_{j} b x_{0}(t_{0})| \\ &\geqslant |a + \gamma_{j} b| - |x(t_{0}) - a| - |\gamma_{j}| |b| |x_{0}(t_{0}) - 1|
\geqslant 1 - \frac{5 \eps}{9}.
\end{align*}
Moreover, $x_{j}(t_{0}) \geqslant 0$ and for every $\mu\in [0,1]$ we have that
$$
\mu x_1(t_0)+(1-\mu)x_2(t_0)  \geqslant \frac{1 -  5\eps/9}{1 + \eps/3} =1 - \frac{8\eps}{ 9 + 3 \eps}>1-\frac{8\eps}{9}\,.
$$
Therefore, we can estimate as follows
\begin{align*}
\left\| y + \theta_0^{-1}J\big(\mu x_{1} + (1-\mu)x_{2}\big)\right\| & \geqslant \left| \theta_0 y(t_{0}) + \mu x_{1}(t_{0}) + (1-\mu)x_{2}(t_{0})\big)\right| \\ & =  |y(t_{0})| + \mu x_1(t_0)+(1-\mu)x_2(t_0)\\ & > 1 - \frac{\eps}{9} + 1-\frac{8\eps}{9} = 2 - \eps
\end{align*}
 and, moreover,
\begin{align*}
\dist{(x, \aconv{(\{ x_{1}, x_{2}\})})} \leqslant \| x - (\lambda\theta_1^{-1} x_{1} + (1 - \lambda)\theta_2^{-1} x_{2})\| = \left\| x - \frac{x}{1 + \eps/3} \right\| \leqslant \eps. \end{align*}
This shows that $J$ is lush by (ii) of Proposition \ref{Prop:characterization-lushness} with $F = \{x_1, x_2\}$.
\end{proof}

The following notion for subspaces of $C_b(\Omega)$ was introduced in \cite[Definition 2.3]{ChGarMaMa_Pol-Daug}, generalizing the analogous concept for subspaces of $C(K)$ introduced in \cite[Definition 2.3]{NumIndexDuality} (for perfect compact spaces $K$ it was considered earlier in \cite{Ka-Po} and was studied extensively in frames of the Daugavet Property theory).

\begin{Defi}
\label{Defi:Crich}

Let $\Omega$ be a completely regular Hausdorff topological space.
A closed subspace $X \subset C_{b}(\Omega)$ is called \emph{C-rich} if for every $\eps > 0$ and every open subset $U \subset \Omega$ there exists a norm-one function $h:\Omega \longrightarrow [0,1]$ in $C_{b}(\Omega)$ such that $\supp{(h)} \subset U$ and $d(h,X) < \eps$.
\end{Defi}
\index{C-rich}%

It follows from Urysohn's Lemma that $C_b(\Omega)$ is C-rich in itself for every completely regular space $\Omega$.

The main tool in the rest of the section will be the following.

\begin{Theo}
\label{Theo:CrichLush}
Let $\Omega$ be a completely regular Hausdorff topological space.
If $X \subset C_{b}(\Omega)$ is C-rich, then the inclusion $J: X \longrightarrow C_{b}(\Omega)$ is lush. In particular, we have $\|J+\T\,T\|=1 + \|T\|$ for every $T\in L(X,C_b(\Omega))$.
\end{Theo}

\begin{proof}
We just have to check that the hypothesis of Theorem \ref{Theo:supNormLushInclusion} are satisfied for $X \subset C_{b}(\Omega) \subset \ell_{\infty}(\Omega)$ and taking as $\mathcal{A}$ the family of all open subsets of $\Omega$. Hypothesis (i) satisfied just using the continuity, while (ii) is consequence of the C-richness of $X$.
\end{proof}

\begin{Rema}
Let us observe that there are natural inclusions $J: X \longrightarrow C_{b}(\Omega)$ which are lush without $X$ being a C-rich subspace of $C_b(\Omega)$. For instance, using Theorem \ref{Theo:lushnessCharacterizationCb} we deduce that the inclusion $$
J: \ell_{1} \longrightarrow C(\mathbb{T}^{\N}), \hspace{4mm} (a_{n})_{n \in \N} \, \longmapsto \, \Big[ (z_{n})_{n \in \N} \longmapsto \sum_{n=1}^{\infty}{a_{n} z_{n}} \Big]
$$
is lush. However, $J(\ell_1)$ is not C-rich in $C(\mathbb{T}^{\N})$. Indeed, we argue by contradiction. Let $\delta > 0$, consider the open set $U = \{ z\in \mathbb{T}^{\N} \colon |z_{1} - 1| < \delta \}$, and suppose that $h \in C(\mathbb{T}^{\N})$ and $a \in \ell_{1}$ satisfy that
$$
h(\mathbb{T}^\N)\subseteq [0,1],\quad \| h\|_{\infty} = 1, \quad \supp{(h)} \subset U \quad \text{and} \quad \| J(a) - h\|_{\infty} < \delta.
$$
Taking supremum over all $z \in \mathbb{T}^{\N} \setminus U$, we deduce that
\begin{equation} \label{equa:lushNotCRichaux1}
(1 - \frac{\delta}{2}) |a_{1}| + \sum_{n \leqslant 2}{|a_{n}|} \leqslant \sup_{|z_{1}-1| \geqslant \delta}{\Real{(a_{1}z_{1})}} + \sum_{n \leqslant 2}{|a_{n}|} \leqslant \sup_{|z_{1} - 1| \geqslant \delta}{|J(a)(z)|} < \delta.
\end{equation}
While $\| h\|_{\infty} = 1$ implies that
\begin{equation}\label{equa:lushNotCRichaux2}
\| a\|_{\ell_{1}} = \sum_{n \geqslant 1}{|a_{n}|} > 1 - \delta.
\end{equation}
Taking $\delta > 0$ small enough, \eqref{equa:lushNotCRichaux2} and \eqref{equa:lushNotCRichaux1} contradict each other.

\end{Rema}

Let us present some applications of Theorem \ref{Theo:CrichLush}. First, it was shown in \cite[Proposition~1.2]{Ka-Po} that if $K$ is a perfect compact Hausdorff topological space, then every
finite-codimensional subspace of $C(K)$ is C-rich, but this is not always the case when $K$ has isolated points. Actually,  finite-codimensional subspaces of general $C(K)$ spaces were characterized in \cite[Proposition~2.5]{NumIndexDuality} in terms of the supports of the functionals defining the subspace. We recall that the \emph{support} of an element $f\in C(K)^\ast$ (represented by the regular measure $\mu_f$) is
$\supp (f):=\bigcap\left\{C\subset K\ : \ C\ \text{closed},\
|\mu_f|(K\setminus C)=0\right\}$. Then, for $f_1,\ldots,f_n\in C(K)^\ast$, the subspace $Y=\bigcap\nolimits_{i=1}^n  \ker f_i$ is C-rich in $C(K)$ if and only if $\ \bigcup\nolimits_{i=1}^n \supp(f_i)$ does
not intersect the set of isolated points of $K$.
\index{supp}%
\begin{Coro}\label{Coro:C-rich-finite-codimension}
Let $K$ be a compact Hausdorff topological space, consider functionals
$f_1,\ldots,f_n\in C(K)^\ast$ and let $Y=\bigcap\nolimits_{i=1}^n  \ker f_i$. If $\ \bigcup\nolimits_{i=1}^n \supp(f_i)$ does
not intersect the set of isolated points of $K$, then the natural inclusion $J:Y\longrightarrow C(K)$ is lush. In particular, if $K$ is perfect, then for every finite-codimensional subspace $Y$ of $C(K)$, the inclusion $J:Y\longrightarrow C(K)$ is lush.
\end{Coro}

For $C[0,1]$ we may even go to smaller subspaces, using a result of \cite{Ka-Po}: if $X$ is a subspace of $C[0,1]$ such that $C[0,1]/X$ does not contain isomorphic copies of $C[0,1]$, then $X$ is C-rich in $C[0,1]$ \cite[Proposition 1.2 and Definition 2.1]{Ka-Po}.

\begin{Coro}\label{Coro:C[01]/XdoesnotcontainsC[01]}
Let $X$ be a subspace of $C[0,1]$ such that $C[0,1]/X$ does not contain isomorphic copies of $C[0,1]$ (in particular, if $C[0,1]/X$ is reflexive). Then, the inclusion $J:X\longrightarrow C[0,1]$ is lush.
\end{Coro}

Let us now go to present the main part of this section. Recall that a \emph{uniform algebra} (on a compact Hausdorff topological space $K$) is a closed subalgebra $A \subset C(K)$ that separates the points of $K$. We refer to \cite{DalesAlgebras} for background.
\index{uniform algebra}%

\begin{Theo}\label{Theo-uniform-algebras} Let $K$ be a compact Hausdorff topological space and let $A$ be a uniform algebra on $K$. Then, there exists subset $\Omega \subset K$ such that $A \subset C_{b}(\Omega)$ (isometrically) is C-rich, and so the inclusion $J:A\longrightarrow C_b(\Omega)$ is lush. Moreover, if $A$ is unital, then $\Omega$ is just its Choquet boundary.
\end{Theo}

\begin{proof}
If $A$ is a unital uniform algebra, then consider $\Omega \subset K$ being the Choquet boundary of $A$. Given $0 < \eps < 1$ and $U \subset K$ with $U \cap \Omega \neq \emptyset$, take $0 < \eta < \eps/4$ small enough so that every $z \in \mathbb{C}$ with $|z| + (1 - \eta)|1 - z| \leqslant 1$ satisfies that $|\Imag{z}| < \eps/2$. By \cite[Lemma~2.5]{BPBalgebras}, there exists $f \in A$ and $t_{0} \in U \cap \Omega$ such that $f(t_{0}) = \| f\|_{\infty} = 1$, $|f(t)|< \eta$ for each $t \in K \setminus U$ and
\[ |f(t)| + (1 - \eta) |1 - f(t)| \leqslant 1 \mbox{ for each $t \in K$}. \]
Put $C :=K \setminus U$ and $B = \{ t \in U\colon |f(t)| \geqslant 2 \eta \}$. These are disjoint compact subsets of $K$, so there exists $\varphi: K \longrightarrow [0,1]$ continuous such that $\varphi|_{C} \equiv 0$ and $\varphi|_{B} \equiv 1$. The element $h:= |\Real{f}|\cdot \varphi : K \longrightarrow [0,1]$ belongs to $S_{C(K)}$ and satisfies $\supp{(h)} \subset U$. We just have to check that $\| h - f\|_{\infty} < \eps$. Indeed, if $t \in B$ then $|1-f(t)|<1$, so $\Real f(t)>0$ and so $|h(t) - f(t)| = |\Imag{f(t)}| \leqslant \eps$; on the other hand, if $t \in K \setminus B$ then $|h(t) - f(t)| \leqslant 4 \eta < \eps$. The restriction $h|_{\Omega}$ satisfies the definition of C-rich for the given $\varepsilon > 0$ and $U \cap \Omega \subset \Omega$.

If $A$ is not unital, then we can repeat the same argument that above but now using \cite[Lemma~2.7]{BPBalgebras} and taking $\Omega$ as the set $\Gamma_{0} \subset K$ that appears in the referenced lemma.
\end{proof}

The Choquet boundary of the disk algebra $A(\mathbb{D})$ is $\T$ (see \cite[Proposition 4.3.13]{DalesAlgebras}, for instance), so by the previous result we have the following consequence.

\begin{Coro}\label{Coro:disk-algebra}
The natural inclusion $J: A(\mathbb{D}) \longrightarrow C(\T)$ is lush.
\end{Coro}

Another family of interesting C-rich subspaces is the following. Let $H$ be an infinite compact Abelian group, let $\sigma$ be the Haar measure on $H$, let $\Gamma$ be the dual group of $H$, let $M(H)$ be the space of all regular Borel measures on $H$ and, finally, let $\mathcal{F}:M(H)\longrightarrow C_b(\Gamma)$ be the Fourier-Stieltjes transform, which is the natural extension of the classical Fourier transform of section \ref{ssecFourier} (see \cite[\S1.3]{RudinGroups}). For $\Lambda\subset \Gamma$, the space of $\Lambda$-spectral continuous functions is defined by
$$
C_\Lambda(H)=\bigl\{f\in C(H)\colon \bigl[\mathcal{F}(f)\bigr](\gamma)=0\ \forall \gamma\in \Gamma\setminus \Lambda\bigr\},
$$
and similarly it is defined the space of $\Lambda$-spectral measures $M_{\Lambda}(H)$. These spaces are known to be precisely the closed translation invariant subspaces of $C(H)$ and $M(H)$, respectively. A subset $\Lambda$ of $\Gamma$ is said to be a \emph{semi-Riesz set} \cite[p.~126]{Werner} if all elements of $M_\Lambda$ are diffuse (i.e.\ if they map singletons to $0$). Semi-Riesz sets include \emph{Riesz sets}, defined as those $\Lambda\subset \Gamma$ such that $M_\Lambda\subset L_1(\sigma)$; the chief example of a Riesz subset of the dual group $\Gamma=\mathbb{Z}$ of $H=\T$ is $\Lambda=\N$. We refer to \cite[\S IV.4]{HWW} for background.
It is shown in \cite[Theorem 4.13]{Luecking} that $\Lambda\setminus \Gamma^{-1}$ is a semi-Riesz set if and only if $C_\Lambda(H)$ is C-rich in $C(H)$. Therefore, we have the following consequence of Theorem \ref{Theo:CrichLush}.
\index{semi-Riesz set}%
\index{Riesz set}%

\begin{Coro}\label{Coro:semi-Riesz-sets}
Let $H$ be an infinite compact Abelian group and let $\Lambda$ be a subset of the dual group $\Gamma$. If $\Gamma\setminus \Lambda^{-1}$ is a semi-Riesz set, then the inclusion $J:C_\Lambda(H)\longrightarrow C(H)$ is lush.
\end{Coro}

\begin{Rema}
It is proved in \cite[Theorem 3.7]{Werner} that if $\Gamma\setminus \Lambda^{-1}$ is a semi-Riesz set, then $C_\Lambda(H)$ is \emph{nicely embedded} into $C(H)$, that is, the isometric embedding $J:C_\Lambda(H)\longrightarrow C(H)$ satisfies that for every $t\in H$, $\|J^\ast \delta_t\|=1$ and the linear span of $J^\ast \delta_t$ is an $L$-summand in $X^\ast$ (this is actually a straightforward consequence of the definition of semi-Riesz set).
\index{nicely embedded}%
Then, it follows immediately from Example \ref{Exam:basicsSpears}.(g) that $J^\ast\delta_t\in \Spear(C_\Lambda(H)^\ast)$ for every $t\in H$, so $J$ is lush by Proposition \ref{Prop:lushSufficient}.(a). This is thus an alternative elementary proof of Corollary \ref{Coro:semi-Riesz-sets} which does not need the more complicated \cite[Theorem 4.13]{Luecking}.

Let us also comment that it was proved in \cite[Proof of Theorem 3.3]{Werner} that unital function algebras are nicely embedded into $C_b(\Omega)$, where $\Omega$ is the Choquet boundary of the algebra, so Theorem \ref{Theo-uniform-algebras}, in the unital case, and Corollary \ref{Coro:disk-algebra}, can be also proved by using the argument above.
\end{Rema}

We next provide with more applications of Theorem \ref{Theo:supNormLushInclusion}. The following definition appears in \cite[Definition 3.2]{Richsubspaces} for vector-valued spaces of continuous functions.

\begin{Defi}\label{Defi:CKsuperspace}
Let $K$ be a compact and Hausdorff space. We say that a closed subspace $X \subset \ell_{\infty}(K)$ is a \emph{$C(K)$-superspace} if it contains $C(K)$ and for each $x \in X$, every open subset $U \subset K$ and each $\eps > 0$, there are an open subset $V \subset U$ and an element $\theta \in \mathbb{K}$ such that
\index{cksuperspace@$C(K)$-superspace}%
$$|\theta| = \sup_{t \in U}{|x(t)|} \hspace{3mm} \mbox{ and } \hspace{3mm} | x(t) - \theta | < \eps \hspace{3mm} \mbox{ for each $t \in V$.}  $$
\end{Defi}

The result for $C(K)$-superspaces is the following.

\begin{Coro}
If $X$ is a $C(K)$-superspace, then the inclusion $C(K) \longrightarrow X$ is lush.
\end{Coro}

\begin{proof}
Using Theorem \ref{Theo:supNormLushInclusion} for the inclusions $C(K) \subset X \subset \ell_{\infty}(K)$ with $\mathcal{A}$ being the set of open subsets of $K$, we have that (ii) is satisfied by Urysohn's Lemma, while (i) is just the definition of $C(K)$-superspace.
\end{proof}

An interesting application is given by the next example.

\begin{Exam} Let $D[0,1]$ be the  space of bounded functions on $[0,1]$ which are right-continuous, have left limits everywhere and are continuous at $t=1$. It is shown in \cite[Proposition 3.3]{Richsubspaces} that $D[0,1]$ is a $C[0,1]$-superspace (this is because $D[0,1]$ is the closure in $\ell_{\infty}[0,1]$ of the span of the step functions $\mathbbm{1}_{[a,b)}, \, 0 \leqslant a \leqslant b < 1$ and $\mathbbm{1}_{[a,1]}, \, 0 \leqslant a \leqslant 1$). Therefore, the inclusion $J:C[0,1]\longrightarrow D[0,1]$ is lush.
\end{Exam}

\section{Operators acting from spaces of integrable functions}\label{subsect:from_L_1}
Our aim here is to describe operators from $L_1(\mu)$ spaces which have the aDP. For commodity, we only deal with probability spaces, but this is not a mayor restriction as $L_1$-spaces associated to $\sigma$-finite measures are (up to an isometric isomorphism) $L_1$ spaces associated to probability measures (see \cite[Proposition~1.6.1]{CembranosMendoza}, for instance). We introduce some notation. Let $(\Omega, \Sigma, \mu)$ be a probability space and let $Y$ be a Banach space. We write $\Sigma^{+} := \bigl\{ B \in \Sigma \colon \mu(B) > 0 \bigr\}$ and for $A\in \Sigma^+$ we consider
$$
\Sigma_{A} := \bigl\{ B \in \Sigma \colon B \subset A \bigr\}, \qquad  \Sigma_{A}^{+}:= \Sigma_{A} \cap \Sigma^{+},\quad \text{and} \quad
\Gamma_{A} := \left\{ \frac{\mathbbm{1}_{B}}{\mu(B)} \colon B \in \Sigma_{A}^{+} \right\}.
$$
Recall that an operator $T \in L(L_{1}(\mu), Y)$ is \emph{(Riesz) representable} if there exists $g \in L_{\infty}(\mu, Y)$ (i.e.\ a strongly measurable and essentially bounded function) such that
\index{00sigmaA@$\Sigma_A$}%
\index{00sigmaz+@$\Sigma^+$}%
\index{00sigmaz+a@$\Sigma^+_A$}%
\index{00gammaA@$\Gamma_A$}%
\index{Riesz representable}%
$$
T(f) = \int_{\Omega}{f g \: d \mu} \qquad \bigl(f \in L_{1}(\mu)\bigr).
$$
Rank-one operators are representable by the classical Riesz representation theorem assuring that $L_1(\mu)^\ast\equiv L_\infty(\mu)$, and then so are all finite-rank operators. Actually, compact operators \cite[p.~68, Theorem 2]{DiestelUhl} and even weakly compact operators \cite[p.~65, Theorem 12]{DiestelUhl} are representable, but the converse result is not true \cite[p.~79, Example 22]{DiestelUhl}. Finally, let us say that the set of representable operators can be isometrically identified with $L_{\infty}(\mu, Y)$ \cite[p.~62, Lemma 4]{DiestelUhl}. We refer the reader to chapter III of \cite{DiestelUhl} for more information and background on representable operators.

Here is the characterization of aDP operators which is the main result of this section.

\begin{Theo}
\label{Theo:L1->XaDP}
Let $(\Omega, \Sigma, \mu)$ be a probability space, let $Y$ be a Banach space and let $G \in L(L_{1}(\mu), Y)$ be a norm-one operator. The following assertions are equivalent:
\begin{enumerate}
\item[(i)] $G$ has the aDP.
\item[(ii)] $G(\Gamma_{A})$ is a spear set for every $A \in \Sigma^{+}$.
\item[(iii)] $\|G + \T\,T \| = 1 + \| T\|$ for every $T \in L(L_{1}(\mu), Y)$ representable.
\end{enumerate}
\end{Theo}

Let us recall the following exhaustion argument that we briefly prove here.

\begin{Obse}\label{Obse:exhaustion}
{\slshape Let $(\Omega, \Sigma, \mu)$ be a finite measure space. If for each $A \in \Sigma^{+}$ there is $B \in \Sigma_{A}^{+}$ satisfying a certain property $(P)$, then we can find a countable family $\AAA \subset \Sigma^{+}$ of disjoint sets such that every $A \in \AAA$ satisfies property (P) and $\Omega \setminus \bigcup \AAA$ is $\mu$-null.\ }

Indeed, this follows from a simple argument: using Zorn's lemma we can take a maximal family $\AAA$ of disjoint sets in $\Sigma^{+}$ satisfying property (P), which must be countable as $\mu$ is finite. To see the last condition, notice that if $A:=\Omega \setminus \bigcup{\mathcal{A}}$ had positive measure, then we could use the hypothesis to find a subset $B \in \Sigma_{A}^{+}$ satisfying (P), and hence $\mathcal{A} \cup \{ A\}$ would contradict the maximality of $\mathcal{A}$.
\end{Obse}

We need a preliminary result.

\begin{Lemm}
\label{Lemm:normOperatorL1}
Let $(\Omega, \Sigma, \mu)$ be a probability space and let $Y$ be a Banach space. For every $T\in L(L_{1}(\mu),Y)$ one has that
\begin{equation}
\label{EQUA:normOperatorL1}
\| T\| = \sup_{A \in \Sigma_{+}}{\bigl\| T(\mathbbm{1}_{A})/\mu(A) \bigr\|} = \sup_{A \in \Sigma_{+}}{\inf_{B \in \Sigma_{A}^{+}}{\bigl\| T(\mathbbm{1}_{B})/\mu(B)\bigr\|}}.
\end{equation}
\end{Lemm}

\begin{proof}
The inequalities $\geqslant$ are clear in both cases, so we just have to see that
$$
\alpha := \sup_{A \in \Sigma_{+}}{\inf_{B \in \Sigma_{A}^{+}}{\| T(\mathbbm{1}_{B})/\mu(B)\|}}
$$
is greater than or equal to $\|T\|$. Let $h= \sum_{A \in \pi}{c_{A} \mathbbm{1}_{A}}$ be a simple function, where $\pi$ is a finite partition of $\Omega$ into elements of $\Sigma^+$, so $\| h\|_{1} = \sum_{A \in \pi}{|c_{A}| \mu(A)}$. Given $\eps > 0$, we have that for each $A \in \Sigma^{+}$ there is $B \in \Sigma_{A}^{+}$ such that $\| T(\mathbbm{1}_{B})/\mu(B)\| < \alpha + \eps$. Using Observation \ref{Obse:exhaustion} in each set $A \in \pi$, we can find a countable partition $\AAA \subset \Sigma^{+}$ of $\Omega$ such that every $B \in \AAA$ of positive measure is contained in some element of $\pi$ and satisfies $\| T(\mathbbm{1}_{B})/\mu(B)\| < \alpha + \eps$. If we write $c_{B} = c_{A}$ whenever $B \subset A$, then
\[
\| T(h)\| \leqslant \sum_{A \in \pi}{|c_{A}| \| T(\mathbbm{1}_{A})\|} \leqslant \sum_{B \in \mathcal{A}}{\mu(B)|c_{B}| \left\| T(\mathbbm{1}_{B})/\mu(B)\right\|} \leqslant (\alpha+\eps) \| h\|_{1}.
\]
Since $\eps > 0$ was arbitrary, we conclude that $\|T(h)\|\leqslant \alpha \|h\|_1$. As $h$ runs on all simple functions, it follows that $\|T\|\leqslant \alpha$.
\end{proof}

\begin{proof}[Proof of Theorem~\ref{Theo:L1->XaDP}]
(i) $\Rightarrow$ (ii): Fix $A \in \Sigma^{+}$, $y \in S_Y$, and $\eps\in(0,1)$. Consider the rank-one operator $T:L_{1}(\mu) \longrightarrow Y$ given by
$$
T(f) = y \int_{A}{f \: d \mu} \quad \bigl(f\in L_1(\mu)\bigr).
$$
Then, $\|G + \T\, T \| = 1 + \| T\| =2$. By Lemma \ref{Lemm:normOperatorL1} we can find $B \in \Sigma^{+}$ such that
\begin{equation*}
\left\|\frac{G(\mathbbm{1}_{B})}{\mu(B)} + \T \frac{T(\mathbbm{1}_{B})}{\mu(B)} \right\| \geqslant 2 - \eps.
\end{equation*}
It follows that $\left\|\frac{T(\mathbbm{1}_{B})}{\mu(B)}\right\| \geqslant 1 - \eps$, consequently
$\mu(A \cap B) \geqslant (1- \eps)\mu(B)$, so the set $\widehat{B}:= A \cap B \in \Sigma_{A}^{+}$ satisfies that
\begin{equation*}
\left\| \frac{\mathbbm{1}_{\widehat{B}}}{\mu(\widehat{B})} - \frac{\mathbbm{1}_{B}}{\mu(B)} \right\|_{1} = 2 - 2 \frac{\mu(\widehat{B})}{\mu(B)} \leqslant 2 \eps.
\end{equation*}
The combination of the above two inequalities leads to
\[
\left\|G\left(\frac{\mathbbm{1}_{\widehat{B}}}{\mu(\widehat{B})}\right) +  \T\,y\right\| = \left\|  \frac{G(\mathbbm{1}_{\widehat{B}})}{\mu(\widehat{B})} +  \T \frac{T(\mathbbm{1}_{\widehat{B}})}{\mu(\widehat{B})} \right\| \geqslant 2 - \eps - 2 \eps.
\]
(ii) $\Rightarrow$ (iii): Let $T \in L(L_{1}(\mu), Y)$ and $\eps > 0$. By Lemma \ref{Lemm:normOperatorL1}, we can find $A \in \Sigma^{+}$ satisfying $\inf_{B \in \Sigma_{A}^{+}}{\left\| T(\mathbbm{1}_{B})/ \mu(B)\right\|} > \| T\| - \eps$. If $T$ is representable, then there exists $B \in \Sigma_{A}^{+}$ such that $\diam\bigl(T(\Gamma_{B})\bigr) < \eps$ (see \cite[p. 62, Lemma 4 and p. 135, Lemma 6]{DiestelUhl}), so taking any $y \in T(\Gamma_{B})$ and using that $G(\Gamma_{B})$ is a spear set we obtain
\[
\| G +\T\, T\| \geqslant \| G(\Gamma_{B}) + \T \, y \| - \|T(\Gamma_{B}) - y \| \geqslant 1 + \| y\| - \eps \geqslant 1 + \| T\| - 2 \eps.
\]
(iii) $\Rightarrow$ (i) is obvious as rank-one operators are representable.
\end{proof}

\begin{Rema}
A direct way to prove (i) $\Rightarrow$ (iii) in Theorem~\ref{Theo:L1->XaDP} is the following: every representable operator $T:L_{1}(\mu) \longrightarrow X$ factorizes through $\ell_{1}$, i.e.\ there are operators $S: L_{1}(\mu) \longrightarrow \ell_{1}$ and $R: \ell_{1} \longrightarrow X$ such that $T = R \circ S$. But then, $S$ is an SCD operator satisfying $\| Tf\| \leqslant \| Sf\|$ for each $f \in L_{1}(\mu)$, and so Corollary \ref{Coro:HSCD-majorized-Operators} implies that $T$ is a target for $G$. However, item (ii) in Theorem~\ref{Theo:L1->XaDP} gives an intrinsic characterization of aDP operators acting from an $L_1$ space which has its own interest.
\end{Rema}

As an obvious consequence of Theorem~\ref{Theo:L1->XaDP}, if for a Banach space $Y$ all bounded linear operators from $L_1(\mu)$ to $Y$ are representable, then the aDP is equivalent to be spear for every $G\in L(L_1(\mu),Y)$. This is the case when $Y$ has the Radon-Nikod\'{y}m Property with respect to $\mu$ \cite[p.~63, Theorem 5]{DiestelUhl}. Therefore, the following result follows.

\begin{Coro}\label{Coro:fromL1RNP-adp=>spear}
Let $(\Omega, \Sigma, \mu)$ be a probability space, let $Y$ be a Banach space which has the Radon-Nikod\'{y}m Property with respect to $\mu$, and let $G \in L(L_{1}(\mu), Y)$ be a norm-one operator. Then $G$ has the aDP if and only if $G$ is a spear operator.
\end{Coro}

Let us observe that if $Y$ in the above corollary has the Radon-Nikod\'{y}m Property, then the result also follows from Corollary \ref{Coro:aDP+SCD=target}. On the other hand, it is easy to see that every Banach space $Y$ has the Radon-Nikod\'{y}m Property with respect to the counting measure on $\N$ (use \cite[Theorem 8, p.~66]{DiestelUhl}, for instance) and so $G\in L(\ell_1,Y)$ has the aDP if and only if $G$ is a spear operator. But, in this case, we will show that actually $G$ has the aDP if and only if $G$ is lush (see Example \ref{Exam:ell_1-gamma-all-equivalent}).

We next present an example showing that spearness and lushness are not equivalent for an operator having an $L_1(\mu)$ space as domain,

\begin{Exam}\label{Exam:fromL1-spear-no-lush}
There is a Banach space $Y  \subset L_{1}[0,2]$ such that the quotient operator $\pi: L_{1}[0,2] \longrightarrow L_{1}[0,2]/Y $ is a spear operator but it is not lush.

Indeed, let $Y$ be the space introduced in \cite[Theorem 4.1]{LushNumOneDual}: this is defined as $Y=^{\perp}Z  \subset L_1[0,2]$, where $Z  \subset L_{\infty}[0,2]$ is a $w^\ast$-closed subspace with no modulus-one (a.e.) functions. If $\pi$ were lush, then its dual operator $\pi^{\ast}: Z  \longrightarrow L_{\infty}[0,2]$ would send ``many'' extreme points to spear vectors by Theorem \ref{Theo-lush=casi-almost-CL}. However, $\pi^\ast$ is an isometric embedding and $Z$ does not contain any modulus-one function, while the spears of $L_{\infty}[0,2]$ are characterized that way, see Example \ref{Exam:basicsSpears}.(e). This is a contradiction. On the other hand, $Z$ is constructed to be C-rich in $L_\infty[0,2]$ (viewed as a $C(K)$ space), and so $\pi^{\ast}$ is lush by Theorem \ref{Theo:CrichLush}. This, in particular, implies that $\pi$ is a spear operator.
\end{Exam}

We finish this section showing that a representable operator $G:L_1(\mu)\longrightarrow Y$ has the aDP if and only if it is represented by a spear vector of $L_\infty(\mu,Y)$. As a consequence, we will describe the spear vectors of $L_\infty(\mu,Y)$ as those functions which take spear values almost everywhere, extending Example \ref{Exam:basicsSpears}.(e) to the vector-valued case.

\begin{Coro}
\label{Coro:spearLinfinity}
Let $G\in L(L_{1}(\mu),Y)$ be a norm-one operator which is representable by $g \in L_{\infty}(\mu, Y)$. Then, the following are equivalent:
\begin{enumerate}
\item[(i)] $G$ has the aDP.
\item[(ii)] $g(t)\in \Spear(Y)$ for a.e.\ $t$.
\item[(iii)] $g\in \Spear\bigl(L_{\infty}(\mu, Y)\bigr)$.
\end{enumerate}
\end{Coro}

\begin{proof}
The implication (ii) $\Rightarrow$ (iii) is an easy adaptation of the scalar case proved in Example \ref{Exam:basicsSpears}.(e). Indeed, for $f\in L_\infty(\mu,Y)$ and $\eps>0$, there is $A\in \Sigma^+$ such that $\|f(t)\|\geqslant \|f\|_\infty-\eps$ for every $t\in A$. By (ii), there is $A'\in \Sigma^+$ with $A'\subset A$ such that $g(t)\in \Spear(Y)$ for every $t\in A'$. Now, fixed an $\eps$-net $\mathbb{T}_\eps$ of $\T$ we can find an element $\theta_1 \in \mathbb{T}_\eps$ and a subset $A''$ of $A'$ with positive measure such that $\|g(t)+\theta_1 f(t)\|\geqslant 1+\|f(t)\|(1-\eps)$ for every $t\in A''$. Therefore, we can write
    \begin{align*}
    \|g + \T\,f\|_\infty&\geqslant \inf_{t\in A''} \|g(t)+ \T\,f(t)\|\\
    &\geqslant \inf_{t\in A''} 1 + \|f(t)\|(1-\eps) \geqslant 1 + (\|f\|_\infty - \eps)(1-\eps)
    \end{align*}
    and the arbitrarines of $\eps$ gives the result.
(iii) $\Rightarrow$ (i) is consequence of the fact that every rank-one operator is representable, so for every rank-one $T\in L(L_1(\mu),Y)$ there is $f\in L_\infty(\mu,Y)$ which represents $T$ and so
$$
\|G + \T\, T\|=\|g + \T\,f\|_\infty = 1 + \|f\|_\infty = 1 + \|T\|.
$$
Let us prove (i) $\Rightarrow$(ii). Fix $\eps > 0$. Since $g$ is strongly measurable, given any $A \in \Sigma^{+}$ there exists $B \in \Sigma_{A}^{+}$ such that $\diam{(g(B))} < \eps$ (see \cite[Proposition 2.2]{multiFunctions}, for instance). By Observation \ref{Obse:exhaustion} we can take a countable family $\mathcal{A} \subset \Sigma^{+}$ of disjoint sets with the property that $\diam{(g(A))} < \eps$ for each $A \in \mathcal{A}$ and $N_{\eps}:= \Omega \setminus \bigcup{\mathcal{A}}$ is $\mu$-null. Given $t \in \Omega \setminus N_{\eps}$, it must belong to some $A \in \mathcal{A}$, and since $G(\Gamma_{A}) \subset \overline{\conv}{\:g(A)}$ (see \cite[p. 48, Corollary 8]{DiestelUhl}), we deduce that $\diam{(G(\Gamma_{A}))} < \eps$. Using that $G(\Gamma_{A})$ is a spear set by Theorem \ref{Theo:L1->XaDP}, it follows that for every $x \in X$, $$\|g(t) + \T\, x\| \geqslant \| G(\Gamma_{A})+ \mathbb{T} \, x\| - \| G(\Gamma_{A}) - g(t) \| \geqslant  \| G(\Gamma_{A})+ \mathbb{T} \, x\| - \varepsilon = 1 + \|x\| - \eps.$$ Finally, if we take now a decreasing sequence $(\eps_{n})_{n \in \N}$ of positive numbers converging to zero and consider the correspondent $N_{\eps_{n}}$ for each $n \in \N$, then every $t \in \Omega \setminus \bigcup_{n \in \N}{N_{\eps_{n}}}$ satisfies $\|g(t) + \T\, x\| = 1 + \| x\|$ for each $x \in X$, i.e.\ $g(t)$ is a spear.
\end{proof}

\chapter{Further results} \label{sec.examples-appl}
Our goal here is to complement the previous chapter with some interesting results. We characterize lush operators when the domain space has the Radon-Nikod\'{y}m Property or the codomain space is Asplund, and we get better results when the domain or the codomain is finite-dimensional or when the operator has rank one. Further, we study the behaviour of lushness, spearness and the aDP with respect to the operation of taking adjoint operators. Finally, we collect some stability results.

\section{Radon-Nikod\'{y}m Property in the domain or Asplund codomain}

We first provide a result about the relationship of an operator with the aDP and spear vectors of the unit ball of the range space and the spear vectors of the dual ball of the domain space.

\begin{Prop}
\label{Prop:aDPDentingPoints}
Let $X$, $Y$ be Banach spaces and let $G\in L(X,Y)$ be an operator with the aDP, then
\begin{enumerate}
\item[(a)] $Gx\in \Spear(Y)$ for every denting point $x$ of $B_{X}$.
\item[(b)] $G^{\ast} y^{\ast}\in \Spear(X^\ast)$ for every weak$^\ast$-denting point $y^{\ast}$ of $B_{Y^{\ast}}$.
\end{enumerate}
\end{Prop}

\begin{proof}
We only illustrate the proof of (a), since the other one is completely analogous. If $x$ is denting, then we can find a decreasing sequence $(S_{n})_{n\in\N}$ of slices of $B_{X}$ containing $x$ and such that $\diam{S_{n}}$ tends to zero. Since $G$ has the aDP, Theorem \ref{Theo:aDPCharacterization}.(iii) gives that $\bigl(G(S_{n})\bigr)_{n \in \N}$ is a decreasing sequence of spear sets whose diameters tend to zero, so $Gx \in \bigcap_{n}{G(S_{n})}$ is a spear vector by Lemma~\ref{Lemm:convergenceSpears}.
\end{proof}

We now characterize spear operators acting from a Banach space with the Radon-Nikod\'{y}m Property.

\begin{Prop}\label{prop:RNP-equivalence}
Let $X$ be a Banach space with the Radon-Nikod\'{y}m Property, let $Y$ be a Banach space and let $G\in L(X,Y)$ be a norm-one operator. Then, the following assertions are equivalent:
\begin{enumerate}
\item[(i)] $G$ is lush.
\item[(ii)] $G$ is a spear operator.
\item[(iii)] $G$ has the aDP.
\item[(iv)] $|y^\ast(Gx)|=1$ for every $y^\ast\in \ext(B_{Y^\ast})$ and every denting point $x$ of $B_X$.
\item[(v)] $B_{X} = \overline{\conv} \bigl\{ x \in B_{X} \colon Gx\in \Spear(Y)\bigr\}$ or, equivalently,
\begin{align*}
     B_X &=
    \overline{\conv}\bigl\{x\in B_X\colon |y^\ast (Gx)|=1\ \forall y^\ast \in \ext(B_{Y^\ast}) \bigr\} \\ &=
    \overline{\conv}\left( \bigcap_{y^\ast \in \ext (B_{Y^\ast})} \T\,\Face(S_X,G^\ast y^\ast)\right).
    \end{align*}
\end{enumerate}
\end{Prop}

\begin{proof}
(i) $\Rightarrow$ (ii) $\Rightarrow$ (iii) are clear. (iii) $\Rightarrow$ (iv) follows from Proposition \ref{Prop:aDPDentingPoints} and Corollary \ref{Defi:spearVector}.(iv). (iv) $\Rightarrow$ (v) is consequence of the fact that $B_X$ is the closed convex hull of its denting points since $X$ has the Radon-Nikod\'{y}m Property (see \cite[\S 2]{Bourgin} for instance), and the equivalent reformulation is a consequence of Theorem \ref{Theo:spears_of_the_dual}. Finally, (v) $\Rightarrow$ (i) follows from Proposition \ref{Prop:lushSufficient}.(c).
\end{proof}

For Asplund spaces, we have the following characterization.

\begin{Prop}\label{Prop:characLushAsplund}
Let $X$ be a Banach space, let $Y$ be an Asplund space and let $G\in L(X,Y)$ be a norm-one operator. Then, the following assertions are equivalent:
\begin{enumerate}
\item[(i)] $G$ is lush.
\item[(ii)] $G$ is a spear operator.
\item[(iii)] $G$ has the aDP
\item[(iv)] $|x^{\ast\ast}(G^\ast y^\ast)|=1$ for every $x^{\ast\ast}\in \ext (B_{X^{\ast\ast}})$ and every weak$^\ast$-denting point $y^\ast$ of $B_{Y^\ast}$.
\item[(v)] The set
$
\bigl\{ y^{\ast} \in \ext B_{Y^{\ast}} \colon G^{\ast}y^{\ast}\in \Spear(X^\ast)\bigr\}
$
is dense in $(\ext B_{Y^{\ast}},\weakstar)$ or, equivalently,
there is a dense subset $K$ of $(\ext B_{Y^{\ast}},\weakstar)$ such that
$$
B_X = \overline{\aconv}\bigl(\Face(S_X,G^\ast y^\ast)\bigr)
$$
for every $y^\ast\in K$.
\item[(vi)] $B_{Y^\ast} = \overline{\conv}^{\weakstar}{\bigl\{ y^{\ast}  \in B_{Y^{\ast}} \colon G^{\ast}y^{\ast}\in \Spear(X^\ast) \bigr\}}$.
\end{enumerate}
\end{Prop}

\begin{proof}
(i) $\Rightarrow$ (ii) $\Rightarrow$ (iii) are clear. (iii) $\Rightarrow$ (iv) follows from Proposition \ref{Prop:aDPDentingPoints} and Corollary \ref{Defi:spearVector}.(iv). (iv) $\Rightarrow$ (v): the set  contains all weak$^\ast$-denting points of $B_{Y^\ast}$ by Proposition \ref{Prop:aDPDentingPoints}, so it is weak$^\ast$ dense since for Asplund spaces, weak$^\ast$-denting points are weak$^\ast$ dense in the set of extreme points of the dual ball (see \cite[\S 2]{Bourgin} for instance). The equivalent reformulation is consequence of Theorem \ref{Theo:spears_of_the_dual}. Finally, (v) $\Rightarrow$ (vi) is  clear and (vi) $\Rightarrow$ (i) follows from Proposition \ref{Prop:lushSufficient}.(b).
\end{proof}

We do not know whether the above result extends to the case when $Y$ is SCD. What is easily true, using Theorem \ref{Theo:SCDOperators}, is that aDP and spearness are equivalent in this case.

\begin{Rema}\label{Rema:aDP-spear-equivalent-Y-SCD}
Let $X$ be a Banach space, let $Y$ be an SCD Banach space, and let $G\in L(X,Y)$ be a norm-one operator. Then, $G$ has the aDP if and only if $G$ is a spear operator.
\end{Rema}

As a consequence of the results above, we may improve Proposition \ref{Prop:examplesell1gamma-c0gamma}.

\begin{Exam}\label{Exam:ell_1-gamma-all-equivalent}
Let $\Gamma$ be an arbitrary set, let $X$, $Y$ be Banach spaces and let $\{e_\gamma\}_{\gamma \in \Gamma}$ be the canonical basis of $\ell_{1}(\Gamma)$ (as defined in Example \ref{Exam:basicsSpears}.(a)).
\begin{enumerate}
\item[(a)] For $G \in L(\ell_{1}(\Gamma),Y)$ the following are equivalent: $G$ is lush, $G$ is a spear operator, $G$ has the aDP, $G(e_{\gamma})\in \Spear(Y)$ for every $\gamma\in \Gamma$, $|y^\ast(G(e_\gamma))|=1$ for every $y^\ast\in \ext(B_{Y^\ast})$ and every $\gamma\in \Gamma$.
\item[(b)] For $G \in L(X, c_{0}(\Gamma))$ the following are equivalent: $G$ is lush, $G$ is a spear operator, $G$ has the aDP, $G^{\ast}(e_\gamma)\in \Spear(X^\ast)$ for every $\gamma\in \Gamma$, $B_X=\overline{\aconv}\bigl(\Face(S_X,G^\ast e_\gamma)\bigr)$ for every $\gamma\in \Gamma$.
\end{enumerate}
\end{Exam}

Part of assertion (a) above also follows from Corollary \ref{Coro:fromL1RNP-adp=>spear}; the whole assertion (b) also follows from Proposition \ref{Prop-ejem-X-into-CK}.

\section{Finite-dimensional domain or codomain} Our goal now is to discuss the situation about spear operators when the domain or the codomain is finite-dimensional. We start with the case in which the domain is finite-dimensional, where the result is just an improvement of Proposition \ref{prop:RNP-equivalence}. To get it, we only have to recall that for finite-dimensional spaces, the concepts of denting point and extreme point coincide thanks to the compactness of the unit ball and Choquet's Lemma (Lemma \ref{Kreinlemma}.(a)).

\begin{Prop}\label{Prop:characLushfinitedimDOMAIN}
Let $X$ be a finite-dimensional space, let $Y$ be a Banach space and let $G\in L(X,Y)$ be a norm-one operator. Then, the following are equivalent:
\begin{enumerate}
\item[(i)] $G$ is lush.
\item[(ii)] $G$ is a spear operator.
\item[(iii)] $G$ has the aDP.
\item[(iv)] $|y^\ast(Gx)|=1$ for every $y^\ast\in \ext(B_{Y^\ast})$ and every $x\in \ext(B_X)$.
\item[(v)] $Gx\in \Spear(Y)$ for every $x\in \ext(B_X)$.
\item[(vi)] $\displaystyle B_X = \conv\left( \bigcap_{y^\ast \in \ext (B_{Y^\ast})} \T\,\Face(S_X,G^\ast y^\ast)\right)$.
\end{enumerate}
\end{Prop}

The next example shows that even in the finite-dimensional case, bijective lush operators can be very far away from being isometries and that their domain and codomain are not necessarily spaces with numerical index one.

\begin{Exam}\label{ex:ADPoperator-index-less-1-domain-codomain}
There exists a bijective lush operator such that neither its domain nor its codomain has the aDP.

Indeed, let $X_1$ be the real four-dimensional space whose unit ball is given by
$$
B_{X_1}=\conv\bigl\{(\eps_1,\eps_2,\eps_3,\eps_4) \colon \eps_k\in\{-1,1\} \text{ and } \eps_i\neq \eps_j \text{ for some } i,j \bigr\}.
$$
Let $Y_1$ be the real space $\ell_\infty^4$, let $X_2=Y_1^\ast=\ell_1^4$ and, finally, let $Y_2=X_1^\ast$. Consider the operator $G_1\in L(X_1,Y_1)$ given by $G_1(x_1)=x_1$ for every $x_1\in X_1$ and consider $G_2=G_1^\ast\in L(X_2,Y_2)$. Finally, calling $X=X_1\oplus_\infty X_2$ and $Y=Y_1\oplus_\infty Y_2$, the operator we are looking for is $G\in L(X,Y)$ given by
$G(x_1,x_2)=(G_1x_1,G_2x_2) $ for every $(x_1, x_2)\in X$.

We start showing that $G$ is lush. To this end, by Proposition~\ref{Prop:characLushfinitedimDOMAIN}, all we have to do is to check that $G$ carries extreme points of $B_X$ to spear vectors of $Y$. By Example \ref{Exam:basicsSpears}.(i), this is equivalent to show that both $G_1$ and $G_2$ carry extreme points to spear elements. This is evident for $G_1$ and it is also straightforward to show for $G_2$ (alternatively, the first assertion gives that $G_1$ is lush by Proposition~\ref{Prop:characLushfinitedimDOMAIN}, so $G_2=G_1^\ast$ is also lush by Corollary \ref{Coro:duality-fin-dim} in the next section, so $G_2$ carries extreme points of $B_{X_2}$ to spear elements in $Y_2$ by using again Proposition~\ref{Prop:characLushfinitedimDOMAIN}).

Finally, let us show that $X$ does not have the aDP (i.e.\ that $\Id_X$ does not have the aDP). By Proposition~\ref{Prop:characLushfinitedimDOMAIN}, it is enough to find an extreme point of $B_X$ which is not a spear vector of $X$. By Example \ref{Exam:basicsSpears}.(i), it is enough to find an extreme point of $B_{X_1}$ which is not a spear vector of $X_1$. Let us show that this happens for $x_1=(1,1,-1,-1)\in X_1$. On the one hand, $x_1$ is clearly an extreme point of $B_{X_1}$ by construction. On the other hand, if $x_1$ were a spear vector, we would have $|x_1^\ast(x_1)|=1$ for every $x_1^\ast\in \ext(B_{X_1^\ast})$ by Corollary \ref{Defi:spearVector}.(iv), so we would get a contradiction if we show that the functional  $x_1^\ast=(\frac12,\frac12,\frac12,\frac12)\in X_1^\ast$ is an extreme point of $B_{X_1^\ast}$. Let us show this last assertion. First, $x_1^\ast$ belongs to $B_{X^\ast_1}$ since for every $x_1=(\eps_1,\eps_2,\eps_3,\eps_4)\in \ext(B_{X_1})$ we have that
$$
|x^*(x)|=\frac12|\eps_1+\eps_2+\eps_3+\eps_4|\leqslant 1.
$$
Next, consider $y_1^\ast\in X_1^\ast$ such that both $x_1^\ast+y_1^\ast$ and $x_1^\ast - y_1^\ast$ lie in $B_{X^\ast_1}$. This, together with the fact that
$$
x_1^*(-1,1,1,1)=x_1^*(1,-1,1,1)=x_1^*(1,1,-1,1)=x_1^*(1,1,1,-1)=1,
$$
implies that
$$
y_1^*(-1,1,1,1)=y_1^*(1,-1,1,1)=y_1^*(1,1,-1,1)=y_1^*(1,1,1,-1)=0,
$$
so $y_1^\ast=0$ since it vanishes on a basis of $X_1$. This gives that $x_1^\ast$ is an extreme point, as desired.
\end{Exam}

When the codomain is finite-dimensional, we can improve Proposition \ref{Prop:characLushAsplund} as follows, just taking into account that weak$^\ast$-denting points and extreme points of the dual ball are the same for a finite-dimensional space.

\begin{Prop}\label{Prop:characLushfinitedimCOdomain}
Let $X$ be a Banach space, let $Y$ be a finite-dimensional space and let $G\in L(X,Y)$ be a norm-one operator. Then, the following assertions are equivalent:
\begin{enumerate}
\item[(i)] $G$ is lush.
\item[(ii)] $G$ is a spear operator.
\item[(iii)] $G$ has the aDP
\item[(iv)] $|x^{\ast\ast}(G^\ast y^\ast)|=1$ for every $x^{\ast\ast}\in \ext (B_{X^{\ast\ast}})$ and every $y^\ast\in \ext (B_{Y^\ast})$.
\item[(v)] $G^{\ast}y^{\ast}\in \Spear(X^\ast)$ for every $y^\ast\in \ext(B_{Y^\ast})$.
\item[(vi)] $B_X = \overline{\aconv}\bigl(\Face(S_X,G^\ast y^\ast)\bigr)$ for every $y^\ast\in \ext(B_{Y^\ast})$.
\item[(vii)] $B_{Y^\ast} = \conv{\bigl\{ y^{\ast}  \in B_{Y^{\ast}} \colon G^{\ast}y^{\ast}\in \Spear(X^\ast) \bigr\}}$.
\end{enumerate}
\end{Prop}

We do not know whether this result, or part of it, is also true when just the range of the operator $G$ is finite-dimensional. But we can provide with the following result for rank-one operators.

\begin{Coro}\label{coro:charac-lush-rank-one}
Let $X$, $Y$ be Banach spaces and let $G\in L(X,Y)$ be a norm-one rank-one operator, and write $G=x_0^\ast \otimes y_0$ for suitable $x_0^\ast\in S_{X^\ast}$ and $y_0\in S_{Y}$. Then, the following assertions are equivalent:
\begin{enumerate}
\item[(i)] $G$ is lush.
\item[(ii)] $G$ is a spear operator.
\item[(iii)] $G$ has the aDP
\item[(iv)] $x_0^\ast\in \Spear(X^\ast)$ and $y_0\in \Spear(Y)$.
\end{enumerate}
\end{Coro}

\begin{proof}
(i) $\Rightarrow$ (ii) $\Rightarrow$ (iii) are clear. Let us prove (iii) $\Rightarrow$ (iv). First, for every $y\in Y$, consider the rank-one operator $T=x_0^\ast\otimes y\in L(X,Y)$ and observe that
$$
\|y_0+\T\,y\|=\|x_0^\ast\otimes y_0 +\T\,x_0^\ast\otimes y\|= \|G + \T\,T\|=1 + \|T\|=1 + \|y\|,
$$
so $y_0\in \Spear(Y)$. Next, we have that $G:X\longrightarrow G(X)=\K y_0\equiv \K$, also has the aDP (use Remark  \ref{Rema:aDPOperators}) and we may use Proposition \ref{Prop:aDPDentingPoints} to get that $G^\ast(1)=x_0^\ast \in \Spear(X^\ast)$.

(iv) $\Rightarrow$ (i). Observe that  $G^\ast(y^\ast)=y^\ast(y_0)x_0^\ast$ for every $y^\ast\in Y^\ast$. Now, for each $y^\ast\in \ext(B_{Y^\ast})$ we have that $|y^\ast(y_0)|=1$  by Corollary \ref{Defi:spearVector}.(iv) (as $y_0\in \Spear(Y)$), so $G^\ast(y^\ast)\in \T\,x_0^\ast\subset \Spear(X^\ast)$. Now, Proposition \ref{Prop:lushSufficient}.(b) gives us that $G$ is lush.
\end{proof}

\section{Adjoint Operators}\label{subsec:duality}

We would like to discuss here the relationship of the aDP, spearness and lushness with the operation of taking the adjoint.

As the norm of an operator and the one of its adjoint coincide, the following observation is immediate.

\begin{Rema}\label{Remark:dualspearadp=>spearadp}
Let $X$, $Y$ be Banach spaces and let $G\in L(X,Y)$ be a norm-one operator. If $G^\ast$ is a spear operator, then $G$ is a spear operator. If $G^\ast$ has the aDP, then $G$ has the aDP.
\end{Rema}

With respect to lushness, the above result is not true, even for $G$ equal to the Identity, as the following example shows.

\begin{Exam}[\mbox{\cite[Theorem 4.1]{LushNumOneDual}}]\label{Exam:dual-lush-no-lush}
There is a (separable) Banach space $X$ such that $X^\ast$ is lush, but $X$ is not lush. Therefore, $G:=\Id:X\longrightarrow X$ is not lush, while $G^\ast$ is lush. Actually, $X=L_1[0,2]/Y$ where $Y$ is the space defined in Example \ref{Exam:fromL1-spear-no-lush}.
\end{Exam}

This example has some more properties which are interesting.

\begin{Rema}
Let $X$ be the Banach space of Example \ref{Exam:dual-lush-no-lush} and consider the operator $G:=\Id:X\longrightarrow X$.
\begin{enumerate}
\item[(a)] $G$ is a spear operator but it is not lush (use Remark \ref{Remark:dualspearadp=>spearadp}).
\item[(b)] $\Spear(X^\ast)=\emptyset$ (\cite[Remarks 4.2.(c)]{LushNumOneDual}). Therefore:
    \begin{itemize}
    \item Theorem \ref{Theo-lush=casi-almost-CL} is far from being true for spear operators;
    \item Proposition \ref{prop:RNP-equivalence} is far from being true for spaces without the Radon-Nikod\'{y}m Property;
    \item there is no lush operator whose domain is $X$.
    \end{itemize}
\end{enumerate}
\end{Rema}

We may give a positive result in this line: if the second adjoint of an operator is lush, then the operator itself is lush. This will be given in Corollary \ref{Coro:G**lush=>Glush}, but we need some preliminary work to get the result. We start with a general result which allows to restrict the domain of a lush operator.

\begin{Prop}\label{Prop:G**lushimplicaJ_YG=G**J_X-lush}
Let $X$, $Z$ be Banach spaces and let $H\in L(X^{\ast\ast},Z^\ast)$ be a weak$^\ast$-weak$^\ast$ continuous norm-one operator. If $H$ is lush, then $H\circ J_X:X\longrightarrow Z^\ast$ is lush.
\end{Prop}

For the sake of clearness, we include the most technical part of the proof of this result in the following lemma.

\begin{Lemm}
\label{Lemm:lushG**=>G}
Let $X$, $Y$, $W$ be Banach spaces and let $G_1\in L(X,Y)$ and $G_2\in L(Y,W)$ be norm-one operators. Suppose that there is a subset $A_1\subset B_{Y^\ast}$ such that $G_1$ satisfies the following property
 \begin{equation}\label{equation:P1}\tag{P1}
 \begin{split}
&\text{For every slice $S$ of $B_X$, every $y^{\ast} \in A_1$, and every $\eps>0$,}\\
&\bigl[G_1(S) \cap \conv{\GS(B_{Y},\T\,y^{\ast}, \eps)} \neq \emptyset\bigr]  \Rightarrow  \bigl[G_1(S) \cap \GS(B_{Y},\T\,y^{\ast}, \eps) \neq \emptyset\bigr].
 \end{split}
\end{equation}
Suppose also that $G_2$ is lush and there is a subset $A_2\subset S_{W^\ast}$ with $\overline{\conv}^{w^\ast}\,A_2=B_{W^\ast}$ such that $G_2^\ast(A_2)\subset A_1$. Then $G:=G_2\circ G_1$ is lush.
\end{Lemm}

\begin{proof}
Fix $x_{0} \in S_X$, $w_{0} \in B_{W}$ and $\varepsilon > 0$.
Let $\delta \in (0,\eps/3)$. Since $G_2$ is lush, applying Proposition \ref{Prop:characterization-lushness}.(iii) with  $\mathcal{A}=A_2\subset S_{W^\ast}$, we can find $w^{\ast} \in A_2$ such that
\begin{equation}\label{eq:lemma(P1)-eq1}
\re w^{\ast}(w_{0}) > 1 - \delta \qquad \text{ and } \qquad \dist\bigl(G_1 x_{0}, \conv \GS(B_{Y}, \T\,G_2^\ast w^\ast, \delta)\bigr) < \delta.
\end{equation}
Then, there are $m \in \N$, $\lambda_{j}\in [0,1]$, $y_{j} \in  \GS(B_{Y},\T\,G_2^\ast w^\ast, \delta)$ for each $j=1, \ldots, m$, with $\sum_{j}{\lambda_{j}} = 1$ and
\[
u:= G_1 x_{0} - \sum_{j=1}^{m}{\lambda_{j} y_{j}} \in \delta B_{Y}.
\]
Notice that for each $j$ we have that $\|y_{j} + u\| \leqslant 1 + \delta$, and
\[
\left|G_2^\ast w^\ast\left( \frac{y_{j} + u}{1 + \delta} \right) \right| \geqslant \frac{1 - 2 \delta}{1 + \delta} = 1 - \frac{3 \delta}{1 + \delta}. \]
This implies that
\[ G_1\left(\frac{x_{0}}{1 + \delta}\right) = \sum_{j=1}^{m}{\lambda_{j} \frac{y_{j} + u}{1 + \delta}} \in \conv \GS\left(B_{Y}, \T\,G_2^\ast w^\ast, \frac{3\delta}{1 + \delta}\right)\,.
\]
Hence, every slice $S$ of $B_{X}$ containing $x_{0}/(1 + \delta)$ satisfies that
\[
G_1(S) \cap \conv \GS\left(B_{Y},\T\,G_2^\ast w^\ast, \frac{3\delta}{1 + \delta}\right) \neq \emptyset.
\]
As $G_2^\ast w^\ast \in A_1$, the hypothesis \eqref{equation:P1} yields that
\[
G_1(S) \cap \GS\left(B_{Y},\T\,G_2^\ast w^\ast, \frac{3\delta}{1 + \delta}\right) \neq \emptyset
\]
and, therefore,
\[
S \cap \GS\left(B_{X},\T\,G_1^\ast G_2^\ast w^\ast, \frac{3\delta}{1 + \delta}\right) \neq \emptyset.
\]
Since $S$ was arbitrary, we conclude that
\[
\frac{x_{0}}{1 + \delta} \in \aconv \GS\left(B_{X}, G^\ast w^\ast, \frac{3\delta}{1 + \delta}\right),
\]
and so
\[
\dist \left( x_{0}, \aconv \GS\left(B_{X}, G^\ast w^\ast, \frac{3\delta}{1 + \delta}\right)\right) < \frac{\delta}{1 + \delta}.
\]
As $0<\delta<\eps/3$, we get that
$\dist \left( x_{0}, \aconv \GS\left(B_{X}, G^\ast w^\ast, \eps\right)\right) < \eps$. This, together with the first part of \eqref{eq:lemma(P1)-eq1}, gives that $G$ is lush by using again Proposition \ref{Prop:characterization-lushness}.(iii).
\end{proof}

\begin{proof}[Proof of Proposition \ref{Prop:G**lushimplicaJ_YG=G**J_X-lush}] We will use the above lemma with $X=X$, $Y=X^{\ast\ast}$, $W=Z^\ast$, $G_1=J_X$ and $G_2=H$. To this end, we first show that $G_1=J_X$ satisfies condition \eqref{equation:P1} of Lemma \ref{Lemm:lushG**=>G} with $A_1 = J_{X^\ast}(B_{X^{\ast}})\subset B_{X^{\ast\ast\ast}}$. Indeed, let us fix a slice of the form $\Slice(B_{X},x_{1}^{\ast}, \delta)$ of $B_X$, $J_{X^\ast}(x^\ast)\in A_1$ and $\eps>0$, and suppose that
$$
J_X\bigl(\Slice(B_{X},x_{1}^{\ast}, \delta)\bigr) \cap \conv\bigl( \GS(B_{X^{\ast\ast}},\T\,J_{X^\ast}(x^\ast),\eps)\bigr)\neq \emptyset.
$$
Since
$$
\conv\bigl( \GS(B_{X^{\ast\ast}},\T\,J_{X^\ast}(x^\ast),\eps)\bigr) \subset \overline{\conv}^{\sigma(X^{\ast\ast},X^\ast)}\, J_X\bigl(\GS(B_X,\T\,x^\ast,\eps)\bigr),
$$
we actually have that
$$
J_X\bigl(\Slice(B_{X},x_{1}^{\ast}, \delta)\bigr) \cap \overline{\conv}^{\sigma(X^{\ast\ast},X^\ast)}\, J_X\bigl(\GS(B_X,\T\,x^\ast,\eps)\bigr)\neq \emptyset
$$
and so, a fortiori,
$$
\Slice(B_{X^{\ast\ast}},J_{X^\ast}(x_{1}^{\ast}), \delta) \cap \overline{\conv}^{\sigma(X^{\ast\ast},X^\ast)}\, J_X\bigl(\GS(B_X,\T\,x^\ast,\eps)\bigr)\neq \emptyset.
$$
But it then follows that
$$
\Slice(B_{X^{\ast\ast}},J_{X^\ast}(x_1^\ast), \delta) \cap J_X\bigl(\GS(B_X,\T\,x^\ast,\eps)\bigr)\neq \emptyset.
$$
This clearly implies that $J_X\bigl(\Slice(B_{X},x_{1}^{\ast}, \delta)\bigr) \cap \GS(B_{X^{\ast\ast}},\T\,J_{X^\ast}(x^\ast),\eps) \neq \emptyset$, as desired.

Now, let $G_2=H\in L(X^{\ast\ast},Z^{\ast})$ and $A_2=J_{Z}(S_{Z})\subset Z^{\ast\ast}$ which is norming for $Z^{\ast}$. As $H$ is weak$^\ast$-weak$^\ast$ continuous, we have that $G_2^\ast(A_2)\subset A_1$ (indeed, let $H_\ast\in L(Z,X^\ast)$ such that $[H_\ast]^\ast=H$ and observe that  $G_2^\ast(J_Z z)=[H_\ast]^{\ast\ast}(J_z z)= J_{X^\ast}(Hz)$ for every $z\in S_Z$).

Therefore, all the requirements of Lemma \ref{Lemm:lushG**=>G} are satisfied, so $G_2\circ G_1=H\circ J_X$ is lush.
\end{proof}

We get a couple of corollaries of this result. The first one deals with the natural inclusion of a lush Banach space into its bidual. It is an immediate consequence of the result above applied to $H=\Id_{X^{\ast\ast}}$.

\begin{Coro}\label{Coro:inclusionintobiduallush}
Let $X$ be a Banach space. If $X^{\ast \ast}$ is lush, then the canonical inclusion $J_X: X \longrightarrow X^{\ast \ast}$ is lush.
\end{Coro}

The next consequence is the promised result saying that lushness passes from the biadjoint operator to the operator.

\begin{Coro}\label{Coro:G**lush=>Glush}
Let $X$, $Y$ be Banach spaces and let $G\in L(X,Y)$ be a norm-one operator. If $G^{\ast\ast}$ is lush, then $G$ is lush.
\end{Coro}

\begin{proof}
Apply Proposition \ref{Prop:G**lushimplicaJ_YG=G**J_X-lush} to $H=G^{\ast\ast}\in L(X^{\ast\ast},Y^{\ast\ast})$, which is weak$^\ast$-weak$^\ast$ continuous, to get that $G^{\ast\ast}\circ J_X:X\longrightarrow Y^{\ast\ast}$ is lush. But, clearly, $G^{\ast\ast}\circ J_X = J_Y \circ G$ and then, restricting the codomain and considering that $Y$ and $J_Y(Y)$ are isometrically isomorphic, Remark \ref{Rema:lushOperators-elementary} gives us that $G$ is lush.
\end{proof}

These two corollaries improve \cite[Proposition 4.3]{LushNumOneDual} where it is proved that a Banach space $X$ is lush whenever $X^{\ast\ast}$ is lush.

\begin{Rema}
The technical hypothesis $G_2^\ast(A_2)\subset A_1$ in Lemma \ref{Lemm:lushG**=>G} is fundamental to get the result. Indeed, consider the inclusion $J: c_{0} \longrightarrow \ell_{\infty}$ and the projection $P: \ell_{\infty} \longrightarrow \ell_{\infty}/c_{0}$. Notice that $J = J_{c_{0}}$ satisfies the condition \eqref{equation:P1} of Lemma \ref{Lemm:lushG**=>G} with $A_{1} = J_{\ell_{1}}(B_{\ell_{1}}) \subset B_{\ell_{\infty}^{\ast}}$ (this is shown in the proof of Proposition \ref{Prop:G**lushimplicaJ_YG=G**J_X-lush}). On the other hand $P$ is lush since it carries every spear vector of $\ell_{\infty}$ into a spear vector of $\ell_{\infty}/c_{0}$. This can be easily seen using the canonical (isometric) identifications $\ell_{\infty} \equiv C(\beta \N)$ and $\ell_{\infty}/c_{0} = C(\beta \N \setminus \N)$, so that $P$ is just the restriction operator. On the other hand, $P \circ J = 0$, which clearly is not lush. The technical hypothesis of the lemma is not satisfied, since every $\mu \in C(\beta \N \setminus \N)^{\ast}$ with $P^{\ast}\mu \in B_{\ell_{1}}$ must be zero.

The same example also shows the need of the operator $H$ in Proposition \ref{Prop:G**lushimplicaJ_YG=G**J_X-lush} to be weak$^\ast$-weak$^\ast$-continuous: indeed, just take $H=P:\ell_\infty \longrightarrow \ell_\infty/c_0$ and observe that $H\circ J_X=0$.
\end{Rema}

Let us now discuss the more complicated direction: when lushness, spearness or the aDP passes from an operator to its adjoint. It is easy to provide examples of operators with the aDP whose adjoint do not share the property: for instance this is the case of the Identity operator on the space $C([0,1],\ell_2)$ (indeed, this space has the aDP by \cite[Example in p.~858]{KSSW}, while its dual contains $\ell_2$ as $L$-summand and so it fails the aDP by \cite[Proposition 3.1]{martinOikhberg}). But the same question for spearness of the Identity (that is, whether the numerical index one passes from a Banach space to its dual) was one of the main open questions of the theory of numerical index until 2007, when it was solved in the negative \cite{NumIndexDuality}. The counterexample given there is actually lush, and its dual even fails the aDP. Let us state all of this here.

\begin{Exam}[\mbox{\cite[\S 3]{NumIndexDuality} }]\label{Example:lushdualnotaDP} Consider the Banach space
$$
X=\bigl\{(x,y,z)\in c\oplus_\infty c\oplus_\infty c\colon \lim x + \lim y + \lim z=0\bigr\},
$$
and write $G:=\Id:X\longrightarrow X$. Then, $G$ is lush ($X$ is actually a C-rich subspace of $c\oplus_\infty c\oplus_\infty c$), but $G^\ast$ does not even have the aDP.
\end{Exam}

Some remarks about this example are interesting. Recall that a James boundary for a Banach space $X$ is a subset $C$ of $B_{X^\ast}$ such that $\|x\|=\max_{x^\ast\in C}|x^\ast(x)|$. As a consequence of the Hanh-Banach and the Krein-Milman theorems, the set $\ext(B_{X^\ast})$ is a James boundary for $X$.

\begin{Rema}
Let $X$ be the space of Example \ref{Example:lushdualnotaDP} and $G:=\Id_X$. Then the set $\Spear(X^\ast)$ is norming for $X$, but it is not a James boundary for $X$ so, in particular, it does not coincide with $\ext(B_{X^\ast})$ \cite[Example 3.4]{NumIndexDuality}. Therefore:
\begin{itemize}
\item Theorem \ref{Theo-lush=casi-almost-CL} cannot be improved to get that the set $\Omega$ is the whole set of extreme points, nor a James boundary for $X$;
\item the G$_\delta$ dense set in Proposition \ref{Prop:characLushAsplund}.(v) does not always coincide with the set of all extreme points of the dual ball, nor is always a James boundary for $X$.
\end{itemize}
\end{Rema}

Our next goal is to provide sufficient conditions which allow to pass the properties of an operator to its adjoint. The first of these conditions is that the domain space has the Radon-Nikod\'{y}m Property.

\begin{Prop}\label{Prop:RNP-duality}
Let $X$ be a Banach space with the Radon-Nikod\'{y}m Property, let $Y$ be a Banach space and let $G\in L(X,Y)$ be a norm-one operator. If $G$ has the aDP, then $G^\ast$ is lush. Therefore, the following six assertions are equivalent: $G$ has the aDP, $G$ is a spear operator, $G$ is lush, $G^\ast$ has the aDP, $G^\ast$ is  a spear operator, $G^\ast$ is lush.
\end{Prop}

\begin{proof}
Write $D$ for the set of denting points of $B_X$. If $G$ has the aDP, Proposition \ref{Prop:aDPDentingPoints} gives that $Gx\in \Spear(Y)$ for every $x\in D$, and then Proposition \ref{Prop:spearVectorsProperties}.(c) gives that $J_Y(Gx)\in \Spear(Y^{\ast\ast})$ for every $x\in D$. Therefore, the set
$$
\bigl\{x^{\ast\ast}\in B_{X^{\ast\ast}}\colon [G^\ast]^\ast(x^{\ast\ast})\in \Spear(Y^{\ast\ast})\bigr\}
$$
contains $J_X(D)$ which is norming for $X^\ast$ as $X$ has the Radon-Nikod\'{y}m Property (see \cite[\S 2]{Bourgin} for instance). Then, Proposition \ref{Prop:lushSufficient}.(a) gives that $G^\ast$ is lush.

Finally, let us comment the proof of the last part. The three first assertions are equivalent by Proposition \ref{prop:RNP-equivalence} since $X$ has the Radon-Nikod\'{y}m Property; $G^\ast$ lush $\ \Rightarrow\ $ $G^\ast$ spear $\ \Rightarrow\ $ $G^\ast$ has the aDP $\ \Rightarrow\ $ $G$ has the aDP by Remark \ref{Remark:dualspearadp=>spearadp}. The remaining implication is just what we have proved above.
\end{proof}

Another result in this line is the following. We recall that a Banach space $X$ is \emph{$M$-embedded} if $J_X(X)^\perp$ is an $L$-summand in $X^{\ast\ast\ast}$ (which is actually equivalent to the fact that the Dixmier projection on $X^{\ast\ast\ast}$ is an $L$-projection). We refer the reader to the monograph \cite{HWW} for more information and background. Examples of $M$-embedded spaces are reflexive spaces (trivial), $c_0$ and all of its closed subspaces, $K(H)$ (the space of compact operators on a Hilbert space $H$), $C(\T)/A(\mathbb{D})$, the little Bloch space $B_0$, among others (see \cite[Examples III.1.4]{HWW}).
\index{Membedded@$M$-embedded}%

\begin{Prop}\label{Prop:M-embedded-duality}
Let $X$ be a Banach space, let $Y$ be an $M$-embedded Banach space, and let $G\in L(X,Y)$ be a norm-one operator. If $G$ has the aDP, then $G^\ast$ is lush. Therefore, the following nine assertions are equivalent: $G$ has the aDP, $G$ is a spear operator, $G$ is lush, $G^\ast$ has the aDP, $G^\ast$ is a spear operator, $G^\ast$ is lush, $G^{\ast\ast}$ has the aDP, $G^{\ast\ast}$ is a spear operator, $G^{\ast\ast}$ is lush.
\end{Prop}

\begin{proof}
We use Proposition \ref{Prop:aDPDentingPoints} to get that the set
$
\bigl\{y^\ast\in B_{Y^\ast}\colon G^\ast y^\ast \in \Spear(X^\ast)\bigr\}
$
contains the set $D$ of those weak$^\ast$-denting points of $B_{Y^\ast}$. By \cite[Corollary~III.3.2]{HWW}, we have that $B_{Y^\ast}=\overline{\conv}(D)$ so, a fortiori,
$$
B_{Y^\ast}=\overline{\conv}\,\bigl\{y^\ast\in B_{Y^\ast}\colon G^\ast y^\ast \in \Spear(X^\ast)\bigr\}.
$$
Then, Proposition \ref{Prop:lushSufficient}.(c) gives that $G^\ast$ is lush.

Finally, for the last part, the three first assertions are equivalent by Proposition \ref{Prop:characLushAsplund} since $Y$ is Asplund \cite[Theorem~III.3.2]{HWW}. The middle three assertions are equivalent by Proposition \ref{prop:RNP-equivalence} since $Y^\ast$ has the Radon-Nikod\'{y}m Property \cite[Theorem~III.3.2]{HWW}. If $G^\ast$ has the aDP, so does $G$ (Remark \ref{Remark:dualspearadp=>spearadp}) and this implies that $G^\ast$ is lush by the above. As $Y^\ast$ has the Radon-Nikod\'{y}m Property, if $G^\ast$ has the aDP, then $G^{\ast\ast}$ is lush by Proposition \ref{Prop:RNP-duality}, and this gives the equivalence with the last three assertions.
\end{proof}

Even though part of what we have used in the proof above is Asplundness of $M$-embedded spaces, just this hypothesis on $Y$ is not enough to get the result as Example \ref{Example:lushdualnotaDP} shows.

A consequence of the two results above is that lushness passes from an operator with finite-dimensional domain or codomain to all of its successive adjoint operators.

\begin{Coro}\label{Coro:duality-fin-dim}
Let $X$, $Y$ be Banach spaces such that at least one of them is finite-dimensional space, and let $G\in L(X,Y)$ be a norm-one operator. If $G$ has the aDP, then all the successive adjoint operators of $G$ are lush.

\end{Coro}

\begin{proof}
If $X$ is finite-dimensional, then it has the Radon-Nikod\'{y}m Property, so Proposition \ref{Prop:RNP-duality} gives that $G^{\ast}$ is lush. If $Y$ is finite-dimensional, then it is clearly $M$-embedded, so Proposition \ref{Prop:M-embedded-duality} gives us that $G^\ast$ is lush. For the successive adjoint operators, one of the above two arguments applies.
\end{proof}

We do not know whether the above result can be extended to finite-rank operators. We may do when the operator has actually rank one.

\begin{Prop}\label{Prop:rank-one-GaDP=>G*lush}
Let $X$, $Y$ be Banach spaces, and let $G\in L(X,Y)$ be a rank-one norm-one operator. If $G$ has the aDP, then $G^\ast$ is lush. Therefore, all the successive adjoints of $G$ are lush.
\end{Prop}

\begin{proof}
If $G$ has the aDP, by Corollary \ref{coro:charac-lush-rank-one} we have that $G=x_0^\ast \otimes y_0$ with $x_0^\ast\in \Spear(X^\ast)$ and $y_0\in \Spear(Y)$. Observe that $G^\ast=J_Y(y_0)\otimes x_0^\ast:Y^\ast \longrightarrow X^\ast$. Since $J_Y(y_0)\in \Spear(Y^{\ast\ast})$ by Proposition \ref{Prop:spearVectorsProperties}.(c), we get that $G^\ast$ is lush by using again Corollary \ref{coro:charac-lush-rank-one}.
\end{proof}

The last result deals with $L$-embedded spaces. Recall that a Banach space $Y$ is \emph{$L$-embedded} if $Y^{\ast\ast}=J_Y(Y)\oplus_1 Y_s$ for suitable closed subspace $Y_s$ of $Y^{\ast\ast}$. We refer to the monograph \cite{HWW} for background.
\index{Lembedded@$L$-embedded}%
Examples of $L$-embedded spaces are reflexive spaces (trival), predual of von Neumann algebras so, in particular, $L_1(\mu)$ spaces, the Lorentz spaces $d(w,1)$ and $L^{p,1}$, the Hardy space $H_0^1$, the dual of the disk algebra $A(\mathbb{D})$, among others (see \cite[Examples IV.1.1 and III.1.4]{HWW}).

\begin{Prop}\label{Prop-Lembedded-duality}
Let $X$ be a Banach space, let $Y$ be an $L$-embedded space, and let $G\in L(X,Y)$ be a norm-one operator.
\begin{enumerate}
\item[(a)] If $G$ is a spear operator, then $G^\ast$ is a spear operator.
\item[(b)] If $G$ has the aDP, then $G^\ast$ has the aDP.
\end{enumerate}
\end{Prop}

\begin{proof}
(a). Write $P_Y:Y^{\ast\ast}\longrightarrow J_Y(Y)$ for the  projection associated to the decomposition $Y^{\ast\ast}=J_Y(Y) \oplus_1 Y_s$. We fix $T\in L(Y^\ast,X^\ast)$ and consider the operators
$$
A:=P_Y\circ T^\ast \circ J_X:X\longrightarrow J_Y(Y) \qquad B:=(\Id-P_Y)\circ T^\ast\circ J_X:X\longrightarrow Y_s,
$$
and observe that $T^\ast\circ J_X=A\oplus B$. Given $\eps>0$, since $J_X(B_X)$ is dense in $B_{X^{\ast\ast}}$ by Goldstine's Theorem and $T^\ast$ is weak$^\ast$-weak$^\ast$-continuous, we may find $x_0\in S_X$ such that
$$
\|T^\ast x_0\|=\|A x_0\|+ \|B x_0\|>\|T\|-\eps.
$$
Now, we may find $y_0\in S_Y$ and $y_s^\ast\in S_{Y_s^\ast}$ such that
$$
\|A x_0\|\,y_0=A x_0\qquad \text{and} \qquad y_s^\ast(Bx_0)=\|Bx_0\|.
$$
We define $S:X\longrightarrow Y$ by $Sx=Ax+y_s^\ast(Bx)y_0$ for every $x\in X$, and observe that $\|S\|\geqslant \|Sx_0\|>\|T\|-\eps$. As $G$ is a spear operator, we have that $\|G + \T\,S\|>1+ \|T\|-\eps$, so we may find $x_1\in S_X$, $w\in \T$, and $y_1^\ast\in S_{Y^\ast}$ such that
$$
\bigl|y_1^\ast(Gx_1 + w Ax_1 + w\, y_s^\ast(Bx_1)y_0)\bigr|>1 + \|T\|-\eps.
$$
Finally, consider $\Phi=(J_{Y^\ast}(y_1^\ast),\,y_1^\ast(y_0)y_s^\ast)\in Y^{\ast\ast\ast}=J_{Y^\ast}(Y^\ast)\oplus_\infty Y_s^\ast$ which has norm-one (here we use the $L$-embeddedness hypothesis) and observe that
\begin{align*}
\|G^\ast + \T\,T\|&=\|G^{\ast\ast}+\T\,T^\ast\| \geqslant \bigl|\bigl[\Phi\bigl(G^{\ast\ast}+ w T^\ast\bigr)\bigr](J_X(x_1))\bigr| \\
&= \bigl|y_1^\ast(Gx_1 + w Ax_1) + w\, y_1^\ast(y_0)y_s^\ast(Bx_1)\bigr|\\
& = \bigl|y_1^\ast(Gx_1 + w A x_1 + w\, y_s^\ast(B x_1)y_0)\bigr| \\ &> 1 + \|T\|-\eps.
\end{align*}
Moving $\eps\downarrow 0$, we get that $G^\ast$ is a spear operator, as desired.

(b). If $G$ just has the aDP, we may repeat the above argument for rank-one operators $T\in L(Y^\ast,X^\ast)$, and everything works fine as the operator $S\in L(X,Y)$ constructed there has finite rank, so $\|G+\T\,S\|=1+\|S\|$ by Theorem \ref{Theo:SCDOperators} (as, clearly, finite-rank operators are SCD).
\end{proof}

\chapter{Isometric and isomorphic consequences}\label{sect:consequences}

Our goal here is to present consequences on the Banach spaces $X$ and $Y$ of the fact that there is $G\in L(X,Y)$ which is a spear operator, is lush or has the aDP.

We first start with a deep structural consequence which generalizes \cite[Corollary 4.10]{SCDsets} where it was proved for real infinite-dimensional Banach spaces with the aDP.

\begin{Theo}\label{Theo:G-finite-rank-noell_1}
Let $X$, $Y$ be \textbf{real} Banach spaces and let $G\in L(X,Y)$. If $G$ has the aDP and has infinite rank, then $X^{\ast}$ contains a copy of $\ell_1$.
\end{Theo}

\begin{proof}
Using Proposition \ref{Prop:separablyDetermined}, we can find separable subspaces $X_{\infty} \subset X$ and $Y_{\infty} \subset Y$ such that $G_\infty:=G|_{X_{\infty}}: X_{\infty} \longrightarrow Y_{\infty}$ has the aDP, and still it has infinite rank. By Remark \ref{Rema:aDPOperators}, we may and do suppose that $\overline{G_\infty(X_\infty)}=Y_\infty$. It is enough to show that $X_{\infty}^{\ast}$ contains a copy of $\ell_1$ since, in this case, $X^{\ast}$ also contains such a copy by the lifting property of $\ell_1$ (see \cite[\mbox{Proposition~2.f.7}]{Lin-Tza-I}
or \cite[p. 11]{vanDulst}). We have two possibilities. If $X_\infty$ contains a copy of $\ell_1$, then $X_\infty^\ast$ contains a quotient isomorphic to $\ell_\infty$ and so $X_\infty^\ast$ contains a copy of $\ell_1$ again by the lifting property of $\ell_1$. If $X_\infty$ does not contain copies of $\ell_1$ then  $B_{X_{\infty}}$ is an SCD set by \cite[Theorem 2.22]{SCDsets} (see Example \ref{Exam:SCD-set-spaces-operators}), so Theorem~\ref{Theo:aDP+SCD=lush}
gives that $G_\infty$ is lush. Then, by Theorem \ref{Theo-lush=casi-almost-CL}, the set $\bigl\{ y^{\ast} \in \ext B_{Y_{\infty}^{\ast}}\colon G_\infty^{\ast}(y^{\ast})\in \Spear(X_\infty^\ast)\bigr\}$ is weak$^\ast$-dense in $\ext B_{Y_\infty^{\ast}}$. As $G_\infty$ has dense range, $G_\infty^\ast$ is injective, and since $Y_\infty^\ast$ is infinite-dimensional, it follows that the set $\Spear(X_\infty^\ast)$ must be infinite. Now, Proposition \ref{Prop:spearVectorsProperties}.(i) gives us that $X_\infty^\ast$ contains a copy of $c_0$ or $\ell_1$. But a dual space contains a copy of $\ell_1$ whenever it contains a copy of $c_0$ \cite[\mbox{Proposition 2.e.8}]{Lin-Tza-I}.
\end{proof}

Another result in this line is the following.

\begin{Prop}
Let $X$ be a \textbf{real} Banach space with the Radon-Nikod\'{y}m Property, let $Y$ be a \textbf{real} Banach space, and let $G\in L(X,Y)$. If $G$ has the aDP and has infinite rank, then $Y\supset c_0$ or $Y\supset \ell_1$.
\end{Prop}

\begin{proof}
By Proposition \ref{prop:RNP-equivalence} we have that $B_{X} = \overline{\conv} \bigl\{ x \in B_{X} \colon Gx\in \Spear(Y)\bigr\}$, so
$$
G(B_X)\subseteq \overline{\conv}\{Gx\colon x\in B_X,\, Gx\in \Spear(Y)\} \subseteq \overline{\conv}\, \Spear(Y).
$$
Now, if $G$ has infinite rank, $\Spear(Y)$ has to be infinite and so Proposition \ref{Prop:spearVectorsProperties}.(j) gives the result.
\end{proof}

\begin{Rema}
Let us observe that both possibilities in the result above may happen. On the one hand, $G=\Id_{\ell_1}:\ell_1\longrightarrow \ell_1$ is lush by Example \ref{Exam:ell_1-gamma-all-equivalent}. On the other hand, the operator $G:\ell_1 \longrightarrow c$ given by $[G(e_n)](k)=-1$ if $k=n$ and $[G(e_n)](k)=1$ if $k\neq n$ is also lush by Example \ref{Exam:ell_1-gamma-all-equivalent}.
\end{Rema}

We next deal with isometric consequences of the existence of operators with the aDP. The following result generalizes \cite[Theorem~2.1]{ConvexSmooth} where it was proved for $G=\Id$. Let us remark that the proof given there relied on a non-trivial result of the theory of numerical range: that the set of operators whose adjoint attain its numerical radius is norm dense in the space of operators.

\begin{Prop} \label{Prop:aDPDconvex-smooth}
Let $X$, $Y$ be Banach spaces and let $G\in L(X,Y)$ be an operator with the aDP. Then
\begin{enumerate}
\item[(a)] If $X^{\ast}$ is strictly convex, then $X=\K$.
\item[(b)] If $X^{\ast}$ is smooth, then $X = \mathbb{K}$.
\item[(c)] If $Y^\ast$ is strictly convex, then $Y=\K$.
\end{enumerate}
\end{Prop}

\begin{proof}
(a). We start showing that $G^*$ has rank one. Using Theorem \ref{Theo:aDPCharacterization}.(iv) we can find $y_{0}^{\ast} \in S_{Y^{\ast}}$ with $\| G^{\ast}y_{0}^{\ast}\| = 1$. By the same result, there is a $\weakstar$-dense subset of $\ext B_{Y^{\ast}}$ whose elements $y^{\ast}$ satisfy that
\begin{equation}\label{eq-Xaststrictlyconvex}
\| G^{\ast}y_{0}^{\ast} + \T\, G^{\ast}y^{\ast}\| = 2.
\end{equation}
It follows from the definition of strict convexity that $G^{\ast}y^{\ast} \in \T G^{\ast}y_{0}^{\ast}$ for every such $y^{\ast}$, and we deduce by the Krein-Milman Theorem and the weak$^\ast$ continuity of $G^\ast$, that $G^{\ast}(B_{Y^{\ast}})$ is contained in $\spn{\{ G^{\ast}y_{0}^{\ast}\}}$. Hence, $G^{\ast}$ has rank one. Therefore, $G$ has rank one and $\Spear(X^\ast)$ is non empty by Corollary~\ref{coro:charac-lush-rank-one}. Finally, $X^\ast$ is one-dimensional by Proposition~\ref{Prop:spearVectorsProperties}.(h).

(b). Given arbitrary elements $x_{0}^{\ast}, x_{1}^{\ast} \in S_{X^{\ast}}$ we use use Theorem~\ref{Theo:aDPCharacterization}.(iv) and the fact that $(\ext B_{Y^{\ast}}, \weakstar)$ is a Baire space (see Lemma~\ref{Kreinlemma}.(c)) to deduce the existence of some $y^{\ast} \in \ext B_{Y^{\ast}}$ with $\| G^{\ast}y^{\ast} + \T\, x_{i}^{\ast}\| = 2$ for $i=0,1$. Take now $x_{i}^{\ast \ast} \in S_{X^{\ast \ast}}$ with $x_{i}^{\ast \ast}(G^{\ast}y^{\ast}) + |x_{i}^{\ast \ast}(x_{i}^{\ast})|=2$ for each $i=0,1$. If $X^{\ast}$ is smooth, then $x_{0}^{\ast \ast} = x_{1}^{\ast \ast}$ and hence $\| x_{0}^{\ast} + \T\, x_{1}^{\ast}\| = 2$. So every element of $S_{X^{\ast}}$ is a spear and then Proposition~\ref{Prop:spearVectorsProperties}.(e) tells us that $X^*$ is one-dimensional.

The proof of (c) follows the lines of the one of (a). Indeed, arguing like in (a), we find $y_0^\ast\in S_{Y^\ast}$ and a weak$^\ast$-dense subset of $\ext(B_{Y^\ast})$ whose elements $y^\ast$ satisfy \eqref{eq-Xaststrictlyconvex} so, a fortiori, they satisfy that
$$
\|y_0^\ast + \T\,y^\ast\|=2.
$$
Being $Y^\ast$ strictly convex, we get that $y^\ast\in \T\,y^\ast_0$ for every such $y^\ast$, but this implies that $Y^\ast$, and so $Y$, is one-dimensional by the Krein-Milman Theorem.
\end{proof}

The following result generalizes \cite[Proposition~2.5]{ConvexSmooth}.

\begin{Prop} \label{Prop:aDPD-Frechet}
Let $X$, $Y$ be Banach spaces and let $G\in L(X,Y)$ be an operator with the aDP.
\begin{enumerate}
\item[(a)] If the norm of $Y$ is Fr\'{e}chet smooth, then $Y=\K$.
\item[(b)] If $X$ and $Y$ are real spaces and the norm of $X$ is Fr\'{e}chet smooth, then $X=\R$.
\end{enumerate}
\end{Prop}

\begin{proof}
(a). By Proposition~\ref{Prop:aDPDentingPoints} we have that $G^\ast y^\ast\in \Spear(X^\ast)$ for every weak$^\ast$-strongly exposed point $y^\ast$ of $B_{Y^\ast}$. Since the norm of $Y$ is Fr\'{e}chet smooth, every functional in $S_{Y^\ast}$ attaining its norm is a weak$^\ast$-strongly exposed point of $B_{Y^\ast}$ (see \cite[Corollary~I.1.5]{DGZ} for instance). As norm-one norm attaining functionals are dense in $S_{Y^\ast}$ by the Bishop-Phelps Theorem and $\Spear(X^*)$ is norm closed by Proposition~\ref{Prop:spearVectorsProperties}.(d), we get in fact that $G^\ast y^\ast\in \Spear(X^\ast)$ for every $y^\ast\in S_{Y^\ast}$. So, given arbitrary elements $y_1^\ast,y_2^\ast\in S_{Y^\ast}$ we can write
\begin{align*}
2=\|G^\ast(y_1^\ast)+\T G^\ast(y_2^\ast)\|\leqslant \|y_1^\ast+\T y_2^\ast\|\leqslant 2
\end{align*}
which gives that every element in $S_{Y^\ast}$ is  a spear. Therefore, $Y^*$ is one-dimensional by Proposition~\ref{Prop:spearVectorsProperties}.(e).

(b). Fixed $X_0\subset X$ and $Y_0\subset Y$ arbitrary separable subspaces we can use Proposition~\ref{Prop:separablyDetermined} to find separable subspaces $X_{0} \subset X_{\infty} \subset X$ and $Y_{0} \subset Y_{\infty} \subset Y$ such that $G(X_{\infty}) \subset Y_{\infty}$ and $G_\infty:=G|_{X_{\infty}} : X_{\infty} \longrightarrow Y_{\infty}$ has norm one and the aDP. Next, we fix a countable dense subset $D\subset S_{X_\infty}$ and we consider $D^\ast \subset S_{X_\infty^\ast}$ given by
$$
D^\ast=\{x^\ast\in S_{X_\infty^\ast} \ : \ \exists x\in D \text{ with } x^\ast(x)=1\}
$$
which is countable since $D$ is countable and $X_\infty$ is smooth. Therefore, we can use the fact that $(\ext B_{Y^{\ast}}, \weakstar)$ is a Baire space (see Lemma~\ref{Kreinlemma}.(c)) and Theorem~\ref{Theo:aDPCharacterization}.(iv) to deduce the existence of some $y^{\ast} \in \ext B_{Y_\infty^{\ast}}$ with $\|G_\infty^{\ast}y^{\ast} + \T\, x^{\ast}\| = 2$ for every $x^\ast\in D^\ast$. We will show that $G_\infty^{\ast}y^{\ast}\in \Spear(X_\infty^\ast)$. To do so, fix $x^\ast\in S_{X_\infty^\ast}$ attaining its norm at $x\in S_{X_\infty}$ and recall that $x$ strongly exposes $x^*$ as $X_\infty$ is Fr\'{e}chet smooth. Let $(x_n)$ be a sequence in $D$ converging to $x$ and let $x_n^\ast\in D^\ast$ satisfying $x_n^\ast(x_n)=1$ for every $n\in \N$. Then we have that
$$
|x_n^\ast(x)-1|=|x_n^\ast(x)-x_n^\ast(x_n)|\leqslant \|x-x_n\|\longrightarrow 0
$$
so $(x_n^\ast)$ converges in norm to $x^*$ and, therefore,
$$
2=\| G_\infty^{\ast}y^{\ast} + \T\, x_n^{\ast}\|\longrightarrow \| G_\infty^{\ast}y^{\ast} + \T\, x^{\ast}\|
$$
which gives $\|G_\infty^{\ast}y^{\ast} + \T\, x^{\ast}\|=2$. Since norm-attaining norm-one functionals are dense in $S_{X_\infty^\ast}$ by the Bishop-Phelps Theorem, we deduce that $G_\infty^{\ast}y^{\ast}$ is a spear in $X_\infty^\ast$. Finally, Proposition~\ref{Prop:spearVectorsProperties}.(k) tells us that $X_\infty$, and thus $X_0$, is one-dimensional as it is smooth. The arbitrariness of $X_0$ implies that $X$ is one-dimensional.
\end{proof}

The next result deals with WLUR points. Given a Banach space $X$, a point $x\in S_X$ is said to be \emph{LUR} (respectively \emph{WLUR}) if for every sequence $(x_n)$ in $B_X$ such that $\|x_n + x\|\longrightarrow 2$ one has that $(x_n)\longrightarrow x$ in norm (respectively weakly).  It is clear that LUR points are WLUR, but the converse result is known to be false \cite{Smith}.
\index{LUR}%
\index{WLUR}%

\begin{Prop} \label{Prop:aDPD-WLURPoint}
Let $X$, $Y$ be Banach spaces and let $G\in L(X,Y)$ be an operator with the aDP. Then
\begin{enumerate}
\item[(a)] If $B_X$ contains a WLUR point, then $X=\K$.
\item[(b)] If $B_Y$ contains a WLUR point, then $Y=\K$.
\end{enumerate}
\end{Prop}

\begin{proof}
(a). Let $x_0$ be a WLUR point of $B_X$. We start showing that $\|Gx_0\|=1$. To do so, take $x_0^\ast\in S_{X^\ast}$ with $x_0^\ast(x_0)=1$ and use Theorem~\ref{Theo:aDPCharacterization}.(iv) to find $y^\ast\in \ext B_{Y^{\ast}}$ such that $\|G^\ast y^\ast +\T x_0^\ast\|=2$. Therefore, there is a sequence $(x_n)$ in $B_X$ satisfying
$$
\bigl|[G^\ast y^\ast](x_n) +\T x_0^\ast(x_n)\bigr|\longrightarrow 2
$$
which clearly implies $|y^\ast(Gx_n)|= \bigl|[G^\ast y^\ast](x_n)\bigr|\longrightarrow1$ and $|x_0^\ast(x_n)|\longrightarrow1$. Hence, there is a sequence $(\theta_n)$ in $\T$ such that $\re x_0^\ast(\theta_nx_n)\longrightarrow1$ and so
$$
\|\theta_nx_n+x_0\|\geqslant \re x_0^\ast(\theta_nx_n+x_0)\longrightarrow2.
$$
Now since $x_0$ is a WLUR point we get that $(\theta_nx_n)$ converges weakly to $x_0$. Therefore, $(G\theta_nx_n)$ converges weakly to $Gx_0$, and the fact that $|y^\ast(G\theta_nx_n)|\longrightarrow1$ tells us that $|y^\ast(Gx_0)|=1$.

Suppose that $X$ is not one-dimensional, then there is $x^\ast\in S_{X^\ast}$ with $x^*(x_0)=0$. Consider the operator $T=x^\ast\otimes Gx_0 \in L(X,Y)$ which satisfies $\|T\|=1$. We have that $\|G+\T\, T\|=2$ since $G$ has the aDP, so there are sequences $(z_n)$ in $S_X$ and $(y_n^\ast)$ in $S_{Y^\ast}$ such that
$$
|y^\ast_n(Gz_n)+\T y^\ast_n(Gx_0)x^\ast(z_n)|\longrightarrow 2
$$
which implies $|x^\ast(z_n)|\longrightarrow1$ and $|y^\ast_n(Gz_n)+\T y^\ast_n(Gx_0)|\longrightarrow2$. Hence, we may find a sequence $(\omega_n)$ in $\T$ such that $| y^\ast_n (\omega_nGz_n+Gx_0)|\longrightarrow 2$ and so
$$
\|\omega_nz_n+x_0\|\geqslant \|G(\omega_nz_n+x_0)\|\geqslant | y^\ast_n (\omega_nGz_n+Gx_0)|\longrightarrow 2.
$$
Since $x_0$ is a WLUR point we get that $(\omega_nz_n)$ converges weakly to $x_0$. This, together with $|x^\ast(z_n)|\longrightarrow1$, tells us that $|x^\ast(x_0)|=1$ which is a contradiction.

(b). Let $y_0$ be a WLUR point of $B_Y$. Since $G$ has the aDP, Theorem~\ref{Theo:aDPCharacterization}.(iv) provides us with a dense G$_\delta$ set $A$ in $(\ext B_{Y^{\ast}},w^\ast)$ such that $\|G^\ast y^\ast\|=1$ for every $y^{\ast}\in A$. We claim that $|y^\ast(y_0)|=1$ for every $y^{\ast}\in A$. Indeed, fixed $y^\ast \in A$, consider the rank-one operator $T=G^\ast y^\ast \otimes y_0$ which satisfies $\|G+\T\, T\|=2$. So there are sequences $(x_n)$ in $S_X$ and $(y^\ast_n)$ in $S_{Y^\ast}$ such that
$$
2\longleftarrow |y_n^\ast(Gx_n)+\T y_n^\ast(Tx_n)|=|y_n^\ast(Gx_n)+\T y^\ast_n(y_0)y^\ast(Gx_n)|.
$$
This implies that $|y^\ast(Gx_n)|\longrightarrow 1$ and that there is a sequence $(\theta_n)$ in $\T$ such that
$$
\|\theta_nGx_n+y_0\|\geqslant |y_n^\ast(\theta_nGx_n+y_0)|\longrightarrow 2.
$$
Being $y_0$ a WLUR point, we deduce that $(\theta_nGx_n)$ converges weakly to $y_0$ and, therefore, we get $|y^\ast(y_0)|=1$, finishing the proof of the claim.

To finish the proof, fix $y\in S_Y$ and observe that
\begin{align*}
\|y_0+\T y\|&\geqslant \sup_{y^\ast\in A}|y^\ast(y_0)+\T y^\ast(y_1)|\\
&=\sup_{y^\ast\in A}|y^\ast(y_0)|+|y^\ast(y_1)|=1+\sup_{y^\ast\in A}|y^\ast(y_1)|=2.
\end{align*}
This,  together with $y_0$ being a WLUR point, gives that $y\in \T y_0$. Therefore, $Y$ is one-dimensional as desired.
\end{proof}

Our next result improves Proposition \ref{Prop:aDPDconvex-smooth} but only for lush operators. We do not know whether it is also true for operators with the aDP.

\begin{Prop} \label{Prop:realLush-strictlyconvex-smooth}
Let $X$, $Y$ be Banach spaces and let $G\in L(X,Y)$ be a norm-one operator which is lush. Then:
\begin{enumerate}
\item[(a)] If $X$ is strictly convex then $X=\K$.
\item[(b)] In the real case, if  $X$ is smooth then $X=\R$.
\item[(c)] If $Y$ is strictly convex then $Y=\K$.
\end{enumerate}
\end{Prop}

\begin{proof}
Given arbitrary separable subspaces $X_0\subset X$ and $Y_0\subset Y$, we can use Proposition~\ref{Prop:characterization-lushness}.(vi) to get the existence of separable subspaces $X_0\subset X_\infty\subset X$ and $Y_0\subset Y_\infty \subset Y$ such that $G(X_\infty)\subset Y_\infty $, $\|G|_{X_\infty}\|=1$, and  $G_\infty :=G|_{X_\infty}: X_\infty\longrightarrow Y_\infty$ is lush. Now Theorem~\ref{Theo-lush=casi-almost-CL} tells us that there exists a G$_\delta$ dense subset $\Omega$ of $(\ext B_{Y_\infty^{\ast}},\weakstar)$ such that $G_\infty^\ast(\Omega)\subset \Spear({X_\infty^\ast})$ or, equivalently, that
\begin{equation}\label{eq:Lush-strictconvex-smooth}
B_{X_\infty} = \overline{\aconv} \bigl(\Face(S_{X_\infty},G_\infty^\ast y^\ast)\bigr)
\end{equation}
for every $y^\ast\in \Omega$.

(a). If $X$ is strictly convex so is $X_\infty$, and then Proposition~\ref{Prop:spearVectorsProperties}.(h) tells us that $X_\infty$ is one-dimensional as $\Spear({X_\infty^\ast})$ is non-empty. Thus, $X_0$ is one-dimensional and its arbitrariness gives that $X$ is one-dimensional.

(b). If $X$ is smooth so is $X_\infty$. Using this time Proposition~\ref{Prop:spearVectorsProperties}.(l) we get that $X_\infty$ is one-dimensional as $\Spear({X_\infty^\ast})$ is non-empty. Therefore, $X_0$ is one-dimensional and its arbitrariness tells us that $X$ is one-dimensional.

(c). In this case we have that $Y_\infty$ is strictly convex. Observe that, fixed $y^\ast\in \Omega$, every element $x$ in the set $\Face(S_{X_\infty},G_\infty^\ast y^\ast)$ satisfies that $y^\ast(G_\infty x)=1$, so by the strict convexity of $Y_\infty$ the set $G_\infty\bigl(\Face(S_{X_\infty},G_\infty^\ast y^\ast)\bigr)$ must consist of one point. This, together with \eqref{eq:Lush-strictconvex-smooth}, implies that $G_\infty$ has rank one. Therefore, $\Spear(Y_\infty)$ is non-empty by Corollary~\ref{coro:charac-lush-rank-one} and so $Y_\infty$ (and thus $Y_0$) is one-dimensional by Proposition~\ref{Prop:spearVectorsProperties}.(h). The arbitrariness of $Y_0$ tells us that $Y$ is one-dimensional.
\end{proof}

Our last result in this chapter is an extension of Theorem \ref{Theo:spears_of_the_dual} to arbitrary lush operators: every lush operator attains its norm (i.e.\ the supremum defining its norm is actually a maximum).

\begin{Prop}\label{Prop:lush-attains-norm}
Let $X$, $Y$ be Banach spaces and let $G\in L(X,Y)$ be a norm-one operator.
If $G$ is lush, then it is norm-attaining. Actually,
$$
B_{X} = \overline{\conv}{\{ x \in S_{X} \colon \|Gx\| = 1 \}}.
$$
\end{Prop}

\begin{proof}
Fix an arbitrary $x_{0} \in B_{X}$. By Proposition \ref{Prop:characterization-lushness}, there are separable Banach spaces  $x_{0} \in X_{\infty} \subset X$ and $Y_{\infty} \subset Y$ satisfying that $G_\infty:=G|_{X_{\infty}}: X_{\infty} \longrightarrow Y_{\infty}$ is lush. Using Theorem \ref{Theo-lush=casi-almost-CL}, there exists $y_{0}^{\ast} \in S_{Y_{\infty}^{\ast}}$ such that $G_\infty^{\ast}y_{0}^{\ast}$ is a spear, so Theorem \ref{Theo:spears_of_the_dual} gives us that
\begin{equation*}
x_{0} \in B_{X_{\infty}} = \overline{\conv}{\{ x \in S_{X_\infty} \colon |G_\infty^{\ast}y_{0}^{\ast}(x)| = 1 \}} \subset  \overline{\conv}{\{ x \in S_{X} \colon \|Gx\| = 1 \}}.\qedhere
\end{equation*}
\end{proof}

We will see in Example \ref{example-aDP-no-NA} that the aDP is not enough to get norm-attainment.

\chapter{Lipschitz spear operators} \label{sec:Lipschitz}
Let $X$, $Y$ be Banach spaces. We denote by $\Lip_{0}{(X, Y)}$ the set of all Lipschitz mappings $F: X \longrightarrow Y$ such that $F(0)=0$. This is a Banach space when endowed with the norm
\[
\| F\|_L = \sup \left\{ \frac{\|F(x) - F(y)\|}{\| x - y\|}\colon x,y\in X,\, x \neq y \right\}.
\]
Observe that, clearly, $L(X,Y)\subset \Lip_0(X,Y)$ with equality of norms.

Our aim in this chapter is to study those elements of $\Lip_0(X,Y)$ which are spears. First, let us give a name for this.

\begin{Defi}
Let $X$, $Y$ be Banach spaces. A norm-one operator $G\in \Lip_0(X,Y)$ is a \emph{Lipschitz spear operator} if $\|G+\T\,F\|_L = 1 + \|F\|_L$ for every $F\in \Lip_0(X,Y)$.
\end{Defi}
\index{Lipschitz spear operator}%

We will prove here that every (linear) lush operator is a Lipschitz spear operator and present similar results for Daugavet centers and for operators with the aDP. To do so, we will use the technique of the Lipschitz-free space. We need some definitions and preliminary results. Let $X$ be a Banach space. Observe that we can associate to each $x \in X$ an element $\delta_{x} \in \Lip_{0}{(X,\K)}^{\ast}$ which is just the evaluation map $\delta_{x}(f) = f(x)$ for every $f\in \Lip_0(X,\K)$. The \emph{Lipschitz-free space} over $X$ is the Banach space
\index{Lipschitz-free space}%
\[
\FF(X):= \overline{\spn}^{\| \cdot \|}{\{ \delta_{x}\colon x \in X\}} \subset \Lip_{0}{(X, \K)}^{\ast}.
\]
It turns out that $\FF(X)$ is an isometric predual of $\Lip_{0}{(X, \K)}$ (which has been very recently shown to be the unique predual \cite{Weaver-Lip}). The map $\delta_X: x \longmapsto \delta_{x}$ establishes an isometric non-linear embedding $X \longrightarrow \FF(X)$ since $\|\delta_{x}-\delta_{y}\|_{\FF(X)}=\|x-y\|_X$ for all $x,y\in X$. The name Lipschitz-free space appeared for the first time in the paper \cite{GodKalt} by G.~Godefroy and N.~Kalton, but the concept was studied much earlier and it is also known as the Arens-Ells space of $X$ (see \cite[\S2.2]{Weaver}). The main features of the Lipschitz-free space which we are going to use here are contained in the following result. The first four assertions are nowadays considered folklore in the the theory of Lipschitz operators, and may be found in the cited paper \cite{GodKalt} (written for the real case, but also working in the complex case), section 2.2 of the book \cite{Weaver} by N.~Weaver, and Lemma 1.1 of \cite{JimenezSepulcreVillegas}. The fifth assertion was proved in \cite[Lemma 2.4]{LipschitzSlices}. For background on Lipschitz-free spaces we refer the reader to the already cited \cite{GodKalt,JimenezSepulcreVillegas,Weaver} and the very recent survey \cite{Godefroy-surveyFF} by G.~Godefroy.

\begin{Lemm}\label{Lemm:elementarypropertiesFF(X)}
Let $X$, $Y$ be Banach spaces.
\begin{enumerate}
\item[(a)] For every $F\in \Lip_0(X,Y)$, there exists a unique linear operator $\widehat{F}: \FF(X) \longrightarrow Y$ such that $\widehat{F}\circ \delta_X = F$ and $\|T_{F}\| = \|F\|_L$. Moreover, the application $F\longmapsto \widehat{F}$ is an isometric isomorphism from $\Lip_0(X,Y)$ onto $L(\FF(X),Y)$.
\item[(b)] There exists a norm-one $\K$-linear quotient map $\beta_X: \FF(X) \longrightarrow X$ which is a left inverse of $\delta_X$, that is, $\beta_X\circ \delta_X=\Id_X$. It is called the \emph{barycenter map} in \cite{GodKalt}, and is given by the formula
\index{barycenter map}%
\[
\beta_X\left(\sum_{x \in X}{a_{x} \delta_{x}}\right) = \sum_{x \in X}{a_{x}  x}.
\]
\item[(c)] From the uniqueness in item (a), it follows that $\widehat{F}=F\circ \beta_X$ for every $F\in L(X,Y)$.
\item[(d)] The set
\begin{equation*}
\mathcal{B}_X = \left\{ \frac{\delta_{x} - \delta_{y}}{\| x - y\|}\colon x,y\in X,\, x \neq y\right\}\subset \FF(X)
\end{equation*}
is norming for $\FF(X)^\ast=\Lip_0(X,\K)$, i.e.\ $B_{\FF(X)} = \overline{\aconv}{(\mathcal{B}_X)}$.
\item[(e)] Given $C\subset S_X$ and a slice $S$ of $\mathcal{B}_X$,
    $$
    \Bigl[\beta_X(S)\cap \overline{\conv}(C)\neq \emptyset\Bigr] \quad \Longrightarrow \quad \Bigr[\beta_X(S)\cap C\neq \emptyset\Bigr].
    $$
\end{enumerate}
\end{Lemm}

A comment on item (e) above could be clarifying. Let $X$ be a Banach space. As $\mathcal{B}_X\subset \FF(X)$ and $\FF(X)^\ast=\Lip_0(X,\K)$, a slice $\mathcal{S}$ of $\mathcal{B}_X$ has the form
$$
\mathcal{S}=\Slice(\mathcal{B}_X,f,\alpha)=\left\{\frac{\delta_x - \delta_y}{\|x-y\|}\colon x,y\in X,\, x \neq y,\, \re \left\langle f, \frac{\delta_x - \delta_y}{\|x-y\|}\right\rangle>1-\alpha\right\},
$$
where $f\in \Lip_0(X,\K)$ has norm one and $\alpha$ is a positive real number. Then, we have that
$$
\beta_X(\mathcal{S})=\left\{\frac{x - y}{\|x-y\|}\colon x,y\in X,\, x \neq y,\, \frac{\re f(x) - \re f(y)}{\|x-y\|}>1-\alpha\right\}
$$
is what is called in \cite{LipschitzSlices} a \emph{Lipschitz slice} of $S_X$. Then, item (e) above means that if a Lipschitz slice of $S_X$ does not intersect a subset $C\subset S_X$, then it does not intersect $\overline{\conv}(C)$ either.
\index{Lipschitz slice}%
This was proved in \cite[Lemma 2.4]{LipschitzSlices} with a completely elemental proof. Let us also say that assertion (e) is equivalent to the following fact \cite[Lemma~2.3]{BecerraLopezRueda-Lipslices}: \emph{given a Lipschitz slice $\beta_X(\mathcal{S})$ of $S_X$ and a point $x_0\in \beta_X(\mathcal{S})$, there is a linear slice $S$ of $S_X$ such that $x_0\in S\subseteq \beta_X(\mathcal{S})$} (indeed, one direction is obvious and for the non trivial one, let $C=S_X\setminus \beta_X(\mathcal{S})$ which clearly satisfies that $\beta_X(\mathcal{S})\cap C=\emptyset$; then, $\beta_X(\mathcal{S})\cap \overline{\conv}\,C=\emptyset$ and so the Hahn-Banach theorem gives the result). The proof of this last result given in \cite{BecerraLopezRueda-Lipslices} is independent of the above one and uses generalized derivatives and the Fundamental Theorem of Calculus for them.

The next one is the main result of this chapter. It is an application of our theory to Lipschitz-free spaces from which we will deduce the commented result about Lipschitz spear operators.

\begin{Theo}\label{Theorem-FreeLipchitz-lush}
Let $X$, $Y$ be Banach spaces and let $G\in L(X,Y)$ be a norm-one operator. If $G$ is lush, then $\widehat{G}:\FF(X)\longrightarrow Y$ is lush.
\end{Theo}

We need the following general technical result.

\begin{Lemm}\label{lemma:transfering-lushness-for-Lipschitz}
Let $X$, $Y$, $Z$ be Banach spaces and let $G_1\in L(Z,X)$ and $G_2\in L(X,Y)$ be norm-one operators. Suppose that there is a  subset $\mathcal{B}\subset B_Z$ norming for $Z^\ast$ (i.e.\ $\overline{\aconv}\,B=B_Z$) such that $G_1$ satisfies the following property
 \begin{equation}\label{equation:P2}\tag{P2}
 \begin{split}
&\text{For every slice $S$ of $\mathcal{B}$, every $x^{\ast} \in S_{X^\ast}$, and every $\eps>0$,}\\
&\bigl[G_1(S) \cap \overline{\conv}{\GS(S_{X},\T\,x^{\ast}, \eps)} \neq \emptyset\bigr] \, \Rightarrow \, \bigl[G_1(S) \cap \GS(S_{X},\T\,x^{\ast}, \eps) \neq \emptyset\bigr].
 \end{split}
\end{equation}
If $G_2$ is lush, then $G:=G_2\circ G_1$ is lush.
\end{Lemm}

\begin{proof}
Fix $z_0\in \mathcal{B}$, $y_0\in S_Y$, and $\eps>0$. As $G_2$ is lush, by Proposition \ref{Prop:characterization-lushness}.(v) we may find $y^\ast \in \ext(B_{Y^\ast})$ such that
$$
y_0\in \Slice(S_Y,y^\ast,\eps) \quad \text{ and } \quad G_1(z_0) \in \overline{\conv}\bigl(\GS(S_{X},\T\,G_2^{\ast}y^{\ast}, \eps)\bigr).
$$
Therefore, for every slice $S$ of $\mathcal{B}$ containing $z_0$ we have that
$$
G_1(S) \cap \overline{\conv}{\: \GS(S_{X},\T\,G_2^{\ast}y^{\ast}, \eps)} \neq \emptyset,
$$
and so \eqref{equation:P2} gives us that
$$
G_1(S) \cap \GS(S_{X},\T\,G_2^{\ast}y^{\ast}, \eps) \neq \emptyset.
$$
Therefore, we have that
$$
S\cap \GS(S_Z,\T\, G_1^\ast G_2^\ast y^\ast,\eps)\neq \emptyset.
$$
This has been proved for every slice $S$ of $\mathcal{B}$ containing $z_0\in \mathcal{B}$, but it is a fortiori also true for every slice $S$ of $S_Z$ containing $z_0$, so it follows that
$$
z_0\in \overline{\conv}{\: \GS(S_{Z},\T\,G^\ast y^{\ast}, \eps)}.
$$
As $\mathcal{B}$ is norming for $Z^\ast$, Proposition \ref{Prop:characterization-lushness}.(v) gives us the result.
\end{proof}

\begin{proof}[Proof of Theorem \ref{Theorem-FreeLipchitz-lush}]
By Lemma \ref{Lemm:elementarypropertiesFF(X)}.(e), it follows that $G_1:=\beta_X:\FF(X)\longrightarrow X$ satisfies condition \eqref{equation:P2} of Lemma \ref{lemma:transfering-lushness-for-Lipschitz}. As $G_2:=G:X\longrightarrow Y$ is lush, it follows from this lemma that $G_2\circ G_1:\FF(X)\longrightarrow Y$ is lush. But $G_2\circ G_1=G\circ \beta_X=\widehat{G}$ by Lemma \ref{Lemm:elementarypropertiesFF(X)}.(c).
\end{proof}

The identification of $L(\FF(X),Y)$ with $\Lip_0(X,Y)$ given in Lemma~\ref{Lemm:elementarypropertiesFF(X)}.(a) allows to deduce the promised result about Lipschitz spear operators from Theorem~\ref{Theorem-FreeLipchitz-lush}.

\begin{Coro}\label{Coro:lush=>Lipschitz-spear}
Let $X$, $Y$ be Banach spaces and let $G\in L(X,Y)$ be a norm-one operator. If $G$ is lush, then $G$ is a Lipschitz spear operator, i.e., $\|G+\T\,F\|_L=1+\|F\|_L$ for every $F\in \Lip_0(X,Y)$.
\end{Coro}

A first particular case of this result follows when we consider a lush Banach space $X$ and $G=\Id_X$. This result appeared previously in \cite{LipschitzSlices,Wang-Huang-Tan}

\begin{Coro}[\mbox{\cite[Theorem 4.1]{LipschitzSlices} and \cite[Theorem 2.6]{Wang-Huang-Tan}}]\label{Coro:lush=>Id--Lipschitz-spear}
Let $X$ be a lush Banach space. Then, $\Id_X$ is a Lipschitz spear operator, i.e.\ $\|\Id_X +\T\,F\|_L=1+\|F\|_L$ for every $F\in \Lip_0(X,Y)$.
\end{Coro}

As we commented, this result is already known, as it is contained in \cite[Theorem 2.6]{Wang-Huang-Tan} and \cite[Theorem 4.1]{LipschitzSlices}. But to get it from those references, the concept of Lipschitz numerical index of a Banach space is needed. Let $X$ be a Banach space. For $F\in \Lip_0(X,X)$, the \emph{Lipschitz numerical range} of $F$ \cite{Wang-Huang-Tan} is
\index{Lipschitz numerical range}%
$$
W_L(F):=\left\{\frac{\xi^\ast\bigl(Fx-Fy)}{\|x-y\|}\colon \xi^\ast\in S_{X^\ast},\, \xi^\ast(x-y)=\|x-y\|,\, x,y\in X,\, x\neq y\right\},
$$
the Lipschitz numerical radius of $F$ is just $w_L(F):=\sup\bigl\{|\lambda|\colon \lambda\in W_L(F)\bigr\}$, and the \emph{Lipschitz numerical index} of $X$ is
\index{Lipschitz numerical index}%
\begin{align*}
n_L(X)&:=\inf\bigl\{w_L(F)\colon F\in \Lip_0(X,X)\, \|F\|_L=1\bigr\} \\
&= \max\bigl\{k\geqslant 0\colon k\|F\|_L \leqslant w_L(F)\bigr\}.
\end{align*}
It is shown in \cite[Corollary 2.3]{Wang-Huang-Tan} that $\Id_X$ is a Lipschitz spear operator if and only if $n_L(X)=1$. With this in mind, Corollary \ref{Coro:lush=>Id--Lipschitz-spear} is just \cite[Theorem 2.6]{Wang-Huang-Tan} in the real case and \cite[Theorem 4.1]{LipschitzSlices} in the complex case. Let us comment that the main difficulty of the proofs in \cite{Wang-Huang-Tan} and \cite{LipschitzSlices} is to deal with Lipschitz operators. With our approach using the Lipschitz-free spaces, we avoid this.

Theorem \ref{Theorem-FreeLipchitz-lush} and Corollary \ref{Coro:lush=>Lipschitz-spear} apply to all the lush operators presented in this manuscript. We would like to emphasise the following two particular ones, which follow from Theorem \ref{Theo:Fourier-transform} and Corollary \ref{Coro:disk-algebra}, respectively.

\begin{Exam}  Let $H$ be a locally compact Abelian group and let $\Gamma$ be its dual group. Then, the Fourier transform $\mathcal{F}: L_{1}(H) \longrightarrow C_{0}(\Gamma)$ is a Lipschitz spear operator, that is,
$$
\|\mathcal{F}+\T\,F\|_L=1 + \|F\|_L
$$
for every $F\in \Lip_0(L_{1}(H), C_{0}(\Gamma))$.
\end{Exam}

\begin{Exam}
The inclusion $J: A(\mathbb{D}) \longrightarrow C(\T)$ is a Lipschitz spear operator, that is,
$$
\|J+\T\,F\|_L=1 + \|F\|_L
$$
for every $F\in \Lip_0(A(\mathbb{D}), C(\T))$.
\end{Exam}

The last consequence of Theorem \ref{Theorem-FreeLipchitz-lush} (actually, of Corollary \ref{Coro:lush=>Lipschitz-spear}) we would like to present here is the following.

\begin{Coro}
Let $X$ be a Banach space. Then,
$$
\Spear(X^\ast)\subset \Spear(\Lip_0(X,\K)).
$$
That is, every spear functional is actually a Lipschitz spear functional.
\end{Coro}

\begin{proof}
Let $x^\ast\in \Spear(X^\ast)$. It follows from Corollary \ref{Prop:characLushfinitedimCOdomain} that $g=x^*\in X^\ast\equiv L(X,\K)$ is lush. Then, Corollary \ref{Coro:lush=>Lipschitz-spear} gives that $\|g+\T\,f\|_L = 1 + \|f\|_L$ for every $f\in \Lip_0(X,\K)$, that is, $g\in \Lip_0(X,\K)$ as desired.
\end{proof}

We would like next to deal with operators with the aDP. The main result here is that we may extend the aDP of a (linear) operator to its linearization to the Lipschitz-free space.

\begin{Theo}\label{Theorem-FreeLipchitz-aDP}
Let $X$, $Y$ be Banach spaces and let $G\in L(X,Y)$ be a norm-one operator. If $G$ has the aDP, then $\widehat{G}:\FF(X)\longrightarrow Y$ has the aDP.
\end{Theo}

\begin{proof}
We fix $y_0\in S_Y$ and $\eps>0$. As $G$ has the aDP, we have that
$$
B_X=\overline{\conv}\,\{x\in S_X\colon \|Gx+\T\,y_0\|>2-\eps\}
$$
by Theorem \ref{Theo:aDPCharacterization}.(iii). Then, if $S$ is an arbitrary slice of $\mathcal{B}_X$, we obviously have that
$$
\beta_X(S)\cap \overline{\conv}\,\{x\in S_X\colon \|Gx+\T\,y_0\|>2-\eps\}\neq \emptyset,
$$
and so Lemma \ref{Lemm:elementarypropertiesFF(X)}.(e) gives us that \begin{equation}\label{eq:Lip-aDP-new}
\beta_X(S)\cap \{x\in S_X\colon \|Gx+\T\,y_0\|>2-\eps\}\neq \emptyset.
\end{equation}
Therefore, $S\cap \bigl\{\xi\in S_{\FF(X)}\colon \bigl\|\bigl[G\circ\beta_X\bigr](\xi) + \T\,y_0\bigr\|>2-\eps\bigr\}\neq \emptyset$, that is, using that $G\circ\beta_X=\widehat{G}$ by Lemma \ref{Lemm:elementarypropertiesFF(X)}.(c),
$$
S\cap \bigl\{\xi\in S_{\FF(X)}\colon \bigl\|\widehat{G}(\xi) + \T\,y_0\bigr\|>2-\eps\bigr\}\neq \emptyset.
$$
The arbitrariness of $S$ gives then that
\begin{equation}\label{eq:LIP-aDP-eq1}
\begin{aligned}
\mathcal{B}_X & \subseteq \overline{\conv} \bigl\{\xi\in S_{\FF(X)}\colon \|\widehat{G}(\xi) + \T\,y_0\|>2-\eps\bigr\} \\  & = \overline{\aconv} \bigl\{\xi\in S_{\FF(X)}\colon \|\widehat{G}(\xi) + \T\,y_0\|>2-\eps\bigr\}.
\end{aligned}
\end{equation}
As $B_{\FF(X)}=\overline{\aconv}\, \mathcal{B}_X$, we actually have that
\begin{equation}\label{eq:LIP-aDP-eq2}
\begin{aligned}
B_{\FF(X)}  =\overline{\aconv} \bigl\{\xi\in S_{\FF(X)}\colon \|\widehat{G}(\xi) + \T\,y_0\|>2-\eps\bigr\} \\  = \overline{\conv} \bigl\{\xi\in S_{\FF(X)}\colon \|\widehat{G}(\xi) + \T\,y_0\|>2-\eps\bigr\},
\end{aligned}
\end{equation}
and Theorem \ref{Theo:aDPCharacterization}.(iii) gives that $\widehat{G}$ has the aDP, as desired.
\end{proof}

The identification of $L(\FF(X),Y)$ with $\Lip_0(X,Y)$ given in Lemma~\ref{Lemm:elementarypropertiesFF(X)}.(a) allows to write Theorem \ref{Theorem-FreeLipchitz-aDP} in terms of the Lipschitz norm of Lipschitz operators. We need some preliminary work to write the results only in terms of the Lipschitz operators. Let $X$, $Y$ be Banach spaces. For  $F\in \Lip_0(X, Y)$, we define the \emph{slope} of $F$ \cite{LipschitzSlices} as the set
\index{slope}%
$$
\slope(F):=\left\{ \frac{F(x_1)-F(x_2) }{ \|x_1 - x_2\| }\colon
x_1\neq x_2 \in X\right\}.
$$
Observe that if $T\in L(X,Y)$, then $\slope(T)=T(S_X)$. On the other hand, it is clear that
$\slope(F)=\widehat{F}\bigl(\mathcal{B}_X\bigr)$ and, in particular,
$$
\overline{\aconv}\,\slope(F) =\overline{\aconv}\,\widehat{F}\bigl(\mathcal{B}_X\bigr) =\overline{\widehat{F}\bigl(B_{\FF(X)}\bigr)}.
$$
With this in mind, we get that if $G\in L(X,Y)$ has the aDP and $F\in \Lip_0(X,Y)$ satisfies that $\overline{\aconv}\,\slope(F)$ is SCD, then $\|\Id + \T\,F\|_L = 1 + \|F\|_L$ by Theorems \ref{Theorem-FreeLipchitz-aDP} and \ref{Theo:SCDOperators}. But, actually, we can go further and avoid to use the absolutely closed convex hull in the assumption.

\begin{Coro}\label{Coro-Lipshitz-aDP-slope-SCD}
Let $X$, $Y$ be Banach spaces and let $G\in \Lip_0(X,Y)$ be a norm-one operator with the aDP. If $F\in \Lip_0(X,Y)$ satisfies that $\slope(F)$ is SCD, then $\|G + \T\, F\|_L=1 + \|F\|_L$.
\end{Coro}

The result will follow from \eqref{eq:Lip-aDP-new} in the proof of Theorem \ref{Theorem-FreeLipchitz-aDP} and the following general result.

\begin{Lemm}
Let $X$, $Y$, $Z$ be Banach spaces, let $\mathcal{B}\subset B_X$ such that $\overline{\aconv}\,\mathcal{B}=B_X$, and let $G\in L(X,Y)$ be a norm-one operator such that $G(S)$ is a spear set for every slice $S$ of $\mathcal{B}$. Then, if $T\in L(X,Z)$ satisfies that $T(\mathcal{B})$ is SCD, then $T$ is a target for $G$. In the case of $Z=Y$, we have $\|G + \T\,T\|=1+\|T\|$.
\end{Lemm}

The proof of this lemma is an easy adaptation of the one of Theorem \ref{Theo:SCDOperators}.

\begin{proof}[Proof of Corollary \ref{Coro-Lipshitz-aDP-slope-SCD}] If $G$ has the aDP, it follows from \eqref{eq:Lip-aDP-new} in the proof of Theorem \ref{Theorem-FreeLipchitz-aDP} that $\widehat{G}:\FF(X)\longrightarrow Y$ satisfies the hypothesis of the above lemma with $\mathcal{B}=\mathcal{B}_X$. Now, if $\slope(F)=\widehat{F}(\mathcal{B}_X)$ is SCD, we have that $\|\widehat{G} + \T\,\widehat{F}\|=1 + \|\widehat{F}\|$, that is, $\|G + \T\,F\|_L = 1 + \|F\|_L$ as desired.
\end{proof}

If we apply this result in the particular case when $X=Y$ and $G=\Id_X$, we get the following result which already appeared in \cite{LipschitzSlices}.

\begin{Coro}[\mbox{\cite[Theorem 3.7]{LipschitzSlices}}] Let $X$ be a Banach space with the aDP. Then $\|\Id_X + \T\, F\|_L=1 + \|F\|_L$ for every $F\in \Lip_0(X,X)$ such that $\slope(F)$ is SCD.
\end{Coro}

Our next result extends Corollary \ref{Coro-Lipshitz-aDP-slope-SCD} to the non-separable case.

\begin{Coro}\label{Coro:Lip-aDP-every-separable-SCD}
Let $X$, $Y$ be Banach spaces and let $G\in \Lip_0(X,Y)$ be a norm-one operator with the aDP. If $F\in \Lip_0(X,Y)$ satisfies that for every separable subspace $X_0$ of $X$, $\slope(F|_{X_0})$ is SCD, then $\|G + \T\, F\|_L=1 + \|F\|_L$.
\end{Coro}

The result follows immediately from Corollary \ref{Coro-Lipshitz-aDP-slope-SCD} and the following lemma.

\begin{Lemm}\label{Lemma:Lipschitz-DaugavetCenters-separablecase}
Let $X$, $Y$ be Banach spaces and let $G\in L(X,Y)$ be an operator with the aDP. Given $F\in \Lip_0(X,Y)$, there are separable subspaces $X_\infty$ of $X$ and $Y_\infty$ of $Y$ such that $G|_{X_\infty}:X_\infty \longrightarrow Y_\infty$ has the aDP, $F(X_\infty)\subseteq Y_\infty$ and $\|F|_{X_\infty}\|_L=\|F\|_L$.
\end{Lemm}

\begin{proof}
Consider two sequences $(x_n)_{n\in \N}$ and $(y_n)_{n\in \N}$ such that $x_n\neq y_n$ for every $n\in \N$ and $\lim\frac{\|F(x_n)-F(y_n)\|}{\|x_n - y_n\|}=\|F\|_L$, let $X_0$ be the closed linear span in $X$ of the elements of the two sequences and let $Y_0$ be the closed linear span of $F(X_0)$. By Proposition \ref{Prop:separablyDetermined}, there are separable subspaces $X_1$ of $X$ and $Y_1$ of $Y$ such that $G|_{X_1}:X_1\longrightarrow Y_1$ has the aDP. By construction, we have that $\|F|_{X_1}\|_L=\|F\|_L$. Now, we may apply again Proposition \ref{Prop:separablyDetermined} starting with $X_1$ and the closed linear span of $Y_1\cup F(X_1)$ to get separable subspaces $X_2$, $Y_2$ such that $G|_{X_2}:X_2\longrightarrow Y_2$ has the aDP, $\|F_{X_2}\|_L=\|F\|_L$, and $F(X_1)\subseteq Y_2$. Repeating the process, it is straightforward to check that the separable subspaces $X_{\infty} := \overline{\bigcup_{n \in \N}{X_{n}}}$ of $X$ and $Y_{\infty}:= \overline{\bigcup_{n \in \N}{Y_{n}}}$ of $Y$ work.
\end{proof}

The main particular cases in which Corollary \ref{Coro:Lip-aDP-every-separable-SCD} apply are the following.

\begin{Coro}\label{Coro:Lipschiz-aDP-RNP...}
Let $X$, $Y$ be Banach spaces and let $G\in L(X,Y)$ be an operator with the aDP. If $F\in \Lip_0(X,Y)$ satisfies that $\overline{\conv}\slope(F)$ has the Radon-Nikod\'{y}m property, the convex point of continuity property or it is an Asplund set, or that $\slope(F)$ does not contain $\ell_1$-sequences, then $\|G + \T\, F\|_L=1 + \|F\|_L$.
\end{Coro}

\begin{proof}
If $\overline{\conv}\slope(F)$ has the Radon-Nikod\'{y}m property, or the convex point of continuity property or it is an Asplund set, then  $\overline{\conv}\slope(F|_{X_0})$ is SCD for every separable subspace $X_0$ of $X$ (use Example \ref{Exam:SCD-set-spaces-operators}). By \cite[Lemma 3.1]{LipschitzSlices}, it follows that $\slope(F|_{X_0})$ is SCD for every separable subspace $X_0$ of $X$, so Corollary \ref{Coro:Lip-aDP-every-separable-SCD} gives the result. If $\slope(F)$ does not contain $\ell_1$-sequences, neither does  $\slope(F|_{X_0})$ for every separable subspace $X_0$ of $X$. This gives that $\slope(F|_{X_0})$ is SCD for every separable subspace $X_0$ of $X$ (the proof of \cite[Theorem 2.22]{SCDsets} also works for not convex subsets) and then Corollary \ref{Coro:Lip-aDP-every-separable-SCD} gives the result.
\end{proof}

It is immediate that the above result applies to the recently introduced Lipschitz compact and Lipschitz weakly compact operators \cite[Definition 2.1]{JimenezSepulcreVillegas}: $F$ in $\Lip_0(X,Y)$ is Lipschitz compact (respectively, Lipschitz weakly compact) if $\slope(F)$ is relatively compact (respectively, relatively weakly compact).

The last aim in this chapter is to give for Daugavet centers analogous results to the ones we have for the aDP. In this case, we have to deal with the real version of the Lipschitz-free space. Given a (real or complex) Banach space $X$, we write $\FF_\R(X)$ to denote the Lipschitz-free space of $X$ in the sense of the real scalars, that is, $\FF_\R(X)$ it is the canonical predual of $\Lip_0(X,\R)$. If $X$ and $Y$ are real Banach spaces, nothing changes, but if they are complex spaces, we are considering only their real structure and so Lemma \ref{Lemm:elementarypropertiesFF(X)} is only valid for real scalars.

The main result for Daugavet centers is the following one.

\begin{Theo}\label{Theorem-FreeLipchitz-Daugavetcenters}
Let $X$, $Y$ be Banach spaces and let $G\in L(X,Y)$ be a norm-one operator. If $G$ is a Daugavet center, then $\widehat{G_\R}:\FF_\R(X)\longrightarrow Y_\R$ is a Daugavet center.
\end{Theo}

We need the following characterization of Daugavet centers which follows immediately from \cite{Bosenko-Kadets} and which will play the role of our Theorem \ref{Theo:aDPCharacterization}.(iii). In particular, it follows from it that to be a Daugavet center only depends on the real structure of the Banach spaces involved.

\begin{Lemm}[\mbox{see \cite[Theorem~2.1]{Bosenko-Kadets}}]
Let $X$, $Y$ be Banach spaces and let $G\in L(X,Y)$ be a norm-one operator. Then $G$ is a Daugavet center if and only if
$$
B_{X} = \overline{\conv}{\bigl( \{ x \in S_{X} \colon \|  Gx + y_{0}\| > 2 - \eps \}\bigr)}
$$
for every $y_0\in S_X$ and every $\eps>0$.
\end{Lemm}

\begin{proof}[Proof of Theorem \ref{Theorem-FreeLipchitz-Daugavetcenters}] We just have to adapt mutatis mutandis the proof of Theorem \ref{Theorem-FreeLipchitz-aDP}, using the above lemma instead of Theorem \ref{Theo:aDPCharacterization}.(iii). But observe that the game played in \eqref{eq:LIP-aDP-eq1} and \eqref{eq:LIP-aDP-eq2} is not valid here as the set $\{ \xi \in S_{\FF_\R(X)} \colon \|  \widehat{G}(\xi) + y_{0}\| > 2 - \eps \}$ is not rounded. At this point is where we have to go to the real version of the Lipschitz-free space, since the set $\mathcal{B}_X$ is clearly $\R$-rounded and then we actually have $B_{\FF_\R(X)}=\overline{\conv}\, \mathcal{B}_X$. With this in mind, everything works.
\end{proof}

Our next aim is to get consequences of Theorem \ref{Theorem-FreeLipchitz-Daugavetcenters} just in terms of the Lipschitz norm and the slope of Lipschitz operators.

\begin{Coro}\label{Coro-Lipshitz-DaugavetCenters-slope-SCD}
Let $X$, $Y$ be Banach spaces and let $G\in \Lip_0(X,Y)$ be a Daugavet center. If $F\in \Lip_0(X,Y)$ satisfies that $\slope(F)$ is SCD, then $\|G + F\|_L=1 + \|F\|_L$.
\end{Coro}

This result immediately follows from Theorem \ref{Theorem-FreeLipchitz-Daugavetcenters} and the next proposition which allows to pass the SCD from a set to its convex hull and which may have its own interest.

\begin{Prop}\label{Prop:SCDconvexhull}
Let $X$ be a Banach space and let $A\subset X$ be a bounded set. If $A$ is SCD, then $\overline{\conv}{(A)}$ is SCD.
\end{Prop}

We need a lemma which relates the slices of a set and the slices of its convex hull.

\begin{Lemm}
Let $X$ be a Banach space, let $A\subset X$ be a bounded set, let $x^\ast\in X^\ast$, and let $\eps,\delta\in \R^+$. Then
\begin{equation*}
\Slice(\conv{A},x^{\ast}, \eps \delta) \subset \conv \bigl(\Slice(A,x^{\ast}, \eps)\bigr) + \frac{2\delta - \delta^{2}}{1 - \delta}\|A\|  \cdot B_{X}.
\end{equation*}
\end{Lemm}

\begin{proof}
Every element of $\Slice(\conv{A},x^{\ast}, \eps \delta)$ is a (finite) convex combination of the form $\sum_{n \in \N}{\lambda_{n}a_{n}}$ where $a_n\in A$ for every $n$ and
\[
\sup_A \Real x^{\ast} - \eps \delta < \sum_{n}{\lambda_{n} \Real{x^{\ast}(a_{n})}}.
\]
Put $J := \left\{ n \in \N\colon  \sup_A \Real{x^{\ast}} - \eps < \Real{x^{\ast}(a_{n})} \right\}$. Then
\begin{align*}
\sup_A \Real{x^{\ast}} - \eps \delta & < \sum_{n \in J}{\lambda_{n} \Real{x^{\ast}(a_{n})}} +  \sum_{n \notin J}{\lambda_{n} \Real{x^{\ast}(a_{n})}}\\
& \leqslant \left( \sum_{n \in J}{\lambda_{n}} \right) \sup_A \Real{x^{\ast}} + \left( \sup_A \Real{x^{\ast}} - \eps \right) \sum_{n \notin J}{\lambda_{n}}\\
& = \sup_A \Real{x^{\ast}} - \eps \sum_{n \notin J}{\lambda_{n}}
\end{align*}
which allows to deduce that $\sum_{n \notin J}{\lambda_{n}} < \delta$. Since $a_{n} \in \Slice(A, x^{\ast}, \eps)$ for each $n \in J$, the result follows from the following estimation:
\begin{align*}
\left\| \sum_{n \in J}{\left( \frac{\lambda_{n}}{\sum_{n \in J}{\lambda_{n}}}\right) a_{n}} - \sum_{n \in \N}{\lambda_{n} a_{n}} \right\| & \leqslant \left\| \sum_{n \in J}{\left( \frac{\lambda_{n}}{\sum_{n \in J}{\lambda_{n}}} - 1 \right) a_{n}} \right\| + \left\| \sum_{n \notin J}{\lambda_{n} a_{n}} \right\|\\
 & \leqslant \left| \frac{1}{\sum_{n \in J}{\lambda_{n}}} - 1 \right| \cdot \left\| \sum_{n \in J}{\lambda_{n} a_{n}} \right\|  + \left\| \sum_{n \notin J}{\lambda_{n} a_{n}} \right\|\\
 & \leqslant \frac{\delta}{1 - \delta} \| A\| + \delta \| A\| = \frac{2 \delta - \delta^{2}}{1 - \delta} \| A\|.\qedhere
\end{align*}
\end{proof}

\begin{proof}[Proof of Proposition \ref{Prop:SCDconvexhull}]
First, by \cite[Remark 2.7]{SCDsets}, it is enough to show that $\conv(A)$ is SCD. Suppose that $\{\Slice(A,x_{n}^{\ast}, \eps_{n}) \colon n \in \N \}$ is a family of slices determining for $A$ and observe that we may and do assume that $\|x_n^\ast\|=1$ for every $n\in \N$. We consider for $\conv{(A)}$ the following (countable) family of slices:
\[ \mathcal{M} := \left\{ \Slice(\conv{(A)}, x^{\ast}_{n}, \eps_{n}/k)\colon n,k \in \N \right\}.   \]
Given any $x^\ast\in S_{X^\ast}$ and any $\eps>0$, we will show that $\Slice(\conv{(A)}, x^{\ast}, \eps)$ contains one element of $\mathcal{M}$, showing then that $\mathcal{M}$ is determining. Indeed, for $\Slice(A,x^{\ast},\eps/2)\subset A$ we know that there is $n_{0} \in \N$ such that $\Slice(A, x^{\ast}_{n_{0}}, \eps_{n_{0}}) \subset \Slice(A,x^{\ast},\eps/2)$. Taking $k \in \N$ big enough, we can assure, using the previous lemma, that
\begin{align*}
\Slice(\conv{(A)}, x^{\ast}_{n_{0}}, \eps_{n_{0}}/k) & \subset \conv\bigl(\Slice(A, x^{\ast}_{n_{0}}, \eps_{n_{0}})\bigr) + \frac{\eps}{2} B_{X}\\ & \subset \conv\bigl(\Slice(A, x^{\ast}, \eps/2)\bigr) + \frac{\eps}{2} B_{X}\\
& \subset \Slice(\conv{(A)},x^{\ast} , \eps/2) + \frac{\eps}{2} B_{X}.
\end{align*}
We can then conclude that
\begin{align*}
\Slice(\conv{(A)}, x^{\ast}_{n_{0}}, \eps_{n_{0}}/k) & \subset \Bigl[\Slice(\conv{(A)},x^{\ast} , \eps/2) + \frac{\eps}{2} B_{X}\Bigr]\cap \conv(A)\\ & \subset \Slice(\conv{(A)}, x^{\ast}, \eps).
\end{align*}
This shows that $\conv(A)$ is SCD, as desired.
\end{proof}

\begin{proof}[Proof of Corollary \ref{Coro-Lipshitz-DaugavetCenters-slope-SCD}]
If $G$ is a Daugavet center, then $\widehat{G_\R}:\FF_\R(X)\longrightarrow Y$ is also a Daugavet center by Theorem \ref{Theorem-FreeLipchitz-Daugavetcenters}. Now, if $F\in \Lip_0(X,Y)$ satisfies that $\slope(F)$ is SCD, so is $\widehat{F_\R}(B_{\FF_\R(X)}) =\overline{\conv}(\slope(F))$. Therefore, $\|\widehat{G_\R}+ \widehat{F}\|=1 + \|\widehat{F}\|$ by \cite[Corollary 1]{BosGnar}. Finally, this is equivalent to $\|G+F\|_L = 1 + \|F\|_L$ by Lemma \ref{Lemm:elementarypropertiesFF(X)}.
\end{proof}

We may extend the result to the non-separable case as we did for the aDP.

\begin{Coro}\label{Coro:Lip-DaugavetCenters-every-separable-SCD}
Let $X$, $Y$ be Banach spaces and let $G\in \Lip_0(X,Y)$ be a Daugavet center. If $F\in \Lip_0(X,Y)$ satisfies that for every separable subspace $X_0$ of $X$, $\slope(F|_{X_0})$ is SCD, then $\|G + F\|_L=1 + \|F\|_L$.
\end{Coro}

The result follows immediately from Corollary \ref{Coro-Lipshitz-DaugavetCenters-slope-SCD} and the following lemma which allows a reduction to the separable case and is completely analogous to Lemma \ref{Lemma:Lipschitz-DaugavetCenters-separablecase}. Its proof follows from \cite[Theorem 1]{Ivashyna-Daugavet_centers} in the same manner that the proof of Lemma \ref{Lemma:Lipschitz-DaugavetCenters-separablecase} follows from Proposition \ref{Prop:separablyDetermined}.

\begin{Lemm}
Let $X$, $Y$ be Banach spaces and let $G\in L(X,Y)$ be a Daugavet center. Given $F\in \Lip_0(X,Y)$, there are separable subspaces $X_\infty$ of $X$ and $Y_\infty$ of $Y$ such that $G|_{X_\infty}:X_\infty \longrightarrow Y_\infty$ is a Daugavet center, $F(X_\infty)\subseteq Y_\infty$ and $\|F|_{X_\infty}\|_L=\|F\|_L$.
\end{Lemm}

The most interesting particular cases of Corollary \ref{Coro:Lip-DaugavetCenters-every-separable-SCD} are summarized in the following result, whose proof is completely analogous to the one of Corollary \ref{Coro:Lipschiz-aDP-RNP...}.

\begin{Coro}\label{Coro:Lipschiz-DaugavetCenters-RNP...}
Let $X$, $Y$ be Banach spaces and let $G\in L(X,Y)$ be a Daugavet center. If $F\in \Lip_0(X,Y)$ satisfies that $\overline{\conv}\slope(F)$ has the Radon-Nikod\'{y}m property, the convex point of continuity property or it is an Asplund set, or that $\slope(F)$ does not contain $\ell_1$-sequences, then $\|G + F\|_L=1 + \|F\|_L$.
\end{Coro}

\chapter{Some stability results}\label{sec:stability}

Our aim here is to provide several results on the stability of our properties by several operations like absolute sums, vector-valued function spaces, and ultraproducts.

\section{Elementary results}

The first result shows that we may produce an injective operator with the aDP from any operator with the aDP.

\begin{Prop}\label{Prop:G-X/kerG-aDP}
Let $X$, $Y$ be Banach spaces and let $G\in L(X,Y)$ be an operator with the aDP and let $P:X\longrightarrow X/\ker G$ be the quotient map. Then, the quotient operator $\widetilde{G}\in L\big(X/\ker G, Y\big)$ satisfying $\widetilde{G}\circ P= G$ has the aDP.
\end{Prop}

\begin{proof}
By Theorem \ref{Theo:aDPCharacterization} it suffices to show that $\widetilde{G}(\widetilde{S})$ is a spear set in $Y$ for every slice $\widetilde{S}$ of $B_{X/\ker G}.$ So, we fix an arbitrary slice $\widetilde{S}$ of $B_{X/\ker G}$ and find $z^\ast\in S_{(X/\ker G)^\ast}$ and $\alpha>0$ such that
$$
\widetilde{S}=\bigl\{x+\ker G \in B_{X/\ker G}\colon \re z^\ast(x+\ker G)>1-\alpha\bigr\}.
$$
Since $P^\ast z^\ast\in S_{X^\ast}$, the set $S=\bigl\{x\in B_X \colon \re [P^\ast z^\ast](x)>1-\alpha\bigr\}$ is a slice of $B_X$. Observe that if $x\in S$ then $\re z^\ast(P(x))>1-\alpha$ and $P(x)\in B_{X/\ker G}$ which give that $P(x)\in \widetilde{S}$. Therefore, we have that $P(S)\subset \widetilde{S}$ and so
$$
G(S)=\widetilde{G}\big(P(S)\big)\subset \widetilde{G}(\widetilde{S}).
$$
Now, as $G(S)$ is a spear set by Theorem \ref{Theo:aDPCharacterization}, so is a fortiori $\widetilde{G}(\widetilde{S})$, as desired.
\end{proof}

The reciprocal result is not true: consider $G:\ell_2^2\longrightarrow \K$ given by $G(x,y)=x$ and observe that $\widetilde{G}\equiv \Id_\K$ is clearly lush, while $G$ does not even have the aDP (use Proposition \ref{Prop:characLushfinitedimDOMAIN} for instance).

Our next aim is to provide a way to extend the domain and the codomain keeping the properties of being spear, lush, or the aDP. The first result deals with extending the domain.

\begin{Prop}
Let $X$, $Y$, $Z$ be Banach spaces, let $G\in L(X,Y)$ be a norm-one operator, and consider the norm-one operator $\widetilde{G}:X\oplus_\infty Z\longrightarrow Y$ given by $\widetilde{G}(x,z)=G(x)$ for every $(x,z)\in X\oplus_\infty Z$. Then:
\begin{enumerate}
\item[(a)] if $G$ is a spear operator, so is $\widetilde{G}$;
\item[(b)] if $G$ has the aDP, so does $\widetilde{G}$;
\item[(c)] if $G$ is lush, so is $\widetilde{G}$.
\end{enumerate}
\end{Prop}

\begin{proof}
(a). Fix $T\in L(X\oplus_\infty Z,Y)$ with $\|T\|>0$ and $\|T\|>\eps>0$. Take $x_0 \in S_X$ and $z_0\in S_Z$ satisfying $\|T(x_0,z_0)\|>\|T\|-\eps$. Now pick $x^*\in S_{X^*}$ so that $x^*(x_0)=1$ and define the operator $S\in L(X,Y)$ by
$$
S(x)=T(x,x^*(x)z_0) \qquad (x\in X)
$$
which satisfies $S(x_0)=T(x_0,z_0)$ and so $\|S\|>\|T\|-\eps$. Now we can estimate as follows
\begin{align*}
\|\widetilde{G}+\T\, T\|&=\sup_{x\in B_X}\sup_{z\in B_Z}\|\widetilde G(x,z)+\T\, T(x,z)\|\\
&\geqslant \sup_{x\in B_X} \|\widetilde{G}(x,x^*(x)z_0)+\T\, T(x,x^*(x)z_0)\|\\
&=\sup_{x\in B_X}\|G(x)+\T S(x)\|=  \|G+\T S\|=1+\|S\|>1+\|T\|-\eps.
\end{align*}
The arbitrariness of $\eps$ gives $\|\widetilde{G}+\T\, T\|\geqslant 1+\|T\|$.

(b). Just observe that if $T$ is a rank-one operator in the argument above, then $S$ also is rank-one.

(c). Consider $(x_0,z_0)\in B_{X\oplus_\infty Z}=B_X \times B_Z$, $y\in S_Y$ and $\eps>0$. As $G$ is lush, Proposition \ref{Prop:characterization-lushness}.(iii) allows to find $y^\ast\in S_{Y^\ast}$ such that
$$
\re y^\ast(y)>1-\eps \quad \text{and} \quad \dist\bigl(x_{0}, \aconv(\GS(S_{X},G^{\ast}y^{\ast}, \eps))\bigr) < \eps.
$$
Now, observe that $\widetilde{G}^\ast y^\ast = (G^\ast y^\ast,0)\in [X\oplus_\infty Z]^\ast$ and it is then immediate that
\[
\dist\bigl((x_{0},z_0), \aconv(\GS(S_{X\oplus_\infty Z},\widetilde{G}^{\ast}y^{\ast}, \eps))\bigr) < \eps.
\]
Now, Proposition \ref{Prop:characterization-lushness}.(iii) gives that $\widetilde{G}$ is lush.
\end{proof}

We can get an analogous result to extend the codomain space.

\begin{Prop}\label{Prop:extending-codomain-L-summand}
Let $X$, $Y$, $Z$ be Banach spaces, let $G\in L(X,Y)$ be a norm-one operator, and consider the norm-one operator $\widetilde{G}:X\longrightarrow Y\oplus_1 Z$ given by $\widetilde{G}(x)=(G(x),0)$ for every $x\in X$. Then:
\begin{enumerate}
\item[(a)] if $G$ is a spear operator, so is $\widetilde{G}$;
\item[(b)] if $G$ has the aDP, so does $\widetilde{G}$;
\item[(c)] if $G$ is lush, so is $\widetilde{G}$.
\end{enumerate}
\end{Prop}

\begin{proof}
(a). Fix $T\in L(X,Y\oplus_1 Z)$ with $\|T\|>0$,  $\|T\|>\eps>0$, and $x_0\in S_X$ such that $\|Tx_0\|\geqslant \|T\|-\eps$. Denote by $P_Y$ and $P_Z$ the respective projections from $Y\oplus_1Z$ to $Y$ and $Z$. Take $z^*\in S_{Z^*}$ satisfying $z^*(P_ZTx_0)=\|P_ZTx_0\|$ and pick $y_0\in S_Y$ so that $P_YTx_0=\|P_YTx_0\|y_0$. Now define $S\in L(X,Y)$ by
$$
Sx=P_YTx+z^*(P_ZTx)y_0 \qquad (x\in X)
$$
which satisfies
\begin{align*}
\|S\|&\geqslant \|Sx_0\|=\big\|P_YTx_0+\|P_ZTx_0\|y_0\big\|=\|P_YTx_0\|+\|P_ZTx_0\|\geqslant \|T\|-\eps.
\end{align*}
Finally, using the triangle inequality and the fact that $G$ is a spear operator, we can estimate as follows:
\begin{align*}
\|\widetilde{G}+\T\, T\|
&=\sup_{x\in B_X}\|Gx+\T P_YTx\|+\|P_ZTx\|\\
&\geqslant\sup_{x\in B_X}\left\|Gx+\T \big(P_YTx + z^*(P_ZTx)y_0\big)\right\| \\
&=\sup_{x\in B_X}\|Gx+\T Sx\|=\|G+\T S \|=1+\|S\|\geqslant 1+\|T\|-\eps.
\end{align*}
The arbitrariness of $\eps$ finishes the proof.

(b). Observe that if $T$ is a rank-one operator in the argument above, then $S$ also is rank-one.

(c). Fix $x_0\in B_X$, $(y,z)\in S_{Y\oplus_1 Z}$ and $\eps>0$. As $G$ is lush, we may find $y^\ast\in S_{Y^\ast}$ such that $$
\re y^\ast(y)>\|y\|-\eps \quad \text{and} \quad \dist\bigl(x_{0}, \aconv(\GS(S_{X},G^{\ast}y^{\ast}, \eps))\bigr) < \eps.
$$
Indeed, if $y\neq 0$, apply Proposition \ref{Prop:characterization-lushness}.(iii) to $y/\|y\|$; if $y=0$, we may apply that proposition to any vector in $S_Y$.
Now, pick $z^\ast\in S_{Z^\ast}$ such that $z^\ast(z)=\|z\|$ and consider $(y^\ast,z^\ast)\in S_{Y^\ast}\times S_{Z^\ast}\subset S_{[Y\oplus_1 Z]^\ast}$. Observe that, on the one hand,
$$
\re [(y^\ast,z^\ast)](y,z)>\|y\|-\eps + \|z\| = 1 - \eps
$$
and, on the other hand, $\widetilde{G}(y^\ast,z^\ast)=G^\ast(y^\ast)$, so we have that
\[
\dist\bigl(x_{0}, \aconv(\GS(S_{X}, \widetilde{G}^{\ast}(y^{\ast},z^\ast), \eps))\bigr) < \eps.
\]
We then get that $\widetilde{G}$ is lush by Proposition \ref{Prop:characterization-lushness}.(iii).
\end{proof}

\section{Absolute sums}

We show in this section the stability of our properties by $c_0$, $\ell_1$, and $\ell_\infty$ sums of Banach spaces. The following result borrows the ideas from \cite[Proposition~1]{MartinPaya-CLspaces} and \cite[\S5]{LushNumerical}.

\begin{Prop}\label{prop:stability-c0l1ellinfty-sums}
Let $\{X_\lambda\colon \lambda\in \Lambda\}$, $\{Y_\lambda\colon \lambda\in \Lambda\}$ be two families of Banach spaces and let $G_\lambda\in L(X_\lambda,Y_\lambda)$ be a norm-one operator for every $\lambda\in \Lambda$. Let $E$ be one of the Banach spaces $c_0$, $\ell_\infty$, or $\ell_1$, let $X=\left[\bigoplus_{\lambda\in\Lambda}X_\lambda\right]_{E}$ and $Y=\left[\bigoplus_{\lambda\in\Lambda}Y_\lambda\right]_{E}$, and define the operator $G:X \longrightarrow Y$  by $G\bigl[(x_\lambda)_{\lambda\in\Lambda}\bigr]=(G_\lambda x_\lambda)_{\lambda\in\Lambda}$ for every $(x_\lambda)_{\lambda\in\Lambda}\in \left[\bigoplus_{\lambda\in\Lambda}X_\lambda\right]_{E}$. Then
\begin{enumerate}
\item[(a)] $G$ is a spear operator if and only if $G_\lambda$ is a spear operator for every $\lambda\in \Lambda$;
\item[(b)] $G$ has the aDP if and only if $G_\lambda$ has the aDP for every $\lambda\in \Lambda$;
\item[(c)] $G$ is lush if and only if $G_\lambda$ is lush for every $\lambda\in \Lambda$.
\end{enumerate}
\end{Prop}

\begin{proof}
(a). We suppose first that $G$ is a spear operator and, fixed $\kappa\in\Lambda$, we have to show that $G_\kappa$ is  a spear operator. Observe that calling $W=\left[\bigoplus_{\lambda\neq \kappa}X_\lambda\right]_{E}$ and $Y=\left[\bigoplus_{\lambda\neq \kappa}Y_\lambda\right]_{E}$, we can write $X=X_\kappa\oplus_\infty W$ and $Y=Y_\kappa\oplus_\infty Z$ when $E$ is $\ell_\infty$ or $c_0$ and $X=X_\kappa\oplus_1 W$ and $Y=Y_\kappa\oplus_1 Z$ when $E$ is $\ell_1$. Given a non-zero operator $T_\kappa \in L(X_\kappa,Y_\kappa)$, define $T\in L(X,Y)$ by $T(x_\kappa, w)=(T_\kappa x_\kappa, 0)$ which obviously satisfies $\|T\|=\|T_\kappa\|$. Let $P_\kappa$ and $P_Z$ denote  the projections from $Y$ onto $Y_\kappa$ and $Z$ respectively. When $E$ is $\ell_\infty$ or $c_0$ we can write
\begin{align*}
1+\|T_\kappa\|&=1+\|T\|=\|G+\T\, T\|=\sup_{(x_\kappa,w)\in B_X}\| G(x_\kappa,w)+\T \, T(x_\kappa,w)\|\\
&\overset{(\ast)}{=}\sup_{(x_\kappa,w)\in B_X} \max\big \{\|P_\kappa G(x_\kappa,w)+\T \,P_\kappa T(x_\kappa,w)\|,\,\\
&\hspace{4.5cm}\| P_ZG(x_\kappa,w)+\T \, P_ZT(x_\kappa,w)\|\big\}\\
&= \max\left \{\|G_\kappa +\T \,T_\kappa \|,\, \sup_{(x_\kappa,w)\in B_X} \| P_ZG(x_\kappa,w)\|\right\}\\
&\leqslant \max\left \{\|G_\kappa +\T \,T_\kappa \|,\, \|G\|\right\}.
\end{align*}
Since $\|G\|=1$, it follows that $1+\|T_\kappa\|\leqslant \|G_\kappa +\T\, T_\kappa\|$ and so $G_\kappa$ is a spear operator. When $E$ is $\ell_1$ equality $(\ast)$ can be continued as follows
\begin{align*}
1+\|T_\kappa\|&\overset{(\ast)}{=}\sup_{(x_\kappa,w)\in B_X}\|P_\kappa G(x_\kappa,w)+\T \,P_\kappa T(x_\kappa,w)\|\\&\hspace{3.5cm}+\| P_ZG(x_\kappa,w)+\T \, P_ZT(x_\kappa,w)\|\\
&=\sup_{(x_\kappa,w)\in B_X}\|G_\kappa x_\kappa+\T \,T_\kappa x_\kappa\|+\|P_ZG(0,w)\|\\
&\leqslant \sup_{(x_\kappa,w)\in B_X}\|G_\kappa +\T \,T_\kappa \|\|x_\kappa\| +\|G\|\|w\|=\max\left \{\|G_\kappa +\T \,T_\kappa \|,\, \|G\|\right\}.
\end{align*}
We prove now the sufficiency when $E$ is $\ell_\infty$ or $c_0$. Given an operator $T\in L(X,Y)$ and fixed $\eps>0$, we find $\kappa\in\Lambda$ such that $\|P_\kappa T\|> \|T\|-\eps$ and write $X=X_\kappa\oplus_\infty W$ where $W=\left[\bigoplus_{\lambda\neq \kappa}X_\lambda\right]_{E}$. Since $B_X$ is the convex hull of $S_{X_\kappa}\times S_W$ we may find $x_0\in S_{X_\kappa}$ and $w_0\in S_W$ such that
$$
\|P_\kappa T(x_0,w_0)\|>\|T\|-\eps.
$$
Now fix $x^\ast\in S_{X^\ast_\kappa}$ with $x^\ast(x_0)=1$ and define the operator $S\in L(X_\kappa,Y_\kappa)$ given by
$$
S(x)=P_\kappa T(x,x^\ast(x)w_0) \qquad (x\in X_\kappa)
$$
which satisfies $\|S\|\geqslant \|Sx_0\|=\|P_\kappa T(x_0, w_0)\|>\|T\|-\eps$. Observe finally that
\begin{align*}
\|G+\T\, T\|&\geqslant \|P_\kappa G+ \T P_\kappa T\| \geqslant \sup_{x\in X_\kappa}\|[P_\kappa G](x,x^\ast(x)w_0)+ \T [P_\kappa T](x,x^\ast(x)w_0)\|\\
&=\sup_{x\in X_\kappa}\|G_\kappa(x)+\T S(x)\|=\|G_\kappa+ \T S\|=1+\|S\|>1+\|T\|-\eps.
\end{align*}
So, the arbitrariness of $\eps$ gives that $\|G+\T\, T\|\geqslant 1+\|T\|$, finishing the proof for $E=c_0, \ell_\infty$.

Suppose now that $E=\ell_1$. Fix an operator $T\in L(X,Y)$ and observe that it may be seen as a family $(T_\lambda)_{\lambda\in\Lambda}$ of operators where $T_\lambda
\in L(X_\lambda, Y)$ for every $\lambda\in \Lambda$, and
$\|T\|=\sup_\lambda \|T_\lambda\|$. Given $\eps>0$, find
$\kappa\in\Lambda$ such that $\|T_{\kappa}\| > \|T\|-\eps$,
and write $X=X_\kappa\oplus_1 W$, $Y=Y_{\kappa}\oplus_1 Z$, and $T_{\kappa}=(A,B)$
where $A\in L(X_{\kappa},Y_\kappa)$ and $B\in L(X_{\kappa},Z)$. Now
we choose $x_0\in S_{X_{\kappa}}$ such that
$$
\|T_{\kappa} x_0\|= \|A x_0\| + \|B x_0\| > \|T\| -\eps,
$$
we find $a_0\in S_{Y_{\kappa}}$, $z^\ast\in S_{Z^\ast}$ satisfying
$$
\|Ax_0\| a_0 = Ax_0 \qquad \text{and} \qquad z^\ast(Bx_0)=\|Bx_0\|,
$$
and define the operator $S\in L(X_{\kappa},Y_\kappa)$ by
$$
S x = Ax + z^\ast(Bx)a_0 \qquad
(x\in X_{\kappa}).
$$
Then
$$
\|S\|\geqslant \|S x_0\| = \Bigl\|Ax_0 + \|Bx_0\| a_0\Bigr\| = \|Ax_0\|+ \|Bx_0\|> \|T\|-\eps.
$$
Moreover, since $G_\kappa$ is a spear operator, fixed $\eps>0$ we may find $x_\kappa\in S_{X_\kappa}$ and $y_\kappa^\ast\in S_{Y_\kappa^\ast}$ such that
$$
\left|y^\ast_\kappa(G_\kappa x_\kappa +\T S x_\kappa)\right|=\left\|G_\kappa x_\kappa +\T S x_\kappa\right\|\geqslant 1+\|S\|-\eps.
$$
Now take $x=(x_\kappa,0)\in S_X$ and $y^\ast=(y^\ast_\kappa, y^\ast_\kappa(a_0)z^\ast)\in S_{Y^\ast}$, and observe that
\begin{align*}
\|G+\T\, T\|&\geqslant  \left|y^\ast(G x +\T T x)\right|=\bigl|y^\ast_\kappa(G_\kappa x_\kappa) +\T[y^\ast_\kappa(Ax_\kappa)+ y^\ast_\kappa(a_0)z^\ast(Bx_\kappa)]\bigr|\\
&=\left|y^\ast_\kappa(G_\kappa x_\kappa +\T S x_\kappa)\right|=\left\|G_\kappa x_\kappa +\T S x_\kappa\right\|\\
&\geqslant 1+\|S\|-\eps\geqslant 1+\|T\|-2\eps.
\end{align*}
So, the arbitrariness of $\eps$ gives that $\|G+\T\, T\|\geqslant 1+\|T\|$, finishing the proof for $E=\ell_1$.

(b). For the aDP, the arguments above apply just taking into account that when one starts with rank-one operators, the constructed operators are also rank-one.

(c). We assume first that $G$ is lush. Fixed $\kappa\in\Lambda$, $x_\kappa\in S_{X_\kappa}$, $y_\kappa\in S_{Y_\kappa}$, $\eps>0$ we consider the elements $(z_\lambda)_{\lambda\in \Lambda}\in B_X$ and  $(w_\lambda)_{\lambda\in \Lambda}\in S_{Y}$ given by
$$
z_\lambda=0, \quad  w_\lambda=0 \quad \text{for}\quad \lambda\neq \kappa \qquad \text{and}\qquad z_\kappa=x_\kappa, \quad w_\kappa=y_\kappa.
$$
Now Proposition \ref{Prop:characterization-lushness}.(iii) provides with $y^\ast\in \Slice(B_{Y^\ast},(w_\lambda)_{\lambda\in \Lambda},\eps)$ such that
\begin{equation}\label{eq:lushness-ell-1-sums}
\dist\bigl((z_\lambda)_{\lambda\in \Lambda}, \aconv(\GS(S_{X},G^{\ast} y^{\ast}, \eps))\bigr)<\eps^2.
\end{equation}
From this point we have to distinguish two cases depending on the space $E$. Suppose first that $E=c_0$ or $E=\ell_\infty$ and
observe that $y^\ast|_{Y_\kappa}\in \Slice(B_{Y^\ast_\kappa},y_\kappa,\eps)$. In this case, given  $(\widetilde{z}_\lambda)_{\lambda\in \Lambda}\in \GS(S_{X},G^{\ast} y^{\ast}, \eps)$, it follows that $\widetilde{z}_\kappa\in \GS(S_{X_\kappa},G^{\ast}_\kappa (y^\ast|_{Y_\kappa}),2 \eps)$ which, together with \eqref{eq:lushness-ell-1-sums}, allows us to deduce that
$$
\dist\bigl(x_\kappa, \aconv(\GS(S_{X_\kappa},G^{\ast}_\kappa (y^\ast|_{Y_\kappa}),2 \eps))\bigr)<\eps.
$$

We consider now the more bulky case in which $E=\ell_1$. Using \eqref{eq:lushness-ell-1-sums}
we can find scalars $\lambda_{i} \in \mathbb{K}$ with $\sum_{i=1}^{n}{|\lambda_{i}|} = 1$ and elements $x^{i} \in gS(S_{X}, G^{\ast}y^{\ast}, \eps)$ such that
$$
\left\| x_{\kappa} - \sum_{i=1}^{n}{\lambda_{i} x_{\kappa}^{i}} \right\|  \leqslant   \sum_{\lambda \in \Lambda}{\left\| z_{\lambda} - \sum_{i=1}^{n}{\lambda_{i} x^{i}_{\lambda}} \right\|_{X_{\lambda}}} < \eps^2.
$$
Since $\| x_{\kappa}\| = 1$, we deduce that $1-\eps^2\leqslant \sum_{i=1}^{n}{|\lambda_{i}| \|x_{\kappa}^{i}\|}$ so Lemma~\ref{Lemm:stabilityL1Technical} tells us that the set $I:=\{ i  \colon \| x_{\kappa}^{i} \| > 1 - \eps\}$ satisfies that $\sum_{i \in I}{|\lambda_{i}|} > 1- \eps$. Hence
\begin{equation}\label{eq:lushness-ell-1-sums-2}
\left\| x_{\kappa} - \sum_{i \in I}{\lambda_{i} x_{\kappa}^{i}}\right\| < 2 \eps.
\end{equation}
But every $i \in I$ satisfies that
\begin{equation*}
1 - \eps < \sum_{\lambda \in \Lambda}{\Real y^{\ast}|_{Y_{\lambda}}(G_{\lambda} x^{i}_{\lambda})} \leqslant \Real{y^{\ast}|_{Y_{\kappa}}(G_{\kappa} x^{j}_{\kappa})} + \sum_{\lambda \neq \kappa}{\| x_{\lambda}^{i}\|}\\
<  \Real{y^{\ast}|_{Y_{\kappa}}(G_{\kappa} x^{i}_{\kappa})} + \eps,
\end{equation*}
from where it follows that
$$
x_{\kappa}^{i} \in gS(B_{X}, G_{\kappa}^{\ast}(y^{\ast}|_{Y_{\kappa}}), 2 \eps)
$$
for each $i \in I$. This, together with \eqref{eq:lushness-ell-1-sums-2}, tells us that
\[
\dist(x_{\kappa}, \aconv gS(B_{X}, G_{\kappa}^{\ast}(y^{\ast}|_{Y_{\kappa}}), 2 \eps)) < 2 \eps,
\]
finishing the proof of the necessity for $E=\ell_1$.

Let us prove the sufficiency when $E=\ell_\infty$ or $E=c_0$. Fixed $(x_\lambda)_{\lambda\in \Lambda}\in B_X$, $(y_\lambda)_{\lambda\in \Lambda}\in S_Y$, and $\eps>0$, there is $\kappa\in \Lambda$ such that $\|y_\kappa\|>1-\eps$. Using that $G_\kappa$ is lush we may find $y^\ast_\kappa\in \Slice(B_{Y^\ast_\kappa},y_\kappa,\eps)$ satisfying
$$
\dist\bigl(x_\kappa, \aconv(\GS(S_{X_\kappa},G^{\ast}_\kappa y^{\ast}_\kappa, \eps))\bigr)<\eps.
$$
Defining $y^\ast\in S_{Y^\ast}$ by $y^\ast[(z_\lambda)_{\lambda\in \Lambda}]=y^\ast_\kappa(z_\kappa)$ for every $(z_\lambda)_{\lambda\in \Lambda}\in X$, we clearly have $y^\ast\in \Slice(B_{Y^\ast},y,\eps)$. Observe that, fixed $\widetilde{x}_\kappa\in \GS(S_{X_\kappa},G^{\ast}_\kappa y^{\ast}_\kappa,\eps)$ and $\theta \in \T$, the element $(z_\lambda)_{\lambda\in \Lambda}\in B_X$ given by $z_\lambda=\overline{\theta}x_\lambda$ for $\lambda\neq\kappa$ and $z_\kappa=\widetilde{x}_\kappa$ belongs to $\GS(S_{X},G^{\ast} y^{\ast}, \eps)$.
Using this it is easy to deduce that
$$
\dist\bigl((x_\lambda)_{\lambda\in \Lambda}, \aconv(\GS(S_{X},G^{\ast} y^{\ast}, \eps))\bigr)<\eps,
$$
which tells us that $G$ is lush by Proposition \ref{Prop:characterization-lushness}.(iii).

Suppose that $E=\ell_1$. We take the  set $\mathcal{B}=\bigl\{(x_\lambda)_{\lambda\in \Lambda}\in B_X \colon \#\supp (x_\lambda)=1\bigr\}$ which is norming for $X^\ast$. Fixed $(x_\lambda)_{\lambda\in \Lambda}\in \mathcal{B}$, there is $\kappa\in\Lambda$ so that $x_\kappa\in B_{X_\kappa}$
and $x_\lambda=0$ for $\lambda\neq\kappa$. Given $(y_\lambda)_{\lambda\in \Lambda}\in S_Y$ and $\eps>0$, we may and do assume that $y_\kappa\neq0$ and, since $G_\kappa$ is lush, we may use Proposition \ref{Prop:characterization-lushness}.(iii) for $x_\kappa\in B_{X_\kappa}$ and $\frac{y_\kappa}{\|y_\kappa\|}\in S_{Y_\kappa}$ to find $y^\ast_\kappa\in \Slice(B_{Y^\ast_\kappa},\frac{y_\kappa}{\|y_\kappa\|},\eps)$ such that
$$
\dist\bigl(x_\kappa, \aconv(\GS(S_{X_\kappa},G^{\ast}_\kappa y^{\ast}_\kappa, \eps))\bigr)<\eps.
$$
For each $\lambda\in \supp (y_\lambda)\setminus\{\kappa\}$ we take $y^\ast_\lambda\in S_{Y^\ast_\lambda}$ satisfying $y^\ast_\lambda(y_\lambda)=\|y_\lambda\|$ and we define $y^\ast\in S_{Y^\ast}$ by
$$
y^\ast\bigl[(w_\lambda)_{\lambda\in \Lambda}\bigr]=\sum_{\lambda\in \supp{(y_\lambda)}}y^\ast_\lambda(w_\lambda) \qquad \bigl((w_\lambda)\in Y\bigr).
$$
Then it obviously follows that $y^\ast\in  \Slice(B_{Y^\ast},(y_\lambda),\eps) $ and, thanks to the shape of $(x_\lambda)_{\lambda\in \Lambda}$, one can easily deduce that
$$
\dist\bigl((x_\lambda)_{\lambda\in \Lambda}, \aconv(\GS(S_{X},G^{\ast} y^{\ast}, \eps))\bigr)<\eps.
$$
This finishes the proof by using Proposition \ref{Prop:characterization-lushness}.(iii) since $\overline{\aconv}\,\mathcal{B}= B_X$.
\end{proof}

\section{Vector-valued function spaces}
Our next aim is to present several results concerning the behaviour of our properties for vector-valued function spaces. We start analysing the situation for spaces of continuous functions.

\begin{Theo}\label{Theo:stabilitityC(K)}
Let $X,Y$ be Banach spaces, let $K$ be a compact Hausdorff topological space and let $G\in L(X,Y)$ be a norm-one operator. Consider the norm-one composition operator $\widetilde{G}:C(K,X) \longrightarrow C(K,Y)$ given by $\widetilde{G}(f)=G\circ f$ for every $f\in C(K,X)$. Then:
\begin{enumerate}
\item[(a)] $\widetilde{G}$ is a spear operator if and only if $G$ is a spear operator.
\item[(b)] $\widetilde{G}$ is lush if and only if $G$ is lush.
\item[(c)] If $K$ contains isolated points, then $\widetilde{G}$ has the aDP if and only if $G$ does.
\item[(d)] If $K$ is perfect, then $\widetilde{G}$ has the aDP if and only if $G(B_X)$ is a spear set.
\end{enumerate}
\end{Theo}

\begin{Rema}
All the information given in the above result was previously known for the case of the identity (see  \cite{LushNumOneDual,martinOikhberg,M-P}).
\end{Rema}

\begin{proof}[Proof of Theorem~\ref{Theo:stabilitityC(K)}.(a)] This is an easy adaptation of \cite[Theorem~5]{M-P}.
Suppose first that $G$ is a  spear operator. Fixed $T\in L(C(K,X),C(K,Y))$ with $\|T\|=1$ and $\eps>0$, find $f_0\in C(K,X)$ with $\|f_0\|=1$ and $t_0\in K$ such that
\begin{equation}\label{eq1-cdek}
\|[Tf_0](t_0)\|>1-\eps.
\end{equation}
Define $z_0=f_0(t_0)$ and find a
continuous function $\varphi: K
\longrightarrow [0,1]$ such that $\varphi(t_0)=1$ and $\varphi(t)=0$ if
$\|f_0(t)-z_0\|\geqslant \eps$. Now write $z_0=\lambda x_1 +
(1-\lambda) x_2$ with $0\leqslant \lambda \leqslant 1$, $x_1,x_2\in S_X$,
and consider the functions
$$
f_j=(1-\varphi)f_0 + \varphi x_j\in C(K,X) \qquad (j=1,2).
$$
Then $\|\varphi f_0 - \varphi z_0\|<\eps$ meaning that
$$
\|f_0 - (\lambda f_1 + (1-\lambda)f_2)\|<\eps,
$$
and, using (\ref{eq1-cdek}), we must have
\begin{equation*} 
\|[Tf_1](t_0)\|>1-2\eps \qquad \text{or} \qquad
\|[Tf_2](t_0)\|>1-2\eps.
\end{equation*}
By making the right choice of $x_0=x_1$ or $x_0=x_2$ we get
$x_0\in S_X$ such that
\begin{equation}\label{eq3-cdek}
\|\left[T\left( (1-\varphi)f_0+\varphi
x_0\right)\right](t_0)\|>1-2\eps .
\end{equation}
Next, we fix $x_0^\ast\in S_{X^\ast}$ with $x_0^\ast(x_0)=1$, denote
$$
\Phi(x)=x_0^\ast(x)(1-\varphi)f_0 + \varphi x\in C(K,X) \qquad
(x\in X),
$$
and consider the operator $S\in L(X)$ given by
$$
Sx = [T(\Phi(x))](t_0) \qquad (x\in X)
$$
which, by (\ref{eq3-cdek}), obviously satisfies $\|S\|\geqslant \|S x_0\|>1-2\eps$. Now, we use that $G$ is a spear operator to find $x\in S_X$ satisfying $\|Gx+\T Sx\|>1+\|S\|-\eps$, and observe that
\begin{align*}
\|\widetilde{G}+\T\, T \|\geqslant \left\|\big[(\widetilde{G}+\T \,T)(\Phi(x))\big](t_0)\right\|=\|Gx+\T S(x)\|>1+\|S\|-\eps>2-3\eps.
\end{align*}
The arbitrariness of $\eps$ gives that $\|\widetilde{G}+\T\, T \|\geqslant 2$ and so $\widetilde{G}$ is a spear operator.

Suppose conversely that $\widetilde{G}$ is a spear operator. Fix $S\in L(X,Y)$, $\eps>0$ and define the operator $T\in L\big(C(K,X),C(K,Y)\big)$ by
$$
[T(f)](t)=S(f(t)) \qquad (t\in K, f\in C(K,X))
$$
which satisfies $\|T\|=\|S\|$. Since $\widetilde{G}$ is a spear operator we may find $f_0\in C(K,X)$ and $t_0\in K$ such that
$$
\left\|\left[(\widetilde{G}+\T \, T)(f_0)\right](t_0)\right\|>1+\|T\|-\eps
$$
and we can write
\begin{align*}
1+\|S\|-\eps&=1+\|T\|-\eps<\left\|\left[(\widetilde{G}+\T \, T)(f_0)\right](t_0)\right\|\\
&=\|G(f_0(t_0))+\T S(f_0(t_0))\|\leqslant \|G+\T S\|.
\end{align*}
The arbitrariness of $\eps$ tells us that $G$ is a spear operator.
\end{proof}

We next deal with lushness for spaces of vector-valued continuous functions.

\begin{proof}[Proof of Theorem~\ref{Theo:stabilitityC(K)}.(b)] Suppose that $G$ is lush and let us show that $\widetilde{G}$ is lush. This part of the proof is an easy adaptation of \cite[Proposition~5.1]{LushNumOneDual}. Let $f\in S_{C(K,Y)}$, $g\in S_{C(K,X)}$, and $\eps>0$ be fixed. Then, we take $t_0\in K$ with $\|f(t_0)\|=1$ and, using that $G$ is lush together with Proposition~\ref{Prop:characterization-lushness}.(iii), we find $y^\ast\in \Slice(B_{Y^\ast},f(t_0),\eps)$ such that
$$
\dist\bigl(g(t_0), \aconv(\GS(S_{X},G^{\ast}y^{\ast}, \eps))\bigr) < \frac{\eps}{2}.
$$
So, there are $\theta_1,\ldots,\theta_n\in \T$, $\lambda_1,\ldots,\lambda_n\in [0,1]$ with $\sum_{k=1}^n\lambda_k=1$, and $x_1,\ldots,x_n\in \GS(S_{X},G^{\ast}y^{\ast}, \eps)$ such that
$$
\left\|g(t_0)-\sum_{k=1}^n\lambda_k\theta_kx_k\right\|<\frac{\eps}{2}.
$$
Next, we take an open set  $U\subset K$ such that $t_0\in U$ and
$$
\|g(t)-g(t_0)\|<\frac{\eps}{2} \qquad (t\in U),
$$
and we fix a continuous function $\varphi:K\longrightarrow [0,1]$ with
$\varphi(t_0)=0$ and $\varphi|_{K\setminus U}\equiv 1$. Now we consider the functional $\xi^\ast\in B_{C(K,Y)^\ast}$ given by $\xi^\ast(h)=y^\ast(h(t_0))$ for $h\in C(K,Y)$ which clearly satisfies $\xi^\ast\in \Slice(B_{C(K,Y)^\ast},f,\eps)$. Finally, for each $k=1,\ldots,n$, we define $g_k\in C(K,X)$ by
$$
g_k(t)=x_k+\varphi(t)(\theta_k^{-1} g(t)-x_k) \qquad (t\in K)
$$
and we observe that
$$
\widetilde{G}^\ast\xi^\ast(g_k)=\xi^\ast(G\circ g_k)=y^\ast(G(g_k(t_0)))=G^\ast y^\ast(x_k).
$$
Therefore, we deduce that $g_k\in \GS(S_{C(K,X)},\widetilde{G}^{\ast}\xi^{\ast}, \eps)$ for every $k=1,\ldots,n$. On the other hand, for an arbitrary $t\in K$ we have that
$$
\left\|g(t)-\sum_{k=1}^n \lambda_k\theta_k g_k(t)\right\|=\left\|(1-\varphi(t))\left(g(t)-\sum_{k=1}^n \lambda_k\theta_k x_k\right)\right\|.
$$
So, if $t\in U$, then $\|g(t)-g(t_0)\|\leqslant \frac{\eps}{2}$ and, therefore,
$$
\left\|g(t)-\sum_{k=1}^n \lambda_k\theta_k g_k(t)\right\|\leqslant \|g(t)-g(t_0)\|+\left\|g(t_0)-\sum_{k=1}^n \lambda_k\theta_k x_k\right\|\leqslant \eps.
$$
If, otherwise, $t\notin U$, then $\varphi(t)=1$ and thus $g(t)-\sum_{k=1}^n \lambda_k\theta_k g_k(t)=0$. All this tells us that
$$
\dist\bigl(g, \aconv(\GS(S_{C(K,X)},\widetilde{G}^{\ast}\xi^{\ast}, \eps))\bigr) <\eps
$$
and shows that $\widetilde{G}$ is lush.

Suppose now that $\widetilde{G}$ is lush and let us show that $G$ is lush. To do so, we will use Proposition~\ref{Prop:characterization-lushness}.(iii) with the set
$$
\mathcal{A}=\{y^\ast\otimes \delta_t \ : \ y^\ast \in S_{Y^\ast}, t\in K\}
$$
where $(y^\ast\otimes \delta_t)(f)=y^\ast(f(t))$ for $f\in C(K,Y)$. Observe that $\mathcal{A}$ is norming and rounded, so $\overline{\conv}^{\omega^{\ast}}{(\mathcal{A})} = B_{C(K,Y)^{\ast}}$. Fixed $x_0\in S_X$, $y_0\in S_Y$, and $\eps>0$, we consider $f\in S_{C(K,Y)}$ and $g\in S_{C(K,X)}$ given respectively by
$$
f(t)=y_0 \qquad \text{and} \qquad g(t)=x_0 \qquad (t\in K).
$$
Now we use that $\widetilde{G}$ is lush and Proposition~\ref{Prop:characterization-lushness}.(iii) to find $y_0^\ast\otimes \delta_{t_0}\in \Slice(\mathcal{A},f,\eps)$ such that
$$
\dist\bigl(g, \aconv(\GS(S_{C(K,X)},\widetilde{G}^{\ast}(y_0^\ast\otimes \delta_{t_0}), \eps))\bigr) <\eps.
$$
Therefore, as $f(t_0)=y_0$, we clearly get that $y_0^\ast \in \Slice(B_{Y^\ast},y_0,\eps)$. Moreover, using that if $h\in \GS(S_{C(K,X)},\widetilde{G}^{\ast}(y_0^\ast\otimes \delta_{t_0}), \eps)$ then $h(t_0)\in \GS(S_X,G^{\ast}y_0^\ast,\eps)$, we easily deduce that
$$
\dist\bigl(x_0, \aconv(\GS(S_X,G^{\ast}y_0^\ast, \eps))\bigr) <\eps
$$
which gives that $G$ is lush.
\end{proof}

\begin{proof}[Proof of Theorem~\ref{Theo:stabilitityC(K)}.(c)]
We start showing that $\widetilde{G}$ has the aDP when $G$ does. Indeed, observe that in the first part of the proof of Theorem~\ref{Theo:stabilitityC(K)}.(a), if the operator $T$ has rank one then so does the operator $S\in L(X,Y)$ constructed there.

To prove the reversed implication, fix an isolated point $t_0\in K$ and observe that we can identify $C(K,X)\equiv X\oplus_\infty C(K\setminus\{t_0\},X)$ and $C(K,X)\equiv Y\oplus_\infty C(K\setminus\{t_0\},Y)$. Now, if for $g\in C(K\setminus\{t_0\},X)$ we write
$$
\widehat{g}(t)=\begin{cases} 0 & \text{ if } t=t_0 \\
g(t) & \text{ if } t\neq t_0
\end{cases} \qquad (t\in K),
$$
and we consider the operator $\widehat{G}:C(K\setminus\{t_0\},X)\longrightarrow C(K\setminus\{t_0\},Y)$ given by $$
\widehat{G}(g)=[\widetilde{G}(\widehat{g})]|_{K\setminus\{t_0\}} \qquad \big(g\in C(K\setminus\{t_0\},X)\big)
$$
then, we can write
$$
\widetilde{G}(x,g)=(Gx, \widehat{G}(g)) \qquad \big(x\in X, g\in C(K\setminus\{t_0\},X)\big).
$$
Therefore, as $\widetilde{G}$ has the aDP, we may use Proposition \ref{prop:stability-c0l1ellinfty-sums} to deduce that $G$ has aDP.
\end{proof}

\begin{proof}[Proof of Theorem~\ref{Theo:stabilitityC(K)}.(d)]
We prove first the sufficiency. We will use Theorem \ref{Theo:aDPCharacterization}.(iii) to show that $\widetilde{G}$ has the aDP. So, fixed $f\in S_{C(K,Y)}$ and $\eps>0$, we write
$$
\Delta_\eps(f)=\{g\in B_{C(K,X)} \ : \ \|\widetilde{G}(g)+\T f\|>2-\eps\}
$$
and we have to show that $\overline{\conv}(\Delta_\eps(f))=B_X$. The argument follows the lines of \cite[p.~81]{WerSur}: let $U$ be the open set $\{t\in K \ : \ \|f(t)\|>1-\eps/2\}$ and pick, given $n\in \N$, open pairwise disjoint non-void subsets $U_1,\ldots,U_n\subset U$ and points $t_j\in U_j$. Next, we use the hypothesis to find $x_j\in B_X$ and $\theta_j \in \T$ such that $\|G(x_j)+\theta_jf(t_j)\|>2-\eps$.
Now, fixed $h\in B_{C(K,X)}$, we may choose functions $g_j\in B_{C(K,X)}$ such that $g_j\equiv h$ in $K\setminus U_j$ and $g_j(t_j)=x_j$. Indeed, take Urysohn functions $\varphi_j :K\longrightarrow [0,1]$ such that
$$
\varphi_j|_{K\setminus U_j}\equiv 1 \qquad \text{and} \qquad \varphi_j(t_j)=0,
$$
and define
$$
g_j(t)=\varphi_j(t)h(t)+(1-\varphi_j(t))x_j \qquad (t\in K).
$$
On the one hand, observe that $g_j\in \Delta_\eps(f)$:
$$
\|\widetilde{G}(g_j)+\theta_jf\|\geqslant \|G(g_j(t_j))+\theta_jf(t_j)\|=\|G(x_j)+\theta_jf(t_j)\|>2-\eps.
$$
On the other hand, for $t\in U_k$ we have that
\begin{align*}
\left\|h(t)-\frac1n\sum_{j=1}^n g_j(t)\right\|=\left\|h(t)-\frac{n-1}{n}h(t)-\frac{1}{n}g_k(t)\right\|=\frac1n\|h(t)-g_k(t)\|\leqslant \frac{2}{n};
\end{align*}
and, for $t\notin \bigcup_jU_j$, it follows that $h(t)-\frac1n\sum_{j=1}^n g_j(t)=0$. This proves that $h\in \overline{\conv}(\Delta_\eps(f))$ and so $\widetilde{G}$ has the aDP.

Suppose now that $\widetilde{G}$ has the aDP. Fixed $\eps>0$ and a non-zero $y\in B_Y$, we take the constant function $f\in C(K,Y)$ given by $f\equiv\frac{y}{\|y\|}$ and we use Theorem \ref{Theo:aDPCharacterization}.(iii) to find $g\in B_{C(K,X)}$ such that $\|\widetilde{G}(g)+\T f\|>2-\eps$. So, there is $t_0\in K$ satisfying $\|G(g(t_0))+\T \frac{y}{\|y\|}\|>2-\eps$ and, therefore, $g(t_0)$ is the element in $B_X$ we are looking for:
\[
\|G(g(t_0))+\T y\|\geqslant \left\|G(g(t_0))+\T \frac{y}{\|y\|}\right\|-\left\|\frac{y}{\|y\|}-y\right\|>2-\eps-(1-\|y\|)=1+\|y\|-\eps. \qedhere
\]
\end{proof}

We next deal with spaces of essentially bounded measurable functions.

\begin{Theo}\label{Theo:stabilityLinfty}
Let $X,Y$ be a Banach spaces, let $(\Omega, \Sigma, \mu)$ be a $\sigma$-finite measure space and let $G\in L(X,Y)$ be a norm-one operator. Consider the norm-one composition operator $\widetilde{G}:L_\infty(\mu,X) \longrightarrow L_\infty(\mu,Y)$ given by $\widetilde{G}(f)=G\circ f$ for every $f\in L_\infty(\mu,X)$. Then:
\begin{enumerate}
\item[(a)] $\widetilde{G}$ is a spear operator if and only if $G$ is a spear operator.
\item[(b)] $\widetilde{G}$ is lush if and only if $G$ is lush.
\item[(c)] If $\mu$ has an atom, then $\widetilde{G}$ has the aDP if and only if $G$ does.
\item[(d)] If $\mu$ is nonatomic, then $\widetilde{G}$ has the aDP if and only if $G(B_X)$ is a spear set.
\end{enumerate}
\end{Theo}

\begin{Rema}
The results in items (a), (c), and (d) of the above theorem were known for the case of the identity (see \cite{martinOikhberg,M-V}). The content of (b) is completely new even for the identity.
\end{Rema}

\begin{Coro}
Let $X$ be a Banach space and let $(\Omega, \Sigma, \mu)$ be a $\sigma$-finite measure space. Then, $L_\infty(\mu, X)$ is lush if and only if $X$ is lush.
\end{Coro}

We will use the following notation during the proof of the theorem:
\[
\Sigma^{+} :=\{ A \in \Sigma \colon 0  < \mu(A) < + \infty\}.
\]
Observe that when the measure is finite, this notation is consistent with the one given in section \ref{subsect:from_L_1}. Moreover, we will also use the subset of $S_{L_{1}(\mu, Y^{\ast})}$ given by
\[
 \mathcal{A} := \left\{ y^{\ast} \frac{\mathbbm{1}_{A}}{\mu(A)} \colon y^{\ast} \in S_{Y^{\ast}}, A \in \Sigma^{+} \right\}
\]
which clearly satisfies that $$\overline{\conv}^{w^\ast}{(\mathcal{A})} = B_{L_{\infty}(\mu, Y)^{\ast}}.$$

\begin{proof} [Proof of Theorem~\ref{Theo:stabilityLinfty}.(a)]
This is an easy adaptation of \cite[Theorem~2.3]{M-V}.
Suppose first that $G$ is a spear operator. We fix $T\in L\big(L_\infty(\mu, X),L_\infty(\mu,Y)\big)$ with $\|T\|=1$. Given $\eps>0$ we may follow the first part of the proof of \cite[Theorem~2.3]{M-V} to find $f\in S_{L_\infty(\mu,X)}$, $x_0\in S_X$, and $A,B\in \Sigma$ with $0<\mu(B)<\infty$, such that
\begin{equation}\label{eq:Linfty(mu,X)}
B\subset A \qquad \text{and}\qquad \left\|\frac{1}{\mu(B)}\int_B T(x_0 \mathbbm{1}_A+f \mathbbm{1}_{\Omega\setminus A})d\mu\right\|>1-\eps.
\end{equation}
Now we fix $x^\ast\in S_{X^\ast}$ with $x^\ast_0(x_0)=1$, we write
$$
\Phi(x)=x_0 \mathbbm{1}_A+x_0^\ast(x)f \mathbbm{1}_{\Omega\setminus A} \qquad (x\in X)
$$
and we define the operator $S\in L(X,Y)$ given by
$$
S(x)=\frac{1}{\mu(B)}\int_BT(\Phi(x))d\mu \qquad (x\in X)
$$
which, by \eqref{eq:Linfty(mu,X)}, satisfies $\|S\|\geqslant\|Sx_0\|>1-\eps$. Next, we use that $G$ is a spear operator to find $x\in S_X$ such that $\|Gx+\T Sx\|>2-\eps$, so we can take $y^\ast\in S_{Y^\ast}$ satisfying
$$
|y^\ast(Gx+\T Sx)|>2-\eps.
$$
Finally, define the functional $g^\ast\in S_{L_\infty(\mu,Y)^\ast}$ by
$$
g^\ast(h)=y^*\left(\frac{1}{\mu(B)}\int_B hd\mu\right) \qquad (h\in L_\infty(\mu, Y))
$$
and observe that
\begin{align*}
\|\widetilde{G}+\T \,T\|&\geqslant \left|g^\ast\big(\widetilde{G}(\Phi(x))+\T \,T(\Phi(x))\big)\right|\\
&=\left|g^\ast(G(x)\mathbbm{1}_A)+\T \,g^\ast\big(T(\Phi(x))\big)\right|=|y^\ast(G(x))+\T y^\ast(S(x))|>2-\eps.
\end{align*}
The arbitrariness of $\eps$ gives that $\widetilde{G}$ is a spear operator.

Assume now that $\widetilde{G}$ is a spear operator. Fix $S\in L(X,Y)$ and define the operator $T\in L\big(L_\infty(\mu, X),L_\infty(\mu,Y)\big)$ by
$$
[T(f)](t)=S(f(t)) \qquad (t\in \Omega, f\in L_\infty(\mu,X))
$$
which clearly satisfies $\|T\|=\|S\|$. As we mentioned at the beginning of the proof, we can find $f\in S_{L_\infty(\mu,X)}$, $x_0\in S_X$, and $A,B\in \Sigma$ with $0<\mu(B)<\infty$, such that $B\subset A$ and
\begin{equation*}
\left\|\frac{1}{\mu(B)}\int_B (\widetilde{G}+\T\, T)(x_0 \mathbbm{1}_A+f \mathbbm{1}_{\Omega\setminus A})d\mu\right\|>\|\widetilde{G}+\T\, T\|-\eps=1+\|T\|-\eps=1+\|S\|-\eps.
\end{equation*}
Therefore, we can write
\begin{align*}
1+\|S\|-\eps&\leqslant \left\|\frac{1}{\mu(B)}\int_B (\widetilde{G}+\T\, T)(x_0 \mathbbm{1}_A+f \mathbbm{1}_{\Omega\setminus A})d\mu\right\|\\
&=\left\|\frac{1}{\mu(B)}\int_B G(x_0)\mathbbm{1}_A+\T\, T(x_0 \mathbbm{1}_A+f \mathbbm{1}_{\Omega\setminus A})d\mu\right\|\\
&=\left\|G(x_0)+\T \frac{1}{\mu(B)}\int_B T(x_0 \mathbbm{1}_A+f \mathbbm{1}_{\Omega\setminus A})d\mu\right\|\\
&=\left\|G(x_0)+\T S\left(\frac{1}{\mu(B)}\int_B x_0 \mathbbm{1}_A+f \mathbbm{1}_{\Omega\setminus A}d\mu\right)\right\|=\|G(x_0)+\T S(x_0)\|.
\end{align*}
Thus, we get $\|G+\T S\|>1+\|S\|-\eps$ and so $G$ is a spear operator.
\end{proof}

We next deal with lushness for $L_\infty(\mu,X)$.

\begin{proof} [Proof of Theorem~\ref{Theo:stabilityLinfty}.(b)]
Assume first that $G$ is lush. To prove that so is $\widetilde{G}$, we will check that Proposition \ref{Prop:characterization-lushness}.(iii) is satisfied. Let $f_{0} \in S_{L_{\infty}(\mu, X)}$, $g_{0} \in S_{L_{\infty}(\mu, Y)}$ and $\eps > 0$. By density, we can assume that $f_{0},g_{0}$ can be written as
$$
f_{0}= \sum_{A \in \pi}{x_{A} \mathbbm{1}_{A}}, \qquad \, g_{0}= \sum_{A \in \pi}{y_{A} \mathbbm{1}_{A}}
$$
where $\pi \subset \Sigma^{+}$ is a countable partition of $\Omega$ and $x_{A} \in B_{X}$, $y_{A} \in B_{Y}$ for each $A \in \pi$. Since $\| g_{0}\|_{L_{\infty}(\mu, Y)} = 1$, we can assume without loss of generality that there is $A_{0} \in \pi$ with $\| y_{A_{0}}\| = 1$. Using that $G$ is lush, we can find $y_{0}^{\ast} \in S_{Y^{\ast}}$ such that
\begin{equation}\label{equa:stabilityLinftyLushAux1}
\Real{y_{0}^{\ast}(y_{A_{0}})} > 1 - \eps
\end{equation}
and elements $x_{j} \in \GS(S_{X}, G^{\ast}y_{0}^{\ast}, \eps)$, $\theta_{j} \in \T$, $\lambda_{j} > 0$ for $j=1, \ldots, m$ ($m \in \N$) with $\sum_{j}{\lambda_{j}} = 1$ satisfying that
\begin{equation}\label{equa:stabilityLinftyLushAux2}
\left\| x_{A_{0}} - \sum_{j=1}^{m}{\lambda_{j} \theta_{j} x_{j}} \right\| < \eps.
\end{equation}
Consider now $h^{\ast}_{0} := y_{0}^{\ast} \mathbbm{1}_{A_{0}}/\mu(A_{0}) \in \mathcal{A}$, and for each $j =1, \ldots, m$ let
$$ f_{j} := \sum_{A \in \pi, A \neq A_{0}}{\overline{\theta_{j}} x_{A} \mathbbm{1}_{A}} + x_{j} \mathbbm{1}_{A_{0}}. $$
Then, by \eqref{equa:stabilityLinftyLushAux1} we have that
$$\Real{h^{\ast}_{0}(g_{0})} = \Real{y_{0}^{\ast}(y_{A_{0}})} > 1- \eps.$$
A similar argument shows that
$$
f_{j} \in \GS(S_{L_{\infty}(\mu, X)}, \widetilde{G}^{\ast}h^{\ast}, \eps)
$$
for every $j=1, \ldots, m$, since $\Real{G^{\ast}y_{0}^{\ast}(x_{j})} >  1 - \eps$. Moreover, using \eqref{equa:stabilityLinftyLushAux2} we immediately conclude that
$$ \left\| f_0 - \sum_{j=1}^{m}{\lambda_{j} \theta_{j} f_{j}} \right\| = \max{\left\{ \sup_{\substack{A \in \pi\\ A \neq A_{0}}}{\left\| x_{A} -\sum_{j=1}^{m}{\lambda_{j} x_{A}} \right\|}, \left\| x_{A_{0}} - \sum_{j=1}^{m}{\lambda_{j} \theta_{j} x_{j}}\right\| \right\}} < \eps. $$

Let us see the converse: assume that $\widetilde{G}$ is lush, fix an element $B \in \Sigma^{+}$, and let $x \in S_{X}$, $y \in S_{Y}$ and $\eps > 0$. Then, defining $f:=x \mathbbm{1}_{B} \in S_{L_{\infty}(\mu, X)}$ and $g:=y \mathbbm{1}_{B} \in S_{L_{\infty}(\mu, Y)}$, we can use the hypothesis to find $h^{\ast} \in S(\mathcal{A},g, \eps )$ such that
\begin{equation}\label{equa:stabilityLinftyLushAux2.5}
\dist{\left( f, \aconv{\GS(B_{L_\infty}(\mu, X), \widetilde{G}^{\ast}h^{\ast}, \eps)} \right)}<\eps.
\end{equation}
Since we can write $h^{\ast} = y_{0}^{\ast} \mathbbm{1}_{A}/\mu(A)$ for some $A \in \Sigma^{+}$ and $y_{0}^{\ast} \in S_{Y^{\ast}}$, condition $h^{\ast} \in S(\mathcal{A},g, \eps )$ can be rewritten as
\begin{equation}\label{equa:stabilityLinftyLushAux2.75}
\Real{y_{0}^{\ast}(y_{0}) \frac{\mu(A \cap B)}{\mu(A)}} > 1 - \eps.
\end{equation}
By \eqref{equa:stabilityLinftyLushAux2.5} we find elements $f_{j} \in \GS(B_{L_{\infty}(\mu,X)}, \widetilde{G}^{\ast}h^{\ast}, \eps)$, $\theta_{j} \in \mathbb{T}$ and $\lambda_{j} >0$ for $j=1, \ldots, m$ satisfying $\sum_{j=1}^{m}{\lambda_{j}} = 1$ and
\begin{equation}\label{equa:stabilityLinftyLushAux3}
\left\| f - \sum_{j=1}^{m}{\lambda_{j} \theta_{j} f_{j}} \right\| < \eps.
\end{equation}
Next, for each $j=1, \ldots, m$, we consider the element
$$
x_{j} := \frac{1}{\mu(A)}\int_{A}{f_{j} \: d \mu} \in B_{X}
$$
which clearly satisfies that
$$
\Real{G^{\ast}y_{0}^{\ast}(x_{j})} = \frac{1}{\mu(A)}\int_{A}{\Real{G^{\ast}y_{0}^{\ast}(f_{j}) \: d \mu}} = \Real{\widetilde{G}^{\ast}h^{\ast}(f_{j})} > 1 - \eps.
$$
Moreover, making use of \eqref{equa:stabilityLinftyLushAux2.75} we get that
$$ \left\| x - \frac{1}{\mu(A)}{\int_{A}{f \:d \mu}} \right\| = \left\| x - \frac{\mu(A \cap B)}{\mu(A)} x \right\| < \eps $$
and, combining this with \eqref{equa:stabilityLinftyLushAux3}, we conclude that
\begin{align*}
\left\| x - \sum_{j=1}^{m}{\lambda_{j} \theta_{j} x_{j}} \right\| & \leqslant \eps + \left\| \frac{1}{\mu(A)}\int_{A}{f \: d \mu} - \sum_{j=1}^{m}{\lambda_{j} \theta_{j} \frac{1}{\mu(A)} \int_{A}{f_{j} \: d \mu}} \right\|\\
& \leqslant \eps + \frac{1}{\mu(A)} \int_{A}{\left\| f - \sum_{j=1}^{m}{\lambda_{j} \theta_{j} f_{j}} \right\| \: d \mu} < 2 \eps,
\end{align*}
which finishes the proof.
\end{proof}

\begin{proof} [Proof of Theorem~\ref{Theo:stabilityLinfty}.(c)]
We start showing that when $G$ has the aDP so does $\widetilde{G}$. Indeed, observe that in the first part of the proof of Theorem \ref{Theo:stabilityLinfty}.(a), if the operator $T$ has rank one then the operator $S\in L(X,Y)$ constructed there also has rank one.

To prove the reversed implication, fix $A_0\in \Sigma$ which is an atom for $\mu$ and observe that we can identify
\begin{align*}
L_\infty((\Omega,\mu),X)&\equiv X\oplus_\infty L_\infty((\Omega\setminus A_0,\mu),X)
\end{align*}
and
\begin{align*}
 L_\infty((\Omega,\mu),Y)&\equiv Y\oplus_\infty L_\infty((\Omega\setminus A_0,\mu),Y).
\end{align*}
Now, if for $g\in L_\infty((\Omega\setminus A_0,\mu),X)$ we write
$$
\widehat{g}(t)=\begin{cases} 0 & \text{ if } t\in A_0 \\
g(t) & \text{ if } t\notin A_0
\end{cases} \qquad (t\in \Omega),
$$
and we consider the operator $\widehat{G}:L_\infty((\Omega\setminus A_0 ,\mu),X)\longrightarrow L_\infty((\Omega\setminus A_0,\mu),Y)$ given by
$$
\widehat{G}(g)=[\widetilde{G}(\widehat{g})]|_{\Omega\setminus A_0} \qquad \big(g\in L_\infty((\Omega\setminus A_0,\mu),X)\big),
$$
then we can write
$$
\widetilde{G}(x,g)=(Gx, \widehat{G}(g)) \qquad \big(x\in X, g\in L_\infty((\Omega\setminus A_0,\mu),X)\big).
$$
Therefore, as $\widetilde{G}$ has the aDP, we may use Proposition \ref{prop:stability-c0l1ellinfty-sums} to deduce that $G$ has the aDP.
\end{proof}

\begin{proof} [Proof of Theorem~\ref{Theo:stabilityLinfty}.(d)] We prove first the sufficiency. We will use Theorem \ref{Theo:aDPCharacterization}.(iii) to show that $\widetilde{G}$ has the aDP. So, fixed $f\in S_{L_\infty(\mu,Y)}$ and $\eps>0$, we write
$$
\Delta_\eps(f)=\{g\in B_{L_\infty(\mu,X)} \ : \ \|\widetilde{G}(g)+\T f\|>2-\eps\}
$$
and we have to show that $\overline{\conv}(\Delta_\eps(f))=B_{L_\infty(\mu,X)}$. Using Lemma~2.2 in \cite{M-V}, we may find $y\in S_Y$ and $A\in \Sigma^{+}$ such that
$$
\|f-(y\mathbbm{1}_A+f\mathbbm{1}_{\Omega\setminus A})\|<\frac{\eps}{2}\,.
$$
As $\mu$ is atomless, for every $n\in \N$ we may and do pick pairwise disjoint sets $U_1,\ldots,U_n$ with positive measure such that $U_i\subset A$ for every $i=1,\ldots,n$. Next, we use that $G(B_X)$ is a spear set to find $x\in B_X$ and $\theta \in \T$ such that $\|G(x)+\theta y\|>2-\eps/2$.
Now, fixed $h\in B_{L_\infty(\mu,X)}$, we define $g_j= x\mathbbm{1}_{U_j} +h\mathbbm{1}_{K\setminus U_j}\in B_{L_\infty(\mu,Y)}$ for $j=1,\ldots,n$.
On the one hand, observe that for every $t\in U_j$ we get the following estimation
$$
\|G(g_j(t))+\theta f(t)\|\geqslant\|G(x)+\theta y\|-\|f(t)-y\|>2-\eps,
$$
so $\|\widetilde{G}(g_j)+\theta f\|\geqslant 2-\eps$ which implies that $g_j\in \Delta_\eps(f)$. On the other hand, for $t\in U_k$ we have that
\begin{align*}
\left\|h(t)-\frac1n\sum_{j=1}^n g_j(t)\right\|=\left\|h(t)-\frac{n-1}{n}h(t)-\frac{1}{n}g_k(t)\right\|=\frac1n\|h(t)-g_k(t)\|\leqslant \frac{2}{n};
\end{align*}
and, for $t\notin \bigcup_jU_j$, it follows that $h(t)-\frac1n\sum_{j=1}^n g_j(t)=0$. This proves that $h\in \overline{\conv}(\Delta_\eps(f))$ and so $\widetilde{G}$ has the aDP.

Conversely, suppose now that $\widetilde{G}$ has the aDP. Fixed $\eps>0$ and a non-zero $y\in B_Y$, we take the constant function $f\in L_\infty(\mu,Y)$ given by $f\equiv\frac{y}{\|y\|}$ and we use Theorem \ref{Theo:aDPCharacterization}.(iii) to find $g\in B_{L_\infty(\mu,X)}$ such that $\|\widetilde{G}(g)+\T f\|>2-\eps$. So, there is $t_0\in \Omega$ satisfying that $\left\|G(g(t_0))+\T \frac{y}{\|y\|}\right\|>2-\eps$ and, therefore,
$$
\|G(g(t_0))+\T y\|\geqslant \left\|G(g(t_0))+\T \frac{y}{\|y\|}\right\|-\left\|\frac{y}{\|y\|}-y\right\|>2-\eps-(1-\|y\|)=1+\|y\|-\eps.
$$
This shows that $G(B_X)$ is a spear set, concluding thus the proof.
\end{proof}

Finally, we would like to work with spaces of integrable functions.

\begin{Theo}\label{Theo:compositionOperatorL1}
Let $X,Y$ be Banach spaces, let $(\Omega, \Sigma, \mu)$ be a $\sigma$-finite measure space, and let $G \in L(X,Y)$ be a norm-one operator. Consider the norm-one composition operator $\widetilde{G}: L_{1}(\mu, X) \longrightarrow L_{1}(\mu, Y)$ given by $\widetilde{G}(f) = G \circ f$ for every $f\in L_1(\mu,X)$. Then:
\begin{enumerate}
\item[(a)] $\widetilde{G}$ is a spear operator if and only if $G$ is a spear operator.
\item[(b)] $\widetilde{G}$ is lush if and only if $G$ is lush.
\item[(c)] If $\mu$ has an atom, then $\widetilde{G}$ has the aDP if and only if $G$ has the aDP.
\item[(d)] If $\mu$ is nonatomic, then $\widetilde{G}$ has the aDP if and only if
\[
B_{X} = \overline{\aconv}{\{ x \in B_{X} \colon \| Gx\| > 1 - \eps \}} \mbox{ for every $\eps > 0$}.
\]
\end{enumerate}
\end{Theo}

\begin{Rema}
The results in items (a), (c), and (d) of the above theorem were known for the case of the identity (see \cite{martinOikhberg,M-P}). The content of (b) is completely new even for the identity.
\end{Rema}

\begin{Coro}
Let $X$ be a Banach space and let $(\Omega, \Sigma, \mu)$ be a $\sigma$-finite measure space. Then, $L_1(\mu,X)$ is lush if and only if $X$ is lush.
\end{Coro}

We claim that for the proof of Theorem~\ref{Theo:compositionOperatorL1} we can assume without loss of generality that $(\Omega, \Sigma, \mu)$ is a probability space. This is clear if $\mu$ is finite by normalizing the measure. If $\mu$ is $\sigma$-finite but not finite, we can find a countable partition $\Lambda \subset \Sigma^{+}$ of $\Omega$.  For each $A \in \Lambda$ we have a finite measure space $(A, \Sigma_{A}, \mu_{A})$ where
$$\Sigma_{A} = \{ B \in \Sigma \colon B \subset A \} \qquad \text{and} \qquad \mu_{A}(B) = \mu(B) \mbox{ for $B \in \Sigma_{A}$}. $$
Moreover, we have a canonical isometry
$$
L_{1}(\Omega, \Sigma, \mu, X) = \left[ \bigoplus_{A \in \Lambda}{L_{1}(A, \Sigma_{A},\mu_{A},)} \right]_{\ell_{1}}
$$
(see \cite[pp.~501]{DefantFloret}, for instance). Then,
the composition operator $\widetilde{G}$ of Theorem \ref{Theo:compositionOperatorL1}, can be seen as $\widetilde{G}[(f_{A})_{A \in \Lambda}] = (\widetilde{G}_{A} f_{A})_{A \in \Lambda}$ where $\widetilde{G}_{A}$ is the composition operator correspondant to $(A, \Sigma_{A}, \mu_{A})$. An application of Proposition \ref{prop:stability-c0l1ellinfty-sums} gives that $\widetilde{G}$ is a spear operator (resp.\ lush, aDP) if and only if $\widetilde{G}_{A}$ is a spear operator (resp.\ lush, aDP) for each $A \in \Lambda$. Therefore, the claim is valid for the proofs of (a) and (b). In the case of (c) and (d) it is also valid taking into account that $(\Omega, \Sigma, \mu)$ has an atom if and only if some $(A, \Sigma_{A}, \mu_{A})$ does, and that if $G$ has the aDP then it satisfies
\[
B_{X} = \overline{\aconv}{\{ x \in B_{X} \colon \| Gx\| > 1 - \eps \}} \mbox{ for every $\eps > 0$}.
\]
Hence, the claim is established.

In order to prove Theorem~\ref{Theo:compositionOperatorL1} we need to introduce some notation. If $(\Omega, \Sigma, \mu)$ is a probability space, the set
$$ \mathcal{A}:=\left\{ \sum_{A \in \pi}{y_{A}^{\ast} \mathbbm{1}_{A}}  \colon \pi \subset \Sigma^{+} \mbox{ finite partition of $\Omega$ }, \, y_{A}^{\ast} \in S_{Y^{\ast}} \right\} \subset S_{L_{\infty}(\mu, Y^{\ast})} $$
satisfies that
\begin{equation}\label{equa:normingSetDualIntegrable}
B_{L_{1}(\mu, Y)^{\ast}} = \overline{\conv}^{\omega^\ast}{(\mathcal{A})},
\end{equation}
since $\mathcal{A}$ is rounded and it is clearly norming for the simple functions of $L_{1}(\mu, Y)$. On the other hand,  we will write
$$ \mathcal{B} := \left\{ x \, \frac{\mathbbm{1}_{B}}{\mu(B)} \colon x \in S_{X}, B \in \Sigma^{+} \right\} $$
which satisfies that
\begin{equation}\label{equa:normingSetIntegrable}
B_{L_{1}(\mu, X)} = \overline{\conv}{(\mathcal{B})}.\end{equation}
Indeed, it is enough to notice that every simple function $f$ in $S_{L_{1}(\mu, X)}$ belongs to the convex hull of $\mathcal{B}$: such an $f$ can be written as $f = \sum_{B \in \pi}{x_{B} \mathbbm{1}_{B}}$, where $\pi \subset \Sigma^{+}$ is a finite family of pairwise disjoint sets of $\Omega$ and $x_{B} \in X \setminus \{ 0\}$ for each $B \in \pi$. Then
$$ \| f\| = \sum_{B \in \pi}{\|x_{B}\| \mu(B)} = 1, $$
and hence
$$ f = \sum_{B \in \pi}{\| x_{B}\| \mu(B) \, \frac{x_{B}}{\| x_{B}\|} \frac{\mathbbm{1}_{B}}{\mu(B)}} \in \conv{\mathcal{B}}. $$

\begin{proof}[Proof of Theorem \ref{Theo:compositionOperatorL1}.(a)]
Suppose first that $\widetilde{G}$ is a spear operator. Fix $T\in L(X,Y)$ and consider $\widetilde{T} \in L\big(L_1(\mu,X),L_1(\mu,Y)\big)$ given by $\widetilde{T}(f) = T \circ f$, which satisfies $\|\widetilde{T}\|=\|T\|$. Given $\eps>0$, we can find $x \in S_{X}$ and $B \in \Sigma^{+}$ such that
\begin{equation}\label{eq:L1(mu,X)-spear}
\left\|(\widetilde{G}+\T \,\widetilde{T})(x \textstyle{\frac{\mathbbm{1}_{B}}{\mu(B)}})\right\|>1+\|T\|-\eps.
\end{equation}
But notice that
\begin{align*}
\left\|\widetilde{G}(\textstyle{x\frac{\mathbbm{1}_{B}}{\mu(B)}})+\T\,\widetilde{T}(\textstyle{x\frac{\mathbbm{1}_{B}}{\mu(B)}})\right\| &  = \left\| \big(G(x) + \T \, T(x)\big) \textstyle{\frac{\mathbbm{1}_{B}}{\mu(B)}} \right\| = \| G(x) + \T \, T(x)\|.
\end{align*}
This, together with \eqref{eq:L1(mu,X)-spear} and the arbitrariness of $\eps$, tells us that $G$ is a spear operator.

Assume now that $G$ is a spear operator. Fixed $T\in L\big(L_1(\mu,X),L_1(\mu,Y)\big)$ with $\|T\|=1$ and $\eps>0$, we may find by \eqref{equa:normingSetIntegrable} elements $x_{0} \in S_{X}$ and $B \in \Sigma^{+}$ such that
$$ \Big\| T\left(x_{0} \textstyle{\frac{\mathbbm{1}_{B}}{\mu(B)}} \right)\Big\| > 1 - \eps. $$
Using now \eqref{equa:normingSetDualIntegrable}, there exists $f^{\ast} = \sum_{A \in \pi}{y_{A}^{\ast} \mathbbm{1}_{A}}$, where $\pi$ is a finite partition of $\Omega$ into sets of $\Sigma^+$ and $ y_{A}^{\ast} \in S_{Y^{\ast}}$ for each $A \in \pi$, satisfying that
\begin{equation}\label{equa:spearStabilityL1aux1}
\Real{f^{\ast}\left( T\left(\textstyle{x_{0}\frac{\mathbbm{1}_{B}}{\mu(B)}}\right)\right)} = \Real{\sum_{A \in \pi}{y_{A}^{\ast}\left( \int_{A}{T\left(\textstyle{x_{0}\frac{\mathbbm{1}_{B}}{\mu(B)}}\right) \: d \mu}\right)}} > 1- \eps.
\end{equation}
Then, we can write
\begin{equation*}
T\left(x_{0} \textstyle{\frac{\mathbbm{1}_{B}}{\mu(B)}} \right) = \sum_{\substack{A \in \pi\\ \mu(A \cap B) \neq 0}}{\frac{\mu(B \cap A)}{\mu(B)} \,  T\left(x_{0} \textstyle{\frac{ \mathbbm{1}_{B \cap A}}{\mu(B\cap A)}}\right)}
\end{equation*}
so by a standard convexity argument we can assume that there is $A_{0} \in \pi$ such that $B\subset A_0$ and \eqref{equa:spearStabilityL1aux1} is still satisfied. By the density of norm-attaining functionals, we can and do assume that every $y_{A}^{\ast}$ is norm-ataining, so there is $y_{A_{0}} \in S_{Y}$ such that $y_{A_{0}}^{\ast}(y_{A_{0}}) = 1$. Define the operator $S: X \longrightarrow Y$ by
$$
S(x)=\int_{A_0}T\left(x\textstyle{\frac{\mathbbm{1}_{B}}{\mu(B)}}\right)d\mu+
\left[\sum_{A \in \pi \setminus \{ A_{0}\}} {y^\ast_A \left(\int_{A}T\left(x\textstyle{\frac{\mathbbm{1}_{B}}{\mu(B)}}\right)d\mu\right)}\right] y_{A_{0}} \quad (x\in X).
$$
It is easy to check that $\| S\| \leqslant 1$, and moreover $\| S\| > 1 - \eps$ since as a consequence of \eqref{equa:spearStabilityL1aux1} we obtain that
$$ \| S(x_{0})\| \geqslant |y_{A_{0}}^{\ast}(S x_{0})| = \left| f^{\ast}\left( T\left(x_{0} \textstyle{\frac{\mathbbm{1}_{B}}{\mu(B)}}\right)\right) \right| > 1 - \eps. $$
By hypothesis, we can then find $x_{1} \in S_{X}$ and $\theta_{1} \in \mathbb{T}$ such that $\| G (x_{1}) + \theta_{1} S(x_{1})\| > 2 - \eps$. We claim that
\begin{equation}\label{equa:spearStabilityL1aux2}
\left\| \widetilde{G}\left(x_{1} \textstyle{\frac{\mathbbm{1}_{B}}{\mu(B)}}\right) + \theta_{1} \, T\left( x_{1} \textstyle{\frac{\mathbbm{1}_{B}}{\mu(B)}}\right) \right\| > 2- \eps.
\end{equation}
Indeed, the left-hand side expression of \eqref{equa:spearStabilityL1aux2} is equal to
\begin{equation*}
\begin{split}
& \int_{B}{\left\| G(x_{1})  \textstyle{\frac{\mathbbm{1}_{B}}{\mu(B)}} + \theta_{1} T\left( x_{1} \textstyle{\frac{\mathbbm{1}_{B}}{\mu(B)}}\right) \right\| \: d \mu} + \int_{\cup{\pi} \setminus B}{\left\|T(x_{1} \textstyle{\frac{\mathbbm{1}_{B}}{\mu(B)}}) \right\| \: d \mu}\\
& \geqslant \left\| G(x_{1}) + \theta_{1} \int_{B}{T(x_{1} \textstyle{\frac{\mathbbm{1}_{B}}{\mu(B)}}) \: d \mu} \right\| + \int_{A_{0} \setminus B}{ \left\| T(x_{1} \textstyle{\frac{\mathbbm{1}_{B}}{\mu(B)}})\right\| \: d \mu} \\
&\hspace{7cm}+ \sum_{A \in \pi \setminus \{ A_{0}\}}{\int_{A}{\left\| T(x_{1} \textstyle{\frac{\mathbbm{1}_{B}}{\mu(B)}}) \right\| \: d \mu}}\\
& \geqslant \left\| G(x_{1}) + \theta_{1} \int_{A_{0}}{T(x_{1} \textstyle{\frac{\mathbbm{1}_{B}}{\mu(B)}}) \: d \mu} \right\| + \sum_{A \in \pi \setminus \{ A_{0}\}}{\left\|y_{A}^{\ast}\left( \int_{A}{ T(x_{1} \textstyle{\frac{\mathbbm{1}_{B}}{\mu(B)}}) \: d \mu} \right) y_{A_{0}} \right\|}\\
& \geqslant \left\| G(x_{1}) + \theta_{1} S(x_{1}) \right\| > 2 - \eps.\qedhere
\end{split}
\end{equation*}
\end{proof}

Our next aim is to deal with lushness. To this end, we will make use of the following immediate numerical result which we prove for the sake of completeness.

\begin{Lemm}\label{Lemm:stabilityL1Technical}
Let $\eps>0$, $\delta > 0$, and let $\lambda_{i} \geqslant 0$ for $i=1, \ldots, n$. Suppose that $\alpha_i,\beta_i \in \R$ are such that $\alpha_{i} \leqslant \beta_{i}$ for $i=1, \ldots, n$ and satisfy $\left(\sum_{i=1}^{n}{\lambda_{i} \beta_{i}}\right) - \eps \delta \leqslant \sum_{i=1}^{n}{\lambda_{i} \alpha_{i}}$. Then,
$$
\sum\{\lambda_{i} \colon \alpha_{i}\leqslant \beta_{i} - \eps\}< \delta.
$$
In particular, if $\sum_{i=1}^{n}{\lambda_{i}} = 1$, then
$$
\sum\{\lambda_{i} \colon \alpha_{i}> \beta_{i} - \eps\}  > 1- \delta.
$$
\end{Lemm}

\begin{proof}
Calling $I=\{ 1\leqslant i\leqslant n \colon  \alpha_{i}> \beta_{i} - \eps\}$ it suffices to observe that
$$
\left(\sum_{i=1}^{n}{\lambda_{i} \beta_{i}}\right) - \eps \delta \leqslant \sum_{i=1}^{n}{\lambda_{i} \alpha_{i}} \leqslant \sum_{i \in I}{\lambda_{i} \beta_{i}} + \sum_{i \notin I}{\lambda_{i} (\beta_{i} - \eps)} = \sum_{i=1}^{n}{\lambda_{i} \beta_{i}} - \eps \sum_{i \notin I}{\lambda_{i}}
$$
from where it easily follows that $\sum_{i \notin I}{\lambda_{i}} < \delta$. The last claim is clear.
\end{proof}

\begin{proof}[Proof of Theorem \ref{Theo:compositionOperatorL1}.(b)]
Assume that $G$ is lush. To check that $\widetilde{G}$ is lush, we just have to show that Proposition \ref{Prop:characterization-lushness}.(iii) is satisfied. Fix $\eps> 0$, $g_{0} \in S_{L_{1}(\mu,Y)}$ and  $f_{0} \in \mathcal{B}$ of the form $f_{0}=x_{0} \mathbbm{1}_{B}/ \mu(B)$ for some $x_{0} \in S_{X}$ and $B \in \Sigma^{+}$. By density, we can assume that $$g_{0}= \sum_{A \in \pi}{y_{A} \frac{\mathbbm{1}_{A}}{\mu(A)}}$$
where $\pi \subset \Sigma^{+}$ is a finite partition of $\Omega$ and $y_{A} \in Y$ satisfy that $\sum_{A \in \pi}{\| y_{A}\|} = 1$. By Proposition \ref{Prop:characterization-lushness}.(iii), the lushness of $G$ lets us find for each $A \in \pi$ an element $y_{A}^{\ast} \in S_{Y^{\ast}}$ such that $\Real{y_{A}^{\ast}(y_{A})} \geqslant (1 - \eps) \| y_{A}\|$ and
\begin{equation}\label{equa:lushL1stabilityAux1}
\dist{\big( x_{0}, \aconv{\big( \GS(S_{X}, G^{\ast}y_{A}^{\ast}, \eps) \big)} \big)}< \eps.
\end{equation}
Let $ h^{\ast} := \sum_{A \in \pi}{y_{A}^{\ast} \mathbbm{1}_{A}}$, which satisfies that $h^{\ast} \in S(S_{L_{\infty}(\mu, Y^{\ast})}, g_{0}, \eps)$ as
$$ \Real{h^{\ast}(g_{0})} = \sum_{A \in \pi}{\Real{y_{A}^{\ast}(y_{A}  )}} >  \sum_{A \in \pi}{(1 - \eps)\| y_{A}\|} = 1 - \eps. $$
Our aim is to prove now that
\begin{equation}\label{equa:lushL1stabilityAux2}
 \dist\Big(f_{0}, \aconv{\Big(\GS(B_{L_{1}(\mu, X)}, \widetilde{G}^{\ast}h^{\ast}, \eps)\Big)} \Big) < \eps,
 \end{equation}
which will finish the proof. First notice that for each $A \in \pi$ with $\mu(B \cap A) \neq 0$ we have that
\begin{equation}\label{equa:lushL1stabilityAux3}
\GS( S_{X}, G^{\ast} y_{A}^{\ast}, \eps) \frac{\mathbbm{1}_{B \cap A}}{\mu(B \cap A)} \subset \GS(B_{L_{1}(\mu, X)}, \widetilde{G}^{\ast}h^{\ast}, \eps),
\end{equation}
since every $x_{A} \in \GS(S_{X}, G^{\ast}y_{A}^{\ast}, \eps)$ satisfies
$$
\Real{\widetilde{G}^{\ast} h^{\ast} \Big( x_{A} \textstyle{\frac{\mathbbm{1}_{B \cap A}}{\mu(B \cap A)}}\Big)} = \Real{G^{\ast}y_{A}^{\ast}(x_{A})} > 1- \eps.
$$
In particular, if for each $A \in \pi$ we take an element $x_{A} \in \aconv{\GS( S_{X}, G^{\ast} y_{A}^{\ast}, \eps)}$ satisfying $\| x_{A} - x_{0}\| < \eps$, which exists by \eqref{equa:lushL1stabilityAux1}, then the inclusion in \eqref{equa:lushL1stabilityAux3} yields that
\begin{align*}
f &:= \sum_{\substack{A \in \pi\\ \mu(A \cap B) \neq 0}}{x_{A} \frac{\mathbbm{1}_{B \cap A}}{\mu(B)}} \\
&= \sum_{\substack{A \in \pi\\ \mu(A \cap B) \neq 0}}{\frac{\mu(B \cap A)}{\mu(B)} x_{A} \frac{ \mathbbm{1}_{B \cap A}}{\mu(B \cap A)}}  \in \aconv{\GS(B_{L_{1}(\mu, X)}, \widetilde{G}^{\ast}h^{\ast}, \eps)}.
\end{align*}
Moreover we have that
$$
\left\| f - f_{0} \right\| = \left\| f - x_{0} \frac{\mathbbm{1}_{B}}{\mu(B)} \right\| = \sum_{\substack{A \in \pi\\ \mu(A \cap B) \neq 0}}{\left\| x_{A} - x_{0} \right\| \frac{\mu(B \cap A)}{\mu(B)}} < \eps
$$
which shows that \eqref{equa:lushL1stabilityAux2} holds.

Let us see the converse: To check that $G$ is lush we will show that condition (iii) of Proposition \ref{Prop:characterization-lushness} is satisfied, for which let $0 < \eps < 1/8$, $x_{0} \in S_{X}$ and $y_{0} \in S_{Y}$. The mentioned condition applied to the lush operator $\widetilde{G}$ for $\eps$, $x_{0} \mathbbm{1}_{\Omega} \in S_{L_{1}(\mu, X)}$ and $y_{0} \mathbbm{1}_{\Omega} \in S_{L_{1}(\mu, Y)}$ provides $n \in \N$, functions
\begin{equation}\label{equa:lushL1stabilityAux4} g^{\ast} \in S(\mathcal{A}, y_{0} \mathbbm{1}_{\Omega}, \eps^{3}) \qquad \text{and} \qquad f_{1}, \ldots, f_{n} \in \GS(S_{L_{1}(\mu, X)}, \widetilde{G}^{\ast}g^{\ast}, \eps^{3})
\end{equation}
and scalars $\theta_{1}, \ldots, \theta_{n} \in \mathbb{T}$, $\lambda_{1}, \ldots, \lambda_{n} \in [0,1]$ with $\sum_{i=1}^{n}{\lambda_{i}} = 1$ satisfying that
\begin{equation}\label{equa:lushL1stabilityAux5}
 \left\| x_{0} \mathbbm{1}_{\Omega} - \sum_{i=1}^{n}{\lambda_{i} \theta_{i} f_{i}} \right\| < \eps^{3}.
 \end{equation}
By density, we can assume that the functions $f_{i}$ are simple and moreover that there is a finite partition $\{A_{1}, \ldots, A_{m}\} \subset \Sigma^{+}$ of $\Omega$ such that
$$
g^{\ast} = \sum_{j=1}^{m}{y_{j}^{\ast} \mathbbm{1}_{A_{j}}} \qquad  \text{and} \qquad  f_{i} = \sum_{j=1}^{m}{x_{i,j} \mathbbm{1}_{A_{j}}} \qquad (i=1, \ldots, n)
$$
where $y_{j}^{\ast} \in S_{Y^{\ast}}$ and $x_{i,j} \in X$ for every $i,j$. Then, conditions \eqref{equa:lushL1stabilityAux4} and \eqref{equa:lushL1stabilityAux5} can be rewritten as
\begin{equation}\label{equa:lushL1stabilityAux6}
1 - \eps^{3} < \Real{g^{\ast}(y_{0} \mathbbm{1}_{\Omega})} = \sum_{j=1}^{m}{\Real{y_{j}^{\ast}(y_{0})} \mu(A_{j})},
\end{equation}
\begin{equation}\label{equa:lushL1stabilityAux7}
\sum_{j=1}^{m}{\mu(A_{j})  \left\| x_{0} - \sum_{i=1}^{n}{\lambda_{i}\theta_{i} x_{i,j}}\right\|}  < \eps^{3},
\end{equation}
and
\begin{align}\label{equa:lushL1stabilityAux8}
& 1 - \eps^{3} < \Real{\widetilde{G}^{\ast}g^{\ast}(f_{i})} =  \sum_{j=1}^{m}{\Real G^{\ast}{y_{j}^{\ast}(x_{i,j})} \mu(A_{j})} \\
&\phantom{ 1 - \eps^{3} < \Real{\widetilde{G}^{\ast}g^{\ast}(f_{i})}}\leqslant \sum_{j=1}^{m}{\| x_{i,j}\| \mu(A_{j})} =1 \qquad \qquad  (i=1, \ldots, n)\notag.
\end{align}
Applying Lemma \ref{Lemm:stabilityL1Technical} to \eqref{equa:lushL1stabilityAux6} and \eqref{equa:lushL1stabilityAux7}, we obtain respectively that
\begin{equation} \label{equa:lushL1stabilityAux9}
\sum{\left\{ \mu(A_{j}) \colon \Real{y_{j}^{\ast}(y_{0})} > 1 - \eps \right\}} > 1 - \eps^{2}
\end{equation}
and
\begin{equation} \label{equa:lushL1stabilityAux10}
\sum{\left\{ \mu(A_{j}) \colon \left\| x_{0} - \sum_{i=1}^{n}{\lambda_{i}\theta_{i} x_{i,j}} \right\| < \eps \right\}} > 1-\eps^{2}.
\end{equation}
Using that $\| x_{0}\| = 1$, the last inequality yields in particular that
\begin{equation}\label{equa:lushL1stabilityAux11}
\sum{\left\{ \mu(A_{j}) \colon  \sum_{i=1}^{n}{\lambda_{i} \|x_{i,j}\|} > 1- \eps \right\}}  > 1-\eps^{2}.
\end{equation}
Combining the relations of \eqref{equa:lushL1stabilityAux8} in a convex sum with the $\lambda_{i}$'s as coefficients we obtain that
\begin{equation}\label{equa:lushL1stabilityAux12}
1 - \eps^{3} < \sum_{j=1}^{m}{\mu(A_{j}) \sum_{i=1}^{n}{\lambda_{i} \Real G^{\ast}{y_{j}^{\ast}(x_{i,j})}}} \leqslant \sum_{j=1}^{m}{\mu(A_{j}) \sum_{i=1}^{n}{\lambda_{i} \| x_{i,j}\|}} = 1.
\end{equation}
Actually, from the right-hand equality of the previous expression we get that
\begin{align*}
1 & = \sum_{i=1}^{n}{\lambda_{i} \, \sum_{j=1}^{m}{\| x_{i,j}\| \mu(A_{j})} } = \sum_{j=1}^{m}{\mu(A_{j}) \, \sum_{i=1}^{n}{\lambda_{i} \|x_{i,j}\|}}\\
& \geqslant (1 + \eps) \sum{\left\{ \mu(A_{j}) \colon  \sum_{i=1}^{n}{\lambda_{i} \|x_{i,j}\|} > 1+ \eps \right\}}\\
& + (1- \eps) \sum{\left\{ \mu(A_{j}) \colon 1 - \eps < \sum_{i=1}^{n}{\lambda_{i} \|x_{i,j}\|} \leqslant 1 + \eps\right\}},
\end{align*}
which, together with \eqref{equa:lushL1stabilityAux11}, implies that the number
$$ \alpha:=  \sum{\left\{ \mu(A_{j}) \colon 1 - \eps < \sum_{i=1}^{n}{\lambda_{i} \| x_{i,j}\| \leqslant 1 + \eps} \right\}}$$
satisfies the relation $$(1 + \eps) (1 - \eps^{2} - \alpha) + (1 - \eps) \alpha \leqslant 1.$$ A simple computation shows that necessarily $1 - \eps - \eps^{2} \leqslant 2 \alpha$, and so
\begin{equation}\label{equa:lushL1stabilityAux13}
\sum{\left\{ \mu(A_{j}) \colon 1 - \eps < \sum_{i=1}^{n}{\lambda_{i} \| x_{i,j}\| \leqslant 1 + \eps} \right\}} = \alpha \geqslant \frac{1}{2} - \eps.
\end{equation}
On the other hand, an application of Lemma \ref{Lemm:stabilityL1Technical} to the left-hand side part of \eqref{equa:lushL1stabilityAux12} gives that
\begin{equation}\label{equa:lushL1stabilityAux14}
\sum{\left\{ \mu(A_{j}) \colon \left(\sum_{i=1}^{n}{\lambda_{i} \| x_{i,j}\|} \right) - \eps^{2} <  \sum_{i=1}^{n}{\lambda_{i} \Real G^{\ast}{y_{j}^{\ast}(x_{i,j})}}  \right\}} >1-\eps.
\end{equation}
Using that $\eps < 1/8$ combined with \eqref{equa:lushL1stabilityAux9}, \eqref{equa:lushL1stabilityAux10}, \eqref{equa:lushL1stabilityAux13}, and \eqref{equa:lushL1stabilityAux14}, we deduce the existence of some $j_{0} \in \N$ satisfying simultaneaously
\begin{equation*}
 1 - \eps < \Real{y_{j_{0}}^{\ast}(y_{0})},
\end{equation*}
\begin{equation}\label{equa:lushL1stabilityAux12.3}
\left\| x_{0} - \sum_{i=1}^{n}{\lambda_{i}\theta_{i} x_{i,j_{0}}} \right\| < \eps,
\end{equation}
\begin{equation}\label{equa:lushL1stabilityAux12.5}
1 - \eps < \sum_{i=1}^{n}{\lambda_{i} \| x_{i, j_{0}}\| \leqslant 1 + \eps},
\end{equation}
and
\begin{equation}\label{equa:lushL1stabilityAux15}
\left(\sum_{i=1}^{n}{\lambda_{i} \| x_{i,j_{0}}\|}\right) - \eps^{2} <  \sum_{i=1}^{n}{\lambda_{i} \Real{G^{\ast} y_{j_{0}}^{\ast}(x_{i,j_{0}})}}.
\end{equation}
If we denote
$$
I := \{ 1 \leqslant i \leqslant n \colon \Real{G^{\ast} y^{\ast}_{j_{0}}(x_{i, j_{0}})} > \| x_{i, j_{0}}\| (1 - \eps) \},
$$
then again we can apply Lemma \ref{Lemm:stabilityL1Technical} to \eqref{equa:lushL1stabilityAux15} with
$$
\beta_{i} = 1, \quad \alpha_{i} = \frac{\Real{G^{\ast}y_{j_0}^{\ast}(x_{i,j_{0}})}}{\| x_{i,j_{0}}\|},\, \quad \text{and} \quad \left(  \sum_{i=1}^{n}{\lambda_{i} \| x_{i,j_{0}}\| \beta_{i}}\right) - \eps^{2} < \sum_{i=1}^{n}{\lambda_{i} \| x_{i, j_{0}}\| \alpha_{i}}
$$
to get that
$$
\sum_{i \notin I}{\lambda_{i} \| x_{i,j_{0}}\|} < \eps.
$$
This, together with \eqref{equa:lushL1stabilityAux12.5}, yields that
\begin{equation}\label{equa:lushL1stabilityAux16}
1 - 2\eps <  \sum_{i \in I}{\lambda_{i} \| x_{i, j_{0}}\|} \leqslant 1 + \eps.
\end{equation}
Consider now the elements
$$
\tilde{x}_{i} := \frac{x_{i, j_{0}}}{\| x_{i, j_{0}}\|} \in S_{X} \qquad  \text{and} \qquad \tilde{\lambda}_{i} :=  \frac{\lambda_{i} \| x_{i, j_{0}}\|}{\sum_{k \in I}{\lambda_{k} \| x_{k, j_{0}}\|}} \geqslant 0
$$
which satisfy
$$
\sum_{i \in I}{\tilde{\lambda}_{i}} = 1 \qquad  \text{and} \qquad \tilde{x}_{i} \in \GS(S_{X}, G^{\ast}y_{j_{0}}^{\ast}, \eps)\qquad  \mbox{ for each $i \in I$}.
$$
Finally, using \eqref{equa:lushL1stabilityAux12.3} and \eqref{equa:lushL1stabilityAux16} we conclude that
\begin{align*}
\left\| x_{0} - \sum_{i \in I}{\tilde{\lambda}_{i}\theta_i \tilde{x}_{i}} \right\| & \leqslant \eps + \left\| \sum_{i=1}^{n}{\lambda_{i}\theta_{i} x_{i, j_{0}}} - \sum_{i \in I}{\tilde{\lambda}_{i}\theta_i \tilde{x}_{i}} \right\|\\
& \leqslant \eps + \sum_{i \notin I}{\lambda_{i} \| x_{i, j_{0}}\|} +  \sum_{i \in I}{\left\|\lambda_{i}  x_{i, j_{0}} - \tilde{\lambda}_{i} \tilde{x}_{i}\right\|  }\\
& \leqslant \eps + \sum_{i \notin I}{\lambda_{i} \| x_{i, j_{0}}\|} + \left| 1 - \frac{1}{\sum_{k \in I}{\lambda_{k} \| x_{k, j_{0}}\|}} \right| \sum_{i \in I}{\lambda_{i} \| x_{i, j_{0}}\|}\\
& = \eps + \sum_{i \notin I}{\lambda_{i} \| x_{i, j_{0}}\|} + \left| 1 - \sum_{i \in I}{\lambda_{i} \| x_{i, j_{0}}\|} \right| \leqslant 4 \eps
\end{align*}
which finishes the proof.
\end{proof}

\begin{proof}[Proof of Theorem \ref{Theo:compositionOperatorL1}.(c)]
Let us fix an atom $A_{0} \in \Sigma^{+}$. Assume that $\widetilde{G}$ has the aDP. We will show that $G$ satisfies condition (iii) of Theorem \ref{Theo:aDPCharacterization}: Let $x_{0} \in S_{X}$, $y_{0} \in S_{Y}$ and $\eps > 0$. By hypothesis, we have that
$$ x_{0} \frac{\mathbbm{1}_{A_{0}}}{\mu(A_{0})} \in \overline{\aconv}{\left( \left\{ f \in \mathcal{B} \colon \left\| \widetilde{G}(f) + y_{0} \frac{\mathbbm{1}_{A_{0}}}{\mu(A_{0})} \right\|> 2 - \eps \, \mu(A_{0}) \right\} \right)}. $$
Then, for each $\eta \in (0,1)$ we can find a finite family $\FF \subset \Sigma^{+}$, elements $x_{B} \in S_{X}$ satisfying
\begin{equation}\label{equa:L1StabilityAux1}
2 - \eps \,  \mu(A_{0})  < \left\| G(x_{B}) \frac{\mathbbm{1}_{B}}{\mu(B)} + y_{0} \frac{\mathbbm{1}_{A_{0}}}{\mu(A_{0})} \right\| \qquad \mbox{ for every $B \in \FF$};
\end{equation}
and scalars $\lambda_{B} \in \mathbb{K}$ with $\sum_{B \in \FF}{|\lambda_{B}|} = 1$ such that
\begin{equation} \label{equa:L1StabilityAux1.5}
\left\| x_{0} \frac{\mathbbm{1}_{A_{0}}}{\mu(A_{0})} - \sum_{B \in \FF}{\lambda_{B}  x_{B} \frac{\mathbbm{1}_{B}}{\mu(B)}}\right\| < \eta.
\end{equation}
But $\mu(A_{0} \cap B)$ is either $0$ or $\mu(A_{0})$ for each $B \in \FF$ as $A_{0}$ is an atom. Then, if we just integrate in  \eqref{equa:L1StabilityAux1.5} over the atom $A_{0}$ we will get that
\begin{equation}\label{equa:L1StabilityAux2}
\left\| x_{0} - \sum_{\substack{B \in \FF\\ \mu(A_{0} \setminus B) = 0}}{\lambda_{B} x_{B} \frac{\mu(A_{0})}{\mu(B)}} \right\| <\eta.
\end{equation}
In particular, this yields that
\begin{equation}\label{equa:L1StabilityAux3} \alpha := \sum_{\substack{B \in \FF\\ \mu(A_{0} \setminus B) = 0}}{|\lambda_{B}| \frac{\mu(A_{0})}{\mu(B)}} > 1 - \eta. \end{equation}
Combining \eqref{equa:L1StabilityAux2} and \eqref{equa:L1StabilityAux3}, we deduce that
\begin{equation}\label{equa:L1StabilityAux3.5}
\begin{split}
\left\| x_{0} - \sum_{\substack{B \in \FF\\ \mu(A_{0} \setminus B) = 0}}{\frac{\lambda_{B} \mu(A_{0})}{\alpha \mu(B)}  x_{B}} \right\| & \leqslant \eta + \left\| \sum_{\substack{B \in \FF\\ \mu(A_{0} \setminus B) = 0}}{\Big(\frac{1}{\alpha} - 1\Big) \lambda_{B} x_{B} \frac{\mu(A_{0})}{\mu(B)}} \right\|\\
& \leqslant \eta + \frac{1 - \alpha}{\alpha} \leqslant \eta + \frac{\eta}{1 - \eta}.
\end{split}
\end{equation}
Since $\eta  \in (0,1)$ was arbitrary, we deduce from \eqref{equa:L1StabilityAux3.5} that $x_{0}$ belongs to the closed absolute convex hull of the set of all $x \in S_{X}$ for which there is $B \in \Sigma^{+}$ satisfying  $\mu(A_{0} \setminus B) = 0$ and
\begin{equation*} 
\left\| G(x) \frac{\mathbbm{1}_{B}}{\mu(B)} + y_{0} \frac{\mathbbm{1}_{A_{0}}}{\mu(A_{0})} \right\| > 2 - \eps \, \mu(A_{0}).
\end{equation*}
But then such elements $x$ and $B$ satisfy in particular that
\begin{equation*}
\begin{split}
2 - \eps\,  \mu(A_{0}) & < \| G(x)\| \frac{\mu(B \setminus A_{0})}{\mu(B)} + \left\| G(x) \frac{\mu(A_{0})}{\mu(B)} + y_{0} \right\|\\
& \leqslant \| G(x) \| \frac{\mu(B \setminus A_{0})}{\mu(B)} + \| y_{0} \| \frac{\mu(B \setminus A_{0})}{\mu(B)} + \frac{\mu(A_{0})}{\mu(B)} \| G(x) + y_{0} \|\\
& \leqslant  \frac{ \mu(B \setminus A_{0})}{\mu(B)} \, 2 + \frac{ \mu(A_{0})}{\mu(B)} \, \| G(x) + y_{0} \|,
\end{split}
\end{equation*}
and hence
$$ 2 - \eps \leqslant 2 - \mu(B) \, \eps \leqslant  \| G(x) + y_{0} \|. $$
We then conclude that
$$x_{0} \in \overline{\aconv}{(\{ x \in B_{X} \colon \| G(x) + y_{0}\| > 2 - \eps \})}.$$
Let us prove now the converse of (c). We remark here that this implication does not use that $\mu$ is atomic. Assuming that $G$ has the aDP, we will now check that $\widetilde{G}$ satisfies Theorem~\ref{Theo:aDPCharacterization}.(ii). For this, it is enough to prove that given a simple function $g_{0} \in S_{L_{1}(\mu, Y)}$ of the form
$$
g_{0} = \sum_{A \in \pi}{y_{A} \frac{\mathbbm{1}_{A}}{\mu(A)}},
$$
where $\pi \subset \Sigma^{+}$ is a finite partition of $\Omega$ and $y_{A} \in Y$ ($A \in \pi$), we have that
$$
\mathcal{B} \subset \overline{\aconv}{(\{ f \in \mathcal{B} \colon \| \widetilde{G}(f) +  g_{0} \|>2-\eps \})}.
$$
Let $x_{0} \in S_{X}$ and $B \in \Sigma^{+}$. Then
$$
x_{0} \frac{\mathbbm{1}_{B}}{\mu(B)} = \sum_{\substack{A \in \pi\\ \mu(A \cap B) \neq 0 }}{\frac{\mu(B \cap A_{j})}{\mu(B)} \, x_{0} \frac{\mathbbm{1}_{B \cap A}}{\mu(B \cap A)}},
$$
and so in order to show that
$$
x_{0} \frac{\mathbbm{1}_{B}}{\mu(B)} \in  \overline{\aconv}{(\{ f \in \mathcal{B} \colon \| \widetilde{G}(f) +  g_{0} \|>2-\eps \})}
$$
we can assume without loss of generality, by using a standard convexity argument, that $B$ is contained in some $A_{0} \in \pi$. Since $G$ has the aDP, using Theorem \ref{Theo:aDPCharacterization}.(iii), for each $\delta > 0$ there is a finite set
$$
F \subset \{ x \in S_{X} \colon \left\| G(x) + y_{A_{0}} \right\| > 1 + \| y_{A_{0}}\| - \eps  \}
$$
such that $\dist{(x_{0}, \aconv{F})} < \delta$.  In particular, this implies that
\begin{equation}\label{equa:L1StabilityAux5}
\dist\left(x_{0} \textstyle{\frac{\mathbbm{1}_{B}}{\mu(B)}}, \aconv{\left\{ x \frac{\mathbbm{1}_{B}}{\mu(B)} \colon x \in F \right\}}\right) < \delta.
\end{equation}
Finally, notice that each $x \in F$ satisfies that
\begin{equation*}
\begin{split}
\| \widetilde{G}(x \textstyle{\frac{\mathbbm{1}_{B}}{\mu(B)}}) +g_{0} \| & = \int_{B}{\| G(x) \textstyle{\frac{\mathbbm{1}_{B}}{\mu(B)}} +y_{A_{0}} \mathbbm{1}_{A_{0}} \| \: d \mu} + \int_{\cup{\pi} \setminus B}{\| g_{0}\| \: d\mu} \\
& = \| G(x) + y_{A_{0}} \mu(B)\| + \| y_{A_{0}}\| \mu(A_{0} \setminus B) +\sum_{A \in \pi \setminus \{ A_{0}\}}{\| y_{A}\| \mu(A)} \\
& \geqslant 1 + \| y_{A_{0}}\| \mu(B) - \eps + \| y_{A_{0}}\| \mu(A_{0} \setminus B) +\sum_{A \in \pi \setminus \{ A_{0}\}}{\| y_{A}\| \mu(A)}\\
& = 1 + \| g_{0}\| - \eps = 2 - \eps.
\end{split}
\end{equation*}
Therefore, \eqref{equa:L1StabilityAux5} leads to
$$
\dist\left( x_{0} \textstyle{\frac{\mathbbm{1}_{B}}{\mu(B)}},  \aconv{(\{ f \in \mathcal{B} \colon \| \widetilde{G}(f) +  g_{0} \|>2-\eps \})}\right) < \delta
$$
for arbitrary $\delta > 0$.
\end{proof}

\begin{proof}[Proof of Theorem \ref{Theo:compositionOperatorL1}.(d)]
Assuming hat $\mu$ has no atoms, we claim that given a simple function $g_{0} \in S_{L_{1}(\mu, Y)}$, for every $\delta > 0$ we can write $g_{0}$ as
\begin{equation}\label{equa:L1StabilityAux6}
g_{0} = \sum_{A \in \pi}{y_{A} \frac{\mathbbm{1}_{A}}{\mu(A)}}
\end{equation}
where $\pi \subset \Sigma^{+}$ is a finite partition of $\Omega$ and the coefficients $y_{A} \in \delta B_{Y}$ for each $A \in \pi$. Let us check this: of course, we can write $g_{0}$ as in \eqref{equa:L1StabilityAux6} for a partition $\pi \subset \Sigma^{+}$ and elements $y_{A} \in Y$ with $\sum_{A \in \pi}{\| y_{A}\|} = 1$. But since $\mu$ has no atoms, we can find for each $A \in \pi$ a partition of $A$ into elements $C \in \Sigma^{+}$ satisfying $\mu(C) \leqslant \delta \mu(A)$. If $\pi'$ is the collection of all such subsets, then this is a finer partition than $\pi$ and
$$ g_{0} = \sum_{A \in \pi}{y_{A} \frac{\mathbbm{1}_{A}}{\mu(A)}} = \sum_{A \in \pi}{\sum_{\substack{C \in \pi'\\ C \subset A}}{\left(\frac{\mu(C)}{\mu(A)}y_{A}\right) \frac{\mathbbm{1}_{C}}{\mu(C)}}} $$
and the proof of the claim is over.

Let us then prove now that if
\begin{equation}\label{equa:compositionOperatorL1Aux1}
B_{X} = \overline{\aconv}{\{ x \in B_{X} \colon \| Gx\| > 1 - \eps \}}
\end{equation}
for every $\eps > 0$, then $\widetilde{G}$ has the aDP. By Theorem \ref{Theo:aDPCharacterization}.(iii), it is enough to show that given a simple function $g_{0} \in S_{L_{1}(\mu, Y)}$ as in \eqref{equa:L1StabilityAux6} and $\eps > 0$, we have that
$$ \mathcal{B} \subset \aconv{\left\{ f \in \mathcal{B} \colon \| \widetilde{G}(f) + g_{0}\| > 2 - \eps \right\}} + \delta B_{L_{1}(\mu, X)} $$
for every $\delta> 0$. Let $x_{0} \in S_{X}$, $B \in \Sigma^{+}$ and $0 < \delta < \eps /3$. Let $\pi \subset \Sigma^{+}$ and  $y_{A} \in \delta B_{Y}$ ($A \in \pi$) as in the claim above for the given $g_{0}$. We can moreover assume that every $A \in \pi$ is either contained in $B$ or in $\Omega \setminus B$. Using \eqref{equa:compositionOperatorL1Aux1}, we can find $m \in \N$, $x_{j} \in B_{X}$ with $\| G (x_{j})\| > 1 - \delta$ and $\lambda_{j}\in \K$ $ (j=1, \ldots,m)$ such that $\sum_{j=1}^{m}|\lambda_{j}| = 1$ and
$$ \left\| x_{0} - \sum_{j=1}^{m}{\lambda_{j} x_{j}} \right\| < \delta. $$
Then, it is easy to check that
\begin{equation}\label{equa:L1StabilityAux7}
\left\| x_{0} \frac{\mathbbm{1}_{B}}{\mu(B)}  - \sum_{\substack{A \in \pi\\ A \subset B}}{\sum_{j=1}^{m}{ \left(\frac{\mu(A)}{\mu(B)} \lambda_{j}\right) x_{j} \frac{\mathbbm{1}_{A}}{\mu(A)}}} \right\| < \delta.
\end{equation}
This shows that $x_{0} \mathbbm{1}_{B}/\mu(B)$ is $\delta$-approximated by an absolutely convex sum of elements of the form $x_{j} \mathbbm{1}_{A}/\mu(A)$ for some $1 \leqslant j \leqslant m$ and $A \in \pi$. Finally, notice that every such element satisfies that
\begin{equation}\label{equa:L1StabilityAux8}
\begin{split}
\left\| \widetilde{G}\left(x_{j} \frac{\mathbbm{1}_{A}}{\mu(A)}\right) + g_{0} \right\| & = \left\| G(x_{j}) \frac{\mathbbm{1}_{A}}{\mu(A)} + g_{0} \right\|\\
& = \| G(x_{j}) + y_{A}\| + \sum_{A' \in \pi, A' \neq A}{\| y_{A'}\|}\\
& \geqslant \| G(x_{j})\| - \delta + 1 - \delta > 2 - 3 \delta
\end{split}
\end{equation}
Therefore,
$$ x_{0} \frac{\mathbbm{1}_{B}}{\mu(B)} \in \aconv{\left\{ f \in B_{L_{1}(\mu, X)} \colon \| \widetilde{G}(f) + g_{0}\| > 2 - 3 \delta \right\}} + \delta B_{L_{1}(\mu, X)}. $$
Using that $0 < \delta < \eps /3$ was arbitrary, we conclude the result.

Conversely, suppose that $\widetilde{G}$ has the aDP and let $\eps > 0$, $x_{0} \in S_{X}$ and $y_{0} \in S_{Y}$. By Theorem~\ref{Theo:aDPCharacterization}.(iii),  we have that
$$ x_{0} \mathbbm{1}_{\Omega} \in \overline{\aconv}{\left\{ f \in \mathcal{B} \colon \left\| \widetilde{G}(f) + y_{0} \mathbbm{1}_{\Omega} \right\| > 2 - \eps  \right\}}. $$
Therefore, given any $\delta> 0$ there is a finite set $\FF \subset \Sigma^{+}$ and elements $x_{A} \in S_{X}$, $\lambda_{A} \in \K$ ($A\in \FF$) such that $\sum_{A \in \FF}{|\lambda_{A}|} = 1$ satisfying that
\begin{equation}\label{equa:L1StabilityAux9}
\left\| \widetilde{G}\left(x_{A} \frac{\mathbbm{1}_{A}}{\mu(A)}\right) + y_{0} \frac{\mathbbm{1}_{B}}{\mu(B)} \right\| > 2 - \eps
\end{equation}
and
\begin{equation}\label{equa:L1StabilityAux10}
\left\| x_{0} \frac{\mathbbm{1}_{B}}{\mu(B)} - \sum_{A \in \FF}{\lambda_{A} x_{A} \frac{\mathbbm{1}_{A}}{\mu(A)}} \right\| < \delta.
\end{equation}
It easily follows from \eqref{equa:L1StabilityAux9} that
\begin{equation*} 
1 - \eps < \left\| G(x_{A}) \frac{\mathbbm{1}_{A}}{\mu(A)} + y_{0} \frac{\mathbbm{1}_{B}}{\mu(B)} \right\| -1 \leqslant \| G(x_{A})\|.
\end{equation*}
On the other hand, \eqref{equa:L1StabilityAux10} yields that
$$ \left\| x_{0} - \sum_{A \in \FF}{\lambda_{A} x_{A}} \right\| = \left\| \int_{\Omega}{ \left( x_{0} \frac{\mathbbm{1}_{B}}{\mu(B)} - \sum_{A \in \FF}{\lambda_{A} x_{A} \frac{\mathbbm{1}_{A}}{\mu(A)}}\right) \: d \mu } \right\| < \delta. $$
Therefore,
$$ \dist{(x_{0}, \aconv{\{ x \in B_{X} \colon \| G x\| > 1 - \eps \}})} < \delta, $$
and since $\delta > 0$ and $x_{0} \in S_{X}$ were arbitrary, we conclude that \eqref{equa:compositionOperatorL1Aux1} holds.
\end{proof}

We may use Theorem~\ref{Theo:compositionOperatorL1}.(d) to produce an example of an operator with the aDP which does not attain its norm. Recall that it was proved in Proposition \ref{Prop:lush-attains-norm} that this cannot happen if the operator is actually lush.

\begin{Exam}\label{example-aDP-no-NA}
Let $G\in L(c_0)$ be the operator given by
$$
(Gx)(n)=\frac{n-1}{n}x(n) \qquad (x\in c_0,n\in \N).
$$
The operator $\widetilde{G}\in L(L_1([0,1],c_0))$ defined by $\widetilde{G}(f)=G\circ f$ for $f\in L_1([0,1],c_0)$ has the aDP and it does not attain its norm.
\end{Exam}

\begin{proof}
It is clear that $\|G\|=1$ and that it does not attain its norm. Moreover, for each $\eps>0$ it easy to check that
\begin{equation}\label{eq:Exam-aDP-operator-not-norm-attaining}
B_{c_0}=\overline{\aconv}\{x\in B_{c_0} \ : \ \|Gx\|>1-\eps\}.
\end{equation}
Indeed, fixed $\eps>0$, $x\in B_{c_0}$, and $\delta>0$, we take $N\in\N$ large enough to satisfy
$$
|x(N)|<\delta \qquad \text{and} \qquad \frac{1}{N}<\eps
$$
and we consider $y,z\in B_{c_0}$ given by
$$
y(n)=z(n)=x(n) \qquad \text{if } n\neq N \qquad \text{and} \qquad y(N)=1, \ z(N)=-1.
$$
Then, it is clear that $\|Gy\|=\|Gz\|\geqslant1-\frac{1}{N}>1-\eps$ and
$$
\left\|x-\frac12(y+z)\right\|=|x(N)|<\delta.
$$
Now the arbitrariness of $\delta>0$ gives that $x\in \overline{\aconv}\{x\in B_{c_0} \ : \ \|Gx\|>1-\eps\}$.
	
Finally, $\widetilde{G}$ has the aDP thanks to Theorem~\ref{Theo:compositionOperatorL1}.(d) as  \eqref{eq:Exam-aDP-operator-not-norm-attaining} holds for every $\eps>0$. Besides, for each non-zero $f\in L_1([0,1],c_0)$ we have that
$$
\|G(f(t))\|<\|f(t)\| \qquad \forall t\in[0,1] \text{ with } f(t)\neq 0
$$
since $G$ does not attain its norm. So we deduce that $\|\widetilde{G}(f)\|<\|f\|$ and, therefore, $\widetilde{G}$ does not attain its norm.
\end{proof}

\section{Target operators, lushness and ultraproducts}

Now, we will prove the stability of target operators and lush operators with respect to the operation of taking ultraproducts. These results extend Corollaries 4.4 and 4.5 of \cite{LushNumerical} about stability of lush spaces with respect to ultraproducts.

Let us recall the basic definitions, taken from \cite{He}.  Let $\mathcal{U}$
be a free ultrafilter on $\N$. The limit of a sequence with respect to the ultrafilter $\mathcal{U}$ is denoted by $\lim_{\mathcal{U}}a_n$, or $\lim_{n, \mathcal{U}}a_n$, if it is necessary to stress that the limit is taken with respect to the variable $n$. Let $\{ X_n \}_{n\in \N}$ be
a sequence of Banach spaces. We can consider the
$\ell_\infty$-sum of the family, $\left[\oplus_{n\in \N}
X_n\right]_{\ell_\infty}$, together with its closed subspace
$$
N({\mathcal{U}}) = \left\{ \{x_n\}_{n\in \N} \in \left[\oplus_{n\in \N}
X_n\right]_{\ell_\infty} \ :\  \lim_{\mathcal{U}} \Vert x_n \Vert=0
\right\}.
$$
The quotient space $(X_n)_{\mathcal{U}}=\left[\oplus_{n\in \N}
X_n\right]_{\ell_\infty}/ N({\mathcal{U}})$ is called the
\emph{ultraproduct} of the family $\{ X_n \}_{n\in \N}$ relative to the
ultrafilter $\mathcal{U}$.
\index{ultraproduct}%
Let $(x_n)_\mathcal{U}$ stand for the element
of $(X_n)_{\mathcal{U}}$ containing a given representative $( x_n ) \in
\left[\oplus_{n\in \N} X_n\right]_{\ell_\infty}$. It is easy to check that
$$
\Vert (x_n)_\mathcal{U}\Vert = \lim_{\mathcal{U}} \Vert x_n \Vert.
$$
Moreover, every  $\tilde x \in (X_n)_{\mathcal{U}}$ can be represented as $\tilde x = (x_n)_{\mathcal{U}}$ in such a way that $\|x_n\| = \|\tilde x\|$ for all $n \in \N$.

If all the $X_n$ are equal to the same Banach space $X$, the ultraproduct
of the family is called the $\mathcal{U}$-\emph{ultrapower} of $X$.
\index{ultrapower}%
We denote this ultrapower by $X_\mathcal{U}$.

Let $\{ X_n \}_{n\in \N}$, $\{ Y_n \}_{n\in \N}$ be two sequences of Banach spaces and let $\{ T_n \}_{n\in \N}$ be a norm-bounded sequence of operators where $T_n \in L(X_n, Y_n)$ for every $n\in \N$. We denote  $(T_n)_{\mathcal{U}}$ the operator that acts from  $(X_n)_{\mathcal{U}}$ to  $(Y_n)_{\mathcal{U}}$ as follows:
 $(T_n)_{\mathcal{U}}(x_n)_\mathcal{U} = (T_nx_n)_{\mathcal{U}}$. Evidently,
 $$
\Vert (T_n)_\mathcal{U}\Vert = \lim_{\mathcal{U}} \Vert T_n \Vert.
$$

Now, we state our main result about ultraproducts.

\begin{Theo}\label{Theor:ultra-target}
Let $\mathcal{U}$ be a free ultrafilter on $\N$, $\{ X_n \}_{n\in \N}$, $\{ Y_n \}_{n\in \N}$, $\{ Z_n \}_{n\in \N}$ be sequences of Banach spaces and let $\{G_n\}_{n\in \N}$, $\{T_n\}_{n\in \N}$ be norm bounded sequences of operators such that $G_n \in S_{L(X_n, Y_n)}$ and $T_n \in L(X_n, Z_n)$ for every $n\in \N$. If each $T_n$ is a target for the corresponding $G_n$ for every $n\in \N$, then $T = (T_n)_{\mathcal{U}} \in L((X_n)_{\mathcal{U}}, (Z_n)_{\mathcal{U}})$  is a target for $G = (G_n)_{\mathcal{U}} \in L((X_n)_{\mathcal{U}}, (Y_n)_{\mathcal{U}})$.
\end{Theo}

We need the following easy remark about the absolutely convex hull of a convex set. In fact, this idea already appeared implicitly in the proof of implication  (i) $\Rightarrow$ (iii) of Corollary \ref{Defi:spearVector}.

\begin{Prop} \label{Prop:conv-aconv}
Let $F \subset B_X$ be a convex set. If $X$ is  a real space, then
$$
\aconv F = \{\lambda_1 x_1 - \lambda_2 x_2 : x_1, x_2 \in F, \lambda_1, \lambda_2 \geqslant 0, \lambda_1 + \lambda_2 = 1\}.
$$
If $X$ is a complex space, then for every   $m  \in \N$  and every  $x \in \aconv F$ there are  $\lambda_1, \ldots, \lambda_m \geqslant 0$,
$\sum_{k=1}^m \lambda_k = 1$ and $x_1, \ldots, x_m \, \in \, F$ such that
\begin{equation}\label{eq.conv-aconv-appr1}
 \left\|x - \sum_{k=1}^m \lambda_k \exp\left(\frac{2\pi i k}{m}\right) x_k\right\| \leqslant
  \frac{2\pi}{m}.
\end{equation}
\end{Prop}

\begin{proof}
We demonstrate only the more complicated complex case. As $x \in \aconv F$ there are $\mu_j \in [0,\,1]$, $j=1, \ldots, N$  with $\sum_{j=1}^N \mu_j = 1$, $\theta_j
\in [0,\,2\pi]$ and $y_j \in F$ satisfying
$$
x = \sum_{j=1}^N  \mu_j \exp\left(i\theta_j\right) y_j .
$$
Taking into account that the points $\left\{\frac{2\pi k}{m}\ : \ k = 1,
\ldots, m\right\}$ form a $\frac{2\pi}{m}$-net of $[0,\,2\pi]$ we can
represent the set of indices $\{1, \ldots, N\}$ as a disjoint union
of sets $A_k$, $k = 1, \ldots, m$ in such a way that
\begin{equation*}\label{eq.appr2}
\left|\,\theta_j - \frac{2\pi k}{m}\,\right| \leqslant \frac{2\pi}{m}
\qquad \text{ for every }\qquad j \in A_k.
\end{equation*}
 Let us show that
$$
 \lambda_k= \sum_{j \in
 A_k}\mu_j, \quad \text{and} \quad x_k= \frac{1}{\lambda_k}
 \sum_{j \in A_k}\mu_j y_j \quad \text{if} \quad A_k \neq
 \emptyset
$$
and
$$
 \lambda_k= 0, \quad \text{and arbitrary} \quad x_k \in F
  \quad \text{if} \quad A_k  = \emptyset
$$
fulfill the desired condition \eqref{eq.conv-aconv-appr1}. Indeed, it is clear that $x_k \in
F$ and $\sum_{k=1}^m \lambda_k = 1$. Now,
\begin{align*}
 \left\|x - \sum_{k=1}^m \lambda_k \exp\left(\frac{2\pi i
 k}{m}\right) x_k\right\| & =\left\|x - \sum_{k: A_k \neq \emptyset}
 \sum_{j \in  A_k}\mu_j \exp\left(\frac{2\pi i k}{m}\right)  y_j \right\|
\\ &
 \leqslant \left\|x - \sum_{k: A_k \neq \emptyset}\sum_{j \in
 A_k}\mu_j \exp(i\theta_j) y_j \right\| + \frac{2\pi}{m} \\ & =
 \left\|x-\sum_{j=1}^N \mu_j\exp(i\theta_j)y_j \right\|  + \frac{2\pi}{m} =
 \frac{2\pi}{m}\,.\qedhere
\end{align*}
\end{proof}

\begin{proof}[Proof of Theorem~\ref{Theor:ultra-target}]
We demonstrate the theorem only for the more complicated complex case. Also, we may and do suppose that $\|T\| = \|T_n\| = 1$ for every $n\in \N$.

Let $x_0 = (x_{0,n})_\mathcal{U} \in B_{(X_n)_{\mathcal{U}}}$,  $y = (y_{n})_\mathcal{U} \in S_{(Y_n)_{\mathcal{U}}}$ and $\eps > 0$ be fixed. Evidently, the ``coordinates''  $x_{0,n}$ can be selected in such a way that  $x_{0,n}\in B_{X_n}$ and $y_{n}\in S_{Y_n}$ for every $n\in \N$. For each $n \in \N$ applying \eqref{eq:diamond} in Definition \ref{Defi:target} for $\eps/2$, $x_{0,n}\in B_{X_n}$ and $y_{n}\in S_{Y_n}$ we obtain the corresponding $F_n \subset B_{X_n}$ satisfying
\begin{equation}\label{eq:diamond++}
\begin{split}
 &\conv{F_n} \subset \left\{ x \in B_{X_n} \colon \| G x + y_n\| > 2 - \frac{\eps}{2} \right\} \quad \mbox{ and } \\ &\dist\bigl(T_nx_{0,n}, T_n\bigl(\aconv(F_n)\bigr)\bigr) < \frac{\eps}{2}.
\end{split}
\end{equation}
Without loss of generality, we assume that  $F_n$ is convex, otherwise we just substitute  $F_n$ by its convex hull. Our choice means that there is $x_n \in \aconv(F_n)$ such that
 $$
\bigl\|T_n x_{0,n} - T_n x_n\bigr\| <  \frac{\eps}{2}.
$$
Select $m \in \N$ such that  $\displaystyle \frac{2\pi}{m} <  \frac{\eps}{2}$.
Using Proposition \ref{Prop:conv-aconv}, we can find for each $n \in \N$ corresponding $\lambda_{n,1} \ldots, \lambda_{n,m} \geqslant 0$,
$\sum_{k=1}^m \lambda_{n,k} = 1$ and $x_{n,1}, \ldots, x_{n,m} \, \in \, F_n$ such that
$$
 \left\|x_n - \sum_{k=1}^m \lambda_{n,k} \exp\left(\frac{2\pi i k}{m}\right) x_{n,k}\right\| <
  \frac{\eps}{2}
$$
and, consequently,
\begin{equation}\label{eq.conv-aconv-apprm}
 \left\|T_n x_{0,n} - \sum_{k=1}^m \lambda_{n,k} \exp\left(\frac{2\pi i k}{m}\right) T_nx_{n,k}\right\| <  \eps.
\end{equation}
For each $k = 1, \ldots , m$ denote $\lambda_k = \lim_{n, \mathcal{U}} \lambda_{n,k}$ and  $\tilde x_k = (x_{n, k})_\mathcal{U} \in B_{(X_n)_{\mathcal{U}}}$. Also, denote $F=\{\tilde x_1, \ldots \tilde x_m\}$. Because $x_{n,1}, \ldots, x_{n,m} \, \in \, F_n$, and $F_n$ is convex,  by \eqref{eq:diamond++} we have
$$
 \conv{F} \subset \left\{ x \in B_{(X_n)_{\mathcal{U}}} \colon \| G x + (y_n)_{\mathcal{U}}\| > 2 - \frac{\eps}{2} \right\}.
$$
Also, since $m$ is fixed, \eqref{eq.conv-aconv-apprm} implies
\begin{align*}
 \left\|Tx_0 - \sum_{k=1}^m \lambda_k \right.&\left.\exp\left(\frac{2\pi i k}{m}\right) T\tilde x_k \right\| =\\ &\lim_{n,\mathcal{U}} \left\|T_n x_{0,n} - \sum_{k=1}^m \lambda_{n,k} \exp\left(\frac{2\pi i k}{m}\right)T_n x_{n,k}\right\| <   \eps,
\end{align*}
that is
$$
\dist\bigl(Tx_0, T\bigl(\aconv(F)\bigr)\bigr) < \eps.
$$
Consequently, $F$ satisfies \eqref{eq:diamond} for $\eps$, $x_0$ and $y$, i.e.,  $F$ is the set we are looking for.
\end{proof}

In the case of ultrapowers the converse result is also true.

\begin{Theo}\label{Theor:ultra-target-invese}
Let $\mathcal{U}$ be a free ultrafilter on $\N$, let $X$, $Y$, $Z$ be Banach spaces, and let $G \in S_{L(X, Y)}$ and $T \in L(X, Z)$ be operators. If $T_{\mathcal{U}} = (T, T, \ldots)_{\mathcal{U}} \in L(X_{\mathcal{U}}, Z_{\mathcal{U}})$  is a target for $G_{\mathcal{U}}= (G, G, \ldots)_{\mathcal{U}} \in L(X_{\mathcal{U}}, Y_{\mathcal{U}})$, then $T$  is a target for $G$.
\end{Theo}

\begin{proof}
For given $x_{0} \in B_{X}$, $\eps > 0$ and $y \in S_{Y}$, we apply the definition of target to $\tilde x_{0} = (x_{0}, x_{0}, \ldots)_{\mathcal{U}} \in B_{X_{\mathcal{U}}}$ and $\tilde y = (y, y, \ldots)_{\mathcal{U}} \in S_{Y_{\mathcal{U}}}$. By \eqref{eq:diamond} in Definition \ref{Defi:target}, we can find $F = \{ (x_{1,n})_{\mathcal{U}}, \ldots, (x_{m,n})_{\mathcal{U}}\}\subset B_{X_{\mathcal{U}}}$, $\lambda_{1}, \ldots, \lambda_{m} \geqslant 0$ with $\sum{\lambda_{k}} = 1$, and $\theta_{1}, \ldots, \theta_{m}\in \mathbb{T}$ such that

\begin{equation}\label{eq:propDiamondAux1++}
\conv{F} \subset \{ \tilde x \in B_{X_{\mathcal{U}}} \colon \| G_{\mathcal{U}} \tilde x + \tilde y\| > 2 - \eps/2 \}.
\end{equation}
and
\begin{equation}\label{eq:propDiamondAux1++++}
\lim_{n,\mathcal{U}} \left\| T x_{0} - \sum_{k=1}^{m}{\lambda_{k} \theta_{k} T (x_{k,n})} \right\| < \eps. \end{equation}
Write $F_n =  \{ x_{1,n}, \ldots, x_{m,n}\}\subset B_X$ and denote by $E$ the set of those $n \in \N$ for which
$$
\conv{F_n} \subset \{  x \in B_{X} \colon \| G  x +  y\| > 2 - \eps \}.
$$
We claim that $E \in \mathcal{U}$. Indeed, if this is not so, then $\N \setminus E \in \mathcal{U}$. For every $n \in \N \setminus E$ choose $\mu_{1,n}, \ldots, \mu_{m,n} \geqslant 0$ with $\sum_{k=1}^m{\mu_{k,n}} = 1$ such that
$$
 \left\| G\left(\sum_{k=1}^{m} \mu_{k,n} x_{k,n}\right) + y \right\| \leqslant 2 - \eps.
$$
Then, for $\mu_k = \lim_{n,\mathcal{U}} \mu_{k,n}$ we have
$$
 \left\| G_{\mathcal{U}}\left(\sum_{k=1}^{m} \mu_{k} (x_{k,n})_{\mathcal{U}}\right) + \tilde y \right\| \leqslant 2 - \eps,
$$
 which contradicts \eqref {eq:propDiamondAux1++}.

 Now, since  $E \in \mathcal{U}$,  according to  \eqref{eq:propDiamondAux1++++} there is an $n_0 \in E$ such that
$$
 \left\| T x_{0} - \sum_{k=1}^{m}{\lambda_{k} \theta_{k} T (x_{k,n_0})} \right\| < \eps.
$$
The corresponding $F_{n_0}$ fulfills \eqref{eq:diamond} in Definition \ref{Defi:target}.
\end{proof}

Since lushness of an operator reduces to the fact that the identity operator is a target for it, we obtain the following two corollaries.

\begin{Coro}\label{Cor:ultra-lush}
Let $\mathcal{U}$ be a free ultrafilter on $\N$, $\{ X_n \}_{n\in \N}$, $\{ Y_n \}_{n\in \N}$ be sequences of Banach spaces, $\{G_n\}_{n\in \N}$ be a sequence of lush operators where $G_n \in S_{L(X_n, Y_n)}$ for every $n\in \N$. Then $G = (G_n)_{\mathcal{U}} \in L((X_n)_{\mathcal{U}}, (Y_n)_{\mathcal{U}})$ is lush.
\end{Coro}

\begin{Coro}\label{cor:ultra-lush-invese}
Let $\mathcal{U}$ be a free ultrafilter on $\N$, $X$, $Y$ be Banach spaces, $G \in S_{L(X, Y)}$. If  $G_{\mathcal{U}}= (G, G, \ldots)_{\mathcal{U}} \in L(X_{\mathcal{U}}, Y_{\mathcal{U}})$ is lush, then $G$ is lush.
\end{Coro}

\chapter{Open problems}\label{sec:OpenProblems}

Corresponding to \emph{Spear sets and spear vectors}:

\begin{Prob}
Let $X$ be a complex Banach space. If $\Spear(X)$ is not compact, does $X$ contain a copy of $c_{0}$ or $\ell_{1}$?
\end{Prob}

\begin{Prob}
If $X$ is a complex smooth Banach space and $\Spear(X^{\ast}) \neq \emptyset$, can we deduce that $X \cong \mathbb{C}$?
\end{Prob}

Corresponding to \emph{Lush operators}:

\begin{Prob} Are items (a) and (b) in Proposition \ref{Prop:lushSufficient} necessary for $G$ to be lush? If there is a counterexample, notice that the domain must be non separable.
\end{Prob}

Corresponding to \emph{Examples in classical Banach spaces}:

\begin{Prob}
Is lush the dual of the Fourier transform on $L_1$? Is lush the Fourier-Stieltjes transform?
\end{Prob}

\begin{Prob}
Is Proposition \ref{Prop-ejem-X-into-CK} always an equivalence?
\end{Prob}

\begin{Prob}
Is there a characterization of lush or spear operators acting from an $L_1(\mu)$ space analogous to the one given in Theorem \ref{Theo:L1->XaDP} for the aDP?
\end{Prob}

Corresponding to \emph{Further results}:

\begin{Prob}
Are spearness and lushness equivalent when the codomain space is SCD? The aDP and spearness are, see Remark \ref{Rema:aDP-spear-equivalent-Y-SCD}.
\end{Prob}

\begin{Prob}
Are the aDP and lushness equivalent when the image of the operator is Asplund? They are equivalent when the codomain is Asplund, see Proposition \ref{Prop:characLushAsplund}.
\end{Prob}

\begin{Prob}
Can the results about rank-one operators be extended to finite-rank operators? They are Corollary \ref{coro:charac-lush-rank-one} and Proposition \ref{Prop:rank-one-GaDP=>G*lush}.
\end{Prob}

\begin{Prob}
If $G:X\longrightarrow Y$ is lush and $Y$ is $L$-embedded, is $G^\ast$ lush? This is true for spearness and the aDP (see Proposition \ref{Prop-Lembedded-duality}).
\end{Prob}

Corresponding to \emph{Isomorphic and isometric consequences}:

\begin{Prob}
Is Theorem \ref{Theo:G-finite-rank-noell_1} valid in the complex case? That is, does the dual of the domain of a complex operator with the aDP always contain $\ell_1$?
\end{Prob}

\begin{Prob} Let
$G:X\longrightarrow Y$ be an operator with the aDP. Does $X=\K$ if $X$ is strictly convex or smooth? Does $Y=\K$ if $Y$ is strictly convex or smooth?
\end{Prob}

\begin{Prob}
Does every spear operator attain its norm? This is true for lush operators, see Proposition \ref{Prop:lush-attains-norm}, but it is not true for operators with the aDP, see Example \ref{example-aDP-no-NA}.
\end{Prob}

Corresponding to \emph{Stability results}:

\begin{Prob}
Are there results about the relationship between spear and lush operators with quotients by the kernel of the operator analogous to the one given in Proposition \ref{Prop:G-X/kerG-aDP} for the aDP?
\end{Prob}

\cleardoublepage

\printindex
\cleardoublepage
\end{document}